\documentclass[12pt]{article}
\usepackage{amsgen,amsmath,amstext,amsbsy,amsopn,amssymb,latexsym}
\usepackage[usenames]{color}
\usepackage[colorlinks=true,linkcolor=blue,citecolor=blue,urlcolor=blue,breaklinks]{hyperref}
\usepackage{array}
\usepackage{cite}
\usepackage{multirow}% http://ctan.org/pkg/multirow
\usepackage{graphicx}% http://ctan.org/pkg/graphicx
\usepackage{amssymb}
\usepackage{tikz}
\usepackage{mathtools}
\usepackage{float}
\usepackage{multicol}
\usepackage[font={small}]{caption}
\usetikzlibrary{shapes.geometric, arrows}
\usetikzlibrary{er,positioning}
\usepackage{verbatim}
\usepackage{algorithm, algorithmicx, algpseudocode}

\DeclareMathOperator*{\argmin}{arg\,min}
\newtheorem{Corollary}{Corollary}
\newtheorem{Proposition}{Proposition}
\newtheorem{Lemma}{Lemma}

\newtheorem{Theorem}{Theorem}
 
\newtheorem{Remark}{Remark}

\newcommand{\R}{\mathbb{R}}
\newcommand{\E}{\mathbf{E}}

\newcommand{\cip}{\overset{p}{\to}}
\newcommand{\cid}{\overset{d}{\to}}

\newcommand{\dif}{\Delta(X_{i})}

\newcommand{\eval}{a_0}

\newcommand{\bel}{\begin{eqnarray}\label}
\newcommand{\eel}{\end{eqnarray}}
\newcommand{\bes}{\begin{eqnarray*}}
\newcommand{\ees}{\end{eqnarray*}}
\newcommand{\bei}{\begin{itemize}}
\newcommand{\eei}{\end{itemize}}
\newcommand{\beiftnt}{\begin{itemize}\footnotesize}

\newcommand{\aveg}{\frac{1}{n_{a}}\sum_{j\in \mathcal{I}_a}}
\newcommand{\I}{{\bf I}}
\newcommand{\1}{{\bf 1}}

\title{Decorrelated Local Linear Estimator: Inference for Non-linear Effects in High-dimensional Additive Models}
\author{Zijian Guo \quad Wei Yuan  \quad Cun-Hui Zhang}

\begin{document}

\maketitle

%\begin{frontmatter}
%\thanksref{T1}}
%\runtitle{Local Inference in Additive Models}
%\thankstext{T1}{The research of Z Guo was supported in part by the NSF DMS 1811857; The research of C Zhang was supported in part by the NSF Grants DMS-1513378, DMS-1721495 and IIS-1741390.}

%\begin{aug}

%\author{\fnms{Zijian} \snm{Guo}\ead[label=e1]{zijguo@stat.rutgers.edu}}
%\and
%\author{\fnms{Cun-Hui} \snm{Zhang}\ead[label=e2]{czhang@stat.rutgers.edu}}
%\runauthor{Z. Guo and C. Zhang}

%\affiliation{Rutgers, The State University of New Jersey}
%\address{Department of Statistics\\
%Rutgers, The State University of New Jersey\\
%Piscataway, NJ 08854-8019\\
%USA\\
%\printead{e1}\\
%\phantom{E-mail:\ }\printead*{e2}\\
%\printead*{u1}\phantom{URL:\ }\\
%}

%\end{aug}

\begin{abstract}
Additive models play an essential role in studying non-linear relationships. Despite many recent advances in estimation, there is a lack of methods and theories for inference in high-dimensional additive models, including confidence interval construction and hypothesis testing. Motivated by inference for non-linear treatment effects, we consider the high-dimensional additive model and make inference for the derivative of the function of interest. We propose a novel decorrelated local linear estimator and establish its asymptotic normality. The main novelty is the construction of the decorrelation weights, which is instrumental in reducing the error inherited from estimating the nuisance functions in the high-dimensional additive model. We construct the confidence interval for the function derivative and conduct the related hypothesis testing. We demonstrate our proposed method over large-scale simulation studies and apply it to identify non-linear effects in the motif regression problem. Our proposed method is implemented in the R package \texttt{DLL} available from CRAN. 

%The asymptotic variance of our proposed estimator matches with the optimal rate in the univariate setting. 
\end{abstract}

%
%\begin{keyword}[class=MSC]
%\kwd[Primary ]{62G08}
%\kwd[; secondary ]{62G10}
%\end{keyword}
%
%\begin{keyword}
%\kwd{High dimension; Decorrelation; Double Estimation Accuracy; Derivative; Exposure Effect; Extreme value location.}
%\end{keyword}
%\end{frontmatter}

%%%%%%%%%%%%%%%%%%%%%%%%%%%%%%%%%%%%%%%%%%%%%%%%
\section{Introduction}
%%%%%%%%%%%%%%%%%%%%%%%%%%%%%%%%%%%%%%%%%%%%%%%%
Additive models play an important role in modern data analysis \cite{buja1989linear,wood2017generalized,hastie1986generalized}. The additive model is useful as it relaxes the stringent linearity assumption imposed in the multiple linear models and generalizes the nice interpretation of linear models. 
In the low-dimensional setting, additive models have been carefully investigated  \cite[e.g.]{buja1989linear,wood2017generalized,hastie1986generalized,opsomer2000asymptotic,horowitz2004nonparametric,mammen1999existence}. 
Recently, there has been a growing interest in the high-dimensional additive model, which generalizes the high-dimensional linear regression. Much progress has been made to understand the prediction performance of various proposals, {including \cite{meier2009high,koltchinskii2010sparsity,raskutti2012minimax, suzuki2013fast,tan2017penalized,yang2015minimax,yuan2016minimax,ravikumar2009sparse}.} However, the statistical inference problem in the high-dimensional additive model is far less understood from both methodological and theoretical perspectives.

Statistical inference in high-dimensional additive models is well-motivated from causal inference with observational studies. Causal conclusions from observational studies are invalidated due to unmeasured confounders \cite[e.g.]{imbens2015causal,pearl2009causality}. A commonly used approach is to condition on a large number of measured covariates such that the conditional ignorability condition holds \cite[e.g.]{belloni2017program,hernan2010causal}. This idea has been carefully investigated by utilizing high-dimensional linear models. However, the linear model imposes a stringent assumption that the exposure has a constant effect on the outcome regardless of the exposure value. Such an assumption might not be plausible for various applications; see the application to motif regression in Section \ref{sec: real data} for example. Non-linear effects have been commonly observed in scientific studies, including return to schooling \cite{card2001estimating}, climate on crop yields \cite{schlenker2008estimating}, and the climate change on the economic outcomes \cite{deschenes2012economic,dell2014we}. The additive model significantly relaxes the linearity assumption and better accommodates the possibly non-linear effect.

In this paper, we consider the additive model for the outcome variable $Y_i\in \R$, 
\begin{equation}
Y_i=f(D_i)+g(X_{i})+\epsilon_i \quad \text{with}\quad g(X_{i})=\sum_{j=1}^{p}g_j(X_{i,j}), \quad \text{for}\quad 1\leq i\leq n,
\label{eq: additive model}
\end{equation} 
where $D_i\in \R$ is the variable of interest (e.g. the exposure or treatment variable), $X_i\in \R^{p}$ is the high-dimensional baseline covariates, $\E(\epsilon_i\mid D_i, X_i)=0,$ and $f:\R \rightarrow \R$ and $g_j: \R \rightarrow \R$ for $1\leq j\leq p$ are unknown functions.
The observed data $\{Y_i, D_i, X_{i}\}_{1\leq i\leq n}$ are assumed to be independently and identically distributed. For a pre-specified $\eval\in \R$ and a small $\tau>0,$ the ratio $\left(f(\eval+\tau)-f(\eval)\right)/\tau$ captures the exposure's effect at $\eval$. With $\tau$ approaching zero, the function derivative $f'(\eval)$ measures the exposure's effect on the outcome \cite{belloni2015some}. The current paper is focused on statistical inference for $f'(\eval)$ with $\eval\in \R$ under the high-dimensional additive model \eqref{eq: additive model}. 

%%%%%%%%%%%%%%%%%%%%%%%%%%%%%%%%%%%%%%%%%%%%%%%%%
\subsection{Results and contributions}
%%%%%%%%%%%%%%%%%%%%%%%%%%%%%%%%%%%%%%%%%%%%%%%%%
In the univariate setting, the local linear estimator is the state-of-the-art method to make inference for $f'(\eval)$ \cite[e.g.]{fan1996local,fan1993local}. However, the inference problem under the high-dimensional additive model \eqref{eq: additive model} is much more challenging due to the presence of the unknown high-dimensional function $g$. With an accurate estimator $\widehat{g}$, the estimated variables $\widehat{Y}_i=Y_i-\widehat{g}(X_i)$ for $1\leq i\leq n$ can be used as proxies for $\{f(D_i)\}_{1\leq i \leq n}$. A natural idea is to estimate $f'(\eval)$ by applying the local linear method to $\{D_i, \widehat{Y}_i\}_{1\leq i\leq n}$ with $\{\widehat{Y}_i\}_{1\leq i\leq n}$ as the outcome variables. However, such a plug-in estimator suffers from large estimation bias due to the estimation error of $\widehat{g}$; see Table \ref{table:1} in Section \ref{sec: sim result 1} for illustrations.

We propose a novel Decorrelated Local Linear (\texttt{DLL}) estimator of $f'(\eval).$ The classical local linear estimator can be expressed as a weighted average of the outcome variables, where the local linear kernel induces the weights. As the major novelty, we construct the {\it decorrelation weights} to mitigate the error inherited from the high-dimensional estimator $\widehat{g}$. Meanwhile, the constructed decorrelation weights ensure that the standard error of our proposed \texttt{DLL} estimator is comparable to that of the classical local linear estimator. The decorrelation weights are constructed in a non-parametric way and designed for the bias correction of the local linear estimator. 

In Theorem \ref{thm: limiting swap}, we establish the asymptotic normality of our proposed \texttt{DLL} estimator as long as the estimator $\widehat{g}$ is consistent. We further show that the asymptotic variance of our proposed estimator matches with the optimal rate in the univariate setting \cite{fan1996local}. We construct the confidence interval for $f'(\eval)$ and test for the hypothesis $H_0: f'(\eval)=0$. %the standard error of our proposed \texttt{DLL} estimator matches with the optimal standard error of the function derivative estimation in the univariate setting

In Section \ref{sec: simulation}, we demonstrate the validity of our theoretical results in moderate sample sizes, address practical issues on algorithm implementation, and provide practical recommendations.
Our proposed method is implemented in the R package \texttt{DLL}, which is available from CRAN. The simulation results show that the \texttt{DLL} estimator significantly outperforms the plug-in estimator and the ReSmoothing estimator \cite{gregory2016optimal}, in terms of the bias correction and coverage property. Regarding the empirical coverage and length, the confidence intervals (CIs) based on the \texttt{DLL} estimator are comparable to the oracle CIs, which are constructed with the oracle knowledge of the high-dimensional function $g$. 

In Section \ref{sec: real data}, we conduct a careful analysis of the motif regression problem \cite{yuan2007motif} and observe a highly non-linear relationship between the gene expression level and the motif scores. Our results demonstrate the advantage of our proposed method over the statistical inference method assuming the linear outcome model.

%%%%%%%%%%%%%%%%%%%%%%%%%%%%%%%%%%%%%%%%%%%%%%%%%
\subsection{Literature review and comparison}
%%%%%%%%%%%%%%%%%%%%%%%%%%%%%%%%%%%%%%%%%%%%%%%%%
Two recent works \cite{gregory2016optimal} and \cite{lu2015kernel} studied the inference problems in high-dimensional additive models. Specifically,  \cite{gregory2016optimal} proposed a two-step ReSmoothing (\texttt{RS}) estimator: in the first step, a pre-smoothing estimator was obtained; in the second step, the pre-smoothing estimator was taken as the proxy outcome, and standard univariate non-parametric technique was then applied. In Section \ref{sec: comparison RS}, we compare our proposed \texttt{DLL} estimator with the \texttt{RS} estimator  and observe that the \texttt{RS} estimator suffers from a large bias of estimating the function derivative while our proposed \texttt{DLL} estimator corrects the bias effectively. Consequently, our proposed confidence interval has better empirical coverage than that based on the \texttt{RS} estimator; see Table \ref{table:4} for the detailed comparison. In addition, \cite{lu2015kernel} considered the confidence band construction problem under the high-dimensional sparse additive model, which is a different inference problem from the current paper. %The method proposed in \cite{lu2015kernel} is to approximate the non-parametric function by a set of basis functions and then apply the debiasing method for the corresponding linear model of the basis functions. 

Inference for function derivative has been actively studied in the non-parametric modeling, including local linear estimator \cite{fan1993local}, regression spline \cite{zhou2000derivative}, kernel methods \cite{gasser1984estimating}, empirical likelihood methods \cite{qin2005empirical}, and others cited therein. However, as discussed, the unknown high-dimensional function $g$ in the additive model poses great challenges to statistical inference for the function derivative at a local point. The paper \cite{belloni2015some} studied the inference for the function derivative in additive models without the sparsity structure. The penalty is essential to recovering the high-dimensional sparse model, which creates an additional bias to correct in the following inference step. The statistical inference problem with the sparsity structure requires extra innovation in terms of both method and theory.

A recent line of active research was focused on statistical inference in high-dimensional linear regression. Debiased estimators or Neyman's Orthogonalization were proposed for inference for single regression coefficients \cite{zhang2014confidence,van2014asymptotically,javanmard2014confidence,belloni2014inference,chernozhukov2015valid,farrell2015robust,chernozhukov2018double}. 
The linear model is a special case of the additive model, where the function derivative $f'(\eval)$ is assumed to be a constant for any $\eval\in \R$. Statistical inference for the non-linear effect in the additive model is a much more challenging problem, which requires novel methods and theories to address the non-linearity. Both the rate of convergence and the sufficient conditions for confidence interval construction are different from those established in the high-dimensional linear regression. A more detailed methodological comparison is presented in Remark \ref{rem: comparison}. The real data analysis in Section \ref{sec: real data} shows that a misleading scientific conclusion might be obtained without accounting for the possible non-linear effects.

Beyond the high-dimensional linear regression, \cite{chernozhukov2018double} and \cite{zhu2017high} studied the inference procedure for the partially linear model. However, the focus is still on the inference problem for the linear component instead of the non-linear component addressed here.

\noindent {\bf Organization.} In Section \ref{sec: DLL}, we introduce the decorrelated local linear estimator. In Section \ref{sec: theory}, we establish the theoretical guarantee of the proposed estimator. In Section \ref{sec: simulation}, we conduct a large-scale simulation study to demonstrate the finite-sample performance of the \texttt{DLL} estimator. In Section \ref{sec: real data}, we apply the \texttt{DLL} estimator to the motif regression problem.  In Section \ref{sec: con-discussion}, we provide conclusion and discussion. %In Section \ref{sec: nuisance error}, we provide the technical analysis to illustrate the effect of decorrelation. Additional proofs are presented as supplementary materials. 

\vspace{2mm}

\noindent {\bf Notations.} For a sequence of random variables $X_n$ indexed by $n$, we use $X_n \cip X$ and $X_n \cid X$ to represent that $X_n$ converges to $X$ in probability and in distribution, respectively. For a matrix $X$, we use $X_{i,j}$, $X_{i}$, and $X_{\cdot,j}$ to
denote its $(i,j)$ entry, the $i$-th row and $j$-th column, respectively; for index sets $S_1$ and $S_2$, $X_{S_1,S_2}$ denotes the sub-matrix of $X$ with row and column indices belonging to $S_1$ and $S_2$, respectively. We use $c$ and $C$ to denote generic positive constants that may vary from place to place. For two positive sequences $a_n$ and $b_n$,  $a_n \lesssim b_n$ means $a_n \leq C b_n$ for all $n$ and $a_n \gtrsim b_n $ if $b_n\lesssim  a_n$ and $a_n \asymp b_n $ if $a_n \lesssim b_n$ and $b_n \lesssim a_n$, and $a_n \ll b_n$ if $\limsup_{n\rightarrow\infty} {a_n}/{b_n}=0$ and $a_n \gg b_n$ if $b_n \ll a_n$. 

\section{Decorrelated local linear estimator}
\label{sec: DLL}
%%%%%%%%%%%%%%%%%%%%%%%%%%%%%%
%\Zijian{Stops here.}
We consider the data $\{X_i, D_i, Y_i\}_{1\leq i\leq n}$ being i.i.d. generated, where for the $i$-th subject, $Y_i\in \R$ denotes the outcome variable, $D_i\in \R$ denotes the variable of interest, and $X_i\in \R^{p}$ denotes the high-dimensional baseline covariates.  
We focus on the additive outcome model \eqref{eq: additive model}. Our goal is to make inference for the function derivative $f'(\eval)$ with $\eval\in \R$ denoting a pre-specified value belonging to the range of $D_i.$

%%%%%%%%%%%%%%%%%%%%%%%%%%%%%
%%%%%%%%%%%%%%%%%%%%%%%%%%%%%
\subsection{Decorrelation for the local linear estimator}
\label{sec: decor}
%%%%%%%%%%%%%%%%%%%%%%%%%%%%%%%%%%%%%%%%%%%%%%%%%

The local polynomial estimator \cite{fan1993local,stone1977consistent,cleveland1979robust,fan1992design} has been developed under the univariate non-parametric regression $
Y_i=f(D_i)+\epsilon_i,$ for $1\leq i\leq n,$
with $f: \R \rightarrow \R$ denoting an unknown smooth function. By the Taylor expansion %of $f_1(x)$ near $\eval$,  
$
f(D_i)=f(\eval)+{f'}(\eval) (D_i-\eval)+r(D_i)$ for $\eval-h\leq D_i\leq \eval+h,$ we approximate $f$ by a linear function in a small neighborhood of $\eval$ and estimate $f'(\eval)$ by fitting a linear model within this small neighborhood.  
For a pre-specified bandwidth $h>0$, define the kernel \begin{equation}
K_{h}(d)= \frac{1}{2h}\cdot \mathbf{1}\left(|D_i-\eval |\leq h\right).
\label{eq: kernel}
\end{equation} The local linear estimator of $f'(\eval)$  
has the following explicit form, 
\begin{equation}
\frac{\sum_{i=1}^{n} W^0_i Y_i K_{h}(D_i)}{\sum_{i=1}^{n} W^0_i (D_i-\eval)K_{h}(D_i)}, \quad \text{with}\quad W^0_i=(D_i-\eval)-\frac{\sum_{j=1}^{n} (D_j-\eval)K_{h}(D_j)}{\sum_{j=1}^{n} K_{h}(D_j)}.
\label{eq: 1dim local}
\end{equation}
The local polynomial estimator is the state-of-the-art method for estimating $f'(\eval)$ in the univariate and also low dimensional setting. To illustrate the main idea, we focus on the uniform kernel and the local linear estimator, but our following proposed method is potentially useful for other kernels and higher-order local polynomial estimator.

 For the high-dimensional sparse additive model in \eqref{eq: additive model}, the existing literature \cite[e.g.]{meier2009high,koltchinskii2010sparsity,suzuki2013fast,tan2017penalized} was focused on accurately estimating the unknown functions $f$ and $g$ in \eqref{eq: additive model}. However, there is a lack of inference methods for $f'(\eval)$.

In the following, we introduce the decorrelation idea and use $\widehat{g}$ to denote an initial estimator of ${g}$; see Section \ref{sec: initial estimator} for the detailed construction. With the estimator $\widehat{g}$, we compute 
$\widehat{Y}_i=Y_i-\widehat{g}(X_{i})$ for $1\leq i\leq n,$ which are proxies for the oracle outcome $Y_i^{\rm ora}=f(D_i)+\epsilon_i.$ As a direct extension of the local linear estimator in \eqref{eq: 1dim local}, we replace $Y_i$ by $\widehat{Y}_i$ and have the following plug-in estimator, 
\begin{equation}
\widetilde{f'(\eval)}=\frac{\sum_{i=1}^{n} W^0_i \widehat{Y}_i K_{h}(D_i)}{\sum_{i=1}^{n} W^0_i (D_i-\eval)K_{h}(D_i)}.
\label{eq: plug-in}
\end{equation}
$\widetilde{f'(\eval)}$ is the same as the local linear estimator applied to the data $\{D_i, \widehat{Y}_i\}_{1\leq i\leq n},$ with $\widehat{Y}_i$ as the outcome variable.  
The simulation results in Section \ref{sec: sim result 1} demonstrate that the plug-in estimator $\widetilde{f'(\eval)}$ suffers from a large bias due to the estimation error of $\widehat{g}$. Consequently, the plug-in estimator is not ready for statistical inference; see Table \ref{table:1} for details. 

To correct the bias of the plug-in estimator, we consider the estimator of the form, \begin{equation}
\widehat{f'(\eval)}=\frac{\frac{1}{n}\sum_{i=1}^{n}W_{i}\widehat{Y}_i K_{h}(D_i)}{\frac{1}{n}\sum_{i=1}^{n}W_{i}(D_i-\eval)K_{h}(D_i)},
\label{eq: generic estimator}
\end{equation}
where $\{{W}_{i}\in \R\}_{1\leq i\leq n}$ are the weights to be specified. The error $\widehat{f'(\eval)}-f'(\eval)$ is decomposed as, 
\begin{equation}
\begin{aligned}
\frac{\frac{1}{n}\sum_{i=1}^{n}{W}_{i}[f(\eval)+r(D_i)+\epsilon_i]K_{h}(D_i)}{\frac{1}{n}\sum_{i=1}^{n}{W}_{i}(D_i-\eval)K_{h}(D_i)}+\frac{\frac{1}{n}\sum_{i=1}^{n}{W}_{i}[\widehat{g}(X_{i})-{g}(X_{i})] K_{h}(D_i)}{\frac{1}{n}\sum_{i=1}^{n}{W}_{i}(D_i-\eval)K_{h}(D_i)}. %+\frac{\frac{1}{n}\sum_{i=1}^{n}{W}_{i}\left[f(\eval)+\right]K_{h}(D_i)}{\frac{1}{n}\sum_{i=1}^{n}{W}_{i}(D_i-\eval)K_{h}(D_i)}. 
\end{aligned}
\label{eq: decomposition 0}
\end{equation}
The first term of \eqref{eq: decomposition 0} appears in the univariate case, while the second term is the new addition in the high-dimensional setting, which results from the error of estimating the high-dimensional function $g$. We construct the weights $\{{W}_{i}\}_{1\leq i\leq n}$ such that the second term of \eqref{eq: decomposition 0} is significantly reduced, but the first term in \eqref{eq: decomposition 0} is of the same scale as the univariate case. 

To achieve this, we define the population decorrelation weights $\{W_i\}_{1\leq i\leq n}$ as
\begin{equation} 
{W}_{i}=(D_i-\eval)-l(X_{i}) \quad \text{with}\quad l(X_{i})\coloneqq \frac{\E \left([D_i-\eval]K_{h}(D_i)|X_{i} \right)}{\E \left(K_{h}(D_i)|X_{i} \right)}.
\label{eq: expression}
\end{equation}
The population decorrelation weights ensure that the estimator $\widehat{f'(\eval)}$ in \eqref{eq: generic estimator} is nearly unbiased and asymptotically normal. In the following Section \ref{sec: nonpara decor}, we propose a non-parametric estimator of the decorrelation weight $W_i$ defined in \eqref{eq: expression}.

We provide intuitions on how the population decorrelation weights significantly reduce the second term in \eqref{eq: decomposition 0}.
The weights $\{W_i\}_{1\leq i\leq n}$ constructed in \eqref{eq: expression} guarantee \begin{equation*}
\E \left[{W}_{i}K_{h}(D_i)\mid X_i\right]=0.
\end{equation*}
 If $(D_i,X_{i}^\intercal,Y_i)^{\intercal}$ is not used to construct $\widehat{g}$, then the above equation implies
\begin{equation}
\E \left[{W}_{i}(\widehat{g}(X_{i})-{g}(X_{i})) K_{h}(D_i)\mid X_i\right]=0.
\label{eq: key decorrelation property}
\end{equation}
The zero mean of the estimation error ${W}_{i}(\widehat{g}(X_{i})-{g}(X_{i})) K_{h}(D_i)$ guarantees the second term of the decomposition \eqref{eq: decomposition 0} converges to zero at a fast rate. 

%%%%%%%%%%%%%%%%%%%%%%%%%%%%%%%%
\subsection{Construction of decorrelation weights}
\label{sec: nonpara decor}
%%%%%%%%%%%%%%%%%%%%%%%%%%%%%%%%%%%%%%%%%%%%%%%%%%%%%%%%%%%%%
In the following, we construct non-parametric estimators of $l(X_i)$ and ${W}_{i}$ defined in \eqref{eq: expression}.  
We decouple the relationship between $D_i$ and $X_i$ by considering a high-dimensional sparse linear model,
\begin{equation}
D_i=X_{i}^{\intercal}\gamma+\delta_i, \quad \text{for}\; 1\leq i\leq n, 
\label{eq: decouple}
\end{equation}
where $\gamma$ is a sparse vector and $\delta_i$ is independent of $X_{i}$. Let $\phi(\delta)$ denote the density function of $\delta_i.$ Under the model \eqref{eq: decouple}, we obtain the following expression for ${l}(X_{i})$,
\begin{equation}
{l}(X_{i})=\frac{ \int_{\mu_i-h}^{\mu_i+h} \left(\delta-\mu_i\right)\phi(\delta) d \delta}{ \int_{\mu_i-h}^{\mu_i+h}\phi(\delta) d\delta} \quad \text{with} \quad \mu_i=\eval-X_{i}^{\intercal}\gamma, \quad \text{for} \quad 1\leq i\leq n.
\label{eq: l expression}
\end{equation}
We also write ${l}(X_{i},\gamma)$ for ${l}(X_{i})$ to highlight its dependence on $\gamma$. As a remark, the independence between $\delta_i$ and $X_i$ is important for establishing the simplified expression of ${l}(X_{i})$ in \eqref{eq: l expression}.
In Section \ref{sec: DLL-S}, we consider possibly non-linear relationship between $D_i$ and $X_i $ and generalize our procedure to the setting where the model between $D_i$ and $X_i$ is a sparse additive model. In Section \ref{sec: sim result 1}, we demonstrate the robustness of our proposed method in finite samples when the independence assumption in the model \eqref{eq: decouple} does not hold; {see Settings 3 and 4 in Section \ref{sec: sim result 1} and Table \ref{table:1} for details.}

We utilize the expression \eqref{eq: l expression} and construct a non-parametric estimator of $l(X_{i})$. We will use sample splitting to create independence required for establishing the decorrelation property in \eqref{eq: key decorrelation property}. Our particular way of data splitting will not lead to the efficiency loss; see Remark \ref{rem: data swapping} for details. 

We randomly split the index set $\{1,2,\cdots,n\}$ into two disjoint subsets $\mathcal{I}_a$ and $\mathcal{I}_b$, with $\mathcal{I}_a\cup \mathcal{I}_b=\{1,2,\cdots,n\},$ $|\mathcal{I}_a|=\lfloor n/2\rfloor,$ and $|\mathcal{I}_b|=n-\lfloor n/2\rfloor.$ With the data $\{Y_i, D_i, X_i\}_{i\in \mathcal{I}_a},$ we estimate $\gamma$ by the Lasso estimator $\widehat{\gamma}^{a}\in \R^{p}$, defined as
\begin{equation}
\widehat{\gamma}^{a}=\argmin_{\gamma\in \R^{p}}\frac{1}{2|\mathcal{I}_a|}\sum_{i \in \mathcal{I}_a}\left(D_i-X_{i}^{\intercal}\gamma\right)^2+\lambda_{1} \sum_{j=1}^{p} \frac{\|X_{\mathcal{I}_a, j}\|_2}{\sqrt{n_a}}|\gamma_j|,
\label{eq: SLasso}
\end{equation} 
where $\lambda_1>0$ is a tuning parameter. We estimate $\{\mu_i\}_{i\in \mathcal{I}_b}$ and $\{\delta_i\}_{i\in \mathcal{I}_b}$ by 
\begin{equation*}
\widehat{\mu}_i=\eval-X_{i}^{\intercal}\widehat{\gamma}^{a} \quad \text{and}\quad \widehat{\delta}_i=D_i-X_{i}^{\intercal}\widehat{\gamma}^{a} \quad \text{for}\quad i\in \mathcal{I}_b.
%\label{eq: first stage residue}
\end{equation*}
For $i\in \mathcal{I}_b,$ we respectively estimate $\int_{\mu_i-h}^{\mu_i+h} \left(\delta-\mu_i\right)\phi(\delta) d \delta$ and $\int_{\mu_i-h}^{\mu_i+h} \phi(\delta) d \delta$ by $$\frac{1}{|\mathcal{I}_b|} \sum_{j\in \mathcal{I}_{b}} (\widehat{\delta}_j-\widehat{\mu}_i) {\bf{1}}(|\widehat{\delta}_j-\widehat{\mu}_i|\leq h) \quad \text{and}\quad \frac{1}{|\mathcal{I}_b|} \sum_{j\in \mathcal{I}_{b}} {\bf{1}}(|\widehat{\delta}_j-\widehat{\mu}_i|\leq h).$$
Then we estimate ${l}(X_{i})$ by 
\begin{equation}
\widehat{l}(X_{i},\widehat{\gamma}^{a})=\frac{ \sum_{j\in \mathcal{I}_{b}} (\widehat{\delta}_j-\widehat{\mu}_i) {\bf{1}}(|\widehat{\delta}_j-\widehat{\mu}_i|\leq h)}{ \sum_{j\in \mathcal{I}_{b}} {\bf{1}}(|\widehat{\delta}_j-\widehat{\mu}_i|\leq h)} \quad \text{for}\quad i\in \mathcal{I}_b. 
\label{eq: estimator}
\end{equation}
Our above construction ensures that the estimator $\widehat{\gamma}^{a}$ is independent of the data points $\{D_i,X_i\}_{i\in \mathcal{I}_b}.$
We can construct the estimators of $\{l(X_i)\}_{i\in \mathcal{I}_a}$ in a similar way to \eqref{eq: estimator} by switching the roles of $\mathcal{I}_a$ and $\mathcal{I}_b.$
We construct the estimator $\widehat{\gamma}^{b}\in \R^{p}$ by applying the Lasso algorithm in \eqref{eq: SLasso} to the data $\{Y_i, D_i, X_i\}_{i\in \mathcal{I}_b}$. 
We estimate $\{\mu_i,\delta_i\}_{i\in \mathcal{I}_a}$ by $$\widehat{\mu}_i=\eval-X_{i}^{\intercal}\widehat{\gamma}^{b} \quad \text{and}\quad \widehat{\delta}_i=D_i-X_{i}^{\intercal}\widehat{\gamma}^{b} \quad \text{for}\quad i\in \mathcal{I}_a,$$
and then estimate $\{{l}(X_{i})\}_{i\in \mathcal{I}_a}$ by 
\begin{equation}
\widehat{l}(X_{i},\widehat{\gamma}^{b})=\frac{\sum_{j\in \mathcal{I}_a}(\widehat{\delta}_j-\widehat{\mu}_i) {\bf{1}}(|\widehat{\delta}_j-\widehat{\mu}_i|\leq h)}{\sum_{j\in \mathcal{I}_a} {\bf{1}}(|\widehat{\delta}_j-\widehat{\mu}_i|\leq h)} \quad \text{for}\quad i\in \mathcal{I}_a.
\label{eq: estimator sym}
\end{equation}
Then we apply the definition in \eqref{eq: expression} and define 
\begin{equation}
\widetilde{W}_{i} = (D_i-\eval)-\widehat{l}(X_i) \quad \text{with}\quad
\widehat{l}(X_i)= \begin{cases} 
\widehat{l}(X_{i},\widehat{\gamma}^{b})& \text{for } i \in \mathcal{I}_a \\
\widehat{l}(X_{i},\widehat{\gamma}^{a})& \text{for } i \in \mathcal{I}_b \\
 \end{cases},
 \label{eq: swapping est}
\end{equation}
where $\widehat{l}(X_{i},\widehat{\gamma}^{a})$ and $\widehat{l}(X_{i},\widehat{\gamma}^{b})$ are defined in \eqref{eq: estimator} and \eqref{eq: estimator sym}, respectively. By centering $\{\widetilde{W}_i\}_{1\leq i\leq n}$, we construct the decorrelation weights as
\begin{equation}
\widehat{W}_{i}=\widetilde{W}_{i}-[\sum_{j=1}^{n}\widetilde{W}_{j}K_{h}(D_j)]/[\sum_{j=1}^{n}K_{h}(D_j)] \quad \text{for}\quad 1\leq i\leq n.
\label{eq: final D}
\end{equation}
With the data $\{Y_i, D_i, X_i\}_{i\in \mathcal{I}_a},$ we construct the initial estimator $\widehat{g}^{a}(\cdot)$ of $g(\cdot)$ in the following equations \eqref{eq: SAM} and \eqref{eq: prediction def}; we construct the estimators $\widehat{g}^{b}(\cdot)$ by applying the same algorithm to the data $\{Y_i, D_i, X_i\}_{i\in \mathcal{I}_b}.$ 
For $1\leq i\leq n,$ we compute
\begin{equation}
\widehat{Y}_i=Y_i-\widehat{g}(X_{i}) \quad \text{with}\quad 
\widehat{g}(X_{i}) = 
 \begin{cases} 
\widehat{g}^{b}(X_{i})& \text{for } i \in \mathcal{I}_a \\
\widehat{g}^{a}(X_{i})& \text{for } i \in \mathcal{I}_b \\
 \end{cases}.
 \label{eq: swapping est nuisance}
\end{equation}

\begin{Remark} \rm
The construction in \eqref{eq: swapping est} and \eqref{eq: swapping est nuisance} uses the ``data swapping" idea. That is, we {swap} the data and the initial estimators. Such a procedure is used to create the independence required for the proof but does not lead to loss of efficiency. The data swapping idea dated at least back to \cite{schick1986asymptotically,klaassen1987consistent} and was recently developed under the name of ``cross-fitting" in the double machine learning literature \cite[e.g.]{chernozhukov2018double}. 
\label{rem: data swapping}
\end{Remark}

%%%%%%%%%%%%%%%%%%%%%%%%%%%%%%%%%%%%%%%%%%%%%%%%%%%%%%%%%%%%%%%%%%%%%%%%%%%%%%
\subsection{Decorrelated local linear estimator and inference for $f'(\eval)$}
\label{sec: inference method}
%%%%%%%%%%%%%%%%%%%%%%%%%%%%%%%%%%%%%%%%%%%%%%%%%%%%%%%%%%%%%%%%%%%%%%%%%%%%%%
By combining $\widehat{Y}_i$ defined in \eqref{eq: swapping est nuisance} and the decorrelation weight $\widehat{W}_{i}$ defined in \eqref{eq: final D}, we apply the generic form \eqref{eq: generic estimator} and propose the Decorrelated Local Linear (\texttt{DLL}) estimator as 
\begin{equation}
\widehat{f'(\eval)}={\frac{1}{n\widehat{S}_{n}}\sum_{i=1}^{n}\widehat{W}_{i}\widehat{Y}_i K_{h}(D_i)} \quad \text{where} \quad \widehat{S}_{n}=\frac{1}{n}\sum_{i=1}^{n}\widehat{W}_{i}(D_i-\eval)K_{h}(D_i).
\label{eq: final estimator}
\end{equation}
In Section \ref{sec: theory}, we show that $\widehat{f'(\eval)}-f'(\eval)$ is asymptotically normal if certain reasonably good estimator $\widehat{g}$ is used in our construction. Consequently, we construct the following $1-\alpha$ confidence interval for $f'(\eval),$ 
\begin{equation}
{\rm CI}[f'(\eval)]=\left(\widehat{f'(\eval)}-z_{\alpha/2}\sqrt{\widehat{\rm V}}, \widehat{f'(\eval)}+z_{\alpha/2}\sqrt{\widehat{\rm V}}\right) \quad \text{with}\quad 
\widehat{\rm V}=\frac{\widehat{\sigma}^2}{n^2\widehat{S}^2_{n}}{\sum_{i=1}^{n}\widehat{W}_{i}^2K_{h}^2(D_i)},
\label{eq: ci}
\end{equation}
where $z_{\alpha/2}$ denotes the upper $\alpha/2$ quantile of the standard normal distribution and $\widehat{\sigma}^2$ is the variance level estimator specified in Section \ref{sec: initial estimator}. To test the null hypothesis $H_0: f'(\eval)=0,$ we develop the following procedure,
\begin{equation}
\psi[f'(\eval)]={\bf 1}\left(|\widehat{f'(\eval)}|\geq z_{\alpha/2}\sqrt{\widehat{\rm V}}\right).
\label{eq: testing}
\end{equation}
\begin{Remark}[Comparison to debiasing methods in linear models.]\rm
The debiased inference methods have been proposed in 
\cite{zhang2014confidence,van2014asymptotically,javanmard2014confidence} about inference for the regression coefficients in high-dimensional regression models and extended to other high-dimensional parametric models \cite{ning2017general,van2014asymptotically} or other inference targets in high-dimensional linear regression \cite{cai2017confidence,cai2021optimal,athey2018approximate,zhu2018linear}. These methods ultilize the fact that $\widehat{\delta}$ is nearly orthogonal to $X$ and then correct the bias with a linear function of $\widehat{\delta}.$ In contrast, our proposed \texttt{DLL} estimator uses a non-linear transformation of $\widehat{\delta}$ to remove the high-dimensional error. Specifically, we construct the deccorelation weights based on certain kernel estimates with $\{\widehat{\delta}_i\}_{1\leq i\leq n};$ see \eqref{eq: estimator} and \eqref{eq: estimator sym}. This new decorrelation idea is particularly designed for bias correction of the local linear estimator.
\label{rem: comparison}
\end{Remark}

%%%%%%%%%%%%%%%%%%%%%%%%%%%%
\subsection{Initial estimators}
\label{sec: initial estimator}
%%%%%%%%%%%%%%%%%%%%%%%%%%%%n{}

We now specify the initial estimators $\widehat{g}$ and $\widehat{\sigma}^2$ used in the construction of the \texttt{DLL} estimator in \eqref{eq: final estimator} and the related confidence interval in \eqref{eq: ci}. In the existing literature \cite{meier2009high,ravikumar2009sparse,koltchinskii2010sparsity,suzuki2013fast,tan2017penalized}, different types of penalty terms are imposed to ensure that only a small number of the unknown functions $f$ and $\{g_j\}_{1\leq j\leq p}$ are non-zero and these non-zero functions are smooth. We adopt the basis method in \cite{meier2009high,ravikumar2009sparse} and generate a set of basis functions to approximate the smooth functions. In particular, for a positive integer $M,$ we use $\{\phi_{0,l}\}_{1\leq l\leq M}$ to denote a set of B-spline basis functions for $f$ and for $1\leq j\leq p$, $\{\phi_{j,l}\}_{1\leq l\leq M}$ to denote a set of B-spline basis functions for $g_j.$ We write $\Psi_{i,0}=(\phi_{0,1}(D_{i}),\cdots,\phi_{0,M}(D_i))\in \R^{M}$ and $\Psi_{i,j}=(\phi_{j,1}(X_{i,j}),\cdots,\phi_{j,M}(X_{i,j}))\in \R^{M}$ for $1\leq j\leq p.$ Following \cite{ravikumar2009sparse}, we implement the following convex optimization problem,
\begin{equation}
\{\widehat{\beta}^{a}_j\}_{0\leq j\leq p}=\argmin_{\beta_j \in \R^{M}, \; {0\leq j\leq p}}\frac{1}{2|\mathcal{I}_a|}\sum_{i\in \mathcal{I}_a}(Y_i-\sum_{j=0}^{p}\Psi_{i,j}^{\intercal}\beta_j)^2+\lambda\sum_{j=0}^{p}\sqrt{\beta_j^{\intercal} \left(\frac{1}{|\mathcal{I}_a|}\sum_{i\in \mathcal{I}_a} \Psi_{i,j} \Psi_{i,j}^{\intercal}\right) \beta_j},
\label{eq: SAM}
\end{equation}
where $\lambda>0$ is a tuning parameter to be chosen. The choice of the tuning parameter $\lambda>0$ and the number M of basis functions are discussed at the beginning of Section \ref{sec: simulation}. Define 
\begin{equation}
\widehat{g}^{a}(X_i)=\sum_{j=1}^{p}\Psi_{i,j}^{\intercal}\widehat{\beta}^{a}_j \quad \text{and}\quad \widehat{f}^{a}(D_i)=\Psi_{i,0}^{\intercal}\widehat{\beta}^{a}_0 \quad \text{for}\quad i\in \mathcal{I}_b.
\label{eq: prediction def}
\end{equation}
Similarly, we define $\{\widehat{\beta}^b_j\}_{0\leq j\leq p}$ as in \eqref{eq: SAM} by replacing $\mathcal{I}_a$ with $\mathcal{I}_b$ and  \begin{equation*}
\widehat{g}^{b}(X_i)=\sum_{j=1}^{p}\Psi_{i,j}^{\intercal}\widehat{\beta}^{b}_j \quad \text{and}\quad \widehat{f}^{b}(D_i)=\Psi_{i,0}^{\intercal}\widehat{\beta}^{b}_0 \quad \text{for}\quad i\in \mathcal{I}_a.
\end{equation*} 
Then we construct 
\begin{equation}
\widehat{g}(X_{i}) = 
 \begin{cases} 
\widehat{g}^{b}(X_{i})& \text{for } i \in \mathcal{I}_a \\
\widehat{g}^{a}(X_{i})& \text{for } i \in \mathcal{I}_b \\
 \end{cases} \quad \text{and} \quad \widehat{f}(D_{i}) = 
 \begin{cases} 
\widehat{f}^{b}(D_{i})& \text{for } i \in \mathcal{I}_a \\
\widehat{f}^{a}(D_{i})& \text{for } i \in \mathcal{I}_b \\
 \end{cases},
 \label{eq: initial estimator}
\end{equation}  
and estimate the variance level $\sigma^2$ by 
\begin{equation}
\widehat{\sigma}^2= \frac{1}{n} \sum_{i=1}^{n} [Y_i-\widehat{f}(D_i)-\widehat{g}(X_{i})]^2.
\label{eq: noise estimator}
\end{equation}
% As a remark, other basis functions might be used, and the basis size M is allowed to vary across different functions. \Wei{Do we emphasize here that M is the same for all functions when tuning the parameters?}
In addition to the aforementioned basis method, we can also adopt the double penalization method \cite{tan2017penalized} by generalizing the smoothing spline; see Section \ref{sec: double penalization} in the supplement. We also discuss the construction of additive models by firstly applying the quantile transformation to the observed covariates; see Section \ref{sec: quantile} in the supplement.

%%%%%%%%%%%%%%%%%%%%%%%%%%%%%%%
\subsection{Algorithm}
%%%%%%%%%%%%%%%%%%%%%%%%%%%%%%%%%%%%%%%%%%%%%%%%%%%%%%%%%%%%
We summarize our proposed \texttt{DLL} estimator (with data swapping) in Algorithm \ref{algo: DDL} and will present the \texttt{DLL} estimator without data swapping in Section \ref{sec: no data swapping} in the supplement. We shall discuss the tuning parameter selection at the beginning of Section \ref{sec: simulation}. 
\begin{algorithm}[htp!]
\caption{Decorrelated Local Linear (\texttt{DLL}) Estimator}
\begin{flushleft}
%\hspace*{\algorithmicindent} 
%\hspace*{\algorithmicindent} 
 \textbf{Input:} Data $X\in \R^{n\times p}, D \in \R^{n}, Y\in \R^{n}$; the evaluation point $\eval\in \R$, bandwidth $h$, tuning parameters $\lambda, \lambda_1>0$, the number of basis $M.$ \\
%\hspace*{\algorithmicindent} 
 \textbf{Output:} Point estimator $\widehat{f'(\eval)}$ and confidence interval ${\rm CI}[f'(\eval)]$. \\
\end{flushleft}
\begin{algorithmic}[1]
\State Implement the sparse additive model in \eqref{eq: SAM} with $M\geq 1$ and $\lambda>0$;
  \State Construct the initial estimators $\{\widehat{g}(X_i)\}_{1\leq i\leq n}$ as in \eqref{eq: initial estimator}; 
  \State Construct the noise estimator $\widehat{\sigma}^2$ as in \eqref{eq: noise estimator}; 
  \State Compute $\widehat{Y}_i=Y_i-\widehat{g}(X_i)$ for $1\leq i\leq n$;
  \Comment{Initial estimators}
  \vspace{2.5mm}
  \State Implement the Lasso algorithm as in \eqref{eq: SLasso} with $\lambda_1>0$;
  \State Construct $\{\widehat{l}(X_i)\}_{1\leq i\leq n}$ as in \eqref{eq: estimator} and \eqref{eq: estimator sym};
\State Construct the weights $\{\widetilde{W}_i\}_{1\leq i\leq n}$ as in \eqref{eq: swapping est};
  \State Construct the centered decorrelation weights $\{\widehat{W}_i\}_{1\leq i\leq n}$ in \eqref{eq: final D};
  \State Construct $\widehat{f'(\eval)}$ as \eqref{eq: final estimator} with $\{\widehat{Y}_i, \widehat{W}_i\}_{1\leq i\leq n}$ and $h>0$; \Comment{\texttt{DLL} estimator}
  \vspace{2.5mm}
  \State Compute the variance estimate $\widehat{\rm V}$ as in \eqref{eq: ci};
  \State Construct ${\rm CI}[f'(\eval)]$ as in \eqref{eq: ci}. \Comment{Confidence interval}
  % \State Implement the testing procedure $\psi[f'(\eval)]$ as in \eqref{eq: testing}. 
\end{algorithmic}
\label{algo: DDL}
\end{algorithm}

%%%%%%%%%%%%%%%%%%%%%%%%%%%%%%%%%%%
\subsection{Decorrelaction with the additive treatment model}
\label{sec: DLL-S}
%%%%%%%%%%%%%%%%%%%%%%%%%%%%%%%%%%%
This section generalizes the construction of decorrelation weights by considering non-linear models between $D_i$ and $X_i$. Particularly, we consider the sparse additive model for $D_i,$
\begin{equation*}
D_i=\sum_{j=1}^{p} \tau_j({X}_{i,j})+\delta_i, \quad \text{for}\; 1\leq i\leq n,
\end{equation*}
where $\tau_j: \R\rightarrow \R$ for $1\leq j\leq p$ are unknown smooth functions. Instead of applying the Lasso algorithm \eqref{eq: SLasso}, we implement another sparse additive model as in \eqref{eq: SAM},
\begin{equation*}
\{\widehat{\gamma}^{a}_j\}_{1\leq j\leq p}=\argmin_{\gamma_j \in \R^{M},\; 1\leq j\leq p}\frac{1}{2|\mathcal{I}_a|}\sum_{i\in \mathcal{I}_a}(D_i-\sum_{j=1}^{p}\Psi_{i,j}^{\intercal}\gamma_j)^2+\lambda_1\sum_{j=1}^{p}\sqrt{\gamma_j^{\intercal} \left(\frac{1}{|\mathcal{I}_a|}\sum_{i\in \mathcal{I}_a} \Psi_{i,j} \Psi_{i,j}^{\intercal}\right) \gamma_j}
\label{eq: SAM first stage}
\end{equation*}
where $\lambda_1>0$ is the tuning parameter to be chosen. We estimate $\{\mu_i,\delta_i\}_{i\in \mathcal{I}_b}$ by \begin{equation*}
\widehat{\mu}_i=\eval-\sum_{j=1}^{p}\Psi_{i,j}^{\intercal}\widehat{\gamma}^{a}_j \quad \text{and}\quad \widehat{\delta}_i=D_i-\sum_{j=1}^{p}\Psi_{i,j}^{\intercal}\widehat{\gamma}^{a}_j \quad \text{for}\quad i\in \mathcal{I}_b.
\end{equation*}
By switching the data in $\mathcal{I}_a$ and $\mathcal{I}_b$, we construct $\{\widehat{\mu}_i, \widehat{\delta}_i\}_{i\in \mathcal{I}_a}$. With the estimates $\{\widehat{\mu}_i, \widehat{\delta}_i\}_{1\leq i\leq n}$, we construct the decorrelation weights in the same way as in \eqref{eq: estimator} and \eqref{eq: estimator sym}. In Section \ref{sec: DLL-S sim}, we compare the finite-sample performance of Algorithm \ref{algo: DDL} and the generalized decorrelation method with the sparse additive model.

%%%%%%%%%%%%%%%%%%%%%%%%%%%%%%%%%%%%%%%%%%%%%%%%%%%%%%%%%%%%%%
%\subsection{Comparison to the Debiased Methods}
%\label{sec: comparison}
%%%%%%%%%%%%%%%%%%%%%%%%%%%%%%%%%%%%%%%%%%%%%%%%%%%%%%%%%%%%%
%\Zijian{Cun-Hui, please check.}
%In this section, we provide a detailed discussion on the difference between our proposed \texttt{DLL} estimator and the debiased estimators in high-dimensional generalized linear models. 

%%%%%%%%%%%%%%%%%%%%%%%%%%%%%%%%%%%%%%%%%%%%%%%%%%%%%%%%%%%%%%%%%%%%%%%%%%%%%%
\section{Theoretical justification}
\label{sec: theory}

\subsection{Technical conditions}
Before presenting the main theorems, we present the technical conditions imposed on the outcome model \eqref{eq: additive model} and the treatment model \eqref{eq: decouple}. Let $\pi(\eval)$ denote the probability density function of $D_i$ evaluated at $\eval\in \R.$ The first condition is on the function of interest $f(\cdot)$, the regression error $\epsilon_i$, and the bandwidth $h>0.$ %{sec: simulation}
\begin{enumerate}
\item[(A1)] $f(\cdot)$ is twicely differentiable at a neighborhood of $\eval$ and  $f''(\cdot)$ is continuous at $\eval$. The error ${\epsilon}_i$ in  \eqref{eq: additive model} satisfies $\E(\epsilon_i\mid D_i, X_i)=0,$ $\E(\epsilon_i^2\mid D_i, X_i)=\sigma^2,$ and $\E(\epsilon_i^{2+c}\mid D_i, X_i)\leq C$ for some positive constants $c>0$ and $C>0$. The bandwidth $h$ used in \eqref{eq: kernel} satisfies 
$
nh \pi(\eval)\gg \log n$ and $n h^{5}\pi(\eval)\leq C$ for some positive constant $C>0.$ %and $h\left(\eval+h+\|\gamma\|_2 \sqrt{\log n}\right)\rightarrow 0.$
\end{enumerate}
%Condition $\rm (A1)$ requires the additive model structure along with a mild moment condition on the error term. In addition,  

Condition (A1) is standard for the analysis of the local polynomial estimator in the univariate case \cite[e.g.]{fan1996local,fan1993local}. The smoothness condition on $f$ ensures a sufficiently small approximation error of $f$ by a linear function in a neighborhood near $\eval$. The conditional moment conditions on $\epsilon_i$ are required to establish the asymptotic normality of our proposed \texttt{DLL} estimator.  
Since the expected number of observations $D_i\in [\eval-h,\eval+h]$ is about $2nh\pi(\eval)$, Condition $\rm (A1)$ requires that there are (asymptotically) infinitely many observations in the local neighborhood of $\eval$ with bandwidth $h$. For a twicely differentiable function $f$, the optimal choice of bandwidth for estimating $f'(\eval)$ is $h\asymp n^{-1/5},$ which satisfies  $n h^{5}\pi(\eval)\leq C.$ 

The second model assumption is imposed on the treatment model \eqref{eq: decouple}. Recall that $\phi(\cdot)$ denotes the probability density function of the regression error $\delta_i=D_i-X_i^{\intercal}\gamma$ in \eqref{eq: decouple} and $\mu_i=\eval-X_{i}^{\intercal}\gamma$ for $1\leq i\leq n.$ We use $C_1(n)>0$ and $C_2(n)>0$ to denote some high probability upper bounds, defined as: with probability larger than $1-\min\{n,p\}^{-c}$ for some positive constant $c>0,$
\begin{equation}
\max_{1\leq i\leq n}\max_{\left|\delta-\mu_i\right|\leq r} \frac{|\phi'(\delta)|}{\phi(\mu_i)}\leq C_1(n), 
 \quad \max_{1\leq i\leq n}\max_{\left|\delta-\mu_i\right|\leq r} \frac{|\phi''(\delta)|}{\phi(\mu_i)}\leq  C_2(n), 
 \label{eq: density ratio}
\end{equation}
where $r=C^*\sqrt{{\|\gamma\|_0 \log p \log n}/{n}}+h$ for some positive constant $C^*>0$. The value $C_1(n)$ (or $C_2(n)$) defined in \eqref{eq: density ratio} captures the ratio of  $\phi'$ (or $\phi''$) over $\phi$ near $\mu_i.$ $C_1(n)$ and $C_2(n)$ are allowed to grow with $n$ and $p$, but in general they grow to infinity at a slow rate; see Remark \ref{rem: growth rate} for details. 
We now state the condition for the model \eqref{eq: decouple} and $C_1(n)$ and $C_2(n)$ defined in \eqref{eq: density ratio}.

\begin{enumerate}
\item [(A2)] The model \eqref{eq: decouple} holds with $k\coloneqq \|\gamma\|_0\ll \min\{1,\pi(\eval)\}\cdot \tfrac{n}{\log p\log n}$, $X_{i}$ and $\delta_i$ being Sub-gaussian, and the error $\delta_i$ being independent of $X_{i}$. The variance of $\delta_i$ is a positive constant and $\Sigma=\E X_{i} X_i^{\intercal}$ satisfies $c_0\leq \lambda_{\min}(\Sigma)\leq \lambda_{\max}(\Sigma)\leq C_0$ for some positive constant $c_0>0$ and $C_0>0$. The density function $\phi$ of $\delta_i$ is upper bounded and %, with probability larger than $1-\min\{n,p\}^{-c}$ for some positive constant $c>0,$ 
\begin{equation}
h^2 C_2(n)+(\sqrt{{\|\gamma\|_0 \log p \log n}/{n}}+h)C_1(n)\rightarrow 0,
\label{eq: density condition}
\end{equation}
where $C_1(n)$ and $C_2(n)$ are defined in \eqref{eq: density ratio}.
\end{enumerate}
For a constant $\pi(\eval)>0,$ $k\ll {n}/{[\log p\log n]}$ is almost the weakest sparsity condition to identify $\gamma$ for the high-dimensional linear model. The conditions on ${\rm Var}(\delta_i)$ and the covariance matrix $\Sigma$ are standard for the high-dimensional analysis. The condition \eqref{eq: density condition} is mild with $C_1(n)$ and $C_2(n)$ growing at the polynomial order of $\log n$ and $h\asymp n^{-1/5};$ see Remark \ref{rem: growth rate}. The independence assumption between $\delta_i$ and $X_i$ is stringent but we believe it is mainly imposed for the technical analysis. In numerical studies, we test the performance of our proposed method when the independence assumption between $\delta_i$ and $X_i$ is violated; {see Settings 3 and 4 and Table \ref{table:1} in Section \ref{sec: sim result 1}  for details.} 
\begin{Remark}[Growth rates of $C_1(n)$ and $C_2(n)$.]\rm
We discuss the order of magnitudes for $C_1(n)$ and $C_2(n)$ over two examples. Firstly, consider the setting that there exist positive constant $C_0>0$ such that $\max_{1\leq i\leq n}\left|\mu_i\right|\leq C_0.$ If $\min_{|\delta|\leq C_0}\phi(\delta)\geq c$ for some positive constant $c>0$, and $\phi(\delta)$ is twicely differentiable for $\delta \in (C_0-r, C_0+r)$, then $C_1(n)$ and $C_2(n)$ are of constant orders. Secondly, we consider that $X_{i}$ is sub-gaussian and $\mu_i$ may be unbounded in this case. If $\phi$ is the Gaussian density, and $(\sqrt{{\|\gamma\|_0 \log p \log n}/{n}}+h)\sqrt{\log n}\lesssim 1$, then with probability larger than $1-n^{-c}$ for some positive constant $c>0$,
$$C_1(n)\lesssim \sqrt{\log n} \quad \text{and}\quad C_2(n)\lesssim \log n.$$  
\label{rem: growth rate}
\end{Remark}

Finally, we require that the initial estimator $\widehat{g}$ estimates $g=\sum_{j=1}^{p}g_j$ to certain accuracy.  We use ${\rm Err}(\widehat{g})$ to denote the estimation accuracy of $\widehat{g}$, which is defined as follows: with probability larger than $1-\min\{n,p\}^{-c}$ for some positive constant $c>0$, the initial estimator $\widehat{g}$ defined in \eqref{eq: swapping est nuisance} satisfies  
\begin{equation}
\sqrt{\E_{X_{*}}(\widehat{g}^a(X_{*})-g(X_{*}))^2+\E_{X_{*}}(\widehat{g}^b(X_{*})-g(X_{*}))^2}\lesssim {\rm Err}(\widehat{g}),
\label{eq: g error}
\end{equation}
where the expectation $\E_{X_{*}}$ is taken with respect to the independent copy $X_{*}$ of $\{X_i\}_{1\leq i\leq n}.$ 
The last assumption is on the rate of convergence of ${\rm Err}(\widehat{g})$.
\begin{enumerate}
\item[(A3)] The estimation accuracy ${\rm Err}(\widehat{g})$ of initial estimator $\widehat{g}$ is required to satisfy 
%The initial estimation error ${\rm Err}(\widehat{g})$ defined in \eqref{eq: g error} satisfies
\begin{equation}
{\max\left\{{(C_1^2(n)+C_2(n))}\sqrt{{h^3 k \log p\log n}},1\right\}} \cdot \frac{{\rm Err}(\widehat{g})}{{\sqrt{\pi(\eval)}}} \rightarrow 0, \label{eq: general inference condition}
\end{equation}
where ${\rm Err}(\widehat{g})$ is defined in \eqref{eq: g error}.
\end{enumerate}
We discuss Condition (A3) by focusing on a commonly used regime with $C_1^2(n)+C_2(n)\lesssim \log n$, $\pi(\eval)$ being of a constant order, and $h\asymp n^{-1/5}.$ If $k\log p(\log n)^3/n^{3/5}\leq c$ for some positive constant $c>0,$ then any consistent estimator $\widehat{g}$ with ${\rm Err}(\widehat{g})\rightarrow 0$ will automatically satisfy the condition  \eqref{eq: general inference condition}. Most estimators proposed in the high-dimensional sparse additive model can be shown to satisfy the assumption (A3). More discussion about (A3) can be found in Section \ref{sec: A3 discussion} in the supplement.

%\Zijian{Stops here on Feb 11, 2022.}

%%%%%%%%%%%%%%%%%%%%%%%%%%%%%%%%%%%%%%%%%%%%%%%%%%%%%%%%%%%%%%%%%%%%%%%%%%%%%%%%
\subsection{Asymptotic normality and inference properties}
\label{sec: property}
%%%%%%%%%%%%%%%%%%%%%%%%%%%%%%%%%%%%%%%%%%%%%%%%%%%%%%%%%%%%%%%%%%%%%%%%%%%%%%%%
We establish the asymptotic limiting distribution for our proposed \texttt{DLL} estimator.
\begin{Theorem}
Suppose that Conditions {\rm (A1)}, {\rm (A2)} and {\rm (A3)} hold. Then our proposed estimator $\widehat{f'(\eval)}$ in \eqref{eq: final estimator} satisfies,
\begin{equation}
\frac{1}{\sqrt{\rm V}}\left(\widehat{f'(\eval)}-f'(\eval)\right)\cid N\left(0, 1 \right)\;\;\text{with}\;\;
{\rm V}\coloneqq\frac{\sigma^2}{n^2\widehat{S}^2_{n}}{\sum_{i=1}^{n}\widehat{W}_{i}^2K_{h}^2(D_i)},
\label{eq: limiting distribution}
\end{equation}
where $\widehat{S}_{n}$ is defined in \eqref{eq: final estimator} and ${\rm V}\cip \tfrac{3\sigma^2}{n h^3\cdot\pi(\eval)}$.
\label{thm: limiting swap}
\end{Theorem}

With $h\asymp n^{-1/5}$ and $\sigma$ and $\pi(\eval)$ being of constant orders, our proposed \texttt{DLL} estimator achieves the optimal rate of convergence $n^{-1/5}$ \cite{tsybakov2009introduction}. Furthermore, the \texttt{DLL} estimator is asymptotically normal and the asymptotic variance depends on the value $\eval$ through the density level $\pi(\eval)$. In finite samples, we compare the variance level of our \texttt{DLL} estimator to that of the oracle estimator by assuming the knowledge of $g$; See Table \ref{table:1} for details.

As a Corollary of Theorem \ref{thm: limiting swap}, we establish the properties of our proposed confidence interval ${\rm CI}[f'(\eval)]$ defined in \eqref{eq: ci}.
\begin{Corollary} Suppose that Conditions {\rm (A1)}, {\rm (A2)} and {\rm (A3)} hold and $\widehat{\sigma}^2\cip \sigma^2$. For any $\alpha\in(0,1/2),$ our proposed confidence interval ${\rm CI}[f'(\eval)]$ defined in \eqref{eq: ci} satisfies,
\begin{equation*}
\liminf_{n\rightarrow \infty} \mathbb{P}(f'(\eval)\in {\rm CI}[f'(\eval)])= 1-\alpha, 
\end{equation*}
and
\begin{equation*}
\limsup_{n\rightarrow \infty} \mathbb{P}\left(\mathbf{L}\left({\rm CI}[f'(\eval)]\right)\geq (2+\delta_0)z_{\alpha/2}\sqrt{\frac{3}{2n h^3\cdot\pi(\eval)}}\sigma\right)=0,
\end{equation*}
where $\mathbf{L}\left({\rm CI}[f'(\eval)]\right)$ denotes the length of the interval, $\delta_0>0$ is any positive constant, and $z_{\alpha/2}$ denotes the upper $\alpha/2$ quantile of the standard normal distribution.
\label{cor: ci}
\end{Corollary}
Beyond Conditions (A1)-(A3), the above corollary requires a consistent estimator of $\sigma^2$ such that our proposed variance estimator $\widehat{\rm V}$ consistently estimates ${\rm V}.$ In Proposition \ref{prop: consistent sigma} in Section \ref{sec: consistent sigmasq}
 in the supplement, we show that our proposed estimator $\widehat{\sigma}^2$ in \eqref{eq: noise estimator} satisfies $\widehat{\sigma}^2\cip \sigma^2$ if both $\widehat{f}$ and $\widehat{g}$ are consistent.   Similarly, we can establish the validity of the proposed testing procedure $\psi[f'(\eval)]$ in \eqref{eq: testing}.
\begin{Corollary} Suppose that the conditions of Corollary \ref{cor: ci} hold. If $f'(\eval)=0$, then the proposed testing procedure $\psi[f'(\eval)]$ defined in \eqref{eq: testing} controls the type I error, that is, 
\begin{equation*}
\limsup_{n \rightarrow \infty} \mathbb{P}(\psi[f'(\eval)]=1)= \alpha.
\end{equation*}
\end{Corollary}

%%%%%%%%%%%%%%%%%%%%%%%%%%%%%%%%%%%%%%%%%%%%%%%%%%%%%%%%%%%%%%%%%%%%%%%%%%%%%%
\subsection{Theoretical reasoning of decorrelation}
\label{sec: decomposition}
%%%%%%%%%%%%%%%%%%%%%%%%%%%%%%%%%%%%%%%%%%%%%%%%%%%%%%%%%%%%%%%%%%%%%%%%%%%%%%
%Since $\sum_{i=1}^{n}\widehat{W}_{i}K_{h}(D_i)=0$,  
In the following, we explain why our constructed decorrelated weight is effective. With  $\dif=g(X_{i})-\widehat{g}(X_{i}),$ the estimation error of the \texttt{DLL} estimator is decomposed as 
{\small
\begin{equation*}
\widehat{f'(\eval)}-f'(\eval)=\underbrace{\frac{1}{n\widehat{S}_{n}}\sum_{i=1}^{n}\widehat{W}_{i}\epsilon_iK_{h}(D_i)}_{\rm Stochastic\; Error}+\underbrace{\frac{1}{n\widehat{S}_{n}}\sum_{i=1}^{n}\widehat{W}_{i}r(D_i)K_{h}(D_i)}_{\rm Approximation \; Error}+\underbrace{\frac{1}{n\widehat{S}_{n}}\sum_{i=1}^{n}\widehat{W}_{i}\dif K_{h}(D_i)}_{\rm High-dimensional \; Error}.
\end{equation*}}

\noindent The stochastic error represents a random component with mean zero and, after rescaling, following an asymptotic normal limiting distribution. The approximation error  is the error of approximating the non-linear function $f$ by a linear function at a local neighborhood of $\eval$. The high-dimensional error is due to the estimation of the unknown function $g$ in high dimensions. Both the stochastic and approximation errors appear in the classical non-parametric regression, while the high-dimensional error is the new addition here. 

The following theorem demonstrates that our proposed decorrelation method is effective in reducing the high-dimensional error.

\begin{Theorem} Suppose that Condition $\rm (A1)$ and $\rm (A2)$ hold. For $\dif={g}(X_{i})-\widehat{g}(X_{i})$ where $\widehat{g}$ is defined in \eqref{eq: swapping est nuisance}, then with probability larger than $1-\frac{1}{t}-\min\{n,p\}^{-c}$ for some $t>1$, 
\begin{equation}
\left|{\frac{1}{n\widehat{S}_{n}}\sum_{i=1}^{n}\widehat{W}_{i}\dif K_{h}(D_i)}\right|\leq t^2 \left[1+\sqrt{{h^3k\log p\log n}} \left(C_1^2(n)+C_2(n)\right)\right]\frac{{\rm Err}(\widehat{g})}{\sqrt{nh^3 \pi^2(\eval)}},
\label{eq: high dim error}
\end{equation}
where ${\rm Err}(\widehat{g})$ is defined in \eqref{eq: g error}.
\label{thm: DLL error reduction}
\end{Theorem}

Our constructed decorrelation weights are instrumental in reducing the error due to estimating $g$. This happens mainly due to the fact that the expectation of $W_i\Delta(X_i) K_{h}(D_i)$ is zero. Condition (A3) and the upper bound \eqref{eq: high dim error} imply $$\frac{1}{\sqrt{\rm V}}\left|{\frac{1}{n\widehat{S}_{n}}\sum_{i=1}^{n}\widehat{W}_{i}\dif K_{h}(D_i)}\right|\cip 0.$$ The data swapping step creates the independence between the error function $\Delta$ and the data $\{X_i, D_i, W_i\},$ which is required for the proof of \eqref{eq: high dim error}. We believe that a more refined analysis might remove the data swapping step.

%%%%%%%%%%%%%%%%%%%%%%%%%%%%%%%%%%%%%%%%%%%%%%%%%%%%%%%%%%%%%%%%
\section{Simulation} 
\label{sec: simulation}
%%%%%%%%%%%%%%%%%%%%%%%%%%%%%%%%%%%%%%%%%%%%%%%%%%%%%%%%%%%%%%%%%%%%%%%%%%%%%%
 We provide more details about the tuning parameter selection for Algorithm \ref{algo: DDL}. For the high-dimensional sparse additive model, we compute the initial estimators $\widehat{f}$ and $\widehat{g}$ by applying the R package \texttt{SAM} \cite{sampackage} and choose the tuning parameter $\lambda$ and the number M of basis functions in \eqref{eq: SAM} by cross validation. We construct the Lasso estimator of $\gamma$ in \eqref{eq: SLasso} by applying the R package \texttt{glmnet} \cite{glmnet} and choose the tuning parameter $\lambda_1$ by cross validation. For local linear methods, choosing a good bandwidth is essential for the finite-sample performance. There are many methods for bandwidth selection. After exploration in the simulation study, we observe that the ``Rule of Thumb'' method proposed in \cite{fan1996local} leads to the most stable performance. This bandwidth selection method is implemented in the R package \texttt{locpol} \cite{locpol} with the \texttt{thumbBw()} function. The ``Rule of Thumb'' is used as our default bandwidth selection method. We demonstrate the performance of the \texttt{DLL} estimator with other bandwidth selection methods in Section \ref{sec: other bw} in the supplement. The codes for replicating our proposed method can be found at \url{https://github.com/zijguo/HighDim-Additive-Inference}. %implemented in R package \texttt{locpol} and \texttt{np} \cite{nppackage}; \Zijian{Wei, here is the standard way to refer to the supplement and please change that at all other places.} 
%searches the parameter space from 0.001 to 1 for 200 candidate $\lambda_n$'s and from 3 to 6 for the number of basis functions.\Zijian{I will revise the following.}As stated in (\ref{eq: initial additive}), the two terms $\|f\|_n$ and $\|f\|_F$ are imposed to serve the two purposes respectively and are controlled by the two parameters $\lambda_n$ and $\rho_n$. The term $\|f\|_n = (\sum_{i=1}^{n} f^2(D_{i}))^{1/2}$ denoted the empirical $L_2$ norm and the term $\|f\|_F = (\int [f''(t)]^2 dt)^{1/2}$ is a measure of the function's smoothness.  Since we utilize B-spline basis for smoothing, the penalty of smoothness can be represented by the number of basis functions. Our modified function  %Every combination of the two parameters is performed and we select the one with minimal cross validated MSE. 

Since we believe that the data swapping is introduced for technical analysis, we mainly report the simulation results for the \texttt{DLL} estimator without the data swapping, which is described in Section \ref{sec: no data swapping} in the supplement. We compare our constructed confidence intervals with and without data swapping in Section \ref{sec: swap and trans} in the supplement. Both confidence intervals attain the desired coverage level. When the sample size is relatively large, they have similar performance; for relatively small sample size, the confidence interval without data swapping can be shorter than that with data swapping.

We demonstrate the finite-sample performance of our proposed \texttt{DLL} estimator across various settings and compare it with three other estimators described as follows, 
\begin{itemize}
\item The plug-in estimator (\texttt{Plug}) is implemented in \eqref{eq: plug-in}, where the initial estimator $\widehat{g}$ and the bandwidth $h$ are constructed in the same way as our proposed \texttt{DLL} estimator. For implementation of the local linear estimator and the related confidence interval, we follow the output of the package \texttt{nprobust} \cite{nprobust}.
\item The oracle estimator (\texttt{Orac}) denotes the local linear estimator applied to the data $\{D_i, Y_i^{\rm ora}\}_{1\leq i\leq n}$ with $Y_i^{\rm ora}=Y_i-g(X_i)=f(D_i)+\epsilon_i$. The oracle estimator is used as the benchmark to compare with. For implementation of the local linear estimator and the related confidence interval, we follow the output of the package \texttt{nprobust} \cite{nprobust}.
\item The ReSmoothing (\texttt{RS}) estimator is a two-step estimator proposed in  \cite{gregory2016optimal}. In the first step, we implement the code available at \url{https://github.com/gregorkb/spaddinf} and obtain a pre-smoothing estimator of $f$, denoted as $\widehat{f}^{\rm pre}$; in the second step, we apply the local polynomial estimator to the data $\{D_i, \widehat{f}^{\rm pre}(D_i)\}_{1\leq i\leq n},$ where $\widehat{f}^{\rm pre}(D_i)$ is used as the outcome. We fit the local linear estimator by the package \texttt{nprobust} \cite{nprobust}. 
\end{itemize}

We generate the outcome following the  model \eqref{eq: additive model} and consider both exactly sparse and approximately sparse settings.\\
%We introduce 15 functions that are not zeros to represent the sparsity in high dimension.
\noindent{\bf Exactly sparse.} We set the first six functions as follows and $g_j=0$ for $6 \leq j \leq p$,
\begin{equation}
    \begin{aligned}
         f(d) &= 1.5\sin(d) & g_1(x) &= 2\exp(-x/2) & g_2(x) &= (x-1)^2 - 25/12 \\
         g_3(x) &= x - 1/3 & g_4(x) &= 0.75x & g_5(x) &= 0.5x.
         \end{aligned}
         \label{eq: basic setting}
\end{equation}

More complicated relationships often exist in real life and the additive model might not be exactly sparse. We further introduce an approximately sparse setting.\\
\noindent{\bf Approximately sparse.} We set $f$ and $\{g_j\}_{1\leq j\leq 5}$ as in \eqref{eq: basic setting}, generate $\{g_j\}_{6\leq j\leq 14}$ as
\begin{equation*}
    \begin{aligned}
         g_6(x) &= 0.5x & g_7(x) &= 0.4x & g_8(x) &= 0.3x &
         g_{9}(x) &= 0.2x & g_{10}(x) &= 0.1\sin(2 \pi x)
    \end{aligned}
\end{equation*}
\begin{equation*}
    \begin{aligned}
        g_{11}(x) &= 0.2\cos(2 \pi x) & g_{12}(x) &= 0.3\sin^2(2 \pi x) & g_{13}(x) &= 0.4\cos^3(2 \pi x) & g_{14}(x) &= 0.5\sin^3(2 \pi x),
    \end{aligned}
\end{equation*}
and generate $\{g_j\}_{15\leq j \leq p}$ as linear functions with $g_j(x)=x/(j-1)$. 

In addition, we explore the finite-sample performance for different non-linear functions by switching the role of $f$ and $g_1$ function; see the results in Section \ref{sec: 1to4 and nonlinear} in the supplement.

%%%%%%%%%%%%%%%%%%%%%%%%%%%%%%%%%%%%%%%%%%%%%%%%%%%%%%%%%%%%%%%%%%%%%%%%
\subsection{Comparison with plug-in and oracle estimators}
\label{sec: sim result 1}
%%%%%%%%%%%%%%%%%%%%%%%%%%%%%%%%%%%%%%%%%%%%%%%%%%%%%%%%%%%%%%%%%%%%%%%%

In the following, we compare our proposed \texttt{DLL} estimator with the plug-in(\texttt{Plug}) and oracle(\texttt{Orac}) estimators. We consider four different settings for generating $D_i$ and $X_i,$ where the independence assumption between $X_i$ and $\delta_i$ in (A2) is violated in Settings 3 and 4. 

\noindent \textbf{Setting 1.} We generate $(D_i,X^{\intercal}_i)^{\intercal}$ following the multivariate Normal distribution $N(\mu,\Sigma)$, where $\mu_j=-0.25$ for $1\leq j\leq p+1$ and $\Sigma\in \R^{(p+1)\times(p+1)}$ is a toeplitz covariance matrix with $\Sigma_{jj}=1$ for $1\leq j\leq p+1$ and for $1\leq j\neq l\leq p+1,$
$$\Sigma_{j,l}=0.7\cdot {\1}(|j-l|=1)+0.5\cdot {\1}(|j-l|=2)+0.3\cdot {\1}(|j-l|=3)+\frac{p-|j-l|}{10(p-4)}{\1}(|j-l|\geq 4).$$
For $|j-l| \geq 4$, the correlation gradually decays from 0.1 to 0.

\noindent \textbf{Setting 2.} $(D_i,X^{\intercal}_i)^{\intercal}$ is generated in the same way as in Setting 1. With $G$ denoting the CDF of $N(-0.25,1),$ we generate the outcome model as 
\begin{equation*}
Y_i=f(G(D_i))+\sum_{j=1}^{p}g_j(G({X}_{i,j}))+\epsilon_i,\quad \text{for}\quad 1\leq i\leq n.
\end{equation*} 
The main difference from Setting 1 is to apply a quantile transformation to $D_i$ and $\{{X}_{i,j}\}_{1\leq j\leq p}$ before applying the additive model transformation. The goal is to make inference for $(f^*)'(\eval)$ with $f^{*}=f\circ G.$ 
We generate the additive model following \eqref{eq: basic setting} but set $f(d) = -1.5\sin(\pi d)$, $g_1(x) = 2\exp(-x)$, $g_4(x) = x^3 - 1/2$, and $g_5(x) = x/(1+x).$

\noindent \textbf{Setting 3.} We generate $({D}^0_i,({X}^0_i)^{\intercal})^{\intercal}$ following $N(\mu,\Sigma)$ with the same $\mu$ and $\Sigma$ as in Setting 1. We define $D_i = 5(G({D}^0_i)-0.5)$ and ${X}_{i,j} = 5(G({X}^0_{i,j})-0.5)$ for $1 \leq j \leq p$, with $G$ denoting the CDF of $N(-0.25,1)$. The marginal distributions of $D_i$ and ${X}_{i,j}$ are  $\text{Uniform}(-2.5,2.5)$ and  $D_i$ is correlated with $\{X_{i,j}\}_{1\leq j\leq p}$.

\noindent \textbf{Setting 4.} We generate $(D_i,X_i^{\intercal})^{\intercal}$ following a centered multivariate t distribution with the same covariance matrix $\Sigma$ as in Setting 1. The degree of freedom is varied across $\{10, 15\}.$

We fix the dimension $p=1500$ and vary the sample size $n$ across \{500, 1000, 1500, 2000\}. The evaluation points $\eval$ are \{-1.25, -0.5, 0.1, 0.25, 1\}. For Setting 1, we generate the outcome using both exactly and approximately sparse models; for Settings 2 to 4, we only consider the exactly sparse outcome model. We generate the simulation data 500 times and then use the following metrics to compare these methods: 1. Bias, the absolute difference between the average of the 500 point estimates and the true value; 2. Root Mean Square Error (RMSE); 3. Standard Error (SE), the empirical standard deviation of the 500 point estimates; 4. Coverage, the empirical coverage out of 500 simulations; 5. Length, the average length of the constructed confidence interval (CI). In Table \ref{table:1}, we compare our proposed \texttt{DLL} with \texttt{Plug}, and \texttt{Orac} across four simulation settings and we take an average of the metrics  across different sample sizes, evaluation functions, and evaluation points.

We summarize the results in Table \ref{table:1}. 
 For the \texttt{Plug} estimator, the bias component is a dominating term in RMSE, while our \texttt{DLL} estimator is effective in bias correction. The RMSE of our proposed \texttt{DLL} estimator is similar to that of the oracle estimator, which is uniformly smaller than that of the \texttt{Plug} estimator. The coverage error is computed as the absolute difference between the empirical coverage and 95\%; in most cases, the coverage error results from the undercoverage. The CIs based on the \texttt{Plug} estimator are in general undercoverage while our proposed CIs achieve the desired coverage level. Our proposed CI is of a similar length to the length of the oracle CI.

\begin{table}[H]
\centering
\resizebox{\linewidth}{!}{
\begin{tabular}[t]{|c|ccc|cc|ccc|ccc|cc|}
\hline
 \multicolumn{1}{|c}{} & \multicolumn{3}{|c|}{Bias Percentage} & \multicolumn{2}{c|}{RMSE Ratio} & \multicolumn{3}{c|}{SE} & \multicolumn{3}{c|}{Coverage Error} & \multicolumn{2}{c|}{Length Ratio} \\
\hline
 Setting & \texttt{DLL} & \texttt{Plug} & \texttt{Orac} & \texttt{DLL} & \texttt{Plug} & \texttt{DLL} & \texttt{Plug} & \texttt{Orac} & \texttt{DLL} & \texttt{Plug} & \texttt{Orac} & \texttt{DLL} & \texttt{Plug}\\
 \hline
 1 & 0.128 & 0.407 & 0.090 & 1.152 & 1.248 & 0.293 & 0.286 & 0.269 & 1.05\% & 4.19\% & 0.88\% & 1.157 & 1.103 \\
 2 & 0.133 & 0.563 & 0.045 & 1.057 & 1.267 & 0.350 & 0.344 & 0.337 & 1.06\% & 7.46\% & 0.80\% & 1.045 & 1.010 \\
 3 & 0.065 & 0.280 & 0.037 & 1.049 & 1.080 & 0.502 & 0.497 & 0.479 & 0.77\% & 1.94\% & 0.81\% & 1.050 & 1.031 \\
 4 & 0.151 & 0.520 & 0.049 & 1.034 & 1.213 & 0.320 & 0.316 & 0.316 & 1.15\% & 6.72\% & 0.91\% & 1.037 & 0.986 \\
\hline
\end{tabular}}
\caption{Comparison of \texttt{DLL}, plug-in (\texttt{Plug}), and oracle (\texttt{Orac}) estimators. For each setting, metrics are averaged over the total 40 combinations of $n\in \{500, 1000, 1500, 2000\}$, $f(d) \in \{\sin(d),\exp(d)\}$, and $\eval \in \{-1.25, -0.5, 0.1, 0.25, 1\}.$ The columns indexed with ``Bias Percentage" report the percentage of the bias out of RMSE; the columns indexed with ``RMSE Ratio" report the ratio of RMSE to the oracle estimator's RMSE; the columns indexed with ``SE" report the empirical standard error; the columns indexed with ``Coverage Error" report the absolute difference between the empirical coverage and 95\%; the columns indexed with ``Length Ratio" report the ratio of the CI length to the length of the CI based on the oracle estimator.}
\label{table:1}
\end{table}

In Table \ref{table:2}, we report the detailed simulation results for Settings 1 to 4 with $\eval\in \{0.1, 0.25\}$ and $n\in \{500, 1000, 1500\},$ and the complete simulation results are presented in Section \ref{sec: 1to4 and nonlinear} in the supplement. The results are consistent with the observations reported in Table \ref{table:1}: our proposed CI achieves the desired coverage and has a similar length to the oracle CI. In addition, the coverage improvement of our proposed CI over the \texttt{Plug} estimator can be quite substantial as our \texttt{DLL} estimator effectively corrects the bias. For Settings 3 and 4, our proposed method is still effective even if the independence assumption required in Condition (A2) is violated.

\begin{table}[htb!]
\centering
\resizebox{\linewidth}{!}{
\begin{tabular}[t]{|c|c|c|ccc|ccc|ccc|ccc|ccc|}
\multicolumn{18}{c}{Setting 1: approximately sparse}\\
\hline
\multicolumn{1}{|c}{ } & \multicolumn{1}{c}{ } & \multicolumn{1}{c|}{ } & \multicolumn{3}{c|}{Bias} & \multicolumn{3}{c|}{RMSE} & \multicolumn{3}{c|}{SE} & \multicolumn{3}{c|}{Coverage} & \multicolumn{3}{c|}{CI Length} \\
\hline
$\eval$ & True & $n$ & \texttt{DLL} & \texttt{Plug} & \texttt{Orac} & \texttt{DLL} & \texttt{Plug} & \texttt{Orac} & \texttt{DLL} & \texttt{Plug} & \texttt{Orac} & \texttt{DLL} & \texttt{Plug} & \texttt{Orac} & \texttt{DLL} & \texttt{Plug} & \texttt{Orac}\\
\hline
 &  & 500 & 0.21 & 0.46 & 0.02 & 0.44 & 0.60 & 0.39 & 0.39 & 0.38 & 0.39 & 0.91 & 0.72 & 0.94 & 1.51 & 1.45 & 1.49\\

 &  & 1000 & 0.07 & 0.31 & 0.00 & 0.35 & 0.45 & 0.33 & 0.34 & 0.33 & 0.33 & 0.93 & 0.83 & 0.93 & 1.31 & 1.27 & 1.27\\

\multirow{-3}{*}{\centering 0.10} & \multirow{-3}{*}{\centering 1.49} & 1500 & 0.05 & 0.26 & 0.01 & 0.31 & 0.39 & 0.29 & 0.31 & 0.30 & 0.29 & 0.94 & 0.86 & 0.95 & 1.18 & 1.15 & 1.15\\
\cline{1-18}
 &  & 500 & 0.20 & 0.45 & 0.00 & 0.45 & 0.60 & 0.39 & 0.41 & 0.39 & 0.39 & 0.91 & 0.77 & 0.94 & 1.56 & 1.50 & 1.55\\

 &  & 1000 & 0.07 & 0.31 & 0.01 & 0.36 & 0.46 & 0.35 & 0.35 & 0.34 & 0.35 & 0.94 & 0.83 & 0.94 & 1.35 & 1.32 & 1.32\\

\multirow{-3}{*}{\centering 0.25} & \multirow{-3}{*}{\centering 1.45} & 1500 & 0.07 & 0.27 & 0.03 & 0.32 & 0.41 & 0.30 & 0.31 & 0.31 & 0.30 & 0.96 & 0.84 & 0.96 & 1.22 & 1.19 & 1.18\\
\hline

\multicolumn{18}{c}{Setting 2: exactly sparse}\\
\hline
\multicolumn{1}{|c}{ } & \multicolumn{1}{c}{ } & \multicolumn{1}{c|}{ } & \multicolumn{3}{c|}{Bias} & \multicolumn{3}{c|}{RMSE} & \multicolumn{3}{c|}{SE} & \multicolumn{3}{c|}{Coverage} & \multicolumn{3}{c|}{CI Length} \\
\hline
$\eval$ & True & $n$ & \texttt{DLL} & \texttt{Plug} & \texttt{Orac} & \texttt{DLL} & \texttt{Plug} & \texttt{Orac} & \texttt{DLL} & \texttt{Plug} & \texttt{Orac} & \texttt{DLL} & \texttt{Plug} & \texttt{Orac} & \texttt{DLL} & \texttt{Plug} & \texttt{Orac}\\
\hline
 &  & 500 & 0.19 & 0.29 & 0.03 & 0.42 & 0.47 & 0.39 & 0.38 & 0.37 & 0.39 & 0.92 & 0.86 & 0.96 & 1.51 & 1.44 & 1.53\\

 &  & 1000 & 0.15 & 0.25 & 0.06 & 0.36 & 0.40 & 0.32 & 0.32 & 0.31 & 0.32 & 0.93 & 0.88 & 0.95 & 1.30 & 1.24 & 1.28\\

\multirow{-3}{*}{\centering 0.10} & \multirow{-3}{*}{\centering 0.74}  & 1500 & 0.11 & 0.20 & 0.03 & 0.33 & 0.36 & 0.31 & 0.31 & 0.30 & 0.31 & 0.94 & 0.88 & 0.93 & 1.17 & 1.13 & 1.15\\
\cline{1-18}
 &  & 500 & 0.21 & 0.31 & 0.06 & 0.45 & 0.49 & 0.41 & 0.40 & 0.38 & 0.40 & 0.92 & 0.86 & 0.95 & 1.55 & 1.48 & 1.58\\

 &  & 1000 & 0.18 & 0.28 & 0.09 & 0.37 & 0.42 & 0.33 & 0.32 & 0.31 & 0.32 & 0.92 & 0.87 & 0.95 & 1.34 & 1.28 & 1.33\\

\multirow{-3}{*}{\centering 0.25} & \multirow{-3}{*}{\centering 0.94}  & 1500 & 0.12 & 0.21 & 0.04 & 0.34 & 0.37 & 0.32 & 0.32 & 0.31 & 0.31 & 0.92 & 0.86 & 0.94 & 1.22 & 1.17 & 1.20\\
\hline

\multicolumn{18}{c}{Setting 3: exactly sparse}\\
\hline
\multicolumn{1}{|c}{ } & \multicolumn{1}{c}{ } & \multicolumn{1}{c|}{ } & \multicolumn{3}{c|}{Bias} & \multicolumn{3}{c|}{RMSE} & \multicolumn{3}{c|}{SE} & \multicolumn{3}{c|}{Coverage} & \multicolumn{3}{c|}{CI Length} \\
\hline
$\eval$ & True & $n$ & \texttt{DLL} & \texttt{Plug} & \texttt{Orac} & \texttt{DLL} & \texttt{Plug} & \texttt{Orac} & \texttt{DLL} & \texttt{Plug} & \texttt{Orac} & \texttt{DLL} & \texttt{Plug} & \texttt{Orac} & \texttt{DLL} & \texttt{Plug} & \texttt{Orac}\\
\hline
 &  & 500 & 0.12 & 0.24 & 0.01 & 0.72 & 0.74 & 0.69 & 0.71 & 0.70 & 0.69 & 0.94 & 0.92 & 0.93 & 2.81 & 2.73 & 2.60\\

 &  & 1000 & 0.05 & 0.19 & 0.05 & 0.64 & 0.66 & 0.60 & 0.63 & 0.63 & 0.60 & 0.95 & 0.93 & 0.95 & 2.42 & 2.39 & 2.26\\

\multirow{-3}{*}{\centering 0.10} & \multirow{-3}{*}{\centering 1.49} & 1500 & 0.04 & 0.15 & 0.02 & 0.57 & 0.58 & 0.55 & 0.57 & 0.56 & 0.55 & 0.96 & 0.94 & 0.95 & 2.19 & 2.16 & 2.05\\
\cline{1-18}
 &  & 500 & 0.08 & 0.21 & 0.00 & 0.73 & 0.74 & 0.68 & 0.72 & 0.71 & 0.68 & 0.94 & 0.93 & 0.94 & 2.79 & 2.73 & 2.59\\

 &  & 1000 & 0.04 & 0.17 & 0.03 & 0.62 & 0.64 & 0.58 & 0.62 & 0.62 & 0.58 & 0.95 & 0.94 & 0.95 & 2.41 & 2.38 & 2.25\\

\multirow{-3}{*}{\centering 0.25} & \multirow{-3}{*}{\centering 1.45} & 1500 & 0.03 & 0.15 & 0.02 & 0.56 & 0.57 & 0.52 & 0.56 & 0.55 & 0.52 & 0.95 & 0.93 & 0.94 & 2.18 & 2.14 & 2.04\\
\hline

\multicolumn{18}{c}{Setting 4: exactly sparse with df=10}\\
\hline
\multicolumn{1}{|c}{ } & \multicolumn{1}{c}{ } & \multicolumn{1}{c|}{ } & \multicolumn{3}{c|}{Bias} & \multicolumn{3}{c|}{RMSE} & \multicolumn{3}{c|}{SE} & \multicolumn{3}{c|}{Coverage} & \multicolumn{3}{c|}{CI Length} \\
\hline
$\eval$ & True & $n$ & \texttt{DLL} & \texttt{Plug} & \texttt{Orac} & \texttt{DLL} & \texttt{Plug} & \texttt{Orac} & \texttt{DLL} & \texttt{Plug} & \texttt{Orac} & \texttt{DLL} & \texttt{Plug} & \texttt{Orac} & \texttt{DLL} & \texttt{Plug} & \texttt{Orac}\\
\hline
 &  & 500 & 0.17 & 0.36 & 0.04 & 0.38 & 0.49 & 0.33 & 0.34 & 0.33 & 0.32 & 0.92 & 0.78 & 0.95 & 1.33 & 1.27 & 1.32\\

 &  & 1000 & 0.10 & 0.28 & 0.03 & 0.30 & 0.40 & 0.28 & 0.29 & 0.28 & 0.28 & 0.94 & 0.80 & 0.95 & 1.13 & 1.06 & 1.06\\

\multirow{-3}{*}{\centering 0.10} & \multirow{-3}{*}{\centering 1.49} & 1500 & 0.06 & 0.23 & 0.02 & 0.25 & 0.34 & 0.23 & 0.24 & 0.24 & 0.23 & 0.96 & 0.83 & 0.97 & 1.00 & 0.94 & 0.92\\
\cline{1-18}
 &  & 500 & 0.18 & 0.38 & 0.05 & 0.40 & 0.51 & 0.36 & 0.35 & 0.35 & 0.35 & 0.89 & 0.76 & 0.94 & 1.34 & 1.28 & 1.34\\

 &  & 1000 & 0.10 & 0.29 & 0.01 & 0.30 & 0.41 & 0.28 & 0.29 & 0.28 & 0.28 & 0.93 & 0.78 & 0.93 & 1.14 & 1.07 & 1.08\\

\multirow{-3}{*}{\centering 0.25} & \multirow{-3}{*}{\centering 1.45} & 1500 & 0.06 & 0.24 & 0.01 & 0.25 & 0.35 & 0.24 & 0.24 & 0.24 & 0.24 & 0.96 & 0.83 & 0.96 & 1.01 & 0.95 & 0.94\\
\hline
\end{tabular}
}
\caption{Comparison of \texttt{DLL}, plug-in (\texttt{Plug}), and oracle (\texttt{Orac}) estimators for Settings 1 to 4, across different sample sizes $n$ and evaluation points $\eval$. The column indexed with ``True'' represents the true value of $f{'}(\eval)$. The columns indexed with ``Bias'', ``RMSE'' and ``SE" report the absolute bias, the root mean square error, and the standard error computed by 500 estimates, respectively; the columns indexed with ``Coverage'' report the empirical coverage level and the columns indexed with ``Length'' report the average CI length.}
\label{table:2}
\end{table}

%%%%%%%%%%%%%%%%%%%%%%%%%%%%%%%%%%%%%%%%%%%%%%%%%%%%%%%%%%%%%%%%%%%%%%%%%%%%%%%%%%%
\subsection{Comparison with the ReSmoothing method}
\label{sec: comparison RS}
%%%%%%%%%%%%%%%%%%%%%%%%%%%%%%%%%%%%%%%%%%%%%%%%%%%%%%%%%%%%%%%%%%%%%%%%%%%%%%%%%%%
We compare \texttt{DLL} with the \texttt{RS} estimator  \cite{gregory2016optimal} and generate the data as in Setting 1 with $p=750$ and $ n\in \{500, 750, 1000\}$. We construct two CIs centered at the \texttt{RS} estimator,
\begin{enumerate}
\item[(a)] \texttt{RS} confidence interval: we apply the local linear estimator to $\{D_i, \widehat{f}^{\rm pre}(D_i)\}_{1\leq i\leq n}$ with the presmoothing estimators $\{\widehat{f}^{\rm pre}(D_i)\}_{1\leq i\leq n}$ as the outcome variables  \cite{gregory2016optimal}; we construct the CI by the output of the package \texttt{nprobust} \cite{nprobust}.
\item[(b)] \texttt{OraRS} confidence interval: we estimate the standard error of the \texttt{RS} estimator by the sample standard deviation of 500 \texttt{RS} estimates and then construct the confidence interval by assuming the asymptotic normality of the \texttt{RS} estimator.
\end{enumerate}

The \texttt{RS} confidence interval does not necessarily achieve the correct coverage since $\{\widehat{f}^{\rm pre}(D_i)\}_{1\leq i\leq n}$ are not i.i.d. The \texttt{OraRS} is not a practical inference procedure, but a favorable implementation of the confidence interval based on the \texttt{RS} estimator as its standard error is computed in an oracle way.

\begin{table}[htb!]
\centering
\resizebox{\linewidth}{!}{
\begin{tabular}[t]{|c|c|ccc|ccc|cccc|cccc|}
\multicolumn{16}{c}{Setting 1, exactly sparse: Comparison with ReSmoothing}\\
\hline
\multicolumn{1}{|c}{ } & \multicolumn{1}{c|}{ } & \multicolumn{3}{c|}{Bias} & \multicolumn{3}{c|}{SE} & \multicolumn{4}{c|}{Coverage} & \multicolumn{4}{c|}{Length} \\
\hline
$\eval$ & $n$ & \texttt{DLL} & \texttt{RS} & \texttt{Orac} & \texttt{DLL} &\texttt{RS} & \texttt{Orac} & \texttt{DLL} & \texttt{RS} & \texttt{OraRS} & \texttt{Orac} & \texttt{DLL} & \texttt{RS} & \texttt{OraRS} & \texttt{Orac}\\
\hline

 &  500 & 0.20 & 0.91 & 0.01 & 0.42 & 0.83 & 0.41 & 0.94 & 0.01 & 0.80 & 0.96 & 1.69 & 0.04 & 3.24 & 1.68\\

 &  750 & 0.07 & 0.69 & 0.01 & 0.40 & 0.90 & 0.39 & 0.93 & 0.01 & 0.88 & 0.95 & 1.55 & 0.02 & 3.53 & 1.51\\

\multirow{-3}{*}{\centering -1.0} & 1000 & 0.05 & 0.53 & 0.01 & 0.35 & 0.89 & 0.35 & 0.96 & 0.01 & 0.92 & 0.95 & 1.46 & 0.02 & 3.48 & 1.42\\
\hline

 &  500 & 0.20 & 0.97 & 0.02 & 0.42 & 0.78 & 0.41 & 0.92 & 0.01 & 0.78 & 0.95 & 1.69 & 0.04 & 3.06 & 1.68\\

 &  750 & 0.09 & 0.66 & 0.01 & 0.39 & 0.91 & 0.40 & 0.95 & 0.01 & 0.89 & 0.93 & 1.57 & 0.02 & 3.57 & 1.53\\

\multirow{-3}{*}{\centering 0.5} & 1000 & 0.09 & 0.51 & 0.03 & 0.37 & 0.87 & 0.37 & 0.94 & 0.01 & 0.92 & 0.95 & 1.46 & 0.02 & 3.41 & 1.41\\
\hline

\multicolumn{16}{c}{Setting 1, approximately sparse: Comparison with ReSmoothing}\\
\hline
\multicolumn{1}{|c}{ } & \multicolumn{1}{c|}{ } & \multicolumn{3}{c|}{Bias} & \multicolumn{3}{c|}{SE} & \multicolumn{4}{c|}{Coverage} & \multicolumn{4}{c|}{Length} \\
\hline
$\eval$ & $n$ & \texttt{DLL} & \texttt{RS} & \texttt{Orac} & \texttt{DLL} &\texttt{RS} & \texttt{Orac} & \texttt{DLL} & \texttt{RS} & \texttt{OraRS} & \texttt{Orac} & \texttt{DLL} & \texttt{RS} & \texttt{OraRS} & \texttt{Orac}\\
\hline
 &  500 & 0.26 & 0.80 & 0.02 & 0.39 & 0.46 & 0.46 & 0.86 & 0.00 & 0.65 & 0.91 & 1.46 & 0.02 & 1.81 & 1.67\\

 &  750 & 0.18 & 0.73 & 0.01 & 0.35 & 0.63 & 0.38 & 0.94 & 0.00 & 0.80 & 0.95 & 1.42 & 0.01 & 2.46 & 1.51\\

\multirow{-3}{*}{\centering -1.0} & 1000 & 0.12 & 0.66 & 0.01 & 0.35 & 0.73 & 0.38 & 0.94 & 0.00 & 0.86 & 0.93 & 1.37 & 0.01 & 2.87 & 1.42\\
\hline

 &  500 & 0.29 & 1.08 & 0.02 & 0.40 & 0.44 & 0.46 & 0.86 & 0.00 & 0.26 & 0.93 & 1.47 & 0.02 & 1.72 & 1.68\\

 &  750 & 0.23 & 0.96 & 0.03 & 0.38 & 0.60 & 0.40 & 0.89 & 0.00 & 0.63 & 0.94 & 1.43 & 0.01 & 2.34 & 1.53\\

\multirow{-3}{*}{\centering 0.5} & 1000 & 0.14 & 0.78 & 0.00 & 0.33 & 0.73 & 0.35 & 0.94 & 0.01 & 0.81 & 0.95 & 1.37 & 0.01 & 2.88 & 1.41\\
\hline
\end{tabular}
}
\caption{Comparison of \texttt{DLL} estimator, ReSmoothing estimator (\texttt{RS}), and the oracle estimator (\texttt{Orac}), across different sample sizes $n$ and evaluation points $\eval$. The columns indexed with ``\texttt{OraRS}" stand for the CI centered at the RS estimator with the standard error computed based on 500 point estimates. The columns indexed with ``Bias'', and ``SE" report the absolute bias, and the standard error computed by 500 estimates, respectively; the columns indexed with ``Coverage'' report the empirical coverage level and the columns indexed with ``Length'' report the average CI length.}
\label{table:4}
\end{table}

As reported in Table \ref{table:4}, our proposed \texttt{DLL} has a much smaller bias than the \texttt{RS} estimator. In terms of coverage, the \texttt{OraRS} confidence interval does not achieve the desired coverage level, even if its standard error is computed in an oracle way. In contrast, our proposed CI achieves the desired coverage level in most settings. The undercoverage of the \texttt{OraRS}  method happens mainly because of the large bias of the \texttt{RS} estimator. We shall further point out that our proposed CI has a similar length to the oracle CI (the benchmark), but the \texttt{OraRS} confidence interval is much wider than our proposed CI and the oracle CI. See Section \ref{sec: compare RS} in the supplement for results with exchanging the roles of $f$ and $g_1$.

We further compare \texttt{DLL}, \texttt{RS}, and \texttt{OraRS} in the simulation settings of \cite{gregory2016optimal}. The \texttt{DLL} estimator has a smaller bias than the \texttt{RS} estimator and outperforms the \texttt{OraRS} confidence interval in terms of empirical coverage; see Section \ref{sec: compare RS} in the supplement for the results.
%%%%%%%%%%%%%%%%%%%%%%%%%%%%%%%%%%%%%%%%%%%%%%%%%%%%%%%%%%%%%%%%%%%
\subsection{Non-linear treatment model}    
\label{sec: DLL-S sim}
%%%%%%%%%%%%%%%%%%%%%%%%%%%%%%%%%%%%%%%%%%%%%%%%%%%%%%%%%%%%%%%%%%%
In this section, we explore the performance of the generalized \texttt{DLL} estimator proposed in Section \ref{sec: DLL-S}, which decorrelates with the sparse additive model. The estimator is referred to as \texttt{DLL-S}.
As the main difference, \texttt{DLL} is using the Lasso algorithm to fit the treatment model while \texttt{DLL-S} is using the sparse additive model. 
We generate $\{X_i\}_{1\leq i\leq n}$ following the same distribution as in Setting 1 but generate $\{D_i\}_{1\leq i\leq n}$ as $D_i = -0.5\exp(-X_{i,1}/2) + 0.5\sin(X_{i,2}) + 0.25X_{i,3}^2 - 0.5X_{i,4} - 0.25X_{i,5}^2 + 0.5\cos(X_{i,6}) - 0.25\exp(-X_{i,7}/2) + 0.25X_{i,8} + \delta_i$ with $\delta_i \sim N(0,0.5)$. The outcome model is generated following the exactly sparse model. We also consider the setting with switching the roles of $f$ and $g_1$ here. The results are reported in Table \ref{table:3}.

\begin{table}[H]
\centering
\resizebox{\linewidth}{!}{
\begin{tabular}[t]{|c|c|c|cccc|cccc|cccc|}
\multicolumn{15}{c}{Non-linear Treatment Model, exactly sparse: $f(d)=1.5\sin(d)$}\\
\hline
\multicolumn{1}{|c}{ } & \multicolumn{1}{c}{ } & \multicolumn{1}{c|}{ } & \multicolumn{4}{c|}{Bias} & \multicolumn{4}{c|}{Coverage} & \multicolumn{4}{c|}{Length} \\
\hline
$\eval$ & True & $n$ & \texttt{DLL} & \texttt{DLL-S} & \texttt{Plug} & \texttt{Orac} & \texttt{DLL} & \texttt{DLL-S} & \texttt{Plug} & \texttt{Orac} & \texttt{DLL} & \texttt{DLL-S} & \texttt{Plug} & \texttt{Orac}\\
\hline
 & & 500 & 0.31 & 0.23 & 0.37 & 0.01 & 0.84 & 0.91 & 0.78 & 0.95 & 1.40 & 1.44 & 1.34 & 1.50\\

 & & 1000 & 0.20 & 0.12 & 0.25 & 0.00 & 0.87 & 0.93 & 0.83 & 0.95 & 1.24 & 1.24 & 1.18 & 1.25\\

 & & 1500 & 0.18 & 0.10 & 0.23 & 0.03 & 0.92 & 0.95 & 0.88 & 0.96 & 1.12 & 1.14 & 1.07 & 1.12\\

\multirow{-4}{*}{\centering\arraybackslash 0.10} & \multirow{-4}{*}{\centering\arraybackslash 1.49} & 2000 & 0.13 & 0.08 & 0.17 & 0.00 & 0.91 & 0.94 & 0.87 & 0.94 & 1.05 & 1.05 & 1.01 & 1.05\\
\cline{1-15}
 & & 500 & 0.33 & 0.23 & 0.39 & 0.06 & 0.80 & 0.87 & 0.74 & 0.94 & 1.34 & 1.36 & 1.28 & 1.42\\

 & & 1000 & 0.19 & 0.12 & 0.24 & 0.01 & 0.90 & 0.92 & 0.85 & 0.95 & 1.17 & 1.18 & 1.12 & 1.19\\

 & & 1500 & 0.16 & 0.10 & 0.21 & 0.02 & 0.92 & 0.95 & 0.87 & 0.95 & 1.06 & 1.07 & 1.02 & 1.06\\

\multirow{-4}{*}{\centering\arraybackslash 0.25} & \multirow{-4}{*}{\centering\arraybackslash 1.45} & 2000 & 0.13 & 0.08 & 0.18 & 0.00 & 0.91 & 0.95 & 0.87 & 0.95 & 1.00 & 0.99 & 0.95 & 0.99\\
\hline

\multicolumn{15}{c}{Non-linear Treatment Model, exactly sparse: $f(d) = 2\exp(-d/2)$}\\
\hline
\multicolumn{1}{|c}{ } & \multicolumn{1}{c}{ } & \multicolumn{1}{c|}{ } & \multicolumn{4}{c|}{Bias} & \multicolumn{4}{c|}{Coverage} & \multicolumn{4}{c|}{Length} \\
\hline
$\eval$ & True & $n$ & \texttt{DLL} & \texttt{DLL-S} & \texttt{Plug} & \texttt{Orac} & \texttt{DLL} & \texttt{DLL-S} & \texttt{Plug} & \texttt{Orac} & \texttt{DLL} & \texttt{DLL-S} & \texttt{Plug} & \texttt{Orac}\\
\hline
 & & 500 & 0.10 & 0.08 & 0.14 & 0.01 & 0.94 & 0.92 & 0.92 & 0.94 & 1.00 & 1.02 & 0.96 & 1.00\\

 & & 1000 & 0.08 & 0.06 & 0.12 & 0.00 & 0.92 & 0.95 & 0.90 & 0.94 & 0.85 & 0.87 & 0.83 & 0.84\\

 & & 1500 & 0.05 & 0.05 & 0.10 & 0.00 & 0.95 & 0.95 & 0.93 & 0.94 & 0.78 & 0.79 & 0.76 & 0.77\\

\multirow{-4}{*}{\centering\arraybackslash 0.10} & \multirow{-4}{*}{\centering\arraybackslash -0.95} & 2000 & 0.05 & 0.03 & 0.10 & 0.01 & 0.95 & 0.94 & 0.92 & 0.95 & 0.73 & 0.74 & 0.71 & 0.72\\
\cline{1-15}
 & & 500 & 0.11 & 0.10 & 0.15 & 0.00 & 0.92 & 0.92 & 0.88 & 0.95 & 0.94 & 0.96 & 0.91 & 0.94\\

 & & 1000 & 0.09 & 0.07 & 0.13 & 0.02 & 0.93 & 0.92 & 0.90 & 0.94 & 0.81 & 0.82 & 0.79 & 0.80\\

 & & 1500 & 0.04 & 0.06 & 0.09 & 0.01 & 0.94 & 0.93 & 0.91 & 0.93 & 0.74 & 0.75 & 0.72 & 0.73\\

\multirow{-4}{*}{\centering\arraybackslash 0.25} & \multirow{-4}{*}{\centering\arraybackslash -0.88} & 2000 & 0.04 & 0.05 & 0.09 & 0.01 & 0.95 & 0.93 & 0.93 & 0.95 & 0.69 & 0.70 & 0.68 & 0.68\\
\hline
\end{tabular}}
\caption{Comparison of \texttt{DLL}, \texttt{DLL-S}, plug-in (\texttt{Plug}), and oracle (\texttt{Orac}) for the non-linear treatment model, across different sample sizes $n$ and evaluation points $\eval$. The column indexed with ``True'' represents the true value of $f{'}(\eval)$. The columns indexed with ``Bias'' report the absolute bias; the columns indexed with ``Coverage'' report the empirical coverage level and the columns indexed with ``Length'' report the average CI length.}
\label{table:3}
\end{table}

In Table \ref{table:3}, we observe that \texttt{DLL-S} improves the performance of \texttt{DLL} in terms of bias correction and empirical coverage. However, the regular \texttt{DLL} still corrects the bias of the plug-in estimator and achieves better coverage than the CI by the plug-in estimator. See Section \ref{sec: 1to4 and nonlinear} in the supplement for results with additional evaluation points.

%%%%%%%%%%%%%%%%%%%%%%%%%%%%%%%%%%%%%%%%%%%%%%%%%%%%%%%%%%%%%%%%%%%%%%%%%%%%%%
\section{Real data analysis}
\label{sec: real data}
%%%%%%%%%%%%%%%%%%%%%%%%%%%%%%%%%%%%%%%%%%%%%%%%%%%%%%%%%%%%%%%%%%%%%%%%%%%%%%
The Motif Regression has important applications to biology, which studies the effect of the motif candidates' matching scores on the gene expression level \cite{yuan2007motif, beer2004predicting, conlon2003integrating, das2004interacting}. Motifs are the DNA sequences bound to transcription factors, which control the transcription activities, e.g., gene expressions \cite{yuan2007motif}. The matching score of a motif describes the abundance of occurrence, that is, how well the motif is represented in the upstream regions of the genes. A gene's expression level can be well-predicted by the matching scores of a set of motifs \cite{yuan2007motif, beer2004predicting, conlon2003integrating, das2004interacting}. The data set consists of the expression values of $n=2587$ genes and the scores of $p+1=666$ motifs. For our analysis, the outcome $\{Y_i\}_{1 \leq i \leq 2587}$ denote the gene expression level and $\{(D_i,X^{\intercal}_i)^{\intercal}\}_{1 \leq i \leq 2587}$ are the matching scores of the 666 motifs. 

\begin{figure}[htp!]
\centering
\includegraphics[width=0.8\textwidth]{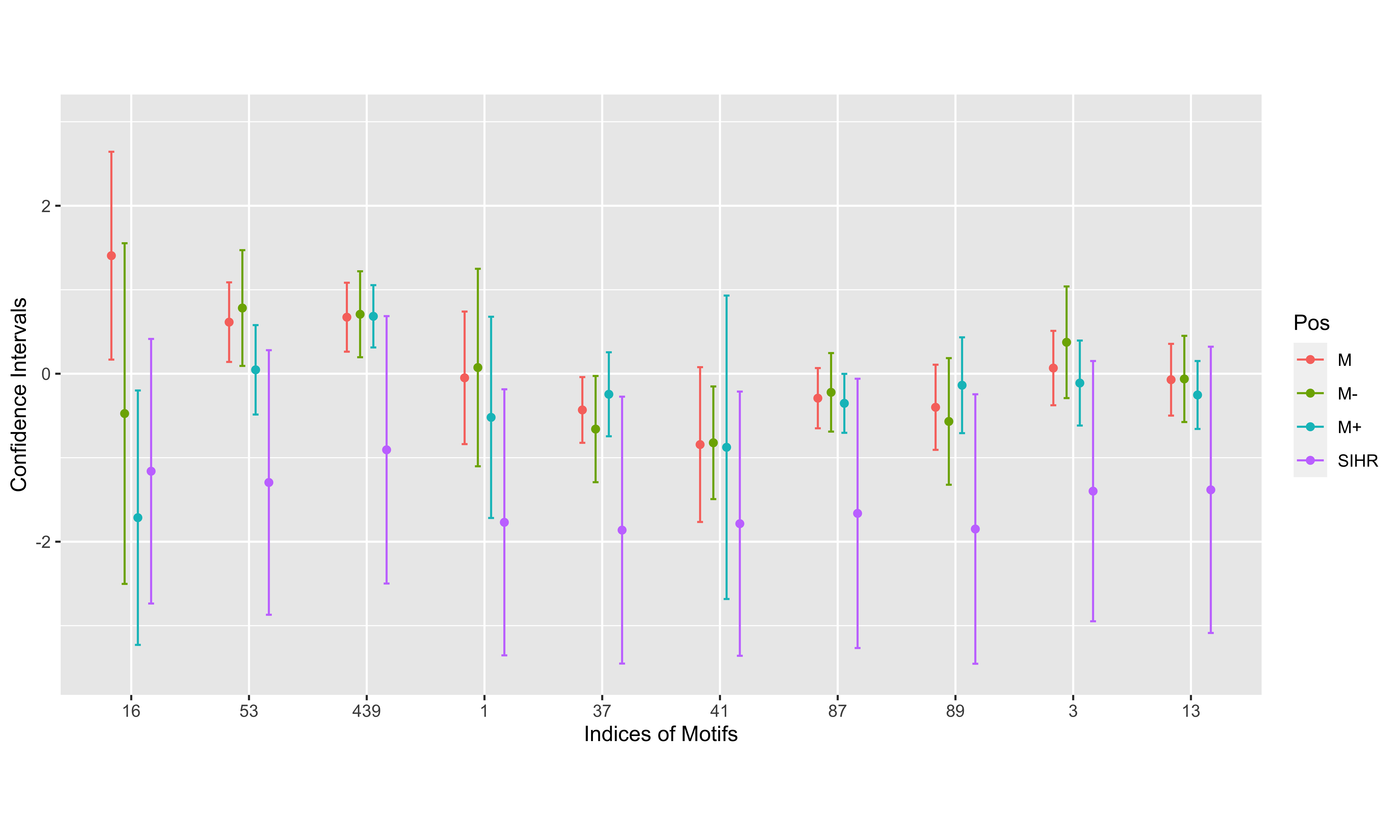}
\caption{Confidence intervals for $f'(\eval)$ by \texttt{DLL} and \texttt{SIHR}. ``M'',``M+'' and ``M-'' represent the \texttt{DLL} estimator for $f'(\eval)$ with $\eval$ set as mean, mean + standard error and mean - standard error, respectively. \texttt{SIHR} represents inference for the constant effect in the high-dimensional linear model.}
\label{fig:1}
\end{figure}

We define an index subset for the motifs $\mathcal{I}=\{1,3,13,16,37,41,53,87,89,439\}\subset\{1,\cdots,666\}.$ To demonstrate our method, we choose one index from $\mathcal{I}$ and set the corresponding motif score as the variable of interest $D$ and the remaining $665$ motif scores as the baseline covariates. We compute its sample mean and standard error for a chosen variable of interest. We choose three different evaluation points $\eval$: mean, mean + standard error, mean - standard error. To demonstrate our method, we compare it with the existing inference method for the high-dimensional linear model, which assumes the linear and constant effect. Specifically, we apply the \texttt{LF()} function in the R package \texttt{SIHR} \cite{rakshit2021sihr} and denote the corresponding estimator as \texttt{SIHR}. We report the comparison in Figure \ref{fig:1}.

Figure \ref{fig:1} demonstrates several interesting observations. First of all, the CI lengths assuming the linear models are in general longer than those of \texttt{DLL}. This happens since the standard deviation of the regression error is about 2.5 by \texttt{SIHR} (assuming the linear model) but 1.45 by \texttt{DLL} (assuming the additive model). This indicates that the relationship between the gene expression levels and the motifs is highly non-linear. For constructing the local linear estimator, the \texttt{DLL} estimator only uses about 40\% of the data  while the \texttt{SIHR} estimator is computed with all data points.

Second, for the motifs with indexes 16, 53, and 439, we observe that those motifs do not have significant effects if we assume the effects to be linear. In contrast, our \texttt{DLL} estimator shows that they might have heterogeneous non-linear effects; for example, for motif 16, the CI at M- location is above zero while the CI at M+ is below zero. Lastly, the motifs with indexes 1, 37, 41, and 87 have significant linear effects, but their non-linear effects vary across different evaluation points.

We design a semi-real simulation study to further compare the finite-sample performance of our proposed {\texttt{DLL}} method and the \texttt{SIHR} method. We keep the data $\{D_i, X_i\}_{1\leq i\leq 2587}$ the same as the real data. After analyzing the original real data, we construct the noise level estimator $\widehat{\sigma}^2$ and $\widehat{f}$ and $\widehat{g}$. We simulate the synthetic response variable $Y^{syn}_i = \widehat{f}(D_i) + \widehat{g}(X_i) + \bar{\epsilon}_i$ for $1\leq i\leq 2587$ with the i.i.d. regression error terms $\{\bar{\epsilon}_i\}_{1\leq i\leq 2587}$ following $N(0,\widehat{\sigma}^2)$. We repeat the simulation 500 times and evaluate \texttt{DLL} on the same three evaluation points as in the real data analysis. 

We compare the results with \texttt{SIHR} and report the comparison in Table \ref{table:5}. The \texttt{SIHR} method, which assumes the linear outcome model, suffers from a large absolute bias and low empirical coverage. Our proposed CIs by the \texttt{DLL} estimators achieve the desired coverage levels in most settings. The lengths of our proposed CIs are, in general, shorter or comparable to the CI output by \texttt{SIHR}. This happens since the regression noise level by \texttt{SIHR} is much larger than the true noise level due to the model misspecification. This matches with the observations for the original real data.

\begin{table}[H]
\centering
\resizebox{\linewidth}{!}{
\begin{tabular}[t]{|c|cccc|cccc|cccc|cccc|}
\hline
\multicolumn{1}{|c|}{ } & \multicolumn{4}{c|}{Bias} & \multicolumn{4}{c|}{SE} & \multicolumn{4}{c|}{Coverage} & \multicolumn{4}{c|}{Length} \\
\hline
Motif & M & M+ & M- & \texttt{SIHR} & M & M+ & M- & \texttt{SIHR} & M & M+ & M- & \texttt{SIHR} & M & M+ & M- & \texttt{SIHR}\\
\hline
1 & 0.08 & 0.15 & 0.14 & 1.37 & 0.23 & 0.45 & 0.32 & 0.41 & 0.93 & 0.87 & 0.94 & 0.45 & 0.94 & 1.73 & 1.19 & 2.66\\
3 & 0.00 & 0.01 & 0.05 & 1.39 & 0.22 & 0.39 & 0.32 & 0.43 & 0.96 & 0.95 & 0.93 & 0.45 & 0.93 & 1.34 & 1.09 & 2.70\\
13 & 0.06 & 0.09 & 0.02 & 1.35 & 0.27 & 0.41 & 0.31 & 0.42 & 0.96 & 0.96 & 0.97 & 0.59 & 1.12 & 1.61 & 1.18 & 2.88\\
16 & 0.24 & 0.26 & 0.00 & 1.30 & 0.36 & 0.71 & 0.45 & 0.40 & 0.87 & 0.94 & 0.98 & 0.53 & 1.03 & 2.28 & 2.07 & 2.66\\
37 & 0.09 & 0.15 & 0.33 & 1.44 & 0.20 & 0.55 & 0.53 & 0.42 & 0.93 & 0.95 & 0.95 & 0.44 & 0.77 & 2.20 & 2.16 & 2.75\\
41 & 0.15 & 0.06 & 0.07 & 1.36 & 0.42 & 0.36 & 0.85 & 0.41 & 0.97 & 0.96 & 0.95 & 0.47 & 1.86 & 1.45 & 3.31 & 2.66\\
53 & 0.22 & 0.12 & 0.08 & 1.35 & 0.25 & 0.36 & 0.27 & 0.41 & 0.89 & 0.94 & 0.96 & 0.49 & 0.93 & 1.30 & 1.02 & 2.66\\
87 & 0.06 & 0.03 & 0.18 & 1.49 & 0.22 & 0.33 & 0.27 & 0.43 & 0.95 & 0.95 & 0.90 & 0.36 & 0.88 & 1.31 & 1.01 & 2.70\\
89 & 0.04 & 0.07 & 0.12 & 1.50 & 0.27 & 0.43 & 0.34 & 0.41 & 0.95 & 0.95 & 0.93 & 0.35 & 1.06 & 1.54 & 1.21 & 2.78\\
439 & 0.01 & 0.05 & 0.05 & 1.46 & 0.29 & 0.44 & 0.29 & 0.43 & 0.92 & 0.92 & 0.96 & 0.39 & 0.99 & 1.51 & 1.19 & 2.69\\
\hline
\end{tabular}}
\caption{Comparison of \texttt{DLL} and \texttt{SIHR} in the semi-real simulation study. The columns indexed with ``M'',``M+'' and ``M-'' report the performance of our proposed \texttt{DLL} inference methods for $f'(\eval)$, with $\eval$ set as mean, mean + standard error, and mean - standard error, respectively. \texttt{SIHR} refers to the high-dimensional inference methods assuming the linear model. The columns indexed with ``Bias'', and ``SE" report the absolute bias, and the standard error computed by 500 estimates, respectively; the columns indexed with ``Coverage'' report the empirical coverage level and the columns indexed with ``Length'' report the average CI length. For \texttt{SIHR}, the bias is taken as the minimal bias across the three evaluation points and the coverage is taken as the maximal coverage across the three evaluation points.}
\label{table:5}
\end{table}
%%%%%%%%%%%%%%%%%%%%%%%%%%%%%%%%%%%%%%%%%%%%%%%%%%%%%%%%%%%%

%%%%%%%%%%%%%%%%%%%%%%%%%%%%%%%%%%%%%%%%%%%%%%%%%%%%%%%%%%%%%%%%%%%%%%%%%%%%%%
\section{Conclusion and discussion}
\label{sec: con-discussion}
%%%%%%%%%%%%%%%%%%%%%%%%%%%%%%%%%%%%%%%%%%%%%%%%%%%%%%%%%%%%%%%%%%%%%%%%%%%%%%
We have proposed the decorrelated local linear estimator to mitigate the error caused by estimating the unknown nuisance functions in the high-dimensional additive model. We have established the asymptotic normality of our proposed estimator. We demonstrate the validity of the theoretical results in moderate samples sizes and provide practical recommendations for the algorithm implementation. Our proposed decorrelation idea is a novel and computationally efficient method designed for bias correction in non-parametric models. An interesting future research direction is extending our proposed method to accommodate other kernel functions and higher-order local polynomials. In addition to inference for $f'(\eval),$ there are other interesting statistical inference problems in the high-dimensional additive model, including confidence interval construction for $f(\eval)$ and the significance test $H_0: f=0$. We leave these problems for future research.

\bibliographystyle{plain}
\bibliography{HDRef}

\newpage
\appendix

%!TEX root = DLL-Draft.tex

\setcounter{page}{1}

\section{Additional Discussions}
%%%%%%%%%%%%%%%%%%%%%%%%%%%%%%%%%%%%%%%%%%%%%%%%%%%%%%%%%%%%%%%%%%%%%%
\subsection{Algorithm without Data Swapping}
\label{sec: no data swapping}
%%%%%%%%%%%%%%%%%%%%%%%%%%%%%%%%%%%
In this section, we present the \texttt{DLL} estimator without data swapping. Different from the \texttt{DLL} estimator with data swapping, we use all the samples, rather than half of them, when fitting the sparse additive model and constructing the decorrelation weights. Following the steps of Algorithm \ref{algo: DDL}, we make a few changes to implement the \texttt{DLL} estimator without data swapping.

In Step 1, implement the sparse additive model as the following optimization problem with $M \geq 1$ and $\lambda > 0$:
\begin{equation}
    \{\widehat{\beta}_j\}_{1\leq j\leq p}=\argmin_{\beta_j \in \R^{M},\; 0\leq j\leq p}\frac{1}{2n}\sum_{i=1}^n(Y_i-\sum_{j=0}^{p}\Psi_{i,j}^{\intercal}\beta_j)^2+\lambda\sum_{j=0}^{p}\sqrt{\beta_j^{\intercal} \left(\frac{1}{n}\sum_{i=1}^n \Psi_{i,j} \Psi_{i,j}^{\intercal}\right) \beta_j}.
    \label{eq: step 1 pred}
\end{equation}

In Step 2, construct the initial estimator $\{\widehat{g}(X_i)\}_{1\leq i\leq n}$ as:
\begin{equation}
    \widehat{g}(X_i)=\sum_{j=1}^{p}\Psi_{i,j}^{\intercal}\widehat{\beta}_j \quad \text{and}\quad \widehat{f}(D_i)=\Psi_{i,0}^{\intercal}\widehat{\beta}_0.
    \label{eq: step 2 pred}
\end{equation}

In Step 5, implement Lasso algorithm as follows with the tuning parameter $\lambda_1 > 0$:
\begin{equation*}
        \widehat{\gamma} = \argmin_{\gamma\in \R^{p}}\frac{1}{2n}\sum_{i=1}^n\left(D_i-X_{i}^{\intercal}\gamma\right)^2+\lambda_{1} \sum_{j=1}^{p} \frac{\|X_j\|_2}{\sqrt{n}}|\gamma_j|,
\end{equation*}
and compute
\begin{equation*}
    \widehat{\mu}_i = \eval-X_{i}^{\intercal}\widehat{\gamma} \quad \text{and}\quad \widehat{\delta}_i=D_i-X_{i}^{\intercal}\widehat{\gamma} \quad \text{for} \quad 1 \leq i \leq n.
\end{equation*}

In Step 6, construct $\{\widehat{l}(X_i)\}_{1\leq i\leq n}$ and the weights $\{\widetilde{W}_i\}_{1\leq i\leq n}$ as:
\begin{equation*}
    \widetilde{W}_{i} = (D_i-\eval)-\widehat{l}(X_i) \quad \text{with}\quad
    \widehat{l}(X_{i},\widehat{\gamma}) = \frac{\frac{1}{n}\sum_{j=1}^n(\widehat{\delta}_j-\widehat{\mu}_i) {\bf{1}}(|\widehat{\delta}_j-\widehat{\mu}_i|\leq h)}{\frac{1}{n}\sum_{j=1}^n {\bf{1}}(|\widehat{\delta}_j-\widehat{\mu}_i|\leq h)}.
\end{equation*}

The other steps are the same as Algorithm \ref{algo: DDL}. We utilize the same parameter tuning procedures stated in Section \ref{sec: simulation}. We compare our constructed confidence intervals with and without data swapping in Section \ref{sec: swap and trans}.

%%%%%%%%%%%%%%%%%%%%%%%%%%%%%%%%%%%
%%%%%%%%%%%%%%%%%%%%%%%%%%%%%%%%%%%%%%%%%%%%%%%%%%%%%%%%%%%%%%%%%%%%%%
\subsection{Double Penalization}
\label{sec: double penalization}
%%%%%%%%%%%%%%%%%%%%%%%%%%%%%%%%%%%%%%%%%%%%%%%%%%%%%%%%%%%%%%%%%%%%%%

In the following, we review the double penalization method \cite{tan2017penalized} to construct initial estimators of $f$ and $g$. This can be viewed as an alternative method to the estimators in \eqref{eq: step 1 pred} and \eqref{eq: step 2 pred}. To construct the penalty term, we define the complexity measure of a univariate function $f$ as 
\begin{equation}
\mathcal{C}(f)=\lambda_n (\|f\|_n+\rho_n\|f\|_{F})
\label{eq: uni complexity}
\end{equation}
where $\lambda_n>0$, $\rho_n>0$, $\|f\|_n=(\sum_{i=1}^{n} f^2(D_{i}))^{1/2}$ denotes the function's empirical $L_2$ norm, and $\|f\|_{F}=(\int [f''(t)]^2 dt)^{1/2}$ is a measure of the function's smoothness. 

%For the additive function $g$, we define 
%\begin{equation}
%\mathcal{C}(g)=\sum_{j=1}^{p} \mathcal{C}(g_j) 
%\end{equation}

 Specifically,  with the positive tuning parameters $\rho_{n}>0$ and $\lambda_{n}>0$, we define the initial estimators as 
\begin{equation*}
\begin{aligned}
\left\{\widehat{f}, \left\{\widehat{g}_j\right\}_{1\leq j\leq p}\right\}=\arg\min\frac{1}{n}\sum_{i=1}^{n}[Y_i-f(D_i)-\sum_{j=1}^{p}g_{j}(X_{ij})]^2+ \mathcal{C}(f)+\sum_{j=1}^{p}\mathcal{C}(g_j),
\end{aligned}
\label{eq: initial additive}
\end{equation*}
where the complexity measure $\mathcal{C}(\cdot)$ is defined in \eqref{eq: uni complexity}. 
%For the complexity term $\mathcal{C}(f),$ the components $\|f\|_n$ and $\|f\|_{F}$ are imposed to encourage sparsity and smoothness, respectively.

%%%%%%%%%%%%%%%%%%%%%%%%%%%%%%%%%%%%%%%%%%%%%%%%%%%%%%%%%%%%%%%%%%%%%%%%%%%%%%%%
\subsection{Initial Estimators with Quantile Transformation}
\label{sec: quantile}
%%%%%%%%%%%%%%%%%%%%%%%%%%%%%%%%%%%%%%%%%%%%%%%%%%%%%%%%%%%%%%%%%%%%%%%%%%%%%%%%
We consider the construction of the initial estimator $\widehat{g}$ by applying the quantile transformation to all variables.  Particularly, we transform $D_i$ to $\widetilde{D}_i$, with $$\widetilde{D}_i=\frac{\text{Ordering of } D_i}{n} \in (0,1];$$ similarly, for $1\leq j\leq p,$ we  transform $X_{i,j}$ to $\widetilde{X}_{i,j}$, with $$\widetilde{X}_{i,j}=\frac{\text{Ordering of } X_{i,j}}{n}\in (0,1].$$ We construct the initial estimators $\widehat{f}$ and $\widehat{g}$ by applying the sparse additive algorithm to $\{Y_i,\widetilde{D}_i,\widetilde{X}_i\}_{1\leq i\leq n}$. Except for constructing the initial estimator differently, the other steps are the same as those in Algorithm \ref{algo: DDL}. This \texttt{DLL} estimator with the extra quantile transformation is refered to as \texttt{Trans} and we do not apply the data swapping for \texttt{Trans} estimator. We compare the performance of \texttt{Trans} with the regular \texttt{DLL} estimator in Section \ref{sec: swap and trans}; see Tables \hyperref[table:S10]{S10} and \hyperref[table:S11]{S11} for details.  

\subsection{Further discussions on the Condition (A3)}
\label{sec: A3 discussion}
In a more general setting, we may plug-in the existing convergence rate of ${\rm Err}(\widehat{g})$ and then the condition \eqref{eq: general inference condition} is reduced to a simultaneous condition on $k\coloneqq\|\gamma\|_0$ and $s\coloneqq\|g\|_0$, where $\|g\|_0$ denotes the number of non-zero functions of $\{g_j\}_{1\leq j\leq p}.$ Particularly, we follow \cite{tan2017penalized} by assuming that the individual functions $\{g_j\}_{1\leq j\leq p}$ belong to the Sobolev space  $\mathcal{W}^{m}_2$ for $m>1/2$ and $f''$ is continuous. We apply Proposition 4 and Theorem 2 in \cite{tan2017penalized} and Corollaries 4 and 5 in \cite{guo2019extreme} and establish that 
$$
{\rm Err}^2(\widehat{g})\lesssim n^{-\frac{4}{5}}+s\cdot n^{-\frac{2m}{2m+1}}+(s+1)\cdot{\log p}/{n}.
$$
 Then the condition  \eqref{eq: general inference condition} is simplified as 
\begin{equation*}
k\cdot s \ll n^{\frac{2m}{2m+1}+\frac{3}{5}}  \; \text{and}\; \max\{k, s\cdot n^{\frac{1}{2m+1}}\}\ll n \quad \text{up to a polynomial order of}\; \log p. 
\end{equation*}
If we set $m=2$, then the above sparsity condition is much weaker than the one in \cite{gregory2016optimal}, which requires $s\ll n^{\frac{3}{10}}$ and $k \ll n^{\frac{4}{15}}$ up to a polynomial order of $\log p.$

\subsection{Consistent estimators of $\sigma^2$}
\label{sec: consistent sigmasq}
Similar to the definition of ${\rm Err}(\widehat{g})$ in \eqref{eq: g error}, we use ${\rm Err}(\widehat{f})$ to denote the accuracy measure of $\widehat{f}$, which is defined as follows: with probability larger than $1-\min$ for some positive constant $c>0,$ 
\begin{equation}
\sqrt{\E_{D_{*}}(\widehat{f}^a(D_{*})-f(D_{*}))^2+\E_{D_{*}}(\widehat{f}^b(D_{*})-f(D_{*}))^2}\lesssim {\rm Err}(\widehat{f}),
\label{eq: f error}
\end{equation}
where $\widehat{f}^a$ and $\widehat{f}^b$ are defined in \eqref{eq: initial estimator} and the expectation is taken with respect to the independent copy $D_{*}$ of $\{D_i\}_{1\leq i\leq n}.$ 
\begin{Proposition} Suppose that Condition {\rm (A1)} holds and $\max\{{\rm Err}(\widehat{f}),{\rm Err}(\widehat{g})\}\rightarrow 0.$ Then the estimator $\widehat{\sigma}^2$ defined in \eqref{eq: noise estimator} satisfies $\widehat{\sigma}^2\cip \sigma^2.$
\label{prop: consistent sigma}
\end{Proposition}
Proposition \ref{prop: consistent sigma} shows that our proposed $\widehat{\sigma}^2$ is consistent if both $\widehat{f}$ and $\widehat{g}$ are consistent. 
The proof of Proposition \ref{prop: consistent sigma} is presented in Section \ref{sec: consistent sigma}.

%%%%%%%%%%%%%%%%%%%%%%%%%%%%%%%%%%%%%%%%%%%%%%%%%%%%%%%%%%%
\section{Notations, Events and Lemmas}
%%%%%%%%%%%%%%%%%%%%%%%%%%%%%%%%%%%%%%%%%%%%%%%%%%%%%%%%%%%
We introduce some notations and events, which will be used throughout the proof. Let ${q}(D_{i}\mid X_{i})$ denote the conditional distribution of $D_{i}$ given $X_{i}$ and $\phi$ denote the density function of the error $\delta_i=D_{i}-X_{i}^{\intercal}\gamma.$ Since $D_{i}-X_{i}^{\intercal}\gamma$ is independent of $X_{i}$, we have 
$${q}(D_{i}=\eval\mid X_{i})=\phi(\eval-X_{i}^{\intercal}{\gamma})=\phi(\mu_i) \quad \text{with}\quad \mu_i=\eval-X_{i}^{\intercal}{\gamma}.$$
We express the density function $\pi$ of the random variable $D_i$ as
\begin{equation*}
\pi(\eval)=\E_{X_i} \left[{q}(D_{i}=\eval\mid X_{i})\right]=\E_{X_i} \left[\phi \left({\eval-X_{i}^{\intercal}{\gamma}}\right)\right].
\end{equation*}
Define the following events,
\begin{equation}
\begin{aligned}
\mathcal{A}_0&=\left\{
\max_{i\in \mathcal{I}_a}\max_{\left|\delta-\mu_i\right|\leq r} \frac{|\phi'(\delta)|}{\phi(\mu_i)}\leq C_1(n),
 \quad \max_{i\in \mathcal{I}_a}\max_{\left|\delta-\mu_i\right|\leq r} \frac{|\phi''(\delta)|}{\phi(\mu_i)} \leq C_2(n)\right\},\\
\mathcal{A}_1&=\left\{\|\widehat{\gamma}^{b}-\gamma\|\lesssim \sqrt{\frac{k \log p}{n}}, \; \max_{i\in \mathcal{I}_a}|X_{i}^{\intercal}(\widehat{\gamma}^{b}-\gamma)|\leq C^* \sqrt{\frac{k \log p \log n}{n}}\right\},\\
\mathcal{A}_2&=\left\{\E_{X_{*}}(\widehat{g}^b(X_{*})-g(X_{*}))^2\lesssim {\rm Err}^2(\widehat{g})\right\}.
%\mathcal{A}_2&=\left\{\left|\widehat{\sigma}_2^2-\sigma_2^2\right|\lesssim \frac{1}{\sqrt{n}}+\frac{k \log p}{n}\right\}\\
%\mathcal{A}_{3,i}&=\left\{\|X_{i}^{\intercal}\gamma\|_2\lesssim \|\gamma\|_2 \sqrt{\log n}\right\}.
%\mathcal{A}_{4,i}&=\left\{\max_{1\leq i\leq n}\left|X_{i}^{\intercal}(\widehat{\gamma}-\gamma)\right|\lesssim \sqrt{\frac{k \log p \log n}{n}}\right\}
\end{aligned}
\label{eq: events}
\end{equation}
%\Zijian{Check whether $\mathcal{A}_2$ and $\mathcal{A}_3$ are used.}
By the definitions of $C_1(n)$ and $C_2(n)$ in \eqref{eq: density ratio}, we have $\mathbf{P}(\mathcal{A}_0^c)\leq \min\{n,p\}^{-c}$ for some positive constant $c>0.$
Throughout the proof, we shall assume $\mathbf{P}(\mathcal{A}_0^c)\ll h\pi(\eval)$ and this will automatically hold in our considered regime.

Define $\mathcal{A}=\mathcal{A}_0\cap \mathcal{A}_1\cap \mathcal{A}_2$. Theorem 7.2 in \cite{bickel2009simultaneous} implies that the Lasso estimator $\widehat{\gamma}^{b}$ satisfies $$\mathbf{P}\left(\|\widehat{\gamma}^{b}-\gamma\|\lesssim \sqrt{\frac{k \log p}{n}}\right)\geq 1-p^{-c},$$
for some positive constant $c>0$. Conditioning on the data in $\mathcal{I}_b$, the random variable ${X_{i}^{\intercal}(\widehat{\gamma}^{b}-\gamma)}/{\|\widehat{\gamma}^{b}-\gamma\|_2}$ is sub-gaussian random variable, which implies $$\mathbf{P}\left(\max_{i\in \mathcal{I}_a}\left|{X_{i}^{\intercal}(\widehat{\gamma}^{b}-\gamma)}/{\|\widehat{\gamma}^{b}-\gamma\|_2}-\E \left[X_{i}^{\intercal}(\widehat{\gamma}^{b}-\gamma)/{\|\widehat{\gamma}^{b}-\gamma\|_2}\right|\mid \mathcal{I}_b\right]\gtrsim \sqrt{\log n}\mid \mathcal{I}_b\right)\leq n^{-c},$$
for some positive constant $c>0.$ Note that
$$\left|\E \left[X_{i}^{\intercal}(\widehat{\gamma}^{b}-\gamma)/{\|\widehat{\gamma}^{b}-\gamma\|_2}\mid \mathcal{I}_b\right]\right|\leq \sqrt{(\widehat{\gamma}-\gamma)^{\intercal}\Sigma(\widehat{\gamma}-\gamma)}/{\|\widehat{\gamma}^{b}-\gamma\|_2}\leq C.$$
The above two inequalities imply that, there exists a constant $C^*>0$ independent of $n$ and $p$ such that 
%for some positive constant $C>0,$ we establish  
$$\mathbf{P}\left(\max_{i\in \mathcal{I}_a}|X_{i}^{\intercal}(\widehat{\gamma}^{b}-\gamma)|\leq C^* \sqrt{\frac{k \log p \log n}{n}}\right)\geq 1-\min\{n,p\}^{-c},$$
for some positive constant $c>0.$ Together with the definition in \eqref{eq: density ratio} and \eqref{eq: g error}, we establish 
\begin{equation}
\mathbf{P}\left(\mathcal{A}\right)\geq 1-\min\{n,p\}^{c} \;\; \text{for some constant}\; c>1.
\label{eq: high prob}
\end{equation}

%\Zijian{Checking and polishing the following two lemmas.}
The following lemma states the expectation of terms involved with $K_h(D_i)$, whose proof can be found in Section \ref{sec: exp lemma}.
\begin{Lemma} Suppose that Condition {\rm (A2)} holds, then we have
\begin{equation}
\left|\frac{\E\left(K_{h}(D_{i})\mid X_{i}\right)}{{q}(\eval\mid X_{i})}-1\right|\cdot {\1}_{\mathcal{A}_0}\leq \frac{h^2}{6}C_2(n);
\label{eq: conditional zero order}
\end{equation}
\begin{equation}
\left|\frac{\E\left(K_{h}(D_{i})\right)}{\pi(\eval)}-1\right|\lesssim {h^2}C_2(n)+ \frac{\mathbf{P}(\mathcal{A}_0^c)}{\pi(\eval)};
\label{eq: expectation 0}
\end{equation}
\begin{equation}
\left|\frac{\E\left[(D_{i}-\eval)K_{h}(D_{i})\right]}{{\pi}(\eval)}\right| \leq \frac{1}{3}h^2 (C_1(n)+\frac{3}{8}hC_2(n))+\frac{\mathbf{P}(\mathcal{A}_0^c)}{{\pi}(\eval)};
\label{eq: first moment bound final}
\end{equation}
\begin{equation}
\left|\frac{\E\left[(D_{i}-\eval)^2K_{h}(D_{i})\right]}{\frac{1}{3}h^2\pi(\eval)}-1\right| \lesssim \frac{1}{10}h^2 C_2(n)+\frac{\mathbf{P}(\mathcal{A}_0^c)}{h\pi(\eval)};
\label{eq: second order bound final}
\end{equation}
\begin{equation}
\left|\frac{\E\left({W}^2_{i}K_{h}(D_{i})\mid X_{i}\right)}{{\frac{1}{3}h^2q(\eval\mid X_i)}}-1\right|\cdot \mathbf{1}_{\mathcal{A}_{0}}\lesssim h^2[C_1^2(n)+C_2(n)];
\label{eq: main part}
\end{equation}
\begin{equation}
\left|\frac{\E W_i^2K_{h}^2(D_{i})}{\frac{1}{3}h\pi(\eval)}-1\right|\lesssim h^2[C_1^2(n)+C_2(n)]+\frac{\mathbf{P}(\mathcal{A}_{0}^c)}{h\pi(\eval)};
\label{eq: expectation 1}
\end{equation}
\begin{equation}
\left|\frac{\E W_i(D_{i}-\eval)K_{h}(D_{i})}{\frac{1}{3}h^2\pi(\eval)}-1\right|\lesssim h^2[C_1^2(n)+C_2(n)]+\frac{\mathbf{P}(\mathcal{A}_{0}^c)}{h\pi(\eval)};
\label{eq: expectation 2}
\end{equation}
\begin{equation}
\left|\E W_{i} \frac{(D_{i}-\eval)^2}{2} K_{h}(D_{i}) \cdot \mathbf{1}_{\mathcal{A}_{0}^{c}}\right|\lesssim {h^4}\left[C_1(n)+hC_2(n)\right]\pi(\eval)+h^2 \mathbf{P}(\mathcal{A}_{0}^{c}).
\label{eq: third moment}
\end{equation}
\label{lem: expectation lemma}
\end{Lemma}

The following lemma is about the concentration results for terms involved with $K_h(D_i)$, whose proof can be found in Section \ref{sec: con lemma}.

\begin{Lemma} Suppose that Condition {\rm (A2)} holds, then for a sufficiently large $n$, with probability $1-\exp(-t^2)$, 
\begin{equation}
c\pi(\eval) \left[ 1-\frac{t}{\sqrt{nh \pi(\eval)}}\right] \leq \left|\frac{1}{n}\sum_{i=1}^{n} K_{h}(D_{i})\right|\leq C \pi(\eval) \left[ 1+ \frac{t}{\sqrt{nh \pi(\eval)}}\right];
\label{eq: prob bound 1}
\end{equation}
\begin{equation}
\left|\frac{1}{n}\sum_{i=1}^{n} W_{i}K_{h}(D_{i})\right|\lesssim t \sqrt{\frac{h}{n}\pi(\eval)};
\label{eq: prob bound 2}
\end{equation}
\begin{equation}
\left|\frac{1}{n}\sum_{i=1}^{n}{(D_{i}-\eval)} K_{h}(D_{i})\right|\lesssim \frac{{\pi}(\eval)}{3}h^2 (C_1(n)+\frac{3}{8}hC_2(n))+{\mathbf{P}(\mathcal{A}_0^c)}+t\sqrt{\frac{h}{n} \pi(\eval)};\label{eq: prob bound 4}
\end{equation}
%\begin{equation}
%\left|\frac{1}{n}\sum_{i=1}^{n}W_i{(D_{i}-\eval)} K_{h}(D_{i})\right|\lesssim 
%{h^4}\left[C_1(n)+hC_2(n)\right]\pi(\eval)+h^2 \mathbf{P}(\mathcal{A}_{0}^{c})+t\sqrt{\frac{h^5}{n} \pi(\eval)}
%\label{eq: prob bound 5}
%\end{equation}
\begin{equation}
\left|\frac{1}{n}\sum_{i=1}^{n}W_i{(D_{i}-\eval)} K_{h}(D_{i})-\E W_i(D_{i}-\eval)K_{h}(D_{i})\right|\leq C h^2\pi(\eval)t\sqrt{\frac{1}{nh\pi(\eval)}};
\label{eq: prob bound 5}
\end{equation}
\begin{equation}
\left|\frac{1}{n}\sum_{i=1}^{n}W_i\frac{(D_{i}-\eval)^2}{2} K_{h}(D_{i})\right|\lesssim {h^4}\left[C_1(n)+hC_2(n)\right]\pi(\eval)+h^2 \mathbf{P}(\mathcal{A}_{0}^{c})+t\sqrt{\frac{h^5}{n} \pi(\eval)};
\label{eq: prob bound 3}
\end{equation}

\begin{equation}
\left|\frac{1}{n}\sum_{i=1}^{n}W_i^2K_{h}^2(D_{i})-\E W_i^2K_{h}^2(D_{i})\right|\lesssim t\sqrt{\frac{h\pi(\eval)}{{n}}}.
\label{eq: prob bound 8}
\end{equation}
%\begin{equation}
%\left|\frac{1}{n}\sum_{i=1}^{n}W_i(D_{i}-\eval)K_{h}(D_{i})-\E W_i(D_{i}-\eval)K_{h}(D_{i})\right|\lesssim\frac{th}{\sqrt{n}}.
%\label{eq: prob bound 7}
%\end{equation}
\label{lem: concentration lemma}
\end{Lemma}

%\Zijian{The following two bounds are not checked yet.}
%In addition, if $X_{i}$ is Sub-gaussian random vector, then with probability larger than $1-\frac{1}{t},$
%\begin{equation}
%\frac{1}{n}\sum_{i=1}^{n}[X_{i}^{\intercal}\gamma]^2K_{h}(D_{i})\lesssim \left(1+\frac{t}{\sqrt{nh}}\right) \|\gamma\|_2^2
%\label{eq: prob bound 6}
%\end{equation}

Combining \eqref{eq: prob bound 1} and \eqref{eq: prob bound 2}, we establish that, with probability $1-\exp(-t^2)$, 
\begin{equation}
\left|\bar{\mu}_{W}\right|\lesssim t \sqrt{\frac{h}{n \pi(\eval)}}.
\label{eq: centering rate}
\end{equation}

\section{Proof of Theorem \ref{thm: DLL error reduction}}
%%%%%%%%%%%%%%%%%%%%%%%%%%%%%%%%%%%%%%%%%%%%%%%%%%%%%%%%%%%
Recall the definition of $\widehat{W}_i$ in \eqref{eq: final D}.
We define an accuracy measure of estimating the decorrelation weights as 
\begin{equation}
{\rm Err}(\widehat{W})=\sqrt{\frac{1}{n}\sum_{i=1}^{n}[\widehat{W}_{i}-\left({W}_{i}-\bar{\mu}_{W}\right)]^2 K_{h}(D_i)} \quad \text{with}\; \bar{\mu}_{W}=\frac{\frac{1}{n}\sum_{i=1}^{n}{W}_{i}K_{h}(D_i)}{\frac{1}{n}\sum_{i=1}^{n}K_{h}(D_i)}.
\label{eq: initial error D}
\end{equation}

We first introduce the following important intermediate results. With probability larger than $1-\frac{C}{t^2}-\min\{n,p\}^{-c}$ for some $t>1$ and positive constants $C>0, c>0$, 
\begin{equation}
\left|\frac{1}{n} \sum_{i=1}^{n}{W}_{i}\dif K_{h}(D_i)\right| \leq C t\sqrt{{h}/{n}} \cdot{\rm Err}(\widehat{g}),
\label{eq: bias b expectation a}
\end{equation}

\begin{equation}
\left|{\frac{1}{n\widehat{S}_{n}}\sum_{i=1}^{n}\widehat{W}_{i}\dif K_{h}(D_i)}\right|\leq t \left(t\sqrt{\frac{1}{nh^3 \pi^2(\eval)}}+\frac{{\rm Err}(\widehat{W})}{h^2\pi(\eval)}\right){\rm Err}(\widehat{g}),%\end{aligned}
\label{eq: bias b ratio}
\end{equation}
and 
\begin{equation}
{\rm Err}(\widehat{W})\lesssim  t\left(\sqrt{\frac{h}{n}}+\sqrt{\frac{k\log p\log n}{n}} h^2 \left(C_1^2(n)+C_2(n)\right)\right).
\label{eq: weight error}
\end{equation}
A combination of \eqref{eq: bias b ratio} and \eqref{eq: weight error} leads to the bound \eqref{eq: high dim error}. We shall prove \eqref{eq: bias b expectation a}, \eqref{eq: bias b ratio}, and \eqref{eq: weight error} in Sections \ref{sec: result 1}, \ref{sec: result 2}, and \ref{sec: result 3}, respectively. 

%%%%%%%%%%%%%%%%%%%%%%%%%%%%%%%%%%%%%%%%%%%%%%%%%%%%%%%%%%%
\subsection{Proof of \eqref{eq: bias b expectation a}} 
\label{sec: result 1}
%%%%%%%%%%%%%%%%%%%%%%%%%%%%%%%%%%%%%%%%%%%%%%%%%%%%%%%%%%%
Recall the following notations, 
\begin{itemize}
\item $\mathcal{I}_a$ and $\mathcal{I}_b$ are two disjoint subsets with approximately equal sample size, with $\mathcal{I}_a\cap \mathcal{I}_b$ empty and  $\mathcal{I}_a\cup \mathcal{I}_b=\{1,2,\cdots,n\}.$
\item $\widehat{g}^{a}$ and $\widehat{g}^{b}$ denote the initial estimator of $g$ based on the data $(X_{i}, D_i, Y_{i})_{i\in\mathcal{I}_a}$ and $(X_{i}, D_i , Y_{i})_{i\in \mathcal{I}_b}$, respectively. 
\end{itemize} The proof relies on the independence created by data swapping. 
Define the estimation error as $\Delta^{a}(X_{i})=\widehat{g}^{a}(X_{i})-g(X_{i})$ and $\Delta^{b}(X_{i})=\widehat{g}^{b}(X_{i})-g(X_{i}).$ We write $\E(\cdot\mid{\mathcal{I}_{a}})$,${\rm Var}(\cdot\mid{\mathcal{I}_{a}})$ and $\mathbf{P}(\cdot\mid{\mathcal{I}_{a}})$ as the expectation, variance and probability conditioning on the sample $(X_{i},D_i,Y_{i})_{i\in \mathcal{I}_a},$ respectively. Similarly, we define $\E(\cdot\mid{\mathcal{I}_{b}}),{\rm Var}(\cdot\mid{\mathcal{I}_{b}})$ and $\mathbf{P}(\cdot\mid{\mathcal{I}_{b}})$ conditioning on $(X_{i},D_i,Y_{i})_{i\in \mathcal{I}_b}$.  For $1\leq i\leq n$, we have the following decomposition
{\small
\begin{equation}
\frac{1}{n}\sum_{i=1}^{n}{W}_{i}\dif K_{h}(D_{i})=\frac{1}{n}\sum_{i\in \mathcal{I}_a}{W}_{i}\Delta^{b}(X_{i})K_{h}(D_{i})+\frac{1}{n}\sum_{i \in \mathcal{I}_b}{W}_{i}\Delta^{a}(X_{i})K_{h}(D_{i}).
\label{eq: swapping decomp1}
\end{equation}
}
We will control the first term $\frac{1}{n}\sum_{i\in \mathcal{I}_a}{W}_{i}\Delta^{b}(X_{i})K_{h}(D_{i})$ in the following and the second term can be controlled by a similar argument. 
Since
$$\mathbf{P}\left(|\frac{1}{n} \sum_{i\in \mathcal{I}_a}{W}_{i}\Delta^{b}(X_{i})K_{h}(D_{i})|\neq  |\frac{1}{n} \sum_{i\in \mathcal{I}_a}{W}_{i}\Delta^{b}(X_{i})K_{h}(D_{i})\cdot \mathbf{1}_{\mathcal{A}_{0}\cap \mathcal{A}_{2}}| \right)\leq \min\{n,p\}^{-c},$$ 
it is sufficient to analyze 
$$\frac{1}{n} \sum_{i=1}^{n}{W}_{i}\Delta^{b}(X_{i})K_{h}(D_{i})\cdot \mathbf{1}_{\mathcal{A}_{0}\cap \mathcal{A}_{2}},$$
where $\mathcal{A}_0$ and $\mathcal{A}_2$ are defined in \eqref{eq: events}.
The above term has two sources of randomness: the initial estimator $\Delta^{b}=\widehat{g}^{b}-g$ and the data $\{X_{i}, D_i\}_{i\in \mathcal{I}_a}$. Since the randomness of $\Delta^{b}$ is induced from the data $(X_{i}, D_i, Y_{i})_{i\in \mathcal{I}_b}$, the estimation error $\Delta^{b}$ is independent of the data $\{X_{i},D_i\}_{i\in \mathcal{I}_a}$. 

Since ${W}_{i}$ is constructed such that $
\E \left[{W}_{i}K_{h}(D_i)\mid X_i\right]=0
$, then we have 
%$\frac{1}{n}\sum_{i\in \mathcal{I}_a}{W}_{i}\Delta^{b}(X_{i})K_{h}(D_{i})\cdot \mathbf{1}_{\mathcal{A}_{0}}$ satisfies
\begin{equation}
\E\left(\frac{1}{n}\sum_{i\in \mathcal{I}_a}{W}_{i}\Delta^{b}(X_{i})K_{h}(D_{i})\cdot \mathbf{1}_{\mathcal{A}_{0}\cap \mathcal{A}_{2}}\mid {\mathcal{I}_b}, \{X_{i}\}_{i\in \mathcal{I}_a}\right)=0.
\label{eq: conditional zero mean}
\end{equation}
%which implies 
%$$\E\left(\frac{1}{n}\sum_{i\in \mathcal{I}_a}{W}_{i}\Delta^{b}(X_{i})K_{h}(D_{i})\cdot \mathbf{1}_{\mathcal{A}_{0}}\mid {\mathcal{I}_b}\right)=0.$$
%By definition of the data swapping estimators defined in \eqref{eq: def of f2est} and \eqref{eq: def of lest}, we have 
We control the second order moment as 
\begin{equation*}
\begin{aligned}
%&{\Var}\left(\frac{1}{n}\sum_{i\in \mathcal{I}_a}{W}_{i}\Delta^{b}(X_{i})K_{h}(D_{i})\cdot \mathbf{1}_{\mathcal{A}_{0}}\mid {\mathcal{I}_b}\right)\\ 
&{\E}\left(\left(\frac{1}{n}\sum_{i\in \mathcal{I}_a}{W}_{i}\Delta^{b}(X_{i})K_{h}(D_{i})\cdot \mathbf{1}_{\mathcal{A}_{0}\cap \mathcal{A}_{2}}\right)^2\right)\\
&=\E\left[ {\E}\left(\left(\frac{1}{n}\sum_{i\in \mathcal{I}_a}{W}_{i}\Delta^{b}(X_{i})K_{h}(D_{i})\cdot \mathbf{1}_{\mathcal{A}_{0}\cap \mathcal{A}_{2}}\right)^2\mid {\mathcal{I}_b}, \{X_{i}\}_{i\in \mathcal{I}_a}\right)\right]\\
&=\frac{1}{n^2}\sum_{i\in \mathcal{I}_a}\E\left[\E\left(W^2_{i}K_{h}^2(D_{i})\mid {\mathcal{I}_b},\{X_{i}\}_{i\in \mathcal{I}_a}\right)(\Delta^{b}(X_{i}))^2\cdot \mathbf{1}_{\mathcal{A}_{0}\cap \mathcal{A}_{2}}\right],
\end{aligned}
\end{equation*}
where the last equality follows from \eqref{eq: conditional zero mean}.
%}
%where $\Delta^{a}(X_{i})=\widehat{f}^a_2(X_{i})-{f}_2(X_{i})$ and $\Delta^{b}(X_{i})=\widehat{f}^b_2(X_{i})-{f}_2(X_{i})$.
%$$
%\left|\frac{(\Delta^{b}(X_{i}))^2\E\left({W}^2_{i}K^2_{h}(D_{i})\mid X_{i}\right)\cdot \mathbf{1}_{\mathcal{A}_{0}}}{(\Delta^{b}(X_{i}))^2\frac{2}{3}h{q}(\eval\mid X_{i})\mathbf{1}_{\mathcal{A}_{0}}}-1\right|\cip 0,
%$$}
%and hence 
By \eqref{eq: main part} and the definition of ${\rm Err}(\widehat{g})$ in \eqref{eq: g error}, we establish 
$$\E\left[\E\left(W^2_{i}K_{h}^2(D_{i})\mid {\mathcal{I}_b},\{X_{i}\}_{i\in \mathcal{I}_a}\right)(\Delta^{b}(X_{i}))^2\cdot \mathbf{1}_{\mathcal{A}_{0}\cap \mathcal{A}_{2}}\right]\lesssim h{\rm Err}(\widehat{g})^2.$$
%We apply the iterative expectation and obtain
%\begin{equation}
%\begin{aligned}
%&\E\left[(\Delta^{b}(X_{i}))^2{W}^2_{i}K^2_{h}(D_{i})\cdot \mathbf{1}_{\mathcal{A}_{0}}\mid {\mathcal{I}_b}\right]\\
%&=\E\left[\left[(\Delta^{b}(X_{i}))^2{W}^2_{i}K^2_{h}(D_{i})\cdot \mathbf{1}_{\mathcal{A}_{0}}\mid {\mathcal{I}_b},\{X_{i}\}_{i\in \mathcal{I}_a}\right]\mid {\mathcal{I}_b} \right]\\
%&=\E\left((\Delta^{b}(X_{i}))^2 \E\left[{W}^2_{i}K^2_{h}(D_{i})\mid X_{i}\right]\cdot \mathbf{1}_{\mathcal{A}_{0}}\right)\lesssim h{\rm Err}(\widehat{g})^2, %\pi(\eval)
%\end{aligned}
%\end{equation}
%where the last inequality follows from . 
%\chzB{Cun-Hui, the part that is not sharp in the above bound is that we use ${q}(\eval\mid X_{i})\lesssim C$ to further bound the following inequality, 
%\begin{equation*}
%\E_{\mathcal{I}_a}\left[(\Delta^{b}(X_{i}))^2{W}^2_{i}K^2_{h}(D_{i})\cdot \mathbf{1}_{\mathcal{A}_{0}}\right]\leq \E_{\mathcal{I}_a}\left[(\Delta^{b}(X_{i}))^2\frac{2}{3}h{q}(\eval\mid X_{i})\cdot \mathbf{1}_{\mathcal{A}_{0}}\right]
%\end{equation*}}
Hence we establish 
{\small
\begin{equation*}
\mathbf{P}\left(\left|\frac{1}{n}\sum_{i\in \mathcal{I}_a}{W}_{i}\Delta^{b}(X_{i})K_{h}(D_{i})\cdot \mathbf{1}_{\mathcal{A}_{0}\cap \mathcal{A}_2}\right|\leq t\sqrt{\frac{h}{n}}\cdot {\rm Err}(\widehat{g})\right) \geq 1-\frac{1}{t^2}-\min\{n,p\}^{-c}.
\end{equation*}
}
%Together with the error bound 
%\begin{equation}
%\mathbf{P}_{\mathcal{I}_b}\left(\sqrt{\E_{\mathcal{I}_a}\left(\widehat{g}^{b}(X_{i})-g(X_{i})\right)^2} \leq {\rm Err}(\widehat{f}^b_2)\right)\geq 1,
%\label{eq: initial condition reform}
%\end{equation}
 
%we have 
%$
%\mathbf{P}\left(\left|\frac{1}{n}\sum_{i\in \mathcal{I}_a}{W}_{i}\Delta^{b}(X_{i})K_{h}(D_{i})\cdot \mathbf{1}_{\mathcal{A}_{0}}\right|\leq t\sqrt{{h}/{n}}\cdot {\rm Err}(\widehat{g})\right) \geq 1-\frac{1}{t}-\frac{1}{n^{c}}.
%$
By symmetry and the decomposition \eqref{eq: swapping decomp1}, we establish \eqref{eq: bias b expectation a}. 

%with taking $t=(nh\pi(\eval))^{\frac{1}{4}}.$
%\begin{equation}
%\mathbf{P}\left(\left|\frac{1}{n}\sum_{i=1}^{n}{W}_{i}\dif K_{h}(D_{i})\cdot \mathbf{1}_{\mathcal{A}_{0}}\right|\leq 2t\sqrt{\frac{h}{n}}\cdot {\rm Err}(\widehat{f}^b_2)\right) \geq 1-\frac{2}{t}-\frac{2}{n^{c}}-2\eta(n).
%\label{eq: control 1}
%\end{equation}

%%%%%%%%%%%%%%%%%%%%%%%%%%%%%%%%%%%%%%%%%%%%%%%%%%%%%%%%%%%%
\subsection{Proof of \eqref{eq: bias b ratio}} 
\label{sec: result 2}
%%%%%%%%%%%%%%%%%%%%%%%%%%%%%%
%%%%%%%%%%%%%%%%%%%%%%%%%%%%%%
We decompose $\frac{1}{n}\sum_{i=1}^{n}\widehat{W}_{i}\dif K_{h}(D_{i}) $ as 
{\small
\begin{equation}
\begin{aligned}
\frac{1}{n}\sum_{i=1}^{n}\left(\widehat{W}_{i}-\left({W}_{i}-\bar{\mu}_{W}\right)\right)\dif K_{h}(D_{i})
+\frac{1}{n}\sum_{i=1}^{n}{W}_{i}\dif K_{h}(D_{i})
-\bar{\mu}_{W}\cdot \frac{1}{n}\sum_{i=1}^{n}\dif K_{h}(D_{i}).
\label{eq: bias b decomposition}
\end{aligned}
\end{equation}
}
By the Cauchy-Schwarz inequality, we have 
\begin{equation}
\left|\frac{1}{n}\sum_{i=1}^{n}\left(\widehat{W}_{i}-\left({W}_{i}-\bar{\mu}_{W}\right)\right)\dif K_{h}(D_{i})\right|\leq {\rm Err}(\widehat{W})\cdot\sqrt{\frac{1}{n}\sum_{i=1}^{n}\dif^2 K_{h}(D_{i})},
\label{eq: cau 1}
\end{equation}
where ${\rm Err}(\widehat{W})$ is defined in \eqref{eq: initial error D}. Hence, it is sufficient to control $\sqrt{\frac{1}{n}\sum_{i=1}^{n}\dif^2 K_{h}(D_{i})}.$
%\begin{equation}
%\lesssim {\rm Err}(\widehat{g}) \sqrt{\pi(\eval)}%+C t\sqrt{\frac{1}{nh} {{\rm Err}^2(\widehat{g})} \pi(\eval)}
%\label{eq: bias b expectation b}
%\end{equation}
%%%%%%%%%%%%%%%%%%%%%%%%%%%%%%%%%%%%%%%%%%%%%%%%%%%%%%%%%%%%%%
%\underline{Proof of \eqref{eq: bias b expectation a} and \eqref{eq: bias b expectation b}}\\ 
%%%%%%%%%%%%%%%%%%%%%%%%%%%%%%%%%%%%%%%%%%%%%%%%%%%%%%%%%%%%%%
 Similar to \eqref{eq: swapping decomp1}, we have
\begin{equation}
\frac{1}{n}\sum_{i=1}^{n}\dif^2 K_{h}(D_{i})=\frac{1}{n}\sum_{i\in \mathcal{I}_a}\left|\Delta^{b}(X_{i})\right|^2K_{h}(D_{i}) +\frac{1}{n}\sum_{i\in \mathcal{I}_b}\left|\Delta^{a}(X_{i})\right|^2K_{h}(D_{i}), 
\label{eq: swapping decomp2}
\end{equation}
and it is sufficient to control $$\frac{1}{n}\sum_{i\in \mathcal{I}_a}\left|\Delta^{b}(X_{i})\right|^2K_{h}(D_{i})\cdot \mathbf{1}_{\mathcal{A}_{0}\cap \mathcal{A}_{2}}.$$ 
%.  Since
%$$\mathbf{P}_{\mathcal{I}_a}\left(\left|\frac{1}{n}\sum_{i\in \mathcal{I}_a}\left|\Delta^{b}(X_{i})\right|^2K_{h}(D_{i})\right|\neq  \left|\frac{1}{n}\sum_{i\in \mathcal{I}_a}\left|\Delta^{b}(X_{i})\right|^2K_{h}(D_{i})\cdot \mathbf{1}_{\mathcal{A}_{0}}\right| \right)\leq n^{-c},$$ %\chzB{Cun-Hui, this is not sharp for the same reasoning as the previous highlight.}
%\chzM{
Note that
{\small
\begin{equation*}
\begin{aligned}
&{\E}\left(\frac{1}{n}\sum_{i\in \mathcal{I}_a}\left|\Delta^{b}(X_{i})\right|^2K_{h}(D_{i})\cdot \mathbf{1}_{\mathcal{A}_{0}\cap \mathcal{A}_{2}}\mid {\mathcal{I}_b},\{X_{i}\}_{i\in \mathcal{I}_a}\right)\\=&\frac{1}{n}\sum_{i\in \mathcal{I}_a}\left|\Delta^{b}(X_{i})\right|^2 \cdot \E\left[K_{h}(D_{i})\mid X_{i}\right]\cdot \mathbf{1}_{\mathcal{A}_{0}\cap \mathcal{A}_{2}}\lesssim  {\rm Err}^2(\widehat{g}),
%&=\frac{2|\mathcal{I}_a|}{n^2h}\E_{\mathcal{I}_a}\left[\left|\Delta^{b}(X_{i})\right|^2 \cdot \E\left[K_{h}(D_{i})\mid X_{i}\right]\cdot \mathbf{1}_{\mathcal{A}_{0}}\right]\\
%&\lesssim \left(1+h^2\left(1+C_{u}^2\right)\exp\left(C_{u}\cdot \frac{h}{\sigma_{\delta}}\right)\right) \frac{1}{h}\E\left[g^2(X_{i}){q}(\eval\mid X_{i})\cdot \mathbf{1}_{A_{3,i}}\right]\lesssim \frac{1}{h}\E\left[g^2(X_{i})\right]
\end{aligned}
\end{equation*}}
where the last inequality follows from \eqref{eq: conditional zero order} and the bounded conditional density ${q}(\eval\mid X_{i}).$ 
%Then we have 
%$$
%{\E}\left(\frac{1}{n}\sum_{i\in \mathcal{I}_a}\left|\Delta^{b}(X_{i})\right|^2 K_{h}(D_{i})\cdot \mathbf{1}_{\mathcal{A}_{0}}\mid {\mathcal{I}_b}\right)\lesssim {\rm Err}^2(\widehat{g})$$
%and 
The above moment bound implies
\begin{equation*}
\mathbf{P}\left(\left|\frac{1}{n}\sum_{i=1}^{n}\left|\Delta^{b}(X_{i})\right|^2 K_{h}(D_{i})\right|\leq C t^2{\rm Err}^2(\widehat{g})\right) \geq 1-\frac{1}{t^2}-\min\{n,p\}^{-c}.
\end{equation*}
By symmetry and the decomposition \eqref{eq: swapping decomp2}, we have 
\begin{equation}
\mathbf{P}\left(\left|\frac{1}{n}\sum_{i=1}^{n}\dif^2 K_{h}(D_{i})\right|\leq C t^2{\rm Err}^2(\widehat{g})\right) \geq 1-\frac{1}{t^2}-\min\{n,p\}^{-c}.
\label{eq: key inter}
\end{equation}
Combined with \eqref{eq: cau 1}, we obtain
{\small
\begin{equation}\mathbf{P}\left(\left|\frac{1}{n}\sum_{i=1}^{n}\left(\widehat{W}_{i}-\left({W}_{i}-\bar{\mu}_{W}\right)\right)\dif K_{h}(D_{i})\right| \lesssim t^2  {\rm Err}(\widehat{W})\cdot{\rm Err}(\widehat{g})\right)\geq 1-\frac{1}{t^2}-\min\{n,p\}^{-c}.
\label{eq: temp bound inter}
\end{equation}
}
Note that 
\begin{equation*}
\frac{1}{n}\sum_{i=1}^{n}\left|\dif \right|K_{h}(D_{i}) \leq \sqrt{\frac{1}{n}\sum_{i=1}^{n}K_{h}(D_{i})}\cdot\sqrt{\frac{1}{n}\sum_{i=1}^{n}\dif^2 K_{h}(D_{i})}.
\label{eq: cau 2}
\end{equation*}
Together with \eqref{eq: prob bound 1}, \eqref{eq: centering rate}, and \eqref{eq: key inter}, we establish
$$
\mathbf{P}\left(\left|\bar{\mu}_{W}\cdot \frac{1}{n}\sum_{i=1}^{n}\dif K_{h}(D_{i})\right|\lesssim t{\rm Err}(\widehat{g})\sqrt{{h}/{n}}\right)\geq 1-\frac{1}{t^2}-\min\{n,p\}^{-c}.
$$
\noindent Together with \eqref{eq: bias b expectation a}, \eqref{eq: temp bound inter} and the decomposition \eqref{eq: bias b decomposition}, we have 
{\small$$
\mathbf{P}\left(\left|\frac{1}{n}\sum_{i=1}^{n}\widehat{W}_{i}\dif K_{h}(D_{i})\right|\lesssim t \left(t\sqrt{{h}/{n}}+{\rm Err}(\widehat{W})\right){\rm Err}(\widehat{g})\right)\geq 1-\frac{1}{t^2}-\min\{n,p\}^{-c}.
$$}
Together with \eqref{eq: second limit}, we establish \eqref{eq: bias b ratio}.

%%%%%%%%%%%%%%%%%%%%%%%%%%%%%%%%%%%%%%%%
\subsection{Proof of \eqref{eq: weight error}}
\label{sec: result 3}
%%%%%%%%%%%%%%%%%%%%%%%%%%%%%%%%%%%%%%%%%%%%%%%%%%%%%%%%
%\Zijian{Stops here.}
In the following, we first establish 
\begin{equation}
{\rm Err}^2(\widehat{W})\lesssim {\frac{1}{n}\sum_{i=1}^{n}\left(l(X_{i},\widehat{\gamma})-l(X_{i},{\gamma})\right)^2 K_{h}(D_{i})}.
\label{eq: simplification}
\end{equation} 
Recall that the uncentered weight $\widetilde{W}_i$ is defined in \eqref{eq: swapping est} and $\widehat{W}_i$ is the corresponding centered weight defined in \eqref{eq: final D}.
Note that 
\begin{equation}
\left(\widehat{W}_{i}-\left({W}_{i}-\bar{\mu}_{W}\right)\right)^2\lesssim \left(\widetilde{W}_{i}-{W}_{i}\right)^2+\left(\frac{\frac{1}{n}\sum_{i=1}^{n}(\widetilde{W}_{i}-{W}_{i})K_{h}(D_{i})}{\frac{1}{n}\sum_{i=1}^{n}K_{h}(D_{i})}\right)^2.
\label{eq: basic upper}
\end{equation}
By Cauchy-Schwarz inequality $$\left(\frac{1}{n}\sum_{i=1}^{n}(\widetilde{W}_{i}-{W}_{i})K_{h}(D_{i})\right)^2\leq \left(\frac{1}{n}\sum_{i=1}^{n}(\widetilde{W}_{i}-{W}_{i})^2K_{h}(D_{i})\right)\cdot\left(\frac{1}{n}\sum_{i=1}^{n}K_{h}(D_{i})\right),$$ 
we have $$\frac{1}{n}\sum_{i=1}^{n}\left(\frac{\frac{1}{n}\sum_{i=1}^{n}(\widetilde{W}_{i}-{W}_{i})K_{h}(D_{i})}{\frac{1}{n}\sum_{i=1}^{n}K_{h}(D_{i})}\right)^2K_{h}(D_{i})\leq \frac{1}{n}\sum_{i=1}^{n}(\widetilde{W}_{i}-{W}_{i})^2K_{h}(D_{i}).$$
By applying the above inequality and \eqref{eq: basic upper}, we obtain 
\begin{equation*}
\begin{aligned}
&\sqrt{\frac{1}{n}\sum_{i=1}^{n}\left(\widehat{W}_{i}-\left({W}_{i}-\bar{\mu}_{W}\right)\right)^2 K_{h}(D_{i})}\lesssim \sqrt{\frac{1}{n}\sum_{i=1}^{n}\left(\widetilde{W}_{i}-{W}_{i}\right)^2 K_{h}(D_{i})}\\
&+\sqrt{\frac{1}{n}\sum_{i=1}^{n}\left(\frac{\frac{1}{n}\sum_{i=1}^{n}(\widetilde{W}_{i}-{W}_{i})K_{h}(D_{i})}{\frac{1}{n}\sum_{i=1}^{n}K_{h}(D_{i})}\right)^2K_{h}(D_{i})}\lesssim \sqrt{\frac{1}{n}\sum_{i=1}^{n}\left(\widetilde{W}_{i}-{W}_{i}\right)^2 K_{h}(D_{i})}.
\end{aligned}
\end{equation*}
By the definitions ${W}_{i}=(D_{i}-\eval)-l(X_{i},\gamma)$ and $\widetilde{W}_{i}=(D_{i}-\eval)-\widehat{l}(X_{i},\widehat{\gamma}),$ the above inequality implies \eqref{eq: simplification}.

Then the proof of \eqref{eq: weight error} is reduced to establishing an upper bound for 
$$ {\frac{2}{n}\sum_{i=1}^{n}\left(\widehat{l}(X_{i},\widehat{\gamma})-l(X_{i},{\gamma})\right)^2 K_{h}(D_{i})}.
$$
We divide the above summation into two parts,
\begin{equation*}
 {\frac{1}{n_a}\sum_{i\in \mathcal{I}_{a}}\left(\widehat{l}(X_{i},\widehat{\gamma}^{b})-l(X_{i},{\gamma})\right)^2 K_{h}(D_{i})}+ {\frac{1}{n_b}\sum_{i\in \mathcal{I}_{b}}\left(\widehat{l}(X_{i},\widehat{\gamma}^{a})-l(X_{i},{\gamma})\right)^2 K_{h}(D_{i})}.
\end{equation*}
By symmetry, we focus on the first summation 
\begin{equation}
{\frac{1}{n_a}\sum_{i\in \mathcal{I}_{a}}\left(\widehat{l}(X_{i},\widehat{\gamma}^{b})-l(X_{i},{\gamma})\right)^2 K_{h}(D_{i})},
\label{eq: main goal}
\end{equation}
 and adopt the notation %s $\mathcal{I}_{a}=\{1,2,\cdots, n_a\}$ and 
 $\widehat{\gamma}^{b}=\widehat{\gamma}.$ %without loss of generality.
We note that 
\begin{equation*}
\widehat{\delta}_j-\widehat{\mu}_i=\delta_j-\mu_i- a_{ij}, \quad \text{with}\quad a_{ij}=(X_{j}-X_{i})^{\intercal}(\widehat{\gamma}-\gamma).
\end{equation*}
On the event $\mathcal{A}_1$ defined in \eqref{eq: events}, we have 
\begin{equation}
|a_{ij}|\leq C^{*}\sqrt{\frac{\|\gamma\|_0 \log p}{n}}\sqrt{\log n}.
\label{eq: inter upper bound}
\end{equation}
%The randomness in ${\frac{1}{n_a}\sum_{i\in \mathcal{I}_{a}}\left(\widehat{l}(X_{i},\widehat{\gamma}^{b})-l(X_{i},{\gamma})\right)^2 K_{h}(D_{i})}$ comes from three sources,

%\begin{equation}
% \mathcal{G}_1=\left\{\max \right\}
% \end{equation}
%and 
%\begin{equation}
%\mathcal{G}_2=\left\{\max_{1\leq i\leq n}\max_{\left|\delta-\mu_i\right|\leq r} \frac{|\phi''(\delta)|}{\phi(\mu_i)}\leq C_1(n) \quad \max_{1\leq i\leq n}\max_{\left|\delta-\mu_i\right|\leq r} \frac{|\phi'(\delta)|}{\phi(\mu_i)}\leq C_2(n)\right\}
%\end{equation}
To facilitate the discussion, we introduce the following notations,
\begin{equation}
\widehat{\I}_{ij}={\bf{1}}(|\delta_j-\mu_i- a_{ij}|\leq h) \quad \text{and} \quad \I_{ij}={\bf{1}}(|\delta_j-\mu_i|\leq h). 
\label{eq: def 1}
\end{equation}
The estimator $\widehat{l}(X_{i},\widehat{\gamma})$ defined in \eqref{eq: estimator sym} can be written as 
\begin{equation}
\widehat{l}(X_{i},\widehat{\gamma})=\frac{\aveg (\delta_j-\mu_i- a_{ij}) \widehat{\I}_{ij}}{\aveg \widehat{\I}_{ij}}.
\label{eq: key simplification}
\end{equation}
The randomness of $\widehat{l}(X_{i},\widehat{\gamma})$ comes from the following three parts,
\begin{itemize} 
\item the noise $\{\delta_j\}_{j\in \mathcal{I}_a}$ 
\item the variable $\mu_i$ with $\mu_i=\eval-X^{\intercal}_{i}\gamma$ for the pre-fixed index $i$
\item the estimation error $a_{ij}$, which depends on $X_i, X_j$ and the initial estimator $\widehat{\gamma}$ computed on the data $\mathcal{I}_b$.  
\end{itemize}
Note that $\{\delta_j\}_{j\in \mathcal{I}_a}$ is independent of the other two random sources. We shall write $\E_{\delta_j}$ as the expectation with respect to $\delta_j$ but condition on the other two components $X_i$ and $\mathcal{I}_b$. We use $\E_{\delta_j}\left[\cdot \mid\widehat{\I}_{ij} \right]$ to denote the conditional expectation by only considering the randomness of $\delta_j$. Specifically, for $j\in \mathcal{I}_a$, $\E_{\delta_j}$ and $\E_{\delta_j}\left[\cdot \mid\widehat{\I}_{ij} \right]$ are shorthanded for 
$$\E_{\delta_j}[\cdot]=\E[\cdot\mid \{X_i\}_{i \in \mathcal{I}_a},\mathcal{I}_b]\quad \text{and}\quad \E_{\delta_j}\left[\cdot \mid\widehat{\I}_{ij} \right]=\E\left[\cdot \mid\widehat{\I}_{ij},\{X_i\}_{i \in \mathcal{I}_a},\mathcal{I}_b \right].$$
We use $\E_{\delta\mid \widehat{\I}_{i\cdot}}$ to denote the conditional expectation of $\{\delta_j\}_{j\in \mathcal{I}_{a}}$ given the events $\widehat{\I}_{i\cdot}=\{\widehat{\I}_{ij}\}_{j\in \mathcal{I}_{a}},$ by only considering the randomness of $\{\delta_j\}_{j\in \mathcal{I}_{a}},$ that is
$$\E_{\delta\mid \widehat{\I}_{i\cdot}}[\cdot]=\E\left[\cdot \mid \widehat{\I}_{i\cdot},\{X_i\}_{i \in \mathcal{I}_a},\mathcal{I}_b \right].$$

We compute the following difference,
\begin{equation}
\begin{aligned}
&\frac{\aveg (\delta_j-\mu_i- a_{ij}) \widehat{\I}_{ij}}{\aveg \widehat{\I}_{ij}}-\frac{\aveg \E_{\delta_j}\left[(\delta_j-\mu_i- a_{ij}) \widehat{\I}_{ij}\mid\widehat{\I}_{ij} \right]}{\aveg \widehat{\I}_{ij} }\\
&=\frac{\aveg \widehat{\I}_{ij}\left((\delta_j-\mu_i- a_{ij}) -\E_{\delta_j}\left(\delta_j-\mu_i- a_{ij} \mid\widehat{\I}_{ij} \right)\right)}{\aveg \widehat{\I}_{ij} }.
\end{aligned}
\label{eq: difference 1}
\end{equation}
Define 
\begin{equation}
F_{ij}(t)=\frac{ \int_{\mu_i+ta_{ij}-h}^{\mu_i+ta_{ij}+h}(\delta-\mu_i-ta_{ij}) \phi(\delta)d \delta}{ \int_{\mu_i+ta_{ij}-h}^{\mu_i+ta_{ij}+h}\phi(\delta)d \delta}.
\label{eq: key def inter}
\end{equation}
Note that 
\begin{equation}
\E_{\delta_j}\left(\delta_j-\mu_i- a_{ij} \mid\widehat{\I}_{ij} \right)=\frac{ \int_{\mu_i+a_{ij}-h}^{\mu_i+a_{ij}+h}(\delta-\mu_i-a_{ij}) \phi(\delta)d \delta}{ \int_{\mu_i+a_{ij}-h}^{\mu_i+a_{ij}+h}\phi(\delta)d \delta}=F_{ij}(1).
\label{eq: key expression}
\end{equation}
By \eqref{eq: key simplification}, \eqref{eq: difference 1}, and \eqref{eq: key expression}, we establish
 \begin{equation}
 \begin{aligned}
 &\widehat{l}(X_{i},\widehat{\gamma})-{l}(X_{i},{\gamma})=\frac{\aveg (\delta_j-\mu_i- a_{ij}) \widehat{\I}_{ij}}{\aveg \widehat{\I}_{ij}}-\frac{ \int_{\mu_i-h}^{\mu_i+h} \left(\delta-\mu_i\right)\phi(\delta) d \delta}{ \int_{\mu_i-h}^{\mu_i+h}\phi(\delta) d\delta} \\
 &=\frac{\aveg \widehat{\I}_{ij}\left[(\delta_j-\mu_i- a_{ij}) -\E_{\delta_j}\left(\delta_j-\mu_i- a_{ij} \mid\widehat{\I}_{ij} \right)\right]}{\aveg \widehat{\I}_{ij} }+\frac{\aveg \widehat{\I}_{ij}\left[F_{ij}(1)-F_{ij}(0)\right]}{\aveg \widehat{\I}_{ij} },
 \end{aligned}
 \label{eq: key decomp inter}
 \end{equation}
where the last component holds since $F_{ij}(0)$ defined in \eqref{eq: key def inter} does not depend on the index $j$. 

The decomposition \eqref{eq: key decomp inter} and the following two inequalities lead to an upper bound for \eqref{eq: main goal}. 
{\small
\begin{equation}
\frac{1}{n_a}\sum_{i\in \mathcal{I}_a} \E\left[\left(\frac{\aveg \widehat{\I}_{ij}\left[(\delta_j-\mu_i- a_{ij}) -\E_{\delta_j}\left(\delta_j-\mu_i- a_{ij} \mid\widehat{\I}_{ij} \right)\right]}{\aveg \widehat{\I}_{ij} }\right)^2 \frac{1}{2h}
{\bf{1}}(|\delta_i-\mu_i|\leq h)\right]\lesssim \frac{h}{n_a}, 
\label{eq: goal 1}
\end{equation}}
and 
\begin{equation} 
\begin{aligned}
&\frac{1}{n_a}\sum_{i\in \mathcal{I}_a} \E \left[\left(\frac{\aveg \widehat{\I}_{ij}\left[F_{ij}(1)-F_{ij}(0)\right]}{\aveg \widehat{\I}_{ij} }\right)^2 \frac{1}{2h}
{\bf{1}}(|\delta_i-\mu_i|\leq h)\cdot {\1}_{\mathcal{A}_0\cap\mathcal{A}_1}\right]\\
&\lesssim h^4 \left(C_1^2(n)+C_2(n)\right)^2 \frac{k\log p\log n}{n}.
\end{aligned}
\label{eq: goal 2}
\end{equation}
We establish \eqref{eq: weight error} by combining \eqref{eq: simplification}, \eqref{eq: key decomp inter}, \eqref{eq: goal 1}, \eqref{eq: goal 2}, and \eqref{eq: high prob}.
%%%%%%%%%%%%%%%%%%%%%%%%%%%%%
\subsubsection{Proof of \eqref{eq: goal 1}}
%%%%%%%%%%%%%%%%%%%%%%%%%%%%%
%For a given $i\in \mathcal{I}_a$, we have 
%\begin{equation*}
%\begin{aligned}
%&\left(\frac{\aveg \widehat{\I}_{ij}\left[(\delta_j-\mu_i- a_{ij}) -\E_{\delta_j}\left(\delta_j-\mu_i- a_{ij} \mid\widehat{\I}_{ij} \right)\right]}{\aveg \widehat{\I}_{ij} }\right)^2 
%\frac{1}{2h}{\bf{1}}(|\delta_i-\mu_i|\leq h) \\
%&=\left(\frac{\aveg \widehat{\I}_{ij}\left[(\delta_j-\mu_i- a_{ij}) -\E_{\delta_j}\left(\delta_j-\mu_i- a_{ij} \mid\widehat{\I}_{ij} \right)\right]}{\frac{1}{n_a}+\frac{1}{n_a}\sum_{j \neq i} \widehat{\I}_{ij} }\right)^2 
%\frac{1}{2h}{\I}_{ii}.\\
%%&=\left(\frac{\frac{1}{n_a}\left[(\delta_i-\mu_i) -\E_{\delta_j}\left(\delta_j-\mu_i \mid\widehat{\I}_{ij} \right)\right]+\frac{1}{n_a}\sum_{j \neq i} \widehat{\I}_{ij}\left[(\delta_j-\mu_i- a_{ij}) -\E_{\delta_j}\left(\delta_j-\mu_i- a_{ij} \mid\widehat{\I}_{ij} \right)\right]}{\frac{1}{n_a}+\frac{1}{n_a}\sum_{j \neq i} \widehat{\I}_{ij} }\right)^2 
%%{\I}_{ii}\\
%%&\lesssim \left(\frac{\frac{1}{n_a}\sum_{j \neq i} \widehat{\I}_{ij}\left[(\delta_j-\mu_i- a_{ij}) -\E_{\delta_j}\left(\delta_j-\mu_i- a_{ij} \mid\widehat{\I}_{ij} \right)\right]}{\frac{1}{n_a}+\frac{1}{n_a}\sum_{j \neq i} \widehat{\I}_{ij} }\right)^2  {\I}_{ii}\\
%%&+\left(\frac{\frac{1}{n_a}\left[(\delta_i-\mu_i) -\E_{\delta_j}\left(\delta_j-\mu_i \mid\widehat{\I}_{ij} \right)\right]}{\frac{1}{n_a}+\frac{1}{n_a}\sum_{j \neq i} \widehat{\I}_{ij} }\right)^2{\I}_{ii}
%\end{aligned}
%\end{equation*}
In the following proof, we fix $i\in \mathcal{I}_a$ and ${\I}_{ii}=\widehat{\I}_{ii}.$ For $j_1,j_2\in \mathcal{I}_a$ with $j_1\neq j_2,$ 
{\small $$\E_{\delta\mid \widehat{\I}_{i\cdot}}\left[(\delta_{j_1}-\mu_i- a_{i{j_1}}) -\E_{\delta_{j_1}}\left(\delta_{j_1}-\mu_i- a_{i{j_1}} \mid\widehat{\I}_{i{j_1}} \right)\right] \left[(\delta_{j_2}-\mu_i- a_{i{j_2}}) -\E_{\delta_{j_2}}\left(\delta_{j_2}-\mu_i- a_{i{j_2}} \mid\widehat{\I}_{i{j_2}} \right)\right]=0,$$}
and 
$$\E_{\delta_j\mid \widehat{\I}_{ij}}\left[(\delta_j-\mu_i- a_{ij}) -\E_{\delta_j}\left(\delta_j-\mu_i- a_{ij} \mid\widehat{\I}_{ij} \right)\right]^2\leq \E_{\delta_j\mid \widehat{\I}_{ij}}\left[(\delta_j-\mu_i- a_{ij})\right]^2\leq h^2.$$
By the above two expressions, we obtain 
\begin{equation}
\begin{aligned}
&\E_{\delta\mid \widehat{\I}_{i\cdot}}\left[\left(\frac{\aveg \widehat{\I}_{ij}\left[(\delta_j-\mu_i- a_{ij}) -\E_{\delta_j}\left(\delta_j-\mu_i- a_{ij} \mid\widehat{\I}_{ij} \right)\right]}{\frac{1}{n_a}+\frac{1}{n_a}\sum_{j \neq i} \widehat{\I}_{ij} }\right)^2 \frac{1}{2h}{\I}_{ii}\right]\\
&=\frac{\frac{1}{n_a^2}\sum_{j\in \mathcal{I}_a} \widehat{\I}_{ij}\E_{\delta_j\mid \widehat{\I}_{ij}}\left[(\delta_j-\mu_i- a_{ij}) -\E_{\delta_j}\left(\delta_j-\mu_i- a_{ij} \mid\widehat{\I}_{ij} \right)\right]^2}{\left(\frac{1}{n_a}+\frac{1}{n_a}\sum_{j \neq i} \widehat{\I}_{ij}\right)^2} \frac{1}{2h}{\I}_{ii}\\
&\leq \frac{h}{2n_a} \frac{ {\I}_{ii}}{\frac{1}{n_a}+\frac{1}{n_a}\sum_{j \neq i} \widehat{\I}_{ij}},
\end{aligned}
\label{eq: step 1}
\end{equation}
%where the first equality holds since 
We further take expectation with respect to $\I$ (but conditioning on $\mu_i$ and $a_{ij}$) in \eqref{eq: step 1} and obtain that 
\begin{equation}
\E_{\delta}\left(\frac{h}{n_a} \frac{ {\I}_{ii}}{\frac{1}{n_a}+\frac{1}{n_a}\sum_{j \neq i} \widehat{\I}_{ij}}\right)=\frac{h}{n_a} \cdot (\E_{\delta} {\I}_{ii}) \cdot \left(\E_{\delta}\frac{ 1}{\frac{1}{n_a}+\frac{1}{n_a}\sum_{j \neq i} \widehat{\I}_{ij}}\right).
\label{eq: step 2}
\end{equation}
Conditioning on $\mu_{i}$ and $a_{ij},$
we define the conditional probability with respect to $\delta_j$ as 
\begin{equation}
e_{ij}= \E_{\delta_j}[\widehat{\I}_{ij}]=  \int_{\mu_i+a_{ij}-h}^{\mu_i+a_{ij}+h} \phi(\delta_j) d\delta_j.
\label{eq: prob expression}
\end{equation}
  By change of variable $\tau_j=\delta_j-(\mu_i+a_{ij})$, we have 
 \begin{equation*}
  \begin{aligned} 
  \int_{\mu_i+a_{ij}-h}^{\mu_i+a_{ij}+h} \phi(\delta_j) d\delta_j&=\int_{-h}^{h}\phi(\tau_{j}+\mu_i+a_{ij}) d\tau_j\\
&=\int_{-h}^{h}[\phi(\mu_i+a_{ij})+\tau_{j}\phi'(\mu_i+a_{ij})+\frac{\tau_j^2}{2}\phi''(\mu_i+a_{ij}+c \tau_j)] d\tau_j,\\
%&=2h\phi(\mu_i+a_{ij})+\\
  \end{aligned}
 \end{equation*}
 for some constant $c\in (0,1).$
 Hence, we have 
 \begin{equation*}
\left|\int_{\mu_i+a_{ij}-h}^{\mu_i+a_{ij}+h} \phi(\delta_j) d\delta_j-2h\phi(\mu_i+a_{ij})\right|\leq \frac{h^3}{3} \max_{|c|\leq 1} \left|\phi''(\mu_i+a_{ij}+ch)\right|,
 \end{equation*}
 and then 
 \begin{equation*}
 \left|   e_{ij}- 2h\phi(\mu_i)\right|\leq  \frac{h^3}{3} \max_{|c|\leq 1} \left|\phi''(\mu_i+a_{ij}+ch)\right|+ 2 h |a_{ij}| \max_{|c|\leq 1} \left|\phi'(\mu_i+c a_{ij})\right|.
 \end{equation*}
 Hence, we establish  
\begin{equation}
\left| \frac{e_{ij}}{2h\phi(\mu_i)}-1\right|\leq \frac{h^2}{3} C_2(n)+ |a_{ij}| C_1(n).
\label{eq: sample version}
\end{equation}
We note that the non-negative random variable $\frac{1}{n_a}\sum_{j \neq i} \widehat{\I}_{ij}$ satisfies 
$$\E\left[\frac{1}{n_a}\sum_{j \neq i} \widehat{\I}_{ij}\mid\{X_i\}_{i \in \mathcal{I}_a},\mathcal{I}_b \right]=\frac{1}{n_a}\sum_{j \neq i} e_{ij},$$ and 
\begin{equation*}
{\rm Var}\left[\frac{1}{n_a}\sum_{j \neq i}  \widehat{\I}_{ij}\mid\{X_i\}_{i \in \mathcal{I}_a},\mathcal{I}_b \right] = \frac{1}{n_a^2}\sum_{j \neq i} e_{ij}(1-e_{ij}).
\end{equation*}
Hence, we apply the equation (5) in \cite{wooff1985bounds} and obtain
\begin{equation}
\begin{aligned}
\E_{\delta}\left[\frac{1}{\frac{1}{n_a}+\frac{1}{n_a}\sum_{j \neq i} \widehat{\I}_{ij}}\right]&=\E\left[\frac{1}{\frac{1}{n_a}+\frac{1}{n_a}\sum_{j \neq i} \widehat{\I}_{ij}}\mid\{X_i\}_{i \in \mathcal{I}_a},\mathcal{I}_b \right]\\&\leq \frac{1}{\frac{1}{n_a}\sum_{j \neq i} e_{ij}+\frac{1}{n_a}}\left(1+\frac{\frac{1}{n_a^2}\sum_{j \neq i} e_{ij}(1-e_{ij})}{\frac{1}{n_a^2}\sum_{j \neq i} e_{ij}}\right)\\
&\leq \frac{2}{\frac{1}{n_a}\sum_{j \neq i} e_{ij}+\frac{1}{n_a}}.\\
\end{aligned}
\label{eq: inverse bound}
\end{equation}
A combination of \eqref{eq: step 2}, \eqref{eq: sample version} and \eqref{eq: inverse bound} leads to 
\begin{equation}
\E_{\delta}\left(\frac{h}{n_a} \frac{ {\I}_{ii}}{\frac{1}{n_a}+\frac{1}{n_a}\sum_{j \neq i} \widehat{\I}_{ij}}\right)\leq \frac{h}{n_a} \frac{2 e_{ii}}{\frac{1}{n_a}\sum_{j \neq i} e_{ij}+\frac{1}{n_a}} \lesssim \frac{h}{n_a}. 
\label{eq: goal 1 achieved}
\end{equation}
By taking expectation with respect to $\{X_i\}_{i \in \mathcal{I}_a}$ and $\mathcal{I}_b,$ we establish \eqref{eq: goal 1}.
\subsubsection{Proof of \eqref{eq: goal 2}}
%\begin{equation} \frac{1}{n_a}\sum_{i=1}^{n_a} \E \left(\frac{\aveg \widehat{\I}_{ij}\left[F_{ij}(1)-F_{ij}(0)\right]}{\aveg \widehat{\I}_{ij} }\right)^2 \frac{1}{2h}
%{\bf{1}}(|\delta_i-\mu_i|\leq h) 
%\label{eq: goal 2}
%\end{equation}
%Define 
%\begin{equation}
%F_{ij}(t)=\frac{ \int_{\mu_i+ta_{ij}-h}^{\mu_i+ta_{ij}+h}(\delta-\mu_i-ta_{ij}) \phi(\delta)d \delta}{ \int_{\mu_i+ta_{ij}-h}^{\mu_i+ta_{ij}+h}\phi(\delta)d \delta}.
%\label{eq: key def inter}
%\end{equation}
%Define 
%\begin{equation}
%F_i(t)=\frac{\aveg \widehat{\I}_{ij}F_{ij}(t)}{\aveg \widehat{\I}_{ij} }
%\end{equation}
%and hence \eqref{eq: goal 2} can be reduced to 
%\begin{equation}
%\frac{1}{n}\sum_{i=1}^{n} \E (F_i(1)-F_i(0))^2\frac{1}{h}
%{\bf{1}}(|\delta_i-\mu_i|\leq h) 
%\end{equation}
%by calculating the difference between $F_i(1)-F_{i}(0)$. 
We calculate the expression of $F_{ij}(t)$ in \eqref{eq: key def inter} as 
%\begin{equation}
%F_i(t)=\frac{\aveg  \int_{\mu_i+ta_{ij}-h}^{\mu_i+ta_{ij}+h}(\delta_j-\mu_i-ta_{ij}) \phi(\delta_j)d \delta_j}{\aveg  \int_{\mu_i+ta_{ij}-h}^{\mu_i+ta_{ij}+h}\phi(\delta_j)d \delta_j}
%\end{equation}
%and shall compute the derivative of $F_i(t)$ with respect to $t$. 
\begin{equation}
\begin{aligned}
\frac{dF_{ij}(t)}{d t}&=\frac{  a_{ij} h \left[\phi(\mu_i+ta_{ij}+h)+\phi(\mu_i+ta_{ij}-h)-\frac{1}{h}\int_{\mu_i+ta_{ij}-h}^{\mu_i+ta_{ij}+h} \phi(\delta_j)d \delta_j\right] }{ \int_{\mu_i+ta_{ij}-h}^{\mu_i+ta_{ij}+h}\phi(\delta_j)d \delta_j}\\
&-\frac{\int_{\mu_i+ta_{ij}-h}^{\mu_i+ta_{ij}+h}(\delta_j-\mu_i-ta_{ij}) \phi(\delta_j)d \delta_j}{\left[ \int_{\mu_i+ta_{ij}-h}^{\mu_i+ta_{ij}+h}\phi(\delta_j)d \delta_j\right]^2}\cdot  a_{ij} \left[\phi(\mu_i+ta_{ij}+h)-\phi(\mu_i+ta_{ij}-h)\right].
\end{aligned}
\label{eq: derivative 1}
\end{equation}
Define 
\begin{equation*}
T_1(a_{ij})=\phi(\mu_i+ta_{ij}+h)+\phi(\mu_i+ta_{ij}-h)-\frac{1}{h}\int_{\mu_i+ta_{ij}-h}^{\mu_i+ta_{ij}+h} \phi(\delta_j)d \delta_j;
\end{equation*}
\begin{equation*}
T_2(a_{ij})= \phi(\mu_i+ta_{ij}+h)-\phi(\mu_i+ta_{ij}-h);
\end{equation*}
and
\begin{equation*}
T_3(a_{ij})= \int_{\mu_i+ta_{ij}-h}^{\mu_i+ta_{ij}+h}[\delta_j-(\mu_i+t a_{ij})]\phi(\delta_j)d \delta_j.
\end{equation*}
Then we can simplify the derivative of $F_{ij}(t)$ in \eqref{eq: derivative 1} as
\begin{equation}
\frac{dF_{ij}(t)}{d t}=\frac{h a_{ij} T_1(a_{ij})}{e_{ij}}-\frac{a_{ij} T_2(a_{ij})}{e_{ij}}\cdot\frac{ T_3(a_{ij})}{e_{ij}},
\label{eq: diff expression}
\end{equation}
where $e_{ij}$ is defined in \eqref{eq: prob expression}.
To bound the above terms, we introduce the following lemma to control all of the above terms.
\begin{Lemma}
Suppose that $\phi(t)$ is twice differentiable for $t\in [\mu-\tau,\mu+\tau],$ there exists some positive constant $C>0$ such that
\begin{equation}
\begin{aligned}
\left|\phi(\mu+\tau)+\phi(\mu-\tau)-\frac{1}{\tau}\int_{\mu-\tau}^{\mu+\tau}\phi(t) dt\right|\leq C\tau^2\cdot \max_{t\in [\mu-\tau,\mu+\tau]} \left|\phi''(t)\right|
%&\leq C\tau^2\cdot \max_{t\in [\mu-\tau,\mu+\tau]} \left|\rho_2(t)\right| \cdot \max_{t\in [\mu-\tau,\mu+\tau]} \left|\phi(t)\right|
\end{aligned}
\label{eq: Gaussian bound 1}
\end{equation}
\begin{equation}
\left|\int_{\mu-\tau}^{\mu+\tau} \left(t-\mu\right)\phi(t) dt\right|\leq C\tau^3\max_{t\in [\mu-\tau,\mu+\tau]} \left|\phi'(t)\right|%\leq C\tau^3\cdot \max_{t\in [\mu-\tau,\mu+\tau]} \left|\rho_1(t)\right| \cdot \max_{t\in [\mu-\tau,\mu+\tau]} \left|\phi(t)\right|
\label{eq: Gaussian bound 2}
\end{equation}
\begin{equation}
\left|\phi(\mu+\tau)-\phi(\mu-\tau)\right|\leq C \tau \max_{t\in [\mu-\tau,\mu+\tau]} \left|\phi'(t)\right|%\leq C\tau \cdot \max_{t\in [\mu-\tau,\mu+\tau]} \left|\rho_1(t)\right| \cdot \max_{t\in [\mu-\tau,\mu+\tau]} \left|\phi(t)\right|
\label{eq: Gaussian bound 3}
\end{equation}
%where 
%$\rho_1(t)=\frac{\phi'(t)}{\phi(t)}$ and $\rho_2(t)=\frac{\phi''(t)}{\phi(t)}.$
\label{lem: Gaussian bound}
\end{Lemma}

It follows from \eqref{eq: Gaussian bound 1} that 
\begin{equation*}
\left|T_1(a_{ij})\right| \lesssim h^2 \cdot \max_{\left|\delta-(\mu_i+ta_{ij})\right|\leq h} |\phi''(\delta)|\leq h^2 \cdot \max_{\left|\delta-\mu_i\right|\leq r} |\phi''(\delta)|.
\end{equation*}
It follows from \eqref{eq: Gaussian bound 2} that 
\begin{equation*}
\left|T_2(a_{ij})\right| \lesssim h \cdot \max_{\left|\delta-(\mu_i+ta_{ij})\right|\leq h} |\phi'(\delta)|\leq h \cdot \max_{\left|\delta-\mu_i\right|\leq r} |\phi'(\delta)|.
\end{equation*}
It follows from \eqref{eq: Gaussian bound 3} that 
\begin{equation*}
\left|T_3(a_{ij})\right| \lesssim h^3 \cdot \max_{\left|\delta-(\mu_i+ta_{ij})\right|\leq h} |\phi'(\delta)|\leq h^3 \cdot \max_{\left|\delta-\mu_i\right|\leq r} |\phi'(\delta)|.
\end{equation*}
Together with \eqref{eq: sample version}, we establish 
\begin{equation*}
\frac{\left|T_1(a_{ij})\right|}{e_{ij}} \lesssim h \cdot  \max_{\left|\delta-\mu_i\right|\leq r} \frac{|\phi''(\delta)|}{\phi(\mu_i)},\quad \frac{\left|T_2(a_{ij})\right|}{e_{ij}} \lesssim \max_{\left|\delta-\mu_i\right|\leq r} \frac{|\phi'(\delta)|}{\phi(\mu_i)},\quad \frac{\left|T_3(a_{ij})\right|}{e_{ij}} \lesssim h^2 \cdot  \max_{\left|\delta-\mu_i\right|\leq r} \frac{|\phi'(\delta)|}{\phi(\mu_i)},
\end{equation*}
where $r=C^*\sqrt{{\|\gamma\|_0 \log p \log n}/{n}}+h.$
Together with the expression \eqref{eq: diff expression} and the upper bound \eqref{eq: inter upper bound}, we establish 
\begin{equation*}
\begin{aligned}
\left|\frac{d F_{ij}(t)}{d t}\right|{\1}_{\mathcal{A}_0\cap\mathcal{A}_1}\lesssim \sqrt{\frac{k \log p\log n}{n}} h^2 \left(C^2_1(n)+C_2(n)\right).
\end{aligned}
\end{equation*}
Hence, we have  
\begin{equation*}
\left(\frac{\aveg \widehat{\I}_{ij}\left[F_{ij}(1)-F_{ij}(0)\right]}{\aveg \widehat{\I}_{ij} }\right)^2 {\1}_{\mathcal{A}_0\cap\mathcal{A}_1}\lesssim \left(\sqrt{\frac{k \log p\log n}{n}} h^2 \left(C_1^2(n)+C_2(n)\right)\right)^2,
\end{equation*}
and 
\begin{equation} 
\begin{aligned}
&\frac{1}{n_a}\sum_{i=1}^{n_a} \E \left(\frac{\aveg \widehat{\I}_{ij}\left[F_{ij}(1)-F_{ij}(0)\right]}{\aveg \widehat{\I}_{ij} }\right)^2 \frac{1}{2h}
{\bf{1}}(|\delta_i-\mu_i|\leq h) \cdot {\1}_{\mathcal{A}_0\cap\mathcal{A}_1} \\
&\lesssim \left(\sqrt{\frac{k \log p\log n}{n}} h^2 \left(C_1^2(n)+C_2(n)\right)\right)^2 \frac{1}{n_a}\sum_{i=1}^{n_a}\E \frac{1}{2h}
{\bf{1}}(|\delta_i-\mu_i|\leq h)\\
&=\left(\sqrt{\frac{k \log p\log n}{n}} h^2 \left(C_1^2(n)+C_2(n)\right)\right)^2 \frac{1}{n_a}\sum_{i=1}^{n_a}\E K_{h}(D_i).
\end{aligned}
\label{eq: goal 2 achieved}
\end{equation}
Together with \eqref{eq: expectation 0}, we establish \eqref{eq: goal 2}.

\section{Proof of Theorem \ref{thm: limiting swap}}

%We recall
%\begin{equation}
%C_{u}=\frac{1}{\sigma_{\delta}}\left(\eval+h+\|\gamma\|_2 \sqrt{\log n}\right).
%\label{eq: second key growing constant}
%\end{equation}

%%%%%%%%%%%%%%%%%%%%%%%%%%%%
%\subsection{Proof of Theorem \ref{thm: limiting swap}}
%%%%%%%%%%%%%%%%%%%%%%%%%%%%
We start with the following error decomposition of
 $\widehat{f'(\eval)}-f'(\eval),$
{
\begin{equation}
\underbrace{\frac{1}{n\widehat{S}_{n}}\sum_{i=1}^{n}\widehat{W}_{i}\epsilon_iK_{h}(D_i)}_{\rm Stochastic\; Error}+\underbrace{\frac{1}{n\widehat{S}_{n}}\sum_{i=1}^{n}\widehat{W}_{i}r(D_i)K_{h}(D_i)}_{\rm Approximation \; Error}+\underbrace{\frac{1}{n\widehat{S}_{n}}\sum_{i=1}^{n}\widehat{W}_{i}\dif K_{h}(D_i)}_{\rm High-dimensional \; Error}.
\label{eq: decomposition}
\end{equation}}
The high-dimensional error is controlled in Theorem \ref{thm: DLL error reduction}. We shall control the stochastic error and the approximation error in Sections \ref{sec: stochastic error} and \ref{sec: approximation error}, respectively.  We present the proof of \eqref{eq: limiting distribution} in Section \ref{sec: limiting proof}. 
%%%%%%%%%%%%%%%%%%%%%%%%%%%%%
%With the above decomposition, we establish Theorem \ref{thm: limiting swap} by applying , together with the following two lemmas. 
%The following lemma establishes that the stochastic error is asymptotically normal, whose proof can be found in Section \ref{sec: stochastic error}.
%%%%%%%%%%%%%%%%%%%%%%%%%%%%
\subsection{Analysis of the Stochastic Error}
\label{sec: stochastic error}
%%%%%%%%%%%%%%%%%%%%%%%%%%%%
%\begin{Lemma} Suppose conditions $\rm (A1)$ hold
%If we further assume  
%\begin{equation}
%\pi(\eval)\gg h^{6} (\log n)^3+\frac{k \log p}{n}\left(h^2+h\frac{k}{n}\max_{i,j}|X_{ij}|^2\right)
%\label{eq: tech condition}
%\end{equation}
% then we have 
We shall establish the following limiting distribution,
\begin{equation}
\frac{1}{\sqrt{\rm V}}\frac{\sum_{i=1}^{n}\widehat{W}_{i}\epsilon_iK_{h}(D_{i})}{n\widehat{S}_{n}} \cid N(0,1),
\label{eq: unconditional distribution}
\end{equation}
with 
\begin{equation}
{\rm V}={\frac{\sigma^2}{n^2\widehat{S}_{n}^2}\sum_{i=1}^{n}\widehat{W}^2_{i}K_{h}^2(D_{i})} \cip \frac{3\sigma^2}{n h^3\cdot\pi(\eval)}.
\label{eq: variance limit}
\end{equation}

In the following, we shall provide proofs for both \eqref{eq: unconditional distribution} and \eqref{eq: variance limit}.

%%%%%%%%%%%%%%%%%%%%%%%%%%%%%%
\subsubsection{Proof of \eqref{eq: variance limit}}
%%%%%%%%%%%%%%%%%%%%%%%%%%%%%%
We decompose the error between $\frac{1}{n}\sum_{i=1}^{n}\widehat{W}_{i}^2K_{h}^2(D_{i})$ and its corresponding estimand,
\begin{equation}
\begin{aligned}
&\left|\frac{1}{n}\sum_{i=1}^{n}\widehat{W}_{i}^2K_{h}^2(D_{i})-\frac{1}{n}\sum_{i=1}^{n}\left(W_i-\bar{\mu}_{W}\right)^2K_{h}^2(D_{i})\right|\\
&=\left|\frac{1}{n}\sum_{i=1}^{n}\left[2\left(W_i-\bar{\mu}_{W}\right)\cdot\left(\widehat{W}_{i}-\left(W_i-\bar{\mu}_{W}\right)\right)+\left(\widehat{W}_{i}-\left(W_i-\bar{\mu}_{W}\right)\right)^2\right]K_{h}^2(D_{i})\right|\\
&\leq \frac{1}{2h}\left(2{\rm Err}(\widehat{W})\cdot\sqrt{\frac{1}{n}\sum_{i=1}^{n}\left(W_i-\bar{\mu}_{W}\right)^2K_{h}(D_{i})}+{\rm Err}^2(\widehat{W})\right),
\end{aligned}
\label{eq: first triangle inequality}
\end{equation}
where the inequality follows from triangle inequality, $|K_{h}(D_{i})|\leq 1/(2h),$ and Cauchy-Schwarz inequality. We bound the difference between the sum of centered variables and that of uncentered variables, 
{\small
\begin{equation*}
%\begin{aligned}
\left|\frac{1}{n}\sum_{i=1}^{n}\left(W_i-\bar{\mu}_{W}\right)^2K_{h}(D_{i})-\frac{1}{n}\sum_{i=1}^{n}W_i^2K_{h}(D_{i})\right|
\leq 2\left|\bar{\mu}_{W}\right|\cdot \left|\frac{1}{n}\sum_{i=1}^{n}W_iK_{h}(D_{i})\right|+2\bar{\mu}_{W}^2 \left|\frac{1}{n}\sum_{i=1}^{n}K_{h}(D_{i})\right|.
%\end{aligned}
\end{equation*}}
By applying \eqref{eq: prob bound 1} and \eqref{eq: prob bound 2} in Lemma \ref{lem: concentration lemma} and \eqref{eq: centering rate}, we establish that, with probability larger than $1-\exp(-t^2)$ for $t\ll \sqrt{nh\pi(\eval)},$
\begin{equation}
\left|\frac{1}{n}\sum_{i=1}^{n}\left(W_i-\bar{\mu}_{W}\right)^2K_{h}(D_{i})-\frac{1}{n}\sum_{i=1}^{n}W_i^2K_{h}(D_{i})\right|\lesssim t^2 \sqrt{\frac{h}{n\pi(\eval)}}  \sqrt{\frac{h}{n}\pi(\eval)}+ \frac{t^2h}{n}\lesssim \frac{ht^2}{n}.
\label{eq: center difference 1}
\end{equation}
Note that 
$$\frac{1}{2h}\cdot \frac{1}{n}\sum_{i=1}^{n}\left(W_i-\bar{\mu}_{W}\right)^2K_{h}(D_{i})=\frac{1}{n}\sum_{i=1}^{n}\left(W_i-\bar{\mu}_{W}\right)^2K_{h}^2(D_{i}).$$
We apply \eqref{eq: expectation 1}, \eqref{eq: prob bound 8}, and \eqref{eq: center difference 1} and establish 
\begin{equation}
\left|\frac{\frac{1}{n}\sum_{i=1}^{n}\left(W_i-\bar{\mu}_{W}\right)^2K_{h}^2(D_{i})}{\frac{1}{3}h\pi(\eval)}-1\right|\lesssim h^2[C_1^2(n)+C_2(n)]+\frac{\mathbf{P}(\mathcal{A}_{0}^c)}{h\pi(\eval)}+\frac{t}{\sqrt{nh\pi(\eval)}}.
\label{eq: inter result 1}
\end{equation}
We shall choose $t=\sqrt{\log n}$ and establish that, with probability larger than $1-n^{-c},$ $$\frac{\frac{1}{n}\sum_{i=1}^{n}\left(W_i-\bar{\mu}_{W}\right)^2K_{h}(D_{i})}{\frac{2}{3}h^2\pi(\eval)}\asymp 1.$$
Combined with \eqref{eq: first triangle inequality} and \eqref{eq: inter result 1}, we establish that, with probability larger than $1-n^{-c},$ 
\begin{equation}
\left|\frac{\frac{1}{n}\sum_{i=1}^{n}\widehat{W}_{i}^2K_{h}^2(D_{i})}{\frac{1}{3}h\pi(\eval)}-1\right|\lesssim
\frac{{\rm Err}^2(\widehat{W})}{h^2\pi(\eval)}+\sqrt{\frac{{\rm Err}^2(\widehat{W})}{h^2\pi(\eval)}}+h^2[C_1^2(n)+C_2(n)]+\frac{\mathbf{P}(\mathcal{A}_{0}^c)}{h\pi(\eval)}+\frac{t}{\sqrt{nh\pi(\eval)}}.
\label{eq: first limit}
\end{equation}
%\begin{equation}
%\left|\frac{1}{n}\sum_{i=1}^{n}\left(W_i-\bar{\mu}_{W}\right)^2K_{h}^2(D_{i})-\frac{1}{n}\sum_{i=1}^{n}W_i^2K_{h}^2(D_{i})\right|\lesssim t \sqrt{\frac{1}{nh\pi(\eval)}}  \sqrt{\frac{h}{n}\pi(\eval)}+ \frac{th}{n}\lesssim \frac{t+ht^2}{n}.
%\end{equation}
%\begin{equation}
%\left|\frac{1}{n}\sum_{i=1}^{n}\widehat{W}_{i}^2K_{h}^2(D_{i})-\E W_i^2K_{h}^2(D_{i})\right|\lesssim \frac{t}{\sqrt{n}}
%\end{equation}

For $\widehat{S}_n$ defined in \eqref{eq: final estimator}, we approximate it by its corresponding estimand, 
\begin{equation}
\begin{aligned}
&\left|\widehat{S}_{n}-\frac{1}{n}\sum_{i=1}^{n}\left(W_i-\bar{\mu}_{W}\right)(D_{i}-\eval)K_{h}(D_{i})\right|\\
&\leq {\rm Err}(\widehat{W}) \cdot \sqrt{\frac{1}{n}\sum_{i=1}^{n}\left(D_{i}-\eval\right)^2K_{h}(D_{i})} \leq h \cdot {\rm Err}(\widehat{W}) \cdot \sqrt{\frac{1}{n}\sum_{i=1}^{n}K_{h}(D_{i})}.
\end{aligned}
\label{eq: inter result 2}
\end{equation}
%Combined with \eqref{eq: first triangle inequality} and \eqref{eq: second triangle inequality}
%\begin{equation}
%\left|\frac{\E W_i(D_{i}-\eval)K_{h}(D_{i})}{\frac{1}{3}h^2\pi(\eval)}-1\right|\lesssim h^2[C_1^2(n)+C_2(n)]+\frac{3\mathbf{P}(\mathcal{A}_{0}^c)}{h\pi(\eval)}.
%\label{eq: expectation 2}
%\end{equation}
%\begin{equation}
%\left|\frac{1}{n}\sum_{i=1}^{n}W_i{(D_{i}-\eval)} K_{h}(D_{i})-\E W_i(D_{i}-\eval)K_{h}(D_{i})\right|\leq C h^2\pi(\eval)t\sqrt{\frac{1}{nh\pi(\eval)}};
%\label{eq: prob bound 5}
%\end{equation}

By applying \eqref{eq: expectation 2} and \eqref{eq: prob bound 5}, we establish that, with probability larger than $1-\exp(-t^2),$
\begin{equation}
\left|\frac{\frac{1}{n}\sum_{i=1}^{n}W_i{(D_{i}-\eval)} K_{h}(D_{i})}{\frac{1}{3}h^2\pi(\eval)}-1\right|\lesssim h^2[C_1^2(n)+C_2(n)]+\frac{\mathbf{P}(\mathcal{A}_{0}^c)}{h\pi(\eval)}+t\sqrt{\frac{1}{nh\pi(\eval)}}.
\label{eq: inter result 3}
\end{equation}

%We apply Law of Large Numbers and Lemma \ref{lem: expectation lemma},
%\begin{equation}
%\frac{\frac{1}{n}\sum_{i=1}^{n}W_i^2K_{h}^2(D_{i})}{\frac{2}{3}h\pi(\eval)} \cip 1 \quad \text{and} \quad \frac{ \frac{1}{n}\sum_{i=1}^{n}W_i(D_{i}-\eval)K_{h}(D_{i})}{\frac{2}{3}h^2\pi(\eval)}\cip 1.
%\label{eq: lln}
%\end{equation}
%For the term
%\begin{equation*}
%%\begin{aligned}
%\left|\frac{1}{n}\sum_{i=1}^{n}\bar{\mu}_{W}(D_{i}-\eval)K_{h}(D_{i})\right|= \left|\bar{\mu}_{W}\right|\cdot \left|\frac{1}{n}\sum_{i=1}^{n}(D_{i}-\eval)K_{h}(D_{i})\right|,
%\end{equation*}
By \eqref{eq: centering rate} and \eqref{eq: prob bound 4} in Lemma \ref{lem: concentration lemma}, we establish that, with probability larger than $1-\exp(-t^2)$,
{\small
\begin{equation}
\left|\frac{1}{n}\sum_{i=1}^{n}\bar{\mu}_{W}(D_{i}-\eval)K_{h}(D_{i})\right|\lesssim  t\sqrt{\frac{h}{n\pi(\eval)}}\cdot\left( \frac{{\pi}(\eval)}{3}h^2 (C_1(n)+hC_2(n))+{\mathbf{P}(\mathcal{A}_0^c)}+t\sqrt{\frac{h}{n} \pi(\eval)}\right). 
\label{eq: center difference 2}
\end{equation}
}
%By taking $t=\sqrt{nh\pi(\eval)}$, we combine \eqref{eq: center difference 1}, \eqref{eq: center difference 2} and \eqref{eq: lln} and establish
%\begin{equation}
%\frac{\frac{1}{n}\sum_{i=1}^{n}\left(W_i-\bar{\mu}_{W}\right)^2K_{h}^2(D_{i})}{\frac{2}{3}h\pi(\eval)} \cip 1 \quad \text{and} \quad \frac{ \frac{1}{n}\sum_{i=1}^{n}\left(W_i-\bar{\mu}_{W}\right)(D_{i}-\eval)K_{h}(D_{i})}{\frac{2}{3}h^2\pi(\eval)}\cip 1.
%\label{eq: lln approximation}
%\end{equation}
%Also note the following fact $$\left|\frac{\frac{1}{h}\left({\rm Err}^2(\widehat{W})+{\rm Err}(\widehat{W})\cdot\sqrt{\frac{1}{n}\sum_{i=1}^{n}\left(W_i-\bar{\mu}_{W}\right)^2K_{h}^2(D_{i})}\right)}{\frac{1}{3}h\pi(\eval)}\right|\lesssim $$ %\leq \left(\frac{h}{\sigma_2^2} \left(1+c_u^2\right) \sqrt{\frac{k \log p\cdot \log n}{n}}\right)^2,$
%and
Together with \eqref{eq: inter result 2} and \eqref{eq: inter result 3}, we establish that, with probability larger than $1-\exp(-t^2),$
\begin{equation}
\left|\frac{\widehat{S}_{n}}{{\frac{1}{3}h^2\pi(\eval)}}-1\right|\lesssim \frac{{\rm Err}(\widehat{W})}{h \sqrt{\pi(\eval)}}+\left(h^2+\sqrt{\frac{h}{n\pi(\eval)}}\right)[C_1^2(n)+C_2(n)]+\frac{\mathbf{P}(\mathcal{A}_{0}^c)}{h\pi(\eval)}+t\sqrt{\frac{1}{nh\pi(\eval)}}.
\label{eq: second limit}
\end{equation}
Under the condition that ${\rm Err}(\widehat{W})\ll h \sqrt{\pi(\eval)},$ $h^2 C_2(n)+hC_1(n)\rightarrow 0,$ $\mathbf{P}(\mathcal{A}_{0}^c)\ll h\pi(\eval),$ and $nh\pi(\eval)\gg \log n,$ we establish \eqref{eq: variance limit} by combining \eqref{eq: first limit} and \eqref{eq: second limit}.

%We introduce the following lemma to facilitate the proof and present the corresponding proof in Section \ref{sec: proof lemma 8}.
%\begin{Lemma} Under the condition that ${\rm Err}(\widehat{W})\ll h\sqrt{\pi(\eval)}$, $hC_{u}\leq 1$  and $nh\pi(\eval)\rightarrow \infty$ where $C_{u}$ is defined in \eqref{eq: second key growing constant}, then
%%\begin{equation} 
%%\frac{1}{n}\sum_{i=1}^{n}\widehat{W}_{i}^2K_{h}^2(D_{i})\rightarrow \frac{2}{3}h \E\phi\left(\frac{\eval-X_{i}^{\intercal}\gamma}{\sigma_2}\right)
%%\label{eq: first limit}
%%\end{equation}
%\begin{equation}
%\frac{\frac{1}{n}\sum_{i=1}^{n}\widehat{W}_{i}^2K_{h}^2(D_{i})}{\frac{2}{3}h\pi(\eval)}\cip 1
%\label{eq: first limit}
%\end{equation}
%%\begin{equation}
%%\widehat{S}_{n}\rightarrow \frac{2}{3}h^2 \E\phi\left(\frac{\eval-X_{i}^{\intercal}\gamma}{\sigma_2}\right)
%%\label{eq: second limit}
%%\end{equation}
%\begin{equation}
%\frac{\widehat{S}_{n}}{\frac{2}{3}h^2 \pi(\eval)}\cip 1 %\frac{2}{3}h^2 \E\phi\left(\frac{\eval-X_{i}^{\intercal}\gamma}{\sigma_2}\right)
%\label{eq: second limit}
%\end{equation}
%In addition, with probability larger than $1-{(nh\pi(\eval))^{-\frac{1}{4}}}$, 
%\begin{equation}
%\left|\frac{\widehat{S}_{n}}{\frac{2}{3}h^2 \pi(\eval)}-1\right| \lesssim \frac{{\rm Err}(\widehat{W})}{h\pi(\eval)}+{(nh\pi(\eval))^{-\frac{1}{4}}}+\frac{(nh\pi(\eval))^{\frac{1}{4}}}{n}
%%\frac{2}{3}h^2 \E\phi\left(\frac{\eval-X_{i}^{\intercal}\gamma}{\sigma_2}\right)
%\label{eq: high prob Sn}
%\end{equation}
% \label{lem: local sum limit}
%\end{Lemma}

\subsubsection{Proof of \eqref{eq: unconditional distribution}}
Define 
$Z_i=\frac{1}{\sqrt{\rm V/\sigma^2}}\frac{\widehat{W}_{i} K_{h}(D_{i})}{n\widehat{S}_{n}}\in \R$. We rewrite the stochastic error as follows, $$\frac{1}{\sqrt{\rm V}}\frac{\sum_{i=1}^{n}\widehat{W}_{i}\epsilon_iK_{h}(D_{i})}{n\widehat{S}_{n}}=\sum_{i=1}^{n} Z_{i} \cdot {\epsilon_i}/{\sigma}.$$ 
We use $\mathcal{O}$ to denote the data $\mathcal{O}=\{D_i,X_i\}_{1\leq i\leq n}.$
Conditioning on $\mathcal{O}$,  $Z_{i} \cdot {\epsilon_i}/{\sigma}$ are independent random variables with $$\E\left(Z_{i} \cdot {\epsilon_i}/{\sigma}\mid \mathcal{O}\right)=0$$ 
and 
$$\sum_{i=1}^{n}{\rm Var}\left(Z_{i} \cdot {\epsilon_i}/{\sigma}\mid \mathcal{O}\right)=1.$$
Define the event 
$$\mathcal{G}_0=\left\{\left|\frac{\frac{1}{n}\sum_{i=1}^{n}\widehat{W}_{i}^2K_{h}^2(D_{i})}{\frac{1}{3}h\pi(\eval)}-1\right|\leq 1/10\right\}.$$
The high probability inequality in \eqref{eq: first limit} implies
\begin{equation}
\mathbf{P}(\mathcal{G}_0)\geq 1-n^{-c}.
\label{eq: high prob clt}
\end{equation}
By applying \eqref{eq: first limit} and the fact that $|\widehat{W}_{i}K_{h}(D_{i})|\leq C$ for a positive constant $C>0$, we obtain that, on the event $\mathcal{G}_0$, \begin{equation}
\left|Z_i\right|=\left|\frac{\widehat{W}_{i}K_{h}(D_{i})}{\sqrt{\sum_{i=1}^{n}\widehat{W}^2_{i}K_{h}^2(D_{i})}}\right|\lesssim \frac{1}{\sqrt{nh\pi(\eval)}}.
\label{eq: clt bound}
\end{equation}

It is sufficient to check the Linderberg condition 
\begin{equation*}
\begin{aligned}
\sum_{i=1}^{n}\E \left[(Z_{i} \cdot {\epsilon_i}/{\sigma})^2 {\1}\left(\left|Z_{i} \cdot {\epsilon_i}/{\sigma}\right|\geq \tau\right)\mid \mathcal{O}\right]
&=\sum_{i=1}^{n}Z_{i}^2 \E \left[\frac{\epsilon_i^2}{\sigma^2} {\1}\left(\left|\frac{\epsilon_i}{\sigma}\right|\geq \frac{\tau}{|Z_{i}|}\right)\mid \mathcal{O}\right]\\
&\leq \sum_{i=1}^{n}Z_{i}^2 \E \left[\frac{\epsilon_i^2}{\sigma^2} {\1}\left(\left|\frac{\epsilon_i}{\sigma}\right|\gtrsim {\tau}\sqrt{nh\pi(\eval)}\right)\mid \mathcal{O}\right]\\
&\leq \sum_{i=1}^{n}Z_{i}^2 \max_{1\leq i\leq n}\E \left[\frac{\epsilon_i^2}{\sigma^2} {\1}\left(\left|\frac{\epsilon_i}{\sigma}\right|\gtrsim {\tau}\sqrt{nh\pi(\eval)}\right)\mid \mathcal{O}\right]\\
&\leq ({\tau}\sqrt{nh\pi(\eval)})^{-c},
\end{aligned}
\end{equation*}
where the first inequality follows from \eqref{eq: clt bound} and the last inequality follows from the condition that $\E(\epsilon_i^{2+c}\mid D_i, X_i)\leq C$ for some positive constant $c>0$ and $C>0.$

Then we apply the Linderberg condition and establish that 
\begin{equation}
\sum_{i=1}^{n}\frac{1}{\sqrt{\rm V}}\frac{\widehat{W}_{i}\epsilon_iK_{h}(D_{i})}{n\widehat{S}_{n}}\mid \mathcal{O} \in \mathcal{G}_0 \cid N(0,1).
\label{eq: conditioning limit}
\end{equation}
By \eqref{eq: high prob clt} and \eqref{eq: conditioning limit}, we have 
\begin{equation*}
\E \left(\E\left[\exp\left(it (\sum_{i=1}^{n} Z_{i}{\epsilon_i}/{\sigma})\right)\mid \mathcal{O}\right] \cdot {\bf 1}_{\mathcal{G}_0}\right)\rightarrow \exp(-t^2/2).
\end{equation*}
Together with \begin{equation*}
\left|\E \exp\left(it (\sum_{i=1}^{n} Z_{i}{\epsilon_i}/{\sigma})\right)-\E \left(\E\left[\exp\left(it (\sum_{i=1}^{n} Z_{i}{\epsilon_i}/{\sigma})\right)\mid \mathcal{O}\right] \cdot {\bf 1}_{\mathcal{G}_0}\right)\right|\leq \mathbf{P}\left(\mathcal{G}_0^{c}\right),
\label{eq: approximation clt}
\end{equation*} 
we establish \eqref{eq: unconditional distribution}.

%Combined with \eqref{eq: first triangle inequality} and \eqref{eq: second triangle inequality}, we establish \eqref{eq: first limit} and \eqref{eq: second limit}.
%In addition, together with \eqref{eq: center difference 2} and \eqref{eq: inter 4}, we apply \eqref{eq: prob bound 5} with $t=(4nh\pi(\eval))^{\frac{1}{4}}$ and establish \eqref{eq: high prob Sn}.

%The following lemma establishes the limiting distribution for the stochastic error and the rate of convergence for the approximation error. %The rate of convergence for the nuisance error is deferred to the next subsection.
%%%%%%%%%%%%%%%%%%%%%%%%%%%%%%%%%%%%%%%%%%
%%%%%%%%%%%%%%%%%%%%%%%%%%%%%%%%%%%%%%%%%%%%%%%%%%%%%%%%%%%%%%%%%

%%%%%%%%%%%%%%%%%%%%%%%%%%%%%%%%%%%%%%%%%%%%%%%%%%%%%%%%%%%%%%%%%
\subsection{Analysis of the Approximation Error}
\label{sec: approximation error}
%%%%%%%%%%%%%%%%%%%%%%%%%%%%%%%%%%%%%%%%%%%%%%%%%%%%%%%%%%%%%%%%%
In the following, we show that, with probability larger than $1-n^{-c},$
{\small
\begin{equation}
\left|\frac{1}{n{\widehat{S}_{n}}}{\sum_{i=1}^{n}\widehat{W}_{i}\left[r(D_{i})-\frac{(D_{i}-\eval)^2}{2} f''(\eval)\right]K_{h}(D_{i})}\right| \lesssim  \max_{|d-\eval|\leq h} \left|f''(d)-f''(\eval)\right|\cdot \left(\frac{{\rm Err}(\widehat{W})}{\sqrt{\pi(\eval)}}+h\right), %o_{p}\left(\frac{{\rm Err}(\widehat{W})}{\sqrt{\pi(\eval)}}+h\right)
\label{eq: first approximation}
\end{equation} }
and 
\begin{equation}
\begin{aligned}
\left|\frac{1}{n{\widehat{S}_{n}}}{\sum_{i=1}^{n}\widehat{W}_{i}\frac{(D_{i}-\eval)^2}{2} f''(\eval)K_{h}(D_{i})}\right|\lesssim \frac{{\rm Err}(\widehat{W})}{\sqrt{\pi(\eval)}}+h c_u,
\end{aligned}
\label{eq: second approximation}
\end{equation} 
with $$c_u=h C_1(n)+h^2 C_2(n)+\sqrt{\frac{\log n}{{n h\pi(\eval)}}}+\frac{\mathbf{P}(\mathcal{A}_0^c)}{h \pi(\eval)}.$$
By the continuity of $f''$ at the point $\eval,$ we combine \eqref{eq: first approximation} and \eqref{eq: second approximation} and establish that, 
\begin{equation}
\left|\frac{1}{n{\widehat{S}_{n}}}{\sum_{i=1}^{n}\widehat{W}_{i}r(D_{i})K_{h}(D_{i})}\right| \lesssim  \frac{{\rm Err}(\widehat{W})}{\sqrt{\pi(\eval)}}+\left(c_u+\max_{|d-\eval|\leq h} \left|f''(d)-f''(\eval)\right|\right) \cdot h.
\label{eq: approximation error}
\end{equation} 
%The following lemma controls the approximation error of the \texttt{DLL} estimator, whose proof can be found in Section \ref{sec: approximation error}. 
%\begin{Lemma}
%Suppose conditions $\rm (A1)$ hold. Then the approximation error satisfies
%
%\label{lem: approximation error}
%\end{Lemma}

%The proof of Theorem \ref{thm: estimation} follows from a combination of Lemma \ref{lem: variance and approximation} and Theorem \ref{thm: nuisance error}. Theorem \ref{thm: known covariance} follows from Theorems \ref{thm: limiting swap} and \ref{thm: estimation} by taking ${\rm Err}(\widehat{W})=0$.
%\subsection{Proof of Theorem \ref{thm: DLL error reduction}}

%\begin{equation}
%\frac{1}{n}\sum_{i=1}^{n} \E (F_i(1)-F_i(0))^2\frac{1}{h}
%{\bf{1}}(|\delta_i-\mu_i|\leq h)  \lesssim \max_{1\leq i, j\leq n}|a_{ij}| h^2 \left(C_1^2(n)+C_2(n)\right) \phi(\mu_i)
%\end{equation}

%%%%%%%%%%%%%%%%%%%%%%%%%%%%%%%%%%%%%%%%%%%%%%%%%%%%%%%%%%%%%%%%%%%%%%%%%%%%%%%%
%\subsection{Proof of Lemma \ref{lem: stochastic error}}
%\label{sec: stochastic error}
%%%%%%%%%%%%%%%%%%%%%%%%%%%%%%%%%%%%%%%%%%%%%%%%%%%%%%%%%%%%%%%%%%%%%%%%%%%%%%%%

%%%%%%%%%%%%%%%%%%%%%%%%%%%%%%%%%%%%%%%%%%%%%%%%%%%%%%%%%%%%%%%%%%%%%%%%%%%%%%%%
\subsubsection{Proof of \eqref{eq: first approximation}}
%\label{sec: approximation error}
%%%%%%%%%%%%%%%%%%%%%%%%%%%%%%%%%%%%%%%%%%%%%%%%%%%%%%%%%%%%%%%%%%%%%%%%%%%%%%%%
%The control of the approximation error in \eqref{eq: first approximation} follows from a combination of \eqref{eq: second limit} in Lemma \ref{lem: local sum limit} and the control of 
%$\frac{1}{n}\sum_{i=1}^{n}\widehat{W}_{i}r(D_{i})K_{h}(D_{i}).$
There exists some $c\in (0,1)$ such that 
\begin{equation*}
\begin{aligned}
r(D_{i})&=f(D_{i})-f(\eval)-(D_{i}-\eval)f'(\eval)\\
&=\frac{(D_{i}-\eval)^2}{2} f''(\eval)+\frac{(D_{i}-\eval)^2}{2}\left[ f''(\eval+c(D_{i}-\eval))-f''(\eval)\right].
\end{aligned}
\end{equation*}
Hence, we have
\begin{equation*}
\frac{2}{h^2}\left|r(D_{i})\mathbf{1}\left(\left|\frac{D_{i}-\eval}{h}\right|\leq 1\right)-\frac{(D_{i}-\eval)^2}{2} f''(\eval)\mathbf{1}\left(\left|\frac{D_{i}-\eval}{h}\right|\leq 1\right)\right|\leq \left|f''(d)-f''(\eval)\right|,
\end{equation*}
for some $d$ satisfying $\eval-h\leq d\leq \eval+h$. 
%Since $h=h(n)\rightarrow 0$ and $f''(d)$ is continuous at $\eval$, then we have
%\begin{equation*}
%\frac{2}{h^2}\left|r(D_{i})\mathbf{1}\left(\left|\frac{D_{i}-\eval}{h}\right|\leq 1\right)-\frac{(D_{i}-\eval)^2}{2} f''(\eval)\mathbf{1}\left(\left|\frac{D_{i}-\eval}{h}\right|\leq 1\right)\right|\rightarrow 0.
%\end{equation*}
The above inequality implies that 
\begin{equation}
\frac{\left|\frac{1}{n}\sum_{i=1}^{n}\widehat{W}_{i}r(D_{i})K_{h}(D_{i})-\frac{1}{n}\sum_{i=1}^{n}\widehat{W}_{i}\frac{(D_{i}-\eval)^2}{2} f''(\eval)K_{h}(D_{i})\right|}{h^2 \pi(\eval)\cdot \frac{1}{n \pi(\eval)}\sum_{i=1}^{n}\left|\widehat{W}_{i}\right|K_{h}(D_{i})}\lesssim \max_{|d-\eval|\leq h}  \left|f''(d)-f''(\eval)\right|.
\label{eq: approximation decomposition}
\end{equation}
%The proof of \eqref{eq: higher order bias} follows from the following inequality, 
We now control the term $\frac{1}{n}\sum_{i=1}^{n}\left|\widehat{W}_{i}\right|K_{h}(D_{i}),$
\begin{equation*}
\begin{aligned}
\frac{1}{n}\sum_{i=1}^{n}\left|\widehat{W}_{i}\right|K_{h}(D_{i})&\leq \frac{1}{n}\sum_{i=1}^{n}\left|\widehat{W}_{i}-\left(W_{i}-\bar{\mu}_{W}\right)\right|K_{h}(D_{i})+\frac{1}{n}\sum_{i=1}^{n}\left(\left|W_{i}\right|+\left|\bar{\mu}_{W}\right|\right)K_{h}(D_{i})\\
&\leq {\rm Err}(\widehat{W})\sqrt{\frac{1}{n}\sum_{i=1}^{n} K_{h}(D_{i})}+\frac{2}{n}\sum_{i=1}^{n}\left|W_{i}\right|K_{h}(D_{i})\\
&\leq {\rm Err}(\widehat{W})\sqrt{\frac{1}{n}\sum_{i=1}^{n} K_{h}(D_{i})}+4h\left(\frac{1}{n}\sum_{i=1}^{n} K_{h}(D_{i})\right),
%&\lesssim \left({\rm Err}(\widehat{W})+2h\right)\cdot\E\phi\left(\frac{\eval-X_{i}^{\intercal}\gamma}{\sigma_2}\right)
\end{aligned}
\label{eq: bound higher order}
\end{equation*}
where the last inequality follows from the fact that $\left|W_{i}\right|K_{h}(D_{i})\leq 2h  K_{h}(D_{i})$. Together with \eqref{eq: prob bound 1} with $t=\sqrt{\log n}$, we establish that, with probability larger than $1-n^{-c},$
\begin{equation}
\frac{1}{n\pi(\eval)}\sum_{i=1}^{n}\left|\widehat{W}_{i}\right|K_{h}(D_{i}) \lesssim \frac{{\rm Err}(\widehat{W})}{\sqrt{\pi(\eval)}}+h.
\label{eq: higher order bias}
\end{equation}
We establish \eqref{eq: first approximation} by combining \eqref{eq: approximation decomposition}, \eqref{eq: higher order bias}, and \eqref{eq: second limit}. 
%With the above calculation, the problem of controlling the approximation error is reduced to the control of the two terms $\frac{1}{n}\sum_{i=1}^{n}\widehat{W}_{i}\frac{(D_{i}-\eval)^2}{2} f''(\eval)K_{h}(D_{i})$ and $\frac{1}{n}\sum_{i=1}^{n}\left|\widehat{W}_{i}\right|K_{h}(D_{i})$, which are established in the following lemma. The proof of the following lemma is present in Section \ref{sec: proof lemma 9}.

%%%%%%%%%%%%%%%%%%%%%%%%%%%%%%%%%%%%%%%%%%%%%%%%%%%%%%%%%%%
\subsubsection{Proof of \eqref{eq: second approximation}}
%%%%%%%%%%%%%%%%%%%%%%%%%%%%%%%%%%%%%%%%%%%%%%%%%%%%%%%%%%%%

%
%\Zijian{Checking the remaining proof.}
%\begin{Lemma}
%Suppose that $h C_{u}\rightarrow 0$ and $nh\pi(\eval) \rightarrow \infty$, 
%then with probability larger than $1-(nh\pi(\eval))^{-\frac{1}{4}}$, 
%\begin{equation}
%\left|\frac{1}{nh^2 \pi(\eval)}\sum_{i=1}^{n}\widehat{W}_{i}\frac{(D_{i}-\eval)^2}{2} f''(\eval)K_{h}(D_{i})\right|\lesssim \frac{{\rm Err}(\widehat{W})}{\sqrt{\pi(\eval)}}+h\left({h C_u}+\frac{1}{(n h\pi(\eval))^{\frac{1}{4}}}\right)
%%({\rm Err}(\widehat{W})+h)\cdot h^2\sqrt{\pi(\eval)}
%%\lesssim \left(1+c_u^3\right)  \sqrt{\frac{k \log p \log n}{n}} h^4+C h^{4} \sqrt{\log n} \pi(\eval)+t\sqrt{\frac{h^{5}}{4n} \pi(\eval)}+h^4\sqrt{\log n}\cdot\pi(\eval)+h^2 \sqrt{\frac{2h}{3n \pi(\eval)}}\left[ \pi(\eval)+ \frac{t}{\sqrt{nh}}\sqrt{\pi(\eval)}\right].
%\label{eq: second order bias}
%\end{equation}
%
%\label{lem: bound bias a}
%\end{Lemma}
%%By applying , we have 
%%\begin{equation}
%%\left|\frac{1}{n}\sum_{i=1}^{n}\widehat{W}_{i}r(D_{i})K_{h}(D_{i})\right|\lesssim ({\rm Err}(\widehat{W})+h)\cdot h^2\sqrt{\pi(\eval)}
%%\label{eq: approximation bias bound}
%%\end{equation}
%Combination of \eqref{eq: second limit} and \eqref{eq: second order bias} leads to \eqref{eq: second approximation}.
%Combining \eqref{eq: second limit}, \eqref{eq: approximation decomposition} and \eqref{eq: higher order bias}, we establish \eqref{eq: first approximation}.

By the expression $\widehat{W}_{i}= \left(W_{i}-\bar{\mu}_{W}\right)+\widehat{W}_{i}-\left(W_{i}-\bar{\mu}_{W}\right)$, we have 
\begin{equation}
\begin{aligned}
&\frac{1}{n}\sum_{i=1}^{n}\widehat{W}_{i}\frac{(D_{i}-\eval)^2}{2} f''(\eval)K_{h}(D_{i})=\frac{1}{n}\sum_{i=1}^{n}\left[\widehat{W}_{i}-\left(W_{i}-\bar{\mu}_{W}\right)\right]\frac{(D_{i}-\eval)^2}{2} f''(\eval)K_{h}(D_{i})\\
&+\frac{1}{n}\sum_{i=1}^{n}W_i\frac{(D_{i}-\eval)^2}{2} f''(\eval)K_{h}(D_{i})-\bar{\mu}_{W}\frac{1}{n}\sum_{i=1}^{n}\frac{(D_{i}-\eval)^2}{2} f''(\eval)K_{h}(D_{i}).
\end{aligned}
\label{eq: approx decomposition}
\end{equation}
By the Cauchy-Schwarz inequality, we have 
\begin{equation}
\begin{aligned}
&\left|\frac{1}{n}\sum_{i=1}^{n}\left[\widehat{W}_{i}-\left(W_{i}-\bar{\mu}_{W}\right)\right]\frac{(D_{i}-\eval)^2}{2} f''(\eval)K_{h}(D_{i})\right|\\
&\lesssim \left|f''(\eval)\right| {\rm Err}(\widehat{W})\cdot \sqrt{\frac{1}{n}\sum_{i=1}^{n}\frac{(D_{i}-\eval)^4}{2} K_{h}(D_{i})}\\
& \leq \left|f''(\eval)\right| \cdot h^2 \cdot {\rm Err}(\widehat{W}) \cdot \sqrt{\frac{1}{n}\sum_{i=1}^{n} K_{h}(D_{i})},
\end{aligned}
\label{eq: approx part 1}
\end{equation}
where the last inequality follows from the fact that ${(D_{i}-\eval)^4} K_{h}(D_{i})\leq h^4  K_{h}(D_{i})$.
In addition, we have
\begin{equation}
\left|\frac{1}{n}\sum_{i=1}^{n}W_i\frac{(D_{i}-\eval)^2}{2} f''(\eval)K_{h}(D_{i})\right|= \left|f''(\eval)\right| \cdot \left| \frac{1}{n}\sum_{i=1}^{n}W_i\frac{(D_{i}-\eval)^2}{2} K_{h}(D_{i})\right|,
%\left(\frac{8}{45}h^{4} \E \left[ \frac{d}{d_{\eval}}\mathbf{P}(\eval\mid X_{i})\right]+O(h^{6})\right).
\label{eq: approx part 2}
\end{equation}
%By \eqref{eq: centering rate}, we have 
and 
\begin{equation}
\begin{aligned}
\left|\bar{\mu}_{W}\frac{1}{n}\sum_{i=1}^{n}\frac{(D_{i}-\eval)^2}{2} f''(\eval)K_{h}(D_{i})\right|&= \left|\bar{\mu}_{W}\right|\cdot \left|f''(\eval)\right|\cdot\left|\frac{1}{n}\sum_{i=1}^{n}\frac{(D_{i}-\eval)^2}{2} K_{h}(D_{i})\right|\\
&\leq \frac{h^2}{2} \left|\bar{\mu}_{W}\right|\cdot \left|f''(\eval)\right|\cdot \frac{1}{n}\sum_{i=1}^{n} K_{h}(D_{i}),\end{aligned}
\label{eq: approx part 3}
\end{equation}
where the last inequality follows from the fact that ${(D_{i}-\eval)^2} K_{h}(D_{i})\leq h^2  K_{h}(D_{i})$.
We now apply \eqref{eq: prob bound 1}, \eqref{eq: prob bound 3}, \eqref{eq: centering rate}, the decomposition \eqref{eq: approx decomposition} with the error bounds in \eqref{eq: approx part 1}, \eqref{eq: approx part 2}, and \eqref{eq: approx part 3}. We establish that, with probability larger than $1-n^{-c},$
\begin{equation*}
\begin{aligned}
&\frac{1}{h^2\pi(\eval)}\left|\frac{1}{n}\sum_{i=1}^{n}\widehat{W}_{i}\frac{(D_{i}-\eval)^2}{2} f''(\eval)K_{h}(D_{i})\right|\\
&\lesssim \frac{{\rm Err}(\widehat{W})}{\sqrt{\pi(\eval)}}+h\left(h C_1(n)+h^2 C_2(n)+\sqrt{\frac{\log n}{{n h\pi(\eval)}}}+\frac{\mathbf{P}(\mathcal{A}_0^c)}{h \pi(\eval)}\right).
%\left(h^4 \left(\log \max\{n,p\}\right)^{\frac{5}{2}}+h^3\sqrt{\frac{1}{n^2h^2 \pi(\eval)}}\right)\sqrt{\pi(\eval)}
\end{aligned}
%\label{eq: second order bias}
\end{equation*}
Together with \eqref{eq: second limit}, we establish \eqref{eq: first approximation}.
%Taking $t=(nh\pi(\eval))^{\frac{1}{4}}$, then we establish \eqref{eq: second order bias}.}
%%%%%%%%%%%%%%%%%%%%%%%%%%%%%
\subsection{Proof of \eqref{eq: limiting distribution}}
\label{sec: limiting proof}
Under the conditions ${\rm Err}(\widehat{W})\ll \min\{\sqrt{nh^3},h \sqrt{\pi(\eval)}\},$ $h^2 C_2(n)+hC_1(n)\rightarrow 0,$ $\mathbf{P}(\mathcal{A}_{0}^c)\ll h\pi(\eval),$ and $n h^5 \pi(\eval)\leq c$, we apply \eqref{eq: unconditional distribution}, \eqref{eq: variance limit}, and \eqref{eq: approximation error} and establish 
\begin{equation}
\frac{1}{\sqrt{\rm V}}\left({\frac{1}{n\widehat{S}_{n}}\sum_{i=1}^{n}\widehat{W}_{i}\epsilon_iK_{h}(D_i)}+{\frac{1}{n\widehat{S}_{n}}\sum_{i=1}^{n}\widehat{W}_{i}r(D_i)K_{h}(D_i)}\right)\cid N(0,1).
\label{eq: part 1}
\end{equation}

% 
%\begin{equation}
%\frac{1}{\sqrt{V}}\frac{1}{n{\widehat{S}_{n}}}{\sum_{i=1}^{n}\widehat{W}_{i}r(D_{i})K_{h}(D_{i})}=o_{p}\left(1\right)
%\label{eq: approximation error control}
%\end{equation}
%where the last equality follows from the condition 
{It follows from \eqref{eq: second limit} and \eqref{eq: high dim error} that, with probability larger than $1-\frac{1}{t}-\min\{n,p\}^{-c}$ for some $t>1$ 
\begin{equation*}
%\begin{aligned}
\frac{1}{\sqrt{V}}\left|\frac{1}{n\widehat{S}_{n}}\sum_{i=1}^{n}\widehat{W}_{i}\dif K_{h}(D_i)\right|\lesssim t^2 \left[1+\sqrt{{h^3k\log p\log n}} \left(C_1^2(n)+C_2(n)\right)\right]\frac{{\rm Err}(\widehat{g})}{\sqrt{\pi(\eval)}}.
%\end{aligned}
\label{eq: nuisance error}
\end{equation*}}
The condition (A3) implies that $$\frac{1}{\sqrt{V}}\left|\frac{1}{n\widehat{S}_{n}}\sum_{i=1}^{n}\widehat{W}_{i}\dif K_{h}(D_i)\right|\cip 0.$$
Combined with \eqref{eq: part 1}, we establish the limiting distribution \eqref{eq: limiting distribution}.

\section{Proofs of Extra Lemmas}
\subsection{Proof of Lemma \ref{lem: Gaussian bound}}
%%%%%%%%%%%%%%%%%%%%%%%%%%%%%%%%%%%%%%%%%%%%%%%%%%%%%%%

By Taylor's expansion, we have for $c_1,c_2\in (0,1)$
\begin{equation*}
\phi(\mu+\tau)+\phi(\mu-\tau)=\phi(\mu)+\tau \cdot \phi'(\mu)+\frac{\tau^2}{2} \cdot \phi''(\mu+c_1 \tau)+\phi(\mu)-\tau \cdot \phi'(\mu)+\frac{\tau^2}{2} \cdot \phi''(\mu+c_2 \tau)
\end{equation*}
and
\begin{equation*}
\frac{1}{\tau}\int_{\mu-\tau}^{\mu+\tau}\phi(t) dt=\frac{1}{\tau}\int_{\mu-\tau}^{\mu+\tau}\left[\phi(\mu)+{(t-\mu)} \phi'(\mu)+\frac{(t-\mu)^2}{2} \phi''(\mu+c_3(t)(t-\mu))\right]dt,
\end{equation*}
where $c_3(t)\in (0,1)$.
Hence, we have 
\begin{equation*}
\begin{aligned}
\phi(\mu+\tau)+\phi(\mu-\tau)-\frac{1}{\tau}\int_{\mu-\tau}^{\mu+\tau}\phi(t) dt&=\frac{\tau^2}{2} \cdot \left(\phi''(\mu+c_1 \tau)+\phi''(\mu+c_2 \tau)\right)\\
&+\frac{1}{\tau}\int_{\mu-\tau}^{\mu+\tau}\frac{(t-\mu)^2}{2} \phi''(\mu+c_3(t)(t-\mu))dt.
\end{aligned}
\end{equation*}
Hence, we establish \eqref{eq: Gaussian bound 1}.
Note that 
\begin{equation*}
\int_{\mu-\tau}^{\mu+\tau} \left(t-\mu\right)\phi(t) dt=\int_{\mu-\tau}^{\mu+\tau}\left[(t-\mu)\phi(\mu)+{(t-\mu)^2}\phi'(\mu+c_4(t)(t-\mu))\right]dt,
\end{equation*}
where $c_4(t) \in (0,1)$. Hence  we establish \eqref{eq: Gaussian bound 2}.
Note that 
\begin{equation*}
\phi(\mu+\tau)-\phi(\mu-\tau)=2\tau \cdot \phi'(\mu+c_5\tau),
\end{equation*}
for $c_5\in (-1,1)$. Hence we establish \eqref{eq: Gaussian bound 3}.

\subsection{Proof of Lemma \ref{lem: expectation lemma}}
\label{sec: exp lemma}
%%%%%%%%%%%%%%%%%%%%%%%%%%%%%%%%%%%%%%%%%%%%%%%%%%%%%%%
\subsubsection{Proof of \eqref{eq: conditional zero order} and \eqref{eq: expectation 0}}
%We focus on the analysis of $\E\left[K_{h}(D_{i})\mid X_{i}\right]$ in the following. 
We start with the expression of $\E\left[K_{h}(D_{i})\mid X_{i}\right],$
\begin{equation*}
\E\left[K_{h}(D_{i})\mid X_{i}\right]=\int_{\left|\frac{D_{i}-\eval}{h}\right|\leq 1}  \frac{1}{2h} {q}(D_{i}\mid X_{i}) dD_{i}.
\end{equation*}
By setting $z=\frac{D_{i}-\eval}{h}$, we simplify the above expression as,
\begin{equation}
%\begin{aligned}
\int_{\left|z\right|\leq 1}  \frac{1}{2}{q}(\eval+hz\mid X_{i}) dz=\frac{1}{2}\int_{\left|z\right|\leq 1} \left[{q}(\eval\mid X_{i})+hz {q}'(\eval\mid X_{i})+\frac{h^2z^2}{2}{q}''(\eval+c(z)hz\mid X_{i})\right]dz\\
%&=\frac{2}{3}h^2{q}(\eval\mid X_{i})+O(h^4)\cdot \max_{\left|z\right|\leq 1}{q}''(x_{h,z}\mid X_{i})
%\end{aligned}
\label{eq: bound constant term}
\end{equation}
for some $c(z)\in (0,1)$. We shall use $c(z)$ as a generic function of $z$ throughout the proof and the specific function $c(z)$ can vary from place to place.
Hence, we have 
\begin{equation}
\left|\E\left[K_{h}(D_{i})\mid X_{i}\right]-{q}(\eval\mid X_{i})\right|\leq \frac{1}{6}h^2\max_{|c|\leq 1}{q}''(\eval+c h\mid X_{i}).
\label{eq: zero order bound}
\end{equation}

%where 
%\begin{equation}
%x_{h,z}=\eval+c hz
%\label{eq: center point}
%\end{equation}
%and 
%\begin{equation*}
%{q}''(x\mid X_{i})=\left(\frac{\left(x-X_{i}^{\intercal}\gamma\right)^2}{\sigma_2^2}-1\right)\phi\left(\frac{x-X_{i}^{\intercal}\gamma}{\sigma_2}\right).
%%\label{eq: second derivative}
%\end{equation*}
By Condition (A2), we have
\begin{equation}
\max_{\left|c\right|\leq 1}\left|\frac{{q}''(\eval+ch\mid X_{i})}{{q}(\eval\mid X_{i})} \right|\cdot {\1}_{\mathcal{A}_0}\leq C_2(n)\label{eq: 2nd derivative bound}
\end{equation}
where $C_2(n)$ is defined in \eqref{eq: density ratio}. Together with \eqref{eq: zero order bound}, we establish \eqref{eq: conditional zero order}.  
%Then we have 
%\begin{equation}
%\left|\frac{\E\left(K_{h}(D_{i})\mid X_{i}\right)}{{q}(\eval\mid X_{i})}-1\right|\cdot {\1}_{\mathcal{A}_0}\leq \frac{h^2}{6}C_2(n)%1_{\mathcal{A}_{3,i}^{c}}+h^2\left(1+C_{u}^2\right)\exp\left(C_{u}\cdot \frac{h}{\sigma_2}\right).
%\label{eq: key bound a}
%\end{equation}

We apply the boundedness of $\phi(\delta)$ and establish $\E\left[K_{h}(D_{i})\mid X_{i}\right]\leq C$ for some positive constant $C>0$. Together with \eqref{eq: conditional zero order}, we establish \eqref{eq: expectation 0}.
\subsubsection{Proof of \eqref{eq: first moment bound final} and \eqref{eq: second order bound final}}

We first prove \eqref{eq: second order bound final} by analyzing the term $\E\left[(D_{i}-\eval)^2K_{h}(D_{i})\mid X_{i}\right]$. Similar to \eqref{eq: bound constant term}, we write down the following explicit expression,
\begin{equation*}
\E\left[(D_{i}-\eval)^2K_{h}(D_{i})\mid X_{i}\right]=\int_{\left|\frac{D_{i}-\eval}{h}\right|\leq 1} \left[D_{i}-\eval\right]^2 \frac{1}{2h} {q}(D_{i}\mid X_{i}) dD_{i}
\end{equation*}
By setting $z=\frac{D_{i}-\eval}{h}$, we have 
\begin{equation}
\begin{aligned}
&\E\left[(D_{i}-\eval)^2K_{h}(D_{i})\mid X_{i}\right]=\int_{\left|z\right|\leq 1} \frac{1}{2}h^2 z^2  {q}(\eval+hz\mid X_{i}) dz\\
&=\int_{\left|z\right|\leq 1} \frac{1}{2}h^2 z^2  \left[{q}(\eval\mid X_{i})+hz {q}'(\eval\mid X_{i})+\frac{h^2z^2}{2}{q}''(\eval+c(z)hz\mid X_{i})\right]dz.
%=\frac{2}{3}h^2{q}(\eval\mid X_{i})+O(h^4)\cdot \max_{\left|z\right|\leq 1}{q}''(x_{h,z}\mid X_{i})
\end{aligned}
\label{eq: bound quadratic term}
\end{equation}
Hence, we have 
\begin{equation}
\left|\E\left[(D_{i}-\eval)^2K_{h}(D_{i})\mid X_{i}\right]-\frac{1}{3}h^2{q}(\eval\mid X_{i})\right|\leq \frac{1}{10}h^4\max_{|c|\leq 1}{q}''(\eval+c h\mid X_{i}).
\label{eq: second order bound}
\end{equation}
Then we have 
\begin{equation}
\left|\frac{\E\left[(D_{i}-\eval)^2K_{h}(D_{i})\mid X_{i}\right]}{\frac{1}{3}h^2{q}(\eval\mid X_{i})}-1\right|\cdot {\1}_{\mathcal{A}_0}\leq \frac{1}{10}h^2 C_2(n).
\label{eq: second order bound}
\end{equation}
Together with $(D_{i}-\eval)^2K_{h}(D_{i})\leq h,$ we establish \eqref{eq: second order bound final}. 

% \begin{equation}
%x_{h,z}=\eval+c(z)hz
%\label{eq: center point}
%\end{equation}
%and 
%\begin{equation}
%{q}''(x_{h,z}\mid X_{i})=\left(\frac{\left(x_{h,z}-X_{i}^{\intercal}\gamma\right)^2}{\sigma_2^2}-1\right)\phi\left(\frac{x_{h,z}-X_{i}^{\intercal}\gamma}{\sigma_2}\right).
%\label{eq: second derivative}
%\end{equation}

We now control \eqref{eq: first moment bound final}.
Similar to \eqref{eq: bound constant term}, we write down the following explicit expression,
\begin{equation*}
\E\left[(D_{i}-\eval)K_{h}(D_{i})\mid X_{i}\right]=\int_{\left|\frac{D_{i}-\eval}{h}\right|\leq 1} \left[D_{i}-\eval\right] \frac{1}{2h} {q}(D_{i}\mid X_{i}) dD_{i}.
\end{equation*}
Then we have 
\begin{equation*}
\begin{aligned}
&\E\left[(D_{i}-\eval)K_{h}(D_{i})\mid X_{i}\right]=\int_{\left|z\right|\leq 1} \frac{h z}{2}  {q}(\eval+hz\mid X_{i}) dz\\
&=\int_{\left|z\right|\leq 1} \frac{h z}{2}   \left[{q}(\eval\mid X_{i})+hz {q}'(\eval\mid X_{i})+\frac{h^2z^2}{2}{q}''(\eval+c(z)hz\mid X_{i})\right]dz. %+ \frac{h^3z^3}{6}{q}'''(\eval+c(z)hz\mid X_{i})\right]dz\\
%&=\frac{2}{3}h^2 {q}'(\eval\mid X_{i})+O(h^4)\cdot \max_{\left|z\right|\leq 1}{q}'''(x_{h,z}\mid X_{i})
\end{aligned}
\label{eq: bound linear term}
\end{equation*}
%Hence, we have 
%\begin{equation}
%\left|\E\left[(D_{i}-\eval)K_{h}(D_{i})\mid X_{i}\right]-\frac{1}{3}h^2 {q}'(\eval\mid X_{i})\right|\leq \frac{1}{8}h^3\max_{|c|\leq 1}{q}''(\eval+c h\mid X_{i})
%\label{eq: first order bound}
%\end{equation}
%where 
%\begin{equation}
%{q}'(\eval\mid X_{i})=-\frac{x_{0}-X_{i}^{\intercal}\gamma}{\sigma_2}{q}(\eval\mid X_{i})
%\label{eq: first derivative}
%\end{equation}
%and 
%\begin{equation}
%{q}'''(x\mid X_{i})=\frac{x-X_{i}^{\intercal}\gamma}{\sigma_2}\left(3-\frac{\left(x-X_{i}^{\intercal}\gamma\right)^2}{\sigma_2^2}\right)\phi\left(\frac{x-X_{i}^{\intercal}\gamma}{\sigma_2}\right).
%\label{eq: third derivative}
%\end{equation}
%%\begin{equation}
%%\E\left[K_{h}(D_{i})\mid X_{i}\right]=\int_{\left|z\right|\leq 1}  {q}(\eval+hz\mid X_{i}) dz
%%\end{equation}
%%
%With a similar argument as \eqref{eq: 2nd derivative bound}, we can establish
%\begin{equation*}
%\left|\frac{{q}'(\eval\mid X_{i}) \cdot \mathbf{1}_{\mathcal{A}_{0}}}{{q}(\eval\mid X_{i})} \right|\leq C_{u}.
%\end{equation*}
%and 
%\begin{equation}
%\left|\frac{\max_{\left|c\right|\leq 1}{q}'''(\eval+ch\mid X_{i}) \cdot \mathbf{1}_{\mathcal{A}_{0}}}{{q}(\eval\mid X_{i})} \right|\lesssim C_{u}\left(1+C_{u}^2\right)\exp\left(C_{u}\cdot \frac{h}{\sigma_2}\right)
%\label{eq: 3rd derivative bound}
%\end{equation}
Hence, we have
\begin{equation}
\frac{\E\left[(D_{i}-\eval)K_{h}(D_{i})\mid X_{i}\right]}{{q}(\eval\mid X_{i})} \cdot {\1}_{\mathcal{A}_{0}} \leq \frac{1}{3}h^2 (C_1(n)+\frac{3}{8}hC_2(n)).
\label{eq: first moment bound}
\end{equation}
Hence, together with  $\left|(D_{i}-\eval)K_{h}(D_{i})\right|\leq 1,$ we have \eqref{eq: first moment bound final}.

\subsubsection{Proof of \eqref{eq: main part} and \eqref{eq: expectation 1}}
By the iterated expectation, we have
\begin{equation*}
\begin{aligned}
&\E\left({W}^2_{i}K^2_{h}(D_{i})\right)=\E\left({W}^2_{i}K^2_{h}(D_{i})\cdot \mathbf{1}_{\mathcal{A}_{0}}\right)+\E\left({W}^2_{i}K^2_{h}(D_{i}) \cdot \mathbf{1}_{\mathcal{A}_{0}^c}\right)\\
&=\E\left[\E\left({W}^2_{i}K^2_{h}(D_{i})\mid X_{i}\right)\mathbf{1}_{\mathcal{A}_{0}}\right]++\E\left({W}^2_{i}K^2_{h}(D_{i}) \cdot \mathbf{1}_{\mathcal{A}_{0}^c}\right).
\end{aligned}
\end{equation*}
We first analyze $\E\left({W}^2_{i}K^2_{h}(D_{i})\mid X_{i}\right)\mathbf{1}_{\mathcal{A}_{0}}$, by noting that 
\begin{equation}
\begin{aligned}
&\E\left({W}^2_{i}K^2_{h}(D_{i})\mid X_{i}\right)=\frac{1}{h}\E\left({W}^2_{i}K_{h}(D_{i})\mid X_{i}\right)\\
&=\frac{1}{h}\left(\E\left[(D_{i}-\eval)^2K_{h}(D_{i})\mid X_{i}\right]-\frac{\left\{\E\left[(D_{i}-\eval)K_{h}(D_{i})\mid X_{i}\right]\right\}^2}{\E\left[K_{h}(D_{i})\mid X_{i}\right]}\right),
\end{aligned}
\label{eq: explicit expression}
\end{equation}
where the last equality follows from the definition of $W_{i}$.
%In the following, we provide upper bounds for $\E\left[(D_{i}-\eval)K_{h}(D_{i})\mid X_{i}\right]$ and $\E\left[(D_{i}-\eval)^2K_{h}(D_{i})\mid X_{i}\right]$. 

Note that 
\begin{equation*}
\begin{aligned}
&\frac{1}{{\frac{1}{3}h^2{q}(\eval\mid X_{i})}}\cdot \frac{\left\{\E\left[(D_{i}-\eval)K_{h}(D_{i})\mid X_{i}\right]\right\}^2}{\E\left[K_{h}(D_{i})\mid X_{i}\right]}\cdot {\1}_{\mathcal{A}_{0}}\\
&=\frac{3}{h^2}\frac{\left\{\E\left[(D_{i}-\eval)K_{h}(D_{i})\mid X_{i}\right]\right\}^2}{{q}^2(\eval\mid X_{i})}\frac{{q}(\eval\mid X_{i})}{\E\left[K_{h}(D_{i})\mid X_{i}\right]}\cdot {\1}_{\mathcal{A}_{0}}
\end{aligned}
\end{equation*}
By the above expression, \eqref{eq: conditional zero order}, and \eqref{eq: first moment bound}, we establish 
\begin{equation}
\frac{1}{{\frac{1}{3}h^2{q}(\eval\mid X_{i})}}\cdot \frac{\left\{\E\left[(D_{i}-\eval)K_{h}(D_{i})\mid X_{i}\right]\right\}^2}{\E\left[K_{h}(D_{i})\mid X_{i}\right]}\cdot {\1}_{\mathcal{A}_{0}}\leq \frac{h^2}{3} \frac{[(C_1(n)+\frac{3}{8}hC_2(n))]^2}{1-\frac{h^2}{6} C_2(n)}.
\label{eq: approximation part 0}
\end{equation} 

We apply \eqref{eq: explicit expression} together with \eqref{eq: second order bound} and \eqref{eq: approximation part 0} and establish \eqref{eq: main part}.
Since $|W_{i}|K_{h}(D_{i})\leq 1$, we have 
\begin{equation}
\left|\frac{\E\left({W}^2_{i}K^2_{h}(D_{i})\mathbf{1}_{\mathcal{A}_{0}^c}\right)}{\frac{1}{3}h\pi(\eval)}\right|\leq \frac{3\mathbf{P}(\mathcal{A}_{0}^c)}{h\pi(\eval)}.
\label{eq: approximation part}
\end{equation}
Combining \eqref{eq: main part} and \eqref{eq: approximation part}, we establish \eqref{eq: expectation 1}.\\
%\begin{equation}
%\E\left({W}^2_{i}K^2_{h}(D_{i})\mid X_{i}\right)\approx\frac{2}{3}h{q}(\eval\mid X_{i})
%\end{equation}
\subsubsection{Proof of \eqref{eq: expectation 2}}
The proof of \eqref{eq: expectation 2} is similar to that of \eqref{eq: expectation 1}.
We first have the following decomposition, 
\begin{equation*}
\E W_i(D_{i}-\eval)K_{h}(D_{i})=\E W_i(D_{i}-\eval)K_{h}(D_{i})\cdot \mathbf{1}_{\mathcal{A}_{0}}+\E W_i(D_{i}-\eval)K_{h}(D_{i}) \cdot \mathbf{1}_{\mathcal{A}_{0}^c}.
\end{equation*}

%By applying  \eqref{eq: explicit expression} and \eqref{eq: key bound a}, we have
%\begin{equation}
%\frac{\E \left[W_i(D_{i}-\eval)K_{h}(D_{i})\mathbf{1}_{\mathcal{A}_{0}}\right]}{{\frac{2}{3}h^2{\pi}(\eval)}}\rightarrow 1
%\label{eq: inter 3}
%W_i(D_{i}-\eval)K_{h}(D_{i})
%\end{equation}
Since $W_i(D_{i}-\eval)K_{h}(D_{i})\leq h$, we have 
\begin{equation}
\left|\frac{\E W_i(D_{i}-\eval)K_{h}(D_{i})\mathbf{1}_{\mathcal{A}^c_0}}{{\frac{2}{3}h^2\pi(\eval)}}\right|\leq \frac{\mathbf{P}(\mathcal{A}_{0}^{c})}{\frac{2}{3}h\pi(\eval)}.
\label{eq: approximation part 2}
\end{equation}
Note that
\begin{equation*}
\begin{aligned}
&\E \left[W_i(D_{i}-\eval)K_{h}(D_{i})\mid X_{i}\right]\mathbf{1}_{\mathcal{A}_{0}}\\
&=\left(\E\left[(D_{i}-\eval)^2K_{h}(D_{i})\mid X_{i}\right]-\frac{\left\{\E\left[(D_{i}-\eval)K_{h}(D_{i})\mid X_{i}\right]\right\}^2}{\E\left[K_{h}(D_{i})\mid X_{i}\right]}\right)\mathbf{1}_{\mathcal{A}_{0}}\\
&=\E\left({W}^2_{i}K_{h}(D_{i})\mid X_{i}\right)\mathbf{1}_{\mathcal{A}_{0}}.
\end{aligned}
\end{equation*}
We apply \eqref{eq: main part}  and \eqref{eq: approximation part 2} to establish \eqref{eq: expectation 2}.

\subsubsection{Proof of \eqref{eq: third moment}}

Note that 
\begin{equation*}
\E W_{i} \frac{(D_{i}-\eval)^2}{2} K_{h}(D_{i})=\E W_{i} \frac{(D_{i}-\eval)^2}{2} K_{h}(D_{i}) \cdot \mathbf{1}_{\mathcal{A}_{0}}+\E W_{i} \frac{(D_{i}-\eval)^2}{2} K_{h}(D_{i}) \cdot \mathbf{1}_{\mathcal{A}_{0}^{c}}.
\end{equation*}
For the first term, we apply the iterated expectation and obtain $$\E W_{i} \frac{(D_{i}-\eval)^2}{2} K_{h}(D_{i})\mathbf{1}_{\mathcal{A}_{0}}=\E\left[\E\left( W_{i} \frac{(D_{i}-\eval)^2}{2} K_{h}(D_{i})\mid X_{i}\right)\mathbf{1}_{\mathcal{A}_{0}}\right],$$
with 
\begin{equation*}
\begin{aligned}
&\E\left( W_{i} \frac{(D_{i}-\eval)^2}{2} K_{h}(D_{i})\mid X_{i}\right)\\
=&\E\left(\frac{(D_{i}-\eval)^3}{2} K_{h}(D_{i})\mid X_{i}\right)-l(X_{i})\E\left( \frac{(D_{i}-\eval)^2}{2} K_{h}(D_{i})\mid X_{i}\right)\\
=&\E\left(\frac{(D_{i}-\eval)^3}{2} K_{h}(D_{i})\mid X_{i}\right)-\frac{\E\left((D_{i}-\eval) K_{h}(D_{i})\mid X_{i}\right)\E\left( \frac{(D_{i}-\eval)^2}{2} K_{h}(D_{i})\mid X_{i}\right)}{\E\left(K_{h}(D_{i})\mid X_{i}\right)}.
\end{aligned}
\end{equation*}
Then it is sufficient to control the terms $$\E\left[(D_{i}-\eval)^3K_{h}(D_{i})\mid X_{i}\right] \mathbf{1}_{\mathcal{A}_{0}},$$ and 
\begin{equation}
\frac{\E\left((D_{i}-\eval) K_{h}(D_{i})\mid X_{i}\right)\E\left( \frac{(D_{i}-\eval)^2}{2} K_{h}(D_{i})\mid X_{i}\right)}{\E\left(K_{h}(D_{i})\mid X_{i}\right)} \mathbf{1}_{\mathcal{A}_{0}}.
\label{eq: bounding moment}
\end{equation}
Since $$\E\left( \frac{(D_{i}-\eval)^2}{2} K_{h}(D_{i})\mid X_{i}\right)\leq \frac{h^2}{2}\E\left(K_{h}(D_{i})\mid X_{i}\right),$$
the term in \eqref{eq: bounding moment} can be upper bounded by $$\frac{h^2}{2}\E\left((D_{i}-\eval) K_{h}(D_{i})\mid X_{i}\right).$$ It follows from \eqref{eq: first moment bound} that 
\begin{equation}
\frac{h^2}{2}\E\left((D_{i}-\eval) K_{h}(D_{i})\mid X_{i}\right) \1_{\mathcal{A}_{0}}\lesssim q(\eval\mid X_i) \cdot h^4 (C_1(n)+\frac{3}{8}hC_2(n)).
\label{eq: approximation part 3}
\end{equation}
 We control the term $\E\left[(D_{i}-\eval)^3K_{h}(D_{i})\mid X_{i}\right]$ in the following.
\begin{equation*}
\begin{aligned}
&\E\left[(D_{i}-\eval)^3K_{h}(D_{i})\mid X_{i}\right]=\int_{\left|z\right|\leq 1} \frac{1}{2}h^3 z^3  {q}(\eval+hz\mid X_{i}) dz\\
&=\int_{\left|z\right|\leq 1} \frac{1}{2}h^3 z^3  \left[{q}(\eval\mid X_{i})+hz {q}'(\eval\mid X_{i})+\frac{h^2z^2}{2}{q}''(\eval+c(z)\mid X_{i})\right]dz,
%&=\frac{2}{5}h^4 {q}'(\eval\mid X_{i})+O(h^6)\cdot \max_{\left|z\right|\leq 1}{q}'''(\eval+chz\mid X_{i}).
\end{aligned}
\end{equation*}
%and then we have 
%\begin{equation*}
%\left|\E\left[(D_{i}-\eval)^3K_{h}(D_{i})\mid X_{i}\right]-\frac{1}{5}h^4 {q}'(\eval\mid X_{i})\right|\leq \frac{1}{12}h^5\max_{\left|c\right|\leq 1}{q}''(\eval+ch\mid X_{i}).
%\end{equation*}
%where 
%\begin{equation}
%{q}'(\eval\mid X_{i})=-\frac{x_{0}-X_{i}^{\intercal}\gamma}{\sigma_2}{q}(\eval\mid X_{i})
%\end{equation}
%\begin{equation}
%{q}'''(\eval+chz\mid X_{i})=\frac{\eval+chz-X_{i}^{\intercal}\gamma}{\sigma_2}\left(3-\frac{\left(\eval+chz-X_{i}^{\intercal}\gamma\right)^2}{\sigma_2^2}\right)\phi\left(\frac{\eval+chz-X_{i}^{\intercal}\gamma}{\sigma_2}\right),
%\end{equation}
and then have 
\begin{equation*}
\left|\frac{\E\left[(D_{i}-\eval)^3K_{h}(D_{i})\mid X_{i}\right]}{q(\eval\mid X_i)}\right|\cdot \1_{\mathcal{A}_{0}}\leq \frac{h^4}{5}\left[C_1(n)+\frac{5}{12}hC_2(n)\right]\pi(\eval).%=\frac{2}{5}h^4 {q}(\eval\mid X_{i})\cdot \frac{-x_{0}+X_{i}^{\intercal}\gamma}{\sigma_2}+O(h^6)\cdot \max_{\left|c\right|\leq 1}{q}'''(\eval+ch\mid X_{i})
\end{equation*}
%By \eqref{eq: third derivative}, we have 
%\begin{equation*}
%\left|\E\left( W_{i} \frac{(D_{i}-\eval)^2}{2} K_{h}(D_{i})\mid X_{i}\right)\mathbf{1}_{\mathcal{A}_{0}}\right|\leq \left(C_{u}\frac{2}{5}h^4+O(h^6) C_{u}\left(3+C_{u}^2\right)\right){q}(\eval\mid X_{i})
%\end{equation*}
%where $C_{u}$ is defined in \eqref{eq: second key growing constant}.
Together with \eqref{eq: approximation part 3} and
$$\left|\E W_{i} \frac{(D_{i}-\eval)^2}{2} K_{h}(D_{i}) \cdot \mathbf{1}_{\mathcal{A}_{0}^{c}}\right|\leq h^2 \mathbf{P}(\mathcal{A}_{0}^{c}),$$
%and $$\left|\E W_{i} \frac{(D_{i}-\eval)^2}{2} K_{h}(D_{i}) \cdot \mathbf{1}_{\mathcal{A}_{0}^{c}}\right|\leq h^2 \mathbf{P}(\mathcal{A}_{0}^{c})=h^2\cdot n^{-c},$$
we establish \eqref{eq: third moment}.

%%%%%%%%%%%%%%%%%%%%%%%%%%%%%%%%%%%%%%%%%%%%%%%%%%%%%%%%%%
\subsection{Proof of Lemma \ref{lem: concentration lemma}}
\label{sec: con lemma}
%%%%%%%%%%%%%%%%%%%%%%%%%%%%%%%%%%%%%%%%%%%%%%%%%%%%%%%%%%
The proofs rely on the Bernstein inequality \cite{bennett1962probability}, which is restated in the following lemma.
\begin{Lemma}
Suppose that $\{H_i\}_{1\leq i\leq n}$ are independent zero mean random variables and $|H_i|\leq M$ almost surely. Then we have 
$$\mathbf{P}\left(\left|\sum_{i=1}^{n} H_i \right|\geq T\right)\leq 2\exp\left(-\frac{T^2/2}{\sum_{i=1}^{n} \E H_i^2+ M T/3}\right).$$
\label{lem: bernstein}
\end{Lemma}
\underline{Proof of \eqref{eq: prob bound 1}}\\
%The term $\frac{1}{n}\sum_{i=1}^{n} K_{h}(D_{i})$ satisfies
%\begin{equation}
%\E\left(\frac{1}{n}\sum_{i=1}^{n} K_{h}(D_{i})\right)=\E \left(K_{h}(D_{i})\right), \; {\rm Var}\left(\frac{1}{n}\sum_{i=1}^{n} K_{h}(D_{i})\right)\leq \frac{1}{nh}\E \left(K_{h}(D_{i})\right)
%\label{eq: mv 1}
%\end{equation}
We shall apply Lemma \ref{lem: bernstein} by taking  $H_i=K_{h}(D_{i})-\E \left(K_{h}(D_{i})\right).$
By \eqref{eq: expectation 0}, there exists $0<c<1/2$ such that 
\begin{equation}
(2-c)\pi(\eval)\leq \E \left(K_{h}(D_{i})\right)\leq (2+c) \pi(\eval). %\quad \text{for} \; c>0, C>2.
\label{eq: upper 1}
\end{equation}
Note that $\left|K_{h}(D_{i})-\E \left(K_{h}(D_{i})\right)\right| \leq 1/h$ and 
$$\E \left(K^2_{h}(D_{i})\right)=\E \left(K_{h}(D_{i})\right)/h\leq (2+c) \pi(\eval)/h.$$
By \eqref{eq: upper 1} and Lemma \ref{lem: bernstein} with $T=t\cdot \max\left\{\sqrt{\frac{n \pi(\eval)}{h}},\frac{1}{h}\right\}$, we establish \eqref{eq: prob bound 1}. \\
\noindent \underline{Proof of \eqref{eq: prob bound 2}}\\
By the definition of $W_i$, the term $\frac{1}{n}\sum_{i=1}^{n} W_{i}K_{h}(D_{i})$ satisfies 
\begin{equation*}
\E\left(\frac{1}{n}\sum_{i=1}^{n} W_{i}K_{h}(D_{i})\right)= 0. 
\label{eq: mv 2}
\end{equation*}
Note that $|W_{i}K_{h}(D_{i})|\leq 2.$ By \eqref{eq: expectation 1}, we apply Lemma \ref{lem: bernstein} with $T=t\cdot \max\left\{\sqrt{nh\pi(\eval)},1\right\}$ and establish \eqref{eq: prob bound 2}.\\
%\begin{equation}
% \E \left(W_{i}^2K_{h}^2(D_{i})\right)= \frac{2}{3}h\pi(\eval) (1+o(1)) %\frac{2}{3n} \E\phi\left(\frac{\eval-X_{i}^{\intercal}\gamma}{\sigma_2}\right)
%\label{eq: upper 2}
%\end{equation}
\underline{Proof of \eqref{eq: prob bound 4}}\\
It follows from \eqref{eq: first moment bound final} that 
\begin{equation*}
\E\left(\frac{1}{n}\sum_{i=1}^{n}{(D_{i}-\eval)} K_{h}(D_{i})\right)=\E{(D_{i}-\eval)} K_{h}(D_{i})\leq \frac{{\pi}(\eval)}{3}h^2 (C_1(n)+\frac{3}{8}hC_2(n))+{\mathbf{P}(\mathcal{A}_0^c)}.
\end{equation*}

We apply \eqref{eq: second order bound final} and establish 
\begin{equation*}
{\rm Var}\left({(D_{i}-\eval)} K_{h}(D_{i})\right)\leq \frac{1}{h}\E{(D_{i}-\eval)}^2 K_{h}(D_{i})\lesssim \frac{1}{3}h\pi(\eval).
\end{equation*}
Note that $\left|{(D_{i}-\eval)} K_{h}(D_{i})\right|\leq 1,$ we apply Lemma \ref{lem: bernstein} with $T=t\cdot \max\left\{\sqrt{nh\pi(\eval)},1\right\}$ and establish \eqref{eq: prob bound 4}.\\ %\begin{equation}
%\left|\frac{\E\left[(D_{i}-\eval)^2K_{h}(D_{i})\right]}{\frac{1}{3}h^2\pi(\eval)}-1\right| \lesssim \frac{1}{10}h^2 C_2(n)+\frac{\mathbf{P}(\mathcal{A}_0^c)}{h\pi(\eval)}.
%\label{eq: second order bound final}
%\end{equation}
%By \eqref{eq: first order bound}, we have 
%\begin{equation}
%\begin{aligned}
%&\left|\E{(D_{i}-\eval)} K_{h}(D_{i})\right|\leq \frac{2}{3}h^2 \E {q}'(\eval\mid X_{i}) \mathbf{1}_{\mathcal{A}_{0}}+\frac{1}{15}h^4\cdot \max_{\left|c\right|\leq 1}{q}'''(\eval+ch\mid X_{i})\mathbf{1}_{\mathcal{A}_{0}}+ \mathbf{P}({A}_{0}^c)\\
%&\lesssim \frac{2}{3}C_{u} h^2 \pi(\eval)+h^4 C_{u}(1+C_{u}^2)\pi(\eval)+\frac{1}{n^{c}}\lesssim C_{u}h^2 \pi(\eval) \\
%\end{aligned}
%\end{equation}
%where the second inequality follows from \eqref{eq: first derivative}, \eqref{eq: third derivative} and \eqref{eq: bounding constant}.
%By \eqref{eq: second order bound}, we have 
%\begin{equation}
%\begin{aligned}
%&\E{(D_{i}-\eval)}^2 K_{h}(D_{i})\leq \frac{2}{3}h^2 \E  {q}(\eval\mid X_{i})\mathbf{1}_{\mathcal{A}_{0}}+ \frac{4}{5}h^4\cdot  \E  \max_{\left|z\right|\leq 1}{q}''(\eval+chz\mid X_{i})\mathbf{1}_{\mathcal{A}_{0}}+ \mathbf{P}({A}_{0}^c)\\
%&\leq \frac{2}{3} h^2 \pi(\eval)+h^4 (1+C_{u}^2)\pi(\eval)+\frac{1}{n^{c}}\lesssim h^2 \pi(\eval) \\
%\end{aligned}
%\end{equation}
%where the second inequality follows from \eqref{eq: third derivative} and \eqref{eq: bounding constant}.
%\\
%%%%%%%%%%%%%%%%%%%%%%%%%%%%%%%%%%%%%%%%%%%%%%%%%%%%%%%%%%%%%%%%%%%%%%%%%%%%%%%%%%%%%
\underline{Proof of \eqref{eq: prob bound 5}}\\
%\begin{equation}
%\begin{aligned}
%&\E \frac{1}{n}\sum_{i=1}^{n}W_i^2K_{h}^2(D_{i})=\E \frac{1}{h}W_i^2K_{h}(D_{i})=\E\left(\E \left[\frac{1}{h}W_i^2K_{h}(D_{i})\mid X_{i}\right]\right)\\
%&=\frac{1}{h}\E\left(\E\left[(D_{i}-\eval)^2K_{h}(D_{i})\mid X_{i}\right]-\frac{\E^2\left[(D_{i}-\eval)K_{h}(D_{i})\mid X_{i}\right]}{\E\left[K_{h}(D_{i})\mid X_{i}\right]}\right)
%\end{aligned}
%\end{equation}
%and hence 
It follows from \eqref{eq: expectation 2} that 
\begin{equation*}
\left|\E \left(W_i(D_{i}-\eval) K_{h}(D_{i})\right)-\frac{2}{3}h^2\pi(\eval)\right|\leq c \frac{2}{3}h^2\pi(\eval),
\label{eq: mv 5a}
\end{equation*}
for some small positive constant $c\in(0,1).$
Note that  
\begin{equation*}
\begin{aligned}
{\rm Var} \left(W_i(D_{i}-\eval)K_{h}(D_{i})\right)\leq \frac{1}{h} \E\left({W}^2_{i}(D_{i}-\eval)^2 K_{h}(D_{i})\right)\lesssim {h^{3}}\E\left(K_{h}(D_{i})\right),
\end{aligned}
\label{eq: mv 5b}
\end{equation*}
and $$\left|W_i(D_{i}-\eval)K_{h}(D_{i})\right|\leq h.$$
We apply Lemma \ref{lem: bernstein} with $T=t\cdot \max\left\{\sqrt{nh^3\pi(\eval)},h\right\}$ and establish \eqref{eq: prob bound 5}.\\
%Combined with \eqref{eq: expectation 2}, we establish that \eqref{eq: prob bound 5}.\\%with probability larger than $1-\frac{1}{t}$, then
%%%%%%%%%%%%%%%%%%%%%%%%%%%%%%%%%%%%%%%%%%%%%%%%%%%%%%%%%%%%%%%%%%%%%%%%%%%%%%%%%%%%%
\underline{Proof of \eqref{eq: prob bound 3}}\\
%\begin{equation}
%\begin{aligned}
%&\E \frac{1}{n}\sum_{i=1}^{n}W_i^2K_{h}^2(D_{i})=\E \frac{1}{h}W_i^2K_{h}(D_{i})=\E\left(\E \left[\frac{1}{h}W_i^2K_{h}(D_{i})\mid X_{i}\right]\right)\\
%&=\frac{1}{h}\E\left(\E\left[(D_{i}-\eval)^2K_{h}(D_{i})\mid X_{i}\right]-\frac{\E^2\left[(D_{i}-\eval)K_{h}(D_{i})\mid X_{i}\right]}{\E\left[K_{h}(D_{i})\mid X_{i}\right]}\right)
%\end{aligned}
%\end{equation}
%and hence 
The term $\frac{1}{n}\sum_{i=1}^{n}W_i\frac{(D_{i}-\eval)^2}{2} K_{h}(D_{i})$ satisfies
\begin{equation*}
\begin{aligned}
\E \left(\frac{1}{n}\sum_{i=1}^{n}W_i\frac{(D_{i}-\eval)^2}{2} K_{h}(D_{i})\right)&= \E\left[ W_{i} \frac{(D_{i}-\eval)^2}{2} K_{h}(D_{i})\right]\\
&\lesssim {h^4}\left[C_1(n)+hC_2(n)\right]+h^2 \mathbf{P}(\mathcal{A}_{0}^{c}),
\end{aligned}
\label{eq: mv 3a}
\end{equation*}
where the last inequality follows from \eqref{eq: third moment}. 
Note that 
\begin{equation*}
\begin{aligned}
{\rm Var} \left(W_i\frac{(D_{i}-\eval)^2}{2} K_{h}(D_{i})\right)\leq \frac{1}{h} \E\left({W}^2_{i}\frac{(D_{i}-\eval)^4}{4} K_{h}(D_{i})\right)\leq \frac{h^{5}}{4}\E\left(K_{h}(D_{i})\right),
\end{aligned}
\label{eq: mv 3b}
\end{equation*}
%In the following, we shall prove that 
%\begin{equation}
%\left|\E W_{i}\frac{(D_{i}-\eval)^2}{2} K_{h}(D_{i})\right|\leq {h^4}\left[C_1(n)+hC_2(n)\right]+h^2 \mathbf{P}(\mathcal{A}_{0}^{c})
%\label{eq: upper 3}
%\end{equation}
and $\left|W_i\frac{(D_{i}-\eval)^2}{2} K_{h}(D_{i})\right|\leq h^2.$ We apply Lemma \ref{lem: bernstein} with $T=t\cdot \max\left\{\sqrt{nh^5\pi(\eval)},h^2\right\}$ and establish \eqref{eq: prob bound 3}.\\
\noindent \underline{Proof of \eqref{eq: prob bound 8}}. Note that both $\left|W_i^2K_{h}^2(D_{i})-\E W_i^2K_{h}^2(D_{i})\right|$ is upper bounded by a constant and $${\rm Var}(W_i^2K_{h}^2(D_{i}))\leq \E W_i^4K_{h}^4(D_{i})\leq \E W_i^2K_{h}^2(D_{i})\lesssim \frac{1}{3}h\pi(\eval),$$ where the second inequality follows from the fact $W_i^2K_{h}^2(D_{i})\leq 1$ and the last inequality follows from \eqref{eq: expectation 1}. We apply Lemma \ref{lem: bernstein} with $T=t\cdot \max\left\{\sqrt{nh\pi(\eval)},1\right\}$ and establish \eqref{eq: prob bound 8}.\\
\subsection{Proof of Lemma \ref{lem: Gaussian bound}}
%%%%%%%%%%%%%%%%%%%%%%%%%%%%%%%%%%%%%%%%%%%%%%%%%%%%%%%

By Taylor's expansion, we have for $c_1,c_2\in (0,1)$
\begin{equation*}
\phi(\mu+\Delta)+\phi(\mu-\Delta)=\phi(\mu)+\Delta \cdot \phi'(\mu)+\frac{\Delta^2}{2} \cdot \phi''(\mu+c_1 \Delta)+\phi(\mu)-\Delta \cdot \phi'(\mu)+\frac{\Delta^2}{2} \cdot \phi''(\mu+c_2 \Delta)
\end{equation*}
and
\begin{equation*}
\frac{1}{\Delta}\int_{\mu-\Delta}^{\mu+\Delta}\phi(t) dt=\frac{1}{\Delta}\int_{\mu-\Delta}^{\mu+\Delta}\left[\phi(\mu)+{(t-\mu)} \phi'(\mu)+\frac{(t-\mu)^2}{2} \phi''(\mu+c_3(t)(t-\mu))\right]dt
\end{equation*}
where $c_3(t)\in (0,1)$.
Hence, we have 
\begin{equation*}
\begin{aligned}
\phi(\mu+\Delta)+\phi(\mu-\Delta)-\frac{1}{\Delta}\int_{\mu-\Delta}^{\mu+\Delta}\phi(t) dt&=\frac{\Delta^2}{2} \cdot \left(\phi''(\mu+c_1 \Delta)+\phi''(\mu+c_2 \Delta)\right)\\
&+\frac{1}{\Delta}\int_{\mu-\Delta}^{\mu+\Delta}\frac{(t-\mu)^2}{2} \phi''(\mu+c_3(t)(t-\mu))dt
\end{aligned}
\end{equation*}
Hence, we establish \eqref{eq: Gaussian bound 1}.
Note that 
\begin{equation}
\int_{\mu-\Delta}^{\mu+\Delta} \left(t-\mu\right)\phi(t) dt=\int_{\mu-\Delta}^{\mu+\Delta}\left[(t-\mu)\phi(\mu)+{(t-\mu)^2}\phi'(\mu+c_4(t)(t-\mu))\right]dt
\end{equation}
where $c_4(t) \in (0,1)$. Hence  we establish \eqref{eq: Gaussian bound 2}.
Note that 
\begin{equation}
\phi(\mu+\Delta)-\phi(\mu-\Delta)=2\Delta \cdot \phi'(\mu+c_5\Delta),
\end{equation}
for $c_5\in (-1,1)$. Hence we establish \eqref{eq: Gaussian bound 3}.

%%%%%%%%%%%%%%%%%%%%%%%%%%%%%%%%%%%%%%%%%%%%%%%%%%%%%%%%%%%%
\subsection{Proof of Proposition \ref{prop: consistent sigma}}
\label{sec: consistent sigma}
%%%%%%%%%%%%%%%%%%%%%%%%%%%%%%%%%%%%%%%%%%%%%%%%%%%%%%%%%%%%
Note that 
\begin{equation*}
\begin{aligned}
\widehat{\sigma}^2-\sigma^2= &\frac{1}{n}\sum_{i \in \mathcal{I}_a} \left[\left(\epsilon_i-[\widehat{f}^b(D_i)-f(D_i)]-[\widehat{g}^b(X_{i})-{g}(X_{i})]\right)^2-\sigma^2\right]\\
+&\frac{1}{n}\sum_{i \in \mathcal{I}_b} \left[\left(\epsilon_i-[\widehat{f}^a(D_i)-f(D_i)]-[\widehat{g}^a(X_{i})-{g}(X_{i})]\right)^2-\sigma^2\right].
\end{aligned}
\end{equation*}
It is sufficient to show that 
\begin{equation}
\frac{1}{|\mathcal{I}_a|}\sum_{i \in \mathcal{I}_a} \left[\left(\epsilon_i-[\widehat{f}^b(D_i)-f(D_i)]-[\widehat{g}^b(X_{i})-{g}(X_{i})]\right)^2-\sigma^2\right]\cip 0.
\label{eq: sufficient}
\end{equation}
For the last hand side of \eqref{eq: sufficient}, we have the decomposition 
\begin{equation}
\begin{aligned}
&\frac{1}{|\mathcal{I}_a|}\sum_{i \in \mathcal{I}_a}  (\epsilon_i^2-\sigma^2)+\frac{1}{|\mathcal{I}_a|}\sum_{i \in \mathcal{I}_a}  \left([\widehat{f}^b(D_i)-f(D_i)]+[\widehat{g}^b(X_{i})-{g}(X_{i})]\right)^2\\
&-\frac{2}{|\mathcal{I}_a|}\sum_{i \in \mathcal{I}_a} \epsilon_i\cdot \left([\widehat{f}^b(D_i)-f(D_i)]+[\widehat{g}^b(X_{i})-{g}(X_{i})]\right) .
\end{aligned}
\label{eq: decomposition noise}
\end{equation}
By the law of large numbers, we have 
\begin{equation}\frac{1}{|\mathcal{I}_a|}\sum_{i \in \mathcal{I}_a}  (\epsilon_i^2-\sigma^2)\cip 0.
\label{eq: law large}
\end{equation}
Define the event 
$$\mathcal{A}'_2=\left\{\E_{X_{*}}(\widehat{g}^b(X_{*})-g(X_{*}))^2\lesssim {\rm Err}^2(\widehat{g}),\;\; \E_{D_{*}}(\widehat{f}^b(D_{*})-f(D_{*}))^2\lesssim {\rm Err}^2(\widehat{f})\right\}
$$
and by the definition of ${\rm Err}(\widehat{f})$ and ${\rm Err}(\widehat{g}),$ we have 
\begin{equation}
\mathbf{P}(\mathcal{A}'_2)\geq 1-\min\{n,p\}^{-c}.
\label{eq: high-prob 2}
\end{equation}
In the following analysis, we condition on the data in $\mathcal{I}_b$ and take the conditional expectation as 
\begin{equation}
\begin{aligned}
\E \left[\frac{1}{|\mathcal{I}_a|}\sum_{i \in \mathcal{I}_a}  \left([\widehat{f}^b(D_i)-f(D_i)]+[\widehat{g}^b(X_{i})-{g}(X_{i})]\right)^2\cdot {\1}_{\mathcal{A}'_2}\mid \mathcal{I}_b\right]\lesssim {\rm Err}^2(\widehat{f})+{\rm Err}^2(\widehat{g}).
\end{aligned}
\label{eq: upper temp 1}
\end{equation}
By the Cauchy inequality, we have 
\begin{equation*}
\begin{aligned}
&\E\left[\left|\frac{1}{|\mathcal{I}_a|}\sum_{i \in \mathcal{I}_a} \epsilon_i\cdot \left([\widehat{f}^b(D_i)-f(D_i)]+[\widehat{g}^b(X_{i})-{g}(X_{i})]\right)\cdot {\1}_{\mathcal{A}'_2} \right|^2\mid \mathcal{I}_b\right]\\
&\leq \E \left[\left(\frac{1}{|\mathcal{I}_a|}\sum_{i \in \mathcal{I}_a} \epsilon^2_i\right)\cdot \left(\frac{1}{|\mathcal{I}_a|}\sum_{i \in \mathcal{I}_a}  \left([\widehat{f}^b(D_i)-f(D_i)]+[\widehat{g}^b(X_{i})-{g}(X_{i})]\right)^2\cdot {\1}_{\mathcal{A}'_2}\right) \mid \mathcal{I}_b\right]\\
&\lesssim \sigma^2 \cdot \left({\rm Err}^2(\widehat{f})+{\rm Err}^2(\widehat{g})\right),
\end{aligned}
\end{equation*}
where the least inequality follows from \eqref{eq: upper temp 1} and $\E(\epsilon_i^2\mid D_i, X_i)=\sigma^2.$ By the Markov inequality, we establish that, with probability larger than $1-\frac{1}{t}$ for some $t>1,$
\begin{equation*}
\begin{aligned}
&\left|\frac{1}{|\mathcal{I}_a|}\sum_{i \in \mathcal{I}_a}  \left([\widehat{f}^b(D_i)-f(D_i)]+[\widehat{g}^b(X_{i})-{g}(X_{i})]\right)^2\right|\cdot {\1}_{\mathcal{A}'_2}\\
&+\left|\frac{2}{|\mathcal{I}_a|}\sum_{i \in \mathcal{I}_a} \epsilon_i\cdot \left([\widehat{f}^b(D_i)-f(D_i)]+[\widehat{g}^b(X_{i})-{g}(X_{i})]\right)\right|\cdot {\1}_{\mathcal{A}'_2}\\
&\lesssim  t\left({\rm Err}^2(\widehat{f})+{\rm Err}^2(\widehat{g})\right)+\sqrt{t\left({\rm Err}^2(\widehat{f})+{\rm Err}^2(\widehat{g})\right)}.
\end{aligned}
\end{equation*}
Combined with \eqref{eq: law large} and \eqref{eq: high-prob 2}, we establish $\widehat{\sigma}^2\cip \sigma^2$ if $\max\{{\rm Err}(\widehat{f}),{\rm Err}(\widehat{g})\}\rightarrow 0.$

%%%%%%%%%%%%%%%%%%%%%%%%%%%%%%
\section{Additional Simulation Results}
%%%%%%%%%%%%%%%%%%%%%%%%%%%%%%
%%%%%%%%%%%%%%%%% Complete Results
\subsection{Setting 1-4 and Nonlinear Treatment Model}
\label{sec: 1to4 and nonlinear}
In this section, we present complete simulation results for Setting 1-4 and the nonlinear treatment model. The sample sizes are varied across \{500,1000,1500,2000\} and $\eval$ is varied across \{-1.25,-0.5,0.1,0.25,1\}. We consider two generating models for $f$ and $g_1$ as follows, 
\begin{itemize}
\item $f(d) = 2\exp(-d/2)$ and $g_1(x) = 1.5\sin(x)$;
\item $f(d) = 1.5\sin(d)$ and $g_1(x) = 2\exp(-x/2)$.
\end{itemize}
The complete results for Setting 1 are summarised in Table \hyperref[table:S1]{S1} and Table \hyperref[table:S2]{S2}. Similar to the results presented in the main paper, our \texttt{DLL} method achieves desired coverage, and the CI length is close to the confidence interval by the oracle estimator. Besides, our \texttt{DLL} method outperforms the plug-in method in terms of coverage since the \texttt{DLL} estimator has a smaller bias. The coverage for the plug-in estimator is relatively good at the boundary \{-1.25,1\} since only a few samples are used with the chosen bandwidth, and the standard error for the plug-in estimator is large, leading to a wide CI. %We can see from the tables that the above statement is true for $f$ being $Exp$ function as well. 

The summarized results for Setting 2 and 3 are in Table \hyperref[table:S3]{S3} and Table \hyperref[table:S4]{S4}, respectively. The results for Setting 2 and Setting 3 are similar to those for Setting 1. For Setting 4, $(D_i,X^{\intercal}_i)^{\intercal}$ follows a t distribution and we vary the degree of freedom in \{10,15\}. The results are reported in Tables \hyperref[table:S5]{S5} and  \hyperref[table:S6]{S6}. In Settings 3 and 4, we test the robustness of our proposed method to the violation of assumption in (A2). The results demonstrate that our proposed \texttt{DLL} method still corrects the bias of the plug-in estimator and attains the desired coverage level, with the CI length similar to the confidence interval by the oracle estimator.

The results of the non-linear treatment model are presented in Table \hyperref[table:S7]{S7}. We see that \texttt{DLL-S} correct more bias than the \texttt{DLL} estimator, and the coverage for \texttt{DLL-S} improves along with this additional bias-correction. However, the bias for \texttt{DLL} is still smaller than the plug-in estimator, and better coverage is obtained.

\begin{table}[htb!]
\centering
\resizebox{\linewidth}{!}{
\begin{tabular}[t]{|c|c|c|ccc|ccc|ccc|ccc|ccc|}
\multicolumn{18}{c}{Setting 1, exactly sparse: $f(d) = 1.5\sin(d)$}\\
\hline
\multicolumn{1}{|c}{ } & \multicolumn{1}{c}{ } & \multicolumn{1}{c|}{ } & \multicolumn{3}{c|}{Bias} & \multicolumn{3}{c|}{RMSE} & \multicolumn{3}{c|}{SE} & \multicolumn{3}{c|}{Coverage} & \multicolumn{3}{c|}{CI Length} \\
\hline
$\eval$ & True & $n$ & \texttt{DLL} & \texttt{Plug} & \texttt{Orac} & \texttt{DLL} & \texttt{Plug} & \texttt{Orac} & \texttt{DLL} & \texttt{Plug} & \texttt{Orac} & \texttt{DLL} & \texttt{Plug} & \texttt{Orac} & \texttt{DLL} & \texttt{Plug} & \texttt{Orac}\\
\hline
 &  & 500 & 0.00 & 0.13 & 0.03 & 0.49 & 0.50 & 0.49 & 0.49 & 0.48 & 0.49 & 0.95 & 0.93 & 0.94 & 2.01 & 1.90 & 1.90\\

 &  & 1000 & 0.02 & 0.10 & 0.03 & 0.41 & 0.42 & 0.40 & 0.41 & 0.40 & 0.40 & 0.96 & 0.94 & 0.95 & 1.68 & 1.60 & 1.59\\

 &  & 1500 & 0.02 & 0.13 & 0.01 & 0.38 & 0.39 & 0.36 & 0.38 & 0.37 & 0.36 & 0.94 & 0.94 & 0.96 & 1.51 & 1.44 & 1.43\\

\multirow{-4}{*}{\centering -1.25} & \multirow{-4}{*}{\centering -1.00} & 2000 & 0.01 & 0.09 & 0.02 & 0.37 & 0.37 & 0.36 & 0.37 & 0.36 & 0.36 & 0.96 & 0.94 & 0.95  & 1.42 & 1.35 & 1.35\\
\cline{1-18}
 &  & 500 & 0.00 & 0.15 & 0.06 & 0.40 & 0.42 & 0.38 & 0.40 & 0.39 & 0.38 & 0.95 & 0.94 & 0.94 & 1.59 & 1.50 & 1.50\\

 &  & 1000 & 0.02 & 0.11 & 0.04 & 0.32 & 0.33 & 0.31 & 0.32 & 0.31 & 0.31 & 0.96 & 0.94 & 0.95 & 1.32 & 1.25 & 1.24\\

 &  & 1500 & 0.01 & 0.13 & 0.01 & 0.29 & 0.31 & 0.29 & 0.29 & 0.29 & 0.29 & 0.97 & 0.95 & 0.96 & 1.19 & 1.14 & 1.13\\

\multirow{-4}{*}{\centering -0.50} & \multirow{-4}{*}{\centering -0.56} & 2000 & 0.03 & 0.08 & 0.03 & 0.28 & 0.29 & 0.27 & 0.28 & 0.28 & 0.27 & 0.95 & 0.94 & 0.94 & 1.12 & 1.07 & 1.07\\
\cline{1-18}
 &  & 500 & 0.17 & 0.30 & 0.06 & 0.45 & 0.51 & 0.42 & 0.42 & 0.42 & 0.41 & 0.92 & 0.85 & 0.92 & 1.61 & 1.53 & 1.53\\

 &  & 1000 & 0.08 & 0.21 & 0.04 & 0.33 & 0.38 & 0.32 & 0.32 & 0.31 & 0.32 & 0.95 & 0.91 & 0.95 & 1.34 & 1.28 & 1.27\\

 &  & 1500 & 0.07 & 0.19 & 0.04 & 0.31 & 0.35 & 0.30 & 0.30 & 0.30 & 0.30 & 0.95 & 0.88 & 0.95 & 1.21 & 1.16 & 1.15\\

\multirow{-4}{*}{\centering 0.10} & \multirow{-4}{*}{\centering 0.74} & 2000 & 0.07 & 0.18 & 0.05 & 0.29 & 0.33 & 0.28 & 0.29 & 0.28 & 0.28 & 0.96 & 0.88 & 0.95  & 1.13 & 1.09 & 1.08\\
\cline{1-18}
 &  & 500 & 0.17 & 0.30 & 0.08 & 0.46 & 0.51 & 0.43 & 0.43 & 0.41 & 0.42 & 0.92 & 0.87 & 0.93 & 1.67 & 1.58 & 1.58\\

 &  & 1000 & 0.08 & 0.22 & 0.04 & 0.36 & 0.40 & 0.35 & 0.35 & 0.34 & 0.35 & 0.93 & 0.89 & 0.94 & 1.38 & 1.33 & 1.32\\

 &  & 1500 & 0.08 & 0.20 & 0.04 & 0.31 & 0.35 & 0.30 & 0.30 & 0.29 & 0.30 & 0.96 & 0.91 & 0.95 & 1.25 & 1.20 & 1.19\\

\multirow{-4}{*}{\centering 0.25} & \multirow{-4}{*}{\centering 0.94} & 2000 & 0.07 & 0.17 & 0.04 & 0.30 & 0.34 & 0.29 & 0.29 & 0.29 & 0.29 & 0.95 & 0.89 & 0.94 & 1.17 & 1.13 & 1.12\\
\cline{1-18}
 &  & 500 & 0.08 & 0.18 & 0.00 & 0.61 & 0.61 & 0.58 & 0.61 & 0.58 & 0.58 & 0.93 & 0.90 & 0.93 & 2.33 & 2.19 & 2.19\\

 &  & 1000 & 0.00 & 0.12 & 0.04 & 0.48 & 0.48 & 0.47 & 0.48 & 0.47 & 0.47 & 0.96 & 0.94 & 0.95 & 1.94 & 1.83 & 1.83\\

 &  & 1500 & 0.03 & 0.14 & 0.01 & 0.44 & 0.45 & 0.42 & 0.44 & 0.43 & 0.42 & 0.95 & 0.93 & 0.95 & 1.75 & 1.66 & 1.66\\

\multirow{-4}{*}{\centering 1.00} & \multirow{-4}{*}{\centering 0.81} & 2000 & 0.01 & 0.10 & 0.02 & 0.40 & 0.40 & 0.39 & 0.40 & 0.39 & 0.39 & 0.95 & 0.93 & 0.95 & 1.63 & 1.55 & 1.55\\
\hline
\multicolumn{18}{c}{Setting 1, approximately sparse: $f(d) = 1.5\sin(d)$}\\
\hline
\multicolumn{1}{|c}{ } & \multicolumn{1}{c}{ } & \multicolumn{1}{c|}{ } & \multicolumn{3}{c|}{Bias} & \multicolumn{3}{c|}{RMSE} & \multicolumn{3}{c|}{SE} & \multicolumn{3}{c|}{Coverage} & \multicolumn{3}{c|}{CI Length} \\
\hline
$\eval$ & True & $n$ & \texttt{DLL} & \texttt{Plug} & \texttt{Orac} & \texttt{DLL} & \texttt{Plug} & \texttt{Orac} & \texttt{DLL} & \texttt{Plug} & \texttt{Orac} & \texttt{DLL} & \texttt{Plug} & \texttt{Orac} & \texttt{DLL} & \texttt{Plug} & \texttt{Orac}\\
\hline
 &  & 500 & 0.02 & 0.08 & 0.01 & 0.47 & 0.47 & 0.50 & 0.47 & 0.47 & 0.50 & 0.96 & 0.94 & 0.94 & 1.89 & 1.78 & 1.92\\

 &  & 1000 & 0.04 & 0.06 & 0.01 & 0.39 & 0.39 & 0.39 & 0.39 & 0.38 & 0.39 & 0.96 & 0.96 & 0.96 &  1.62 & 1.56 & 1.61\\

 &  & 1500 & 0.02 & 0.06 & 0.01 & 0.35 & 0.35 & 0.36 & 0.35 & 0.34 & 0.36 & 0.98 & 0.95 & 0.96 & 1.47 & 1.41 & 1.44\\

\multirow{-4}{*}{\centering -1.25} & \multirow{-4}{*}{\centering -1.00} & 2000 & 0.03 & 0.04 & 0.01 & 0.34 & 0.34 & 0.34 & 0.34 & 0.34 & 0.34 & 0.95 & 0.95 & 0.95 & 1.37 & 1.31 & 1.34\\
\cline{1-18}
 &  & 500 & 0.03 & 0.08 & 0.02 & 0.38 & 0.38 & 0.38 & 0.38 & 0.37 & 0.38 & 0.95 & 0.95 & 0.95 & 1.48 & 1.41 & 1.50\\

 &  & 1000 & 0.03 & 0.07 & 0.02 & 0.33 & 0.33 & 0.33 & 0.33 & 0.32 & 0.32 & 0.94 & 0.93 & 0.94 & 1.28 & 1.23 & 1.26\\

 &  & 1500 & 0.03 & 0.06 & 0.03 & 0.28 & 0.28 & 0.28 & 0.28 & 0.28 & 0.28 & 0.96 & 0.95 & 0.96 & 1.16 & 1.11 & 1.14\\

\multirow{-4}{*}{\centering -0.50} & \multirow{-4}{*}{\centering -0.56} & 2000 & 0.02 & 0.06 & 0.01 & 0.28 & 0.28 & 0.26 & 0.28 & 0.27 & 0.26 & 0.96 & 0.93 & 0.96 & 1.08 & 1.04 & 1.06\\
\cline{1-18}
 &  & 500 & 0.19 & 0.29 & 0.03 & 0.42 & 0.47 & 0.39 & 0.38 & 0.37 & 0.39 & 0.92 & 0.86 & 0.96 & 1.51 & 1.44 & 1.53\\

 &  & 1000 & 0.15 & 0.25 & 0.06 & 0.36 & 0.40 & 0.32 & 0.32 & 0.31 & 0.32 & 0.93 & 0.88 & 0.95 & 1.30 & 1.24 & 1.28\\

 &  & 1500 & 0.11 & 0.20 & 0.03 & 0.33 & 0.36 & 0.31 & 0.31 & 0.30 & 0.31 & 0.94 & 0.88 & 0.93 & 1.17 & 1.13 & 1.15\\

\multirow{-4}{*}{\centering 0.10} & \multirow{-4}{*}{\centering 0.74} & 2000 & 0.10 & 0.18 & 0.04 & 0.29 & 0.32 & 0.28 & 0.27 & 0.27 & 0.27 & 0.95 & 0.90 & 0.94 & 1.10 & 1.05 & 1.08\\
\cline{1-18}
 &  & 500 & 0.21 & 0.31 & 0.06 & 0.45 & 0.49 & 0.41 & 0.40 & 0.38 & 0.40 & 0.92 & 0.86 & 0.95 & 1.55 & 1.48 & 1.58\\

 &  & 1000 & 0.18 & 0.28 & 0.09 & 0.37 & 0.42 & 0.33 & 0.32 & 0.31 & 0.32 & 0.92 & 0.87 & 0.95 & 1.34 & 1.28 & 1.33\\

 &  & 1500 & 0.12 & 0.21 & 0.04 & 0.34 & 0.37 & 0.32 & 0.32 & 0.31 & 0.31 & 0.92 & 0.86 & 0.94 & 1.22 & 1.17 & 1.20\\

\multirow{-4}{*}{\centering 0.25} & \multirow{-4}{*}{\centering 0.94} & 2000 & 0.12 & 0.20 & 0.06 & 0.30 & 0.34 & 0.29 & 0.28 & 0.27 & 0.28 & 0.95 & 0.87 & 0.95 & 1.13 & 1.09 & 1.11\\
\cline{1-18}
 &  & 500 & 0.12 & 0.19 & 0.01 & 0.59 & 0.60 & 0.58 & 0.58 & 0.57 & 0.58 & 0.92 & 0.90 & 0.93 & 2.19 & 2.07 & 2.21\\

 &  & 1000 & 0.04 & 0.12 & 0.02 & 0.46 & 0.47 & 0.46 & 0.46 & 0.45 & 0.46 & 0.96 & 0.94 & 0.95 & 1.88 & 1.79 & 1.86\\

 &  & 1500 & 0.01 & 0.07 & 0.04 & 0.42 & 0.42 & 0.44 & 0.42 & 0.42 & 0.44 & 0.95 & 0.94 & 0.95 & 1.70 & 1.62 & 1.67\\

\multirow{-4}{*}{\centering 1.00} & \multirow{-4}{*}{\centering 0.81} & 2000 & 0.05 & 0.12 & 0.01 & 0.40 & 0.41 & 0.39 & 0.40 & 0.39 & 0.39 & 0.94 & 0.93 & 0.95 & 1.58 & 1.52 & 1.55\\
\hline
\end{tabular}
}
\caption*{Table S1: Comparison of \texttt{DLL}, plug-in (\texttt{Plug}), oracle (\texttt{Orac}) estimators in Setting 1 when $f(d) = 1.5\sin(d)$, across different sample sizes $n$ and evaluation points $\eval$. The column indexed with ``True'' represents the true value of $f'(\eval)$. The columns indexed with ``Bias'', ``RMSE'' and ``SE" report the absolute bias, the root mean square error, and the standard error computed by 500 estimates, respectively; the columns indexed with ``Coverage'' report the empirical coverage level and the columns indexed with ``Length'' report the average CI length.}
\label{table:S1}
\end{table}

\begin{table}[htb!]
\centering
\resizebox{\linewidth}{!}{
\begin{tabular}[t]{|c|c|c|ccc|ccc|ccc|ccc|ccc|}
\multicolumn{18}{c}{Setting 1, exactly sparse: $f(d) = 2\exp(-d/2)$}\\
\hline
\multicolumn{1}{|c}{ } & \multicolumn{1}{c}{ } & \multicolumn{1}{c|}{ } & \multicolumn{3}{c|}{Bias} & \multicolumn{3}{c|}{RMSE} & \multicolumn{3}{c|}{SE} & \multicolumn{3}{c|}{Coverage} & \multicolumn{3}{c|}{CI Length} \\
\hline
$\eval$ & True & $n$ & \texttt{DLL} & \texttt{Plug} & \texttt{Orac} & \texttt{DLL} & \texttt{Plug} & \texttt{Orac} & \texttt{DLL} & \texttt{Plug} & \texttt{Orac} & \texttt{DLL} & \texttt{Plug} & \texttt{Orac} & \texttt{DLL} & \texttt{Plug} & \texttt{Orac}\\
\hline
 &  & 500 & 0.03 & 0.14 & 0.01 & 0.30 & 0.32 & 0.25 & 0.30 & 0.29 & 0.25 & 0.96 & 0.92 & 0.95 & 1.19 & 1.11 & 0.94\\

 &  & 1000 & 0.05 & 0.16 & 0.02 & 0.25 & 0.28 & 0.18 & 0.24 & 0.23 & 0.18 & 0.94 & 0.88 & 0.95  & 0.92 & 0.87 & 0.70\\

 &  & 1500 & 0.02 & 0.13 & 0.01 & 0.21 & 0.24 & 0.15 & 0.21 & 0.21 & 0.15 & 0.94 & 0.90 & 0.95  & 0.82 & 0.78 & 0.62\\

\multirow{-4}{*}{\centering -1.25} & \multirow{-4}{*}{\centering -0.41} & 2000 & 0.02 & 0.12 & 0.01 & 0.20 & 0.23 & 0.15 & 0.20 & 0.19 & 0.15 & 0.95 & 0.89 & 0.95 & 0.74 & 0.71 & 0.56\\
\cline{1-18}
 &  & 500 & 0.00 & 0.01 & 0.02 & 0.24 & 0.23 & 0.20 & 0.24 & 0.23 & 0.20 & 0.95 & 0.94 & 0.92 & 0.93 & 0.87 & 0.73\\

 &  & 1000 & 0.01 & 0.00 & 0.03 & 0.19 & 0.19 & 0.16 & 0.19 & 0.19 & 0.15 & 0.93 & 0.93 & 0.92 & 0.73 & 0.69 & 0.55\\

 &  & 1500 & 0.02 & 0.00 & 0.02 & 0.17 & 0.16 & 0.13 & 0.16 & 0.16 & 0.12 & 0.95 & 0.94 & 0.92 & 0.65 & 0.62 & 0.48\\

\multirow{-4}{*}{\centering -0.50} & \multirow{-4}{*}{\centering -0.52} & 2000 & 0.02 & 0.00 & 0.03 & 0.15 & 0.15 & 0.12 & 0.15 & 0.15 & 0.11 & 0.95 & 0.95 & 0.93 & 0.59 & 0.56 & 0.44\\
\cline{1-18}
 &  & 500 & 0.01 & 0.13 & 0.01 & 0.25 & 0.27 & 0.20 & 0.25 & 0.24 & 0.20 & 0.94 & 0.90 & 0.94 & 0.94 & 0.89 & 0.75\\

 &  & 1000 & 0.01 & 0.10 & 0.00 & 0.19 & 0.21 & 0.15 & 0.19 & 0.18 & 0.15 & 0.95 & 0.90 & 0.94 & 0.74 & 0.70 & 0.56\\

 &  & 1500 & 0.01 & 0.09 & 0.00 & 0.16 & 0.18 & 0.12 & 0.16 & 0.16 & 0.12 & 0.96 & 0.92 & 0.95 & 0.66 & 0.63 & 0.49\\

\multirow{-4}{*}{\centering 0.10} & \multirow{-4}{*}{\centering -0.40} & 2000 & 0.02 & 0.07 & 0.01 & 0.15 & 0.16 & 0.12 & 0.15 & 0.15 & 0.12 & 0.96 & 0.92 & 0.94 & 0.60 & 0.57 & 0.45\\
\cline{1-18}
 &  & 500 & 0.01 & 0.14 & 0.01 & 0.26 & 0.29 & 0.21 & 0.26 & 0.25 & 0.21 & 0.95 & 0.90 & 0.95 & 0.97 & 0.92 & 0.77\\

 &  & 1000 & 0.01 & 0.12 & 0.00 & 0.20 & 0.23 & 0.15 & 0.20 & 0.19 & 0.15 & 0.95 & 0.88 & 0.96 & 0.76 & 0.73 & 0.58\\

 &  & 1500 & 0.02 & 0.11 & 0.01 & 0.17 & 0.19 & 0.13 & 0.17 & 0.16 & 0.13 & 0.95 & 0.91 & 0.95 & 0.68 & 0.65 & 0.51\\

\multirow{-4}{*}{\centering 0.25} & \multirow{-4}{*}{\centering -0.35} & 2000 & 0.02 & 0.09 & 0.02 & 0.16 & 0.18 & 0.12 & 0.16 & 0.15 & 0.12 & 0.95 & 0.92 & 0.95 & 0.62 & 0.59 & 0.46\\
\cline{1-18}
 &  & 500 & 0.01 & 0.20 & 0.03 & 0.37 & 0.41 & 0.29 & 0.37 & 0.35 & 0.29 & 0.94 & 0.91 & 0.93 & 1.37 & 1.29 & 1.07\\

 &  & 1000 & 0.04 & 0.14 & 0.02 & 0.28 & 0.31 & 0.21 & 0.28 & 0.27 & 0.21 & 0.95 & 0.92 & 0.95 & 1.06 & 1.01 & 0.80\\

 &  & 1500 & 0.02 & 0.14 & 0.02 & 0.25 & 0.28 & 0.19 & 0.25 & 0.24 & 0.19 & 0.96 & 0.90 & 0.95 & 0.94 & 0.90 & 0.71\\

\multirow{-4}{*}{\centering 1.00} & \multirow{-4}{*}{\centering -0.15} & 2000 & 0.01 & 0.14 & 0.01 & 0.23 & 0.27 & 0.17 & 0.23 & 0.23 & 0.17 & 0.93 & 0.89 & 0.94 & 0.86 & 0.83 & 0.64\\
\hline

\multicolumn{18}{c}{Setting 1, approximately sparse: $f(d) = 2\exp(-d/2)$}\\
\hline
\multicolumn{1}{|c}{ } & \multicolumn{1}{c}{ } & \multicolumn{1}{c|}{ } & \multicolumn{3}{c|}{Bias} & \multicolumn{3}{c|}{RMSE} & \multicolumn{3}{c|}{SE} & \multicolumn{3}{c|}{Coverage} & \multicolumn{3}{c|}{CI Length} \\
\hline
$\eval$ & True & $n$ & \texttt{DLL} & \texttt{Plug} & \texttt{Orac} & \texttt{DLL} & \texttt{Plug} & \texttt{Orac} & \texttt{DLL} & \texttt{Plug} & \texttt{Orac} & \texttt{DLL} & \texttt{Plug} & \texttt{Orac} & \texttt{DLL} & \texttt{Plug} & \texttt{Orac}\\
\hline
 &  & 500 & 0.01 & 0.08 & 0.01 & 0.29 & 0.29 & 0.24 & 0.29 & 0.28 & 0.24 & 0.95 & 0.94 & 0.96 & 1.11 & 1.03 & 0.92\\

 &  & 1000 & 0.01 & 0.08 & 0.01 & 0.22 & 0.23 & 0.19 & 0.22 & 0.22 & 0.19 & 0.96 & 0.93 & 0.95 & 0.89 & 0.84 & 0.71\\

 &  & 1500 & 0.02 & 0.09 & 0.02 & 0.21 & 0.23 & 0.16 & 0.21 & 0.21 & 0.16 & 0.95 & 0.93 & 0.94  & 0.80 & 0.76 & 0.62\\

\multirow{-4}{*}{\centering -1.25} & \multirow{-4}{*}{\centering -0.41} & 2000 & 0.01 & 0.06 & 0.01 & 0.19 & 0.19 & 0.15 & 0.19 & 0.18 & 0.15 & 0.94 & 0.93 & 0.95  & 0.73 & 0.69 & 0.57\\
\cline{1-18}
 &  & 500 & 0.05 & 0.08 & 0.03 & 0.24 & 0.24 & 0.19 & 0.23 & 0.22 & 0.19 & 0.92 & 0.89 & 0.92 & 0.87 & 0.82 & 0.72\\

 &  & 1000 & 0.03 & 0.04 & 0.03 & 0.19 & 0.19 & 0.15 & 0.19 & 0.18 & 0.15 & 0.93 & 0.92 & 0.92 & 0.70 & 0.67 & 0.56\\

 &  & 1500 & 0.02 & 0.03 & 0.02 & 0.17 & 0.17 & 0.13 & 0.17 & 0.16 & 0.12 & 0.92 & 0.91 & 0.95 & 0.63 & 0.60 & 0.49\\

\multirow{-4}{*}{\centering -0.50} & \multirow{-4}{*}{\centering -0.52} & 2000 & 0.01 & 0.02 & 0.01 & 0.15 & 0.14 & 0.12 & 0.15 & 0.14 & 0.12 & 0.95 & 0.92 & 0.94 & 0.57 & 0.55 & 0.45\\
\cline{1-18}
 &  & 500 & 0.01 & 0.15 & 0.02 & 0.23 & 0.27 & 0.20 & 0.23 & 0.23 & 0.20 & 0.95 & 0.90 & 0.95 & 0.89 & 0.84 & 0.74\\

 &  & 1000 & 0.00 & 0.13 & 0.00 & 0.19 & 0.22 & 0.15 & 0.19 & 0.18 & 0.15 & 0.93 & 0.86 & 0.94 & 0.71 & 0.68 & 0.57\\

 &  & 1500 & 0.01 & 0.12 & 0.00 & 0.16 & 0.19 & 0.13 & 0.16 & 0.15 & 0.13 & 0.95 & 0.86 & 0.94 & 0.64 & 0.61 & 0.50\\

\multirow{-4}{*}{\centering 0.10} & \multirow{-4}{*}{\centering -0.40} & 2000 & 0.01 & 0.10 & 0.01 & 0.15 & 0.18 & 0.11 & 0.15 & 0.15 & 0.11 & 0.95 & 0.86 & 0.96 & 0.58 & 0.56 & 0.46\\
\cline{1-18}
 &  & 500 & 0.00 & 0.16 & 0.03 & 0.23 & 0.27 & 0.21 & 0.23 & 0.22 & 0.21 & 0.96 & 0.90 & 0.93 & 0.91 & 0.86 & 0.76\\

 &  & 1000 & 0.00 & 0.15 & 0.00 & 0.20 & 0.24 & 0.16 & 0.20 & 0.19 & 0.16 & 0.92 & 0.86 & 0.94  & 0.74 & 0.70 & 0.58\\

 &  & 1500 & 0.01 & 0.14 & 0.00 & 0.16 & 0.22 & 0.13 & 0.16 & 0.16 & 0.13 & 0.96 & 0.85 & 0.95 & 0.66 & 0.63 & 0.52\\

\multirow{-4}{*}{\centering 0.25} & \multirow{-4}{*}{\centering -0.35} & 2000 & 0.02 & 0.11 & 0.01 & 0.15 & 0.18 & 0.12 & 0.15 & 0.15 & 0.12 & 0.96 & 0.88 & 0.95 & 0.60 & 0.58 & 0.47\\
\cline{1-18}
 &  & 500 & 0.01 & 0.23 & 0.02 & 0.35 & 0.40 & 0.30 & 0.35 & 0.33 & 0.30 & 0.94 & 0.88 & 0.92 & 1.28 & 1.22 & 1.07\\

 &  & 1000 & 0.02 & 0.17 & 0.02 & 0.29 & 0.32 & 0.23 & 0.29 & 0.27 & 0.23 & 0.94 & 0.89 & 0.93 & 1.03 & 0.98 & 0.82\\

 &  & 1500 & 0.01 & 0.16 & 0.01 & 0.23 & 0.27 & 0.17 & 0.23 & 0.23 & 0.17 & 0.93 & 0.89 & 0.95 & 0.92 & 0.88 & 0.72\\

\multirow{-4}{*}{\centering 1.00} & \multirow{-4}{*}{\centering -0.15} & 2000 & 0.00 & 0.15 & 0.01 & 0.22 & 0.26 & 0.18 & 0.22 & 0.21 & 0.18 & 0.95 & 0.87 & 0.94 & 0.84 & 0.81 & 0.66\\
\hline
\end{tabular}
}
\caption*{Table S2: Comparison of \texttt{DLL}, plug-in (\texttt{Plug}), oracle (\texttt{Orac}) estimators in Setting 1 when $f(d) = 2\exp(-d/2)$, across different sample sizes $n$ and evaluation points $\eval$. The column indexed with ``True'' represents the true value of $f'(\eval)$. The columns indexed with ``Bias'', ``RMSE'' and ``SE" report the absolute bias, the root mean square error, and the standard error computed by 500 estimates, respectively; the columns indexed with ``Coverage'' report the empirical coverage level and the columns indexed with ``Length'' report the average CI length.}
\label{table:S2}
\end{table}

\begin{table}[htb!]
\centering
\resizebox{\linewidth}{!}{
\begin{tabular}[t]{|c|c|c|ccc|ccc|ccc|ccc|ccc|}
\multicolumn{18}{c}{Setting 2, exactly sparse: $f(d) = 1.5\sin(d)$}\\
\hline
\multicolumn{1}{|c}{ } & \multicolumn{1}{c}{ } & \multicolumn{1}{c|}{ } & \multicolumn{3}{c|}{Bias} & \multicolumn{3}{c|}{RMSE} & \multicolumn{3}{c|}{SE} & \multicolumn{3}{c|}{Coverage} & \multicolumn{3}{c|}{CI Length} \\
\hline
$\eval$ & True & $n$ & \texttt{DLL} & \texttt{Plug} & \texttt{Orac} & \texttt{DLL} & \texttt{Plug} & \texttt{Orac} & \texttt{DLL} & \texttt{Plug} & \texttt{Orac} & \texttt{DLL} & \texttt{Plug} & \texttt{Orac} & \texttt{DLL} & \texttt{Plug} & \texttt{Orac}\\
\hline
 &  & 500 & 0.12 & 0.38 & 0.06 & 0.50 & 0.60 & 0.48 & 0.48 & 0.47 & 0.48 & 0.94 & 0.87 & 0.94 & 1.92 & 1.84 & 1.89\\

 &  & 1000 & 0.04 & 0.26 & 0.02 & 0.41 & 0.48 & 0.42 & 0.41 & 0.41 & 0.42 & 0.95 & 0.90 & 0.94 & 1.64 & 1.58 & 1.58\\

 &  & 1500 & 0.03 & 0.21 & 0.02 & 0.39 & 0.44 & 0.38 & 0.39 & 0.39 & 0.38 & 0.96 & 0.90 & 0.95 & 1.49 & 1.44 & 1.44\\

\multirow{-4}{*}{\centering -1.25} & \multirow{-4}{*}{\centering 0.47} & 2000 & 0.03 & 0.19 & 0.00 & 0.35 & 0.39 & 0.33 & 0.34 & 0.34 & 0.33 & 0.94 & 0.92 & 0.94 & 1.39 & 1.36 & 1.35\\
\cline{1-18}
 &  & 500 & 0.19 & 0.44 & 0.03 & 0.41 & 0.58 & 0.38 & 0.37 & 0.37 & 0.38 & 0.91 & 0.77 & 0.94 & 1.49 & 1.43 & 1.47\\

 &  & 1000 & 0.07 & 0.30 & 0.01 & 0.32 & 0.43 & 0.31 & 0.31 & 0.31 & 0.31 & 0.96 & 0.83 & 0.94 & 1.29 & 1.25 & 1.25\\

 &  & 1500 & 0.06 & 0.25 & 0.02 & 0.30 & 0.39 & 0.29 & 0.29 & 0.29 & 0.29 & 0.95 & 0.85 & 0.95 & 1.17 & 1.13 & 1.14\\

\multirow{-4}{*}{\centering -0.50} & \multirow{-4}{*}{\centering 1.32} & 2000 & 0.04 & 0.22 & 0.01 & 0.29 & 0.36 & 0.28 & 0.29 & 0.29 & 0.28 & 0.94 & 0.86 & 0.95 & 1.10 & 1.07 & 1.07\\
\cline{1-18}
 &  & 500 & 0.21 & 0.46 & 0.02 & 0.44 & 0.60 & 0.39 & 0.39 & 0.38 & 0.39 & 0.91 & 0.72 & 0.94 & 1.51 & 1.45 & 1.49\\

 &  & 1000 & 0.07 & 0.31 & 0.00 & 0.35 & 0.45 & 0.33 & 0.34 & 0.33 & 0.33 & 0.93 & 0.83 & 0.93 & 1.31 & 1.27 & 1.27\\

 &  & 1500 & 0.05 & 0.26 & 0.01 & 0.31 & 0.39 & 0.29 & 0.31 & 0.30 & 0.29 & 0.94 & 0.86 & 0.95 & 1.18 & 1.15 & 1.15\\

\multirow{-4}{*}{\centering 0.10} & \multirow{-4}{*}{\centering 1.49} & 2000 & 0.05 & 0.24 & 0.02 & 0.29 & 0.37 & 0.28 & 0.29 & 0.28 & 0.28 & 0.94 & 0.86 & 0.96 & 1.12 & 1.08 & 1.08\\
\cline{1-18}
 &  & 500 & 0.20 & 0.45 & 0.00 & 0.45 & 0.60 & 0.39 & 0.41 & 0.39 & 0.39 & 0.91 & 0.77 & 0.94 & 1.56 & 1.50 & 1.55\\

 &  & 1000 & 0.07 & 0.31 & 0.01 & 0.36 & 0.46 & 0.35 & 0.35 & 0.34 & 0.35 & 0.94 & 0.83 & 0.94 & 1.35 & 1.32 & 1.32\\

 &  & 1500 & 0.07 & 0.27 & 0.03 & 0.32 & 0.41 & 0.30 & 0.31 & 0.31 & 0.30 & 0.96 & 0.84 & 0.96 & 1.22 & 1.19 & 1.18\\

\multirow{-4}{*}{\centering 0.25} & \multirow{-4}{*}{\centering 1.45} & 2000 & 0.05 & 0.24 & 0.02 & 0.29 & 0.37 & 0.28 & 0.28 & 0.28 & 0.28 & 0.95 & 0.87 & 0.96 & 1.15 & 1.12 & 1.12\\
\cline{1-18}
 &  & 500 & 0.17 & 0.38 & 0.00 & 0.55 & 0.64 & 0.54 & 0.52 & 0.51 & 0.54 & 0.96 & 0.90 & 0.94 & 2.19 & 2.06 & 2.15\\

 &  & 1000 & 0.05 & 0.25 & 0.04 & 0.50 & 0.55 & 0.50 & 0.50 & 0.49 & 0.50 & 0.95 & 0.89 & 0.93 & 1.89 & 1.81 & 1.83\\

 &  & 1500 & 0.02 & 0.19 & 0.04 & 0.44 & 0.47 & 0.43 & 0.44 & 0.43 & 0.43 & 0.94 & 0.92 & 0.94 & 1.71 & 1.65 & 1.65\\

\multirow{-4}{*}{\centering 1.00} & \multirow{-4}{*}{\centering 0.81} & 2000 & 0.03 & 0.18 & 0.00 & 0.41 & 0.44 & 0.40 & 0.41 & 0.40 & 0.40 & 0.95 & 0.92 & 0.95 & 1.60 & 1.56 & 1.56\\
\hline

\multicolumn{18}{c}{Setting 2, exactly sparse: $f(d) = 2\exp(-d/2)$}\\
\hline
\multicolumn{1}{|c}{ } & \multicolumn{1}{c}{ } & \multicolumn{1}{c|}{ } & \multicolumn{3}{c|}{Bias} & \multicolumn{3}{c|}{RMSE} & \multicolumn{3}{c|}{SE} & \multicolumn{3}{c|}{Coverage} & \multicolumn{3}{c|}{CI Length} \\
\hline
$\eval$ & True & $n$ & \texttt{DLL} & \texttt{Plug} & \texttt{Orac} & \texttt{DLL} & \texttt{Plug} & \texttt{Orac} & \texttt{DLL} & \texttt{Plug} & \texttt{Orac} & \texttt{DLL} & \texttt{Plug} & \texttt{Orac} & \texttt{DLL} & \texttt{Plug} & \texttt{Orac}\\
\hline
 &  & 500 & 0.04 & 0.23 & 0.03 & 0.44 & 0.49 & 0.41 & 0.44 & 0.43 & 0.41 & 0.95 & 0.90 & 0.94 & 1.70 & 1.62 & 1.57\\

 &  & 1000 & 0.01 & 0.16 & 0.01 & 0.37 & 0.40 & 0.35 & 0.37 & 0.36 & 0.35 & 0.95 & 0.92 & 0.94 & 1.43 & 1.39 & 1.34\\

 &  & 1500 & 0.01 & 0.16 & 0.02 & 0.34 & 0.37 & 0.32 & 0.34 & 0.33 & 0.32 & 0.94 & 0.91 & 0.93 & 1.30 & 1.27 & 1.24\\

\multirow{-4}{*}{\centering -1.25} & \multirow{-4}{*}{\centering -1.87} & 2000 & 0.00 & 0.13 & 0.00 & 0.30 & 0.32 & 0.29 & 0.30 & 0.29 & 0.29 & 0.96 & 0.95 & 0.96 & 1.21 & 1.18 & 1.15\\
\cline{1-18}
 &  & 500 & 0.04 & 0.27 & 0.01 & 0.37 & 0.45 & 0.32 & 0.37 & 0.36 & 0.32 & 0.93 & 0.83 & 0.94 & 1.34 & 1.29 & 1.24\\

 &  & 1000 & 0.01 & 0.19 & 0.01 & 0.30 & 0.34 & 0.27 & 0.30 & 0.29 & 0.27 & 0.94 & 0.89 & 0.94 & 1.13 & 1.09 & 1.06\\

 &  & 1500 & 0.04 & 0.20 & 0.03 & 0.28 & 0.34 & 0.26 & 0.27 & 0.27 & 0.26 & 0.94 & 0.85 & 0.94 & 1.02 & 1.00 & 0.97\\

\multirow{-4}{*}{\centering -0.50} & \multirow{-4}{*}{\centering -1.28} & 2000 & 0.01 & 0.13 & 0.01 & 0.23 & 0.26 & 0.22 & 0.23 & 0.22 & 0.22 & 0.97 & 0.93 & 0.96 & 0.95 & 0.92 & 0.91\\
\cline{1-18}
 &  & 500 & 0.07 & 0.32 & 0.01 & 0.35 & 0.47 & 0.32 & 0.34 & 0.34 & 0.32 & 0.95 & 0.83 & 0.95 & 1.37 & 1.31 & 1.26\\

 &  & 1000 & 0.02 & 0.22 & 0.00 & 0.30 & 0.37 & 0.27 & 0.30 & 0.30 & 0.27 & 0.93 & 0.85 & 0.95 & 1.14 & 1.10 & 1.07\\

 &  & 1500 & 0.00 & 0.17 & 0.01 & 0.25 & 0.30 & 0.24 & 0.25 & 0.25 & 0.24 & 0.96 & 0.90 & 0.96 & 1.04 & 1.01 & 0.99\\

\multirow{-4}{*}{\centering 0.10} & \multirow{-4}{*}{\centering -0.95} & 2000 & 0.01 & 0.14 & 0.01 & 0.24 & 0.28 & 0.23 & 0.24 & 0.24 & 0.23 & 0.96 & 0.92 & 0.95 & 0.97 & 0.94 & 0.92\\
\cline{1-18}
 &  & 500 & 0.06 & 0.31 & 0.01 & 0.37 & 0.47 & 0.34 & 0.36 & 0.35 & 0.34 & 0.94 & 0.86 & 0.94 & 1.41 & 1.35 & 1.30\\

 &  & 1000 & 0.01 & 0.22 & 0.00 & 0.30 & 0.37 & 0.28 & 0.30 & 0.30 & 0.28 & 0.96 & 0.88 & 0.96 & 1.18 & 1.14 & 1.10\\

 &  & 1500 & 0.02 & 0.19 & 0.02 & 0.28 & 0.33 & 0.26 & 0.28 & 0.27 & 0.26 & 0.95 & 0.89 & 0.95 & 1.07 & 1.04 & 1.02\\

\multirow{-4}{*}{\centering 0.25} & \multirow{-4}{*}{\centering -0.88} & 2000 & 0.01 & 0.14 & 0.02 & 0.24 & 0.27 & 0.23 & 0.24 & 0.24 & 0.23 & 0.96 & 0.92 & 0.97 & 1.00 & 0.97 & 0.95\\
\cline{1-18}
 &  & 500 & 0.05 & 0.27 & 0.01 & 0.48 & 0.54 & 0.45 & 0.48 & 0.47 & 0.45 & 0.95 & 0.92 & 0.95 & 1.97 & 1.87 & 1.81\\

 &  & 1000 & 0.02 & 0.17 & 0.02 & 0.41 & 0.44 & 0.40 & 0.41 & 0.40 & 0.40 & 0.95 & 0.92 & 0.95 & 1.65 & 1.58 & 1.54\\

 &  & 1500 & 0.00 & 0.16 & 0.00 & 0.37 & 0.40 & 0.36 & 0.38 & 0.37 & 0.36 & 0.95 & 0.92 & 0.95 & 1.50 & 1.45 & 1.42\\

\multirow{-4}{*}{\centering 1.00} & \multirow{-4}{*}{\centering -0.61} & 2000 & 0.01 & 0.15 & 0.00 & 0.37 & 0.39 & 0.36 & 0.37 & 0.36 & 0.36 & 0.94 & 0.92 & 0.94 & 1.39 & 1.36 & 1.33\\
\hline
\end{tabular}
}
\caption*{Table S3: Comparison of \texttt{DLL}, plug-in (\texttt{Plug}), oracle (\texttt{Orac}) estimators in Setting 2, across different sample sizes $n$ and evaluation points $\eval$. The column indexed with ``True'' represents the true value of $f'(\eval)$. The columns indexed with ``Bias'', ``RMSE'' and ``SE" report the absolute bias, the root mean square error, and the standard error computed by 500 estimates, respectively; the columns indexed with ``Coverage'' report the empirical coverage level and the columns indexed with ``Length'' report the average CI length.}
\label{table:S3}
\end{table}

\begin{table}[htb!]
\centering
\resizebox{\linewidth}{!}{
\begin{tabular}[t]{|c|c|c|ccc|ccc|ccc|ccc|ccc|}
\multicolumn{18}{c}{Setting 3, exactly sparse: $f(d) = 1.5\sin(d)$}\\
\hline
\multicolumn{1}{|c}{ } & \multicolumn{1}{c}{ } & \multicolumn{1}{c|}{ } & \multicolumn{3}{c|}{Bias} & \multicolumn{3}{c|}{RMSE} & \multicolumn{3}{c|}{SE} & \multicolumn{3}{c|}{Coverage} & \multicolumn{3}{c|}{CI Length} \\
\hline
$\eval$ & True & $n$ & \texttt{DLL} & \texttt{Plug} & \texttt{Orac} & \texttt{DLL} & \texttt{Plug} & \texttt{Orac} & \texttt{DLL} & \texttt{Plug} & \texttt{Orac} & \texttt{DLL} & \texttt{Plug} & \texttt{Orac} & \texttt{DLL} & \texttt{Plug} & \texttt{Orac}\\
\hline
 &  & 500 & 0.12 & 0.32 & 0.02 & 0.73 & 0.78 & 0.70 & 0.72 & 0.71 & 0.70 & 0.95 & 0.92 & 0.94 & 2.85 & 2.77 & 2.62\\

 &  & 1000 & 0.02 & 0.19 & 0.02 & 0.62 & 0.64 & 0.57 & 0.62 & 0.62 & 0.57 & 0.96 & 0.92 & 0.96 & 2.42 & 2.36 & 2.23\\

 &  & 1500 & 0.03 & 0.17 & 0.00 & 0.61 & 0.62 & 0.56 & 0.61 & 0.60 & 0.56 & 0.94 & 0.90 & 0.94 & 2.19 & 2.16 & 2.04\\

\multirow{-4}{*}{\centering -1.25} & \multirow{-4}{*}{\centering 0.47} & 2000 & 0.01 & 0.13 & 0.02 & 0.51 & 0.53 & 0.49 & 0.51 & 0.51 & 0.49 & 0.96 & 0.94 & 0.95 & 2.06 & 2.02 & 1.93\\
\cline{1-18}
 &  & 500 & 0.10 & 0.24 & 0.02 & 0.77 & 0.80 & 0.70 & 0.77 & 0.76 & 0.70 & 0.93 & 0.91 & 0.94 & 2.82 & 2.75 & 2.60\\

 &  & 1000 & 0.00 & 0.14 & 0.01 & 0.61 & 0.62 & 0.55 & 0.62 & 0.60 & 0.55 & 0.96 & 0.95 & 0.96 & 2.41 & 2.38 & 2.26\\

 &  & 1500 & 0.04 & 0.17 & 0.04 & 0.57 & 0.59 & 0.52 & 0.57 & 0.57 & 0.52 & 0.95 & 0.95 & 0.96 & 2.19 & 2.18 & 2.06\\

\multirow{-4}{*}{\centering -0.50} & \multirow{-4}{*}{\centering 1.32} & 2000 & 0.02 & 0.11 & 0.01 & 0.48 & 0.49 & 0.45 & 0.48 & 0.48 & 0.45 & 0.98 & 0.96 & 0.96 & 2.05 & 2.04 & 1.94\\
\cline{1-18}
 &  & 500 & 0.12 & 0.24 & 0.01 & 0.72 & 0.74 & 0.69 & 0.71 & 0.70 & 0.69 & 0.94 & 0.92 & 0.93 & 2.81 & 2.73 & 2.60\\

 &  & 1000 & 0.05 & 0.19 & 0.05 & 0.64 & 0.66 & 0.60 & 0.63 & 0.63 & 0.60 & 0.95 & 0.93 & 0.95 & 2.42 & 2.39 & 2.26\\

 &  & 1500 & 0.04 & 0.15 & 0.02 & 0.57 & 0.58 & 0.55 & 0.57 & 0.56 & 0.55 & 0.96 & 0.94 & 0.95 & 2.19 & 2.16 & 2.05\\

\multirow{-4}{*}{\centering 0.10} & \multirow{-4}{*}{\centering 1.49} & 2000 & 0.01 & 0.10 & 0.00 & 0.54 & 0.55 & 0.52 & 0.54 & 0.54 & 0.52 & 0.95 & 0.94 & 0.93 & 2.05 & 2.01 & 1.92\\
\cline{1-18}
 &  & 500 & 0.08 & 0.21 & 0.00 & 0.73 & 0.74 & 0.68 & 0.72 & 0.71 & 0.68 & 0.94 & 0.93 & 0.94 & 2.79 & 2.73 & 2.59\\

 &  & 1000 & 0.04 & 0.17 & 0.03 & 0.62 & 0.64 & 0.58 & 0.62 & 0.62 & 0.58 & 0.95 & 0.94 & 0.95 & 2.41 & 2.38 & 2.25\\

 &  & 1500 & 0.03 & 0.15 & 0.02 & 0.56 & 0.57 & 0.52 & 0.56 & 0.55 & 0.52 & 0.95 & 0.93 & 0.94 & 2.18 & 2.14 & 2.04\\

\multirow{-4}{*}{\centering 0.25} & \multirow{-4}{*}{\centering 1.45} & 2000 & 0.02 & 0.09 & 0.01 & 0.50 & 0.51 & 0.49 & 0.50 & 0.50 & 0.49 & 0.96 & 0.95 & 0.95 & 2.05 & 2.02 & 1.92\\
\cline{1-18}
 &  & 500 & 0.09 & 0.27 & 0.02 & 0.71 & 0.74 & 0.65 & 0.70 & 0.69 & 0.65 & 0.96 & 0.94 & 0.95 & 2.82 & 2.72 & 2.61\\

 &  & 1000 & 0.03 & 0.17 & 0.00 & 0.60 & 0.62 & 0.58 & 0.60 & 0.59 & 0.58 & 0.97 & 0.95 & 0.96 & 2.42 & 2.35 & 2.23\\

 &  & 1500 & 0.04 & 0.17 & 0.03 & 0.59 & 0.61 & 0.56 & 0.58 & 0.58 & 0.56 & 0.95 & 0.94 & 0.94 & 2.19 & 2.14 & 2.04\\

\multirow{-4}{*}{\centering 1.00} & \multirow{-4}{*}{\centering 0.81} & 2000 & 0.02 & 0.14 & 0.03 & 0.52 & 0.53 & 0.48 & 0.52 & 0.51 & 0.48 & 0.95 & 0.94 & 0.95 & 2.05 & 2.02 & 1.93\\
\hline

\multicolumn{18}{c}{Setting 3, exactly sparse: $f(d) = 2\exp(-d/2)$}\\
\hline
\multicolumn{1}{|c}{ } & \multicolumn{1}{c}{ } & \multicolumn{1}{c|}{ } & \multicolumn{3}{c|}{Bias} & \multicolumn{3}{c|}{RMSE} & \multicolumn{3}{c|}{SE} & \multicolumn{3}{c|}{Coverage} & \multicolumn{3}{c|}{CI Length} \\
\hline
$\eval$ & True & $n$ & \texttt{DLL} & \texttt{Plug} & \texttt{Orac} & \texttt{DLL} & \texttt{Plug} & \texttt{Orac} & \texttt{DLL} & \texttt{Plug} & \texttt{Orac} & \texttt{DLL} & \texttt{Plug} & \texttt{Orac} & \texttt{DLL} & \texttt{Plug} & \texttt{Orac}\\
\hline
 &  & 500 & 0.07 & 0.15 & 0.02 & 0.46 & 0.47 & 0.45 & 0.45 & 0.44 & 0.45 & 0.96 & 0.92 & 0.94 & 1.80 & 1.76 & 1.74\\

 &  & 1000 & 0.05 & 0.11 & 0.01 & 0.43 & 0.43 & 0.42 & 0.42 & 0.41 & 0.42 & 0.93 & 0.92 & 0.92 & 1.55 & 1.53 & 1.51\\

 &  & 1500 & 0.02 & 0.11 & 0.01 & 0.38 & 0.39 & 0.37 & 0.38 & 0.37 & 0.37 & 0.95 & 0.92 & 0.93 & 1.41 & 1.39 & 1.37\\

\multirow{-4}{*}{\centering -1.25} & \multirow{-4}{*}{\centering -1.87} & 2000 & 0.04 & 0.08 & 0.02 & 0.36 & 0.36 & 0.34 & 0.35 & 0.36 & 0.34 & 0.94 & 0.93 & 0.96 & 1.33 & 1.32 & 1.29\\
\cline{1-18}
 &  & 500 & 0.00 & 0.13 & 0.03 & 0.46 & 0.47 & 0.45 & 0.46 & 0.46 & 0.45 & 0.95 & 0.93 & 0.94 & 1.78 & 1.75 & 1.75\\

 &  & 1000 & 0.01 & 0.13 & 0.01 & 0.40 & 0.42 & 0.38 & 0.40 & 0.40 & 0.38 & 0.95 & 0.92 & 0.95 & 1.55 & 1.52 & 1.50\\

 &  & 1500 & 0.01 & 0.11 & 0.01 & 0.36 & 0.37 & 0.36 & 0.36 & 0.35 & 0.36 & 0.95 & 0.95 & 0.94 & 1.41 & 1.39 & 1.37\\

\multirow{-4}{*}{\centering -0.50} & \multirow{-4}{*}{\centering -1.28} & 2000 & 0.02 & 0.08 & 0.00 & 0.35 & 0.35 & 0.33 & 0.35 & 0.34 & 0.33 & 0.95 & 0.95 & 0.95 & 1.33 & 1.32 & 1.29\\
\cline{1-18}
 &  & 500 & 0.05 & 0.17 & 0.01 & 0.50 & 0.52 & 0.46 & 0.49 & 0.49 & 0.46 & 0.94 & 0.91 & 0.94 & 1.79 & 1.75 & 1.76\\

 &  & 1000 & 0.01 & 0.12 & 0.01 & 0.41 & 0.42 & 0.40 & 0.41 & 0.40 & 0.40 & 0.94 & 0.94 & 0.93 & 1.54 & 1.52 & 1.51\\

 &  & 1500 & 0.03 & 0.08 & 0.01 & 0.36 & 0.36 & 0.35 & 0.35 & 0.35 & 0.35 & 0.95 & 0.94 & 0.95 & 1.41 & 1.39 & 1.38\\

\multirow{-4}{*}{\centering 0.10} & \multirow{-4}{*}{\centering -0.95} & 2000 & 0.01 & 0.10 & 0.01 & 0.34 & 0.36 & 0.34 & 0.34 & 0.34 & 0.34 & 0.94 & 0.93 & 0.93 & 1.33 & 1.31 & 1.28\\
\cline{1-18}
 &  & 500 & 0.07 & 0.20 & 0.04 & 0.48 & 0.51 & 0.47 & 0.48 & 0.47 & 0.47 & 0.94 & 0.92 & 0.94 & 1.79 & 1.75 & 1.76\\

 &  & 1000 & 0.03 & 0.15 & 0.02 & 0.38 & 0.41 & 0.38 & 0.38 & 0.38 & 0.38 & 0.97 & 0.95 & 0.95 & 1.54 & 1.52 & 1.50\\

 &  & 1500 & 0.02 & 0.09 & 0.01 & 0.36 & 0.37 & 0.34 & 0.36 & 0.36 & 0.34 & 0.95 & 0.93 & 0.95 & 1.41 & 1.39 & 1.38\\

\multirow{-4}{*}{\centering 0.25} & \multirow{-4}{*}{\centering -0.88} & 2000 & 0.00 & 0.11 & 0.02 & 0.34 & 0.35 & 0.32 & 0.34 & 0.33 & 0.32 & 0.95 & 0.94 & 0.95 & 1.33 & 1.31 & 1.29\\
\cline{1-18}
 &  & 500 & 0.04 & 0.14 & 0.03 & 0.47 & 0.48 & 0.46 & 0.47 & 0.46 & 0.46 & 0.94 & 0.91 & 0.95 & 1.78 & 1.73 & 1.77\\

 &  & 1000 & 0.02 & 0.12 & 0.01 & 0.41 & 0.42 & 0.39 & 0.41 & 0.40 & 0.39 & 0.94 & 0.92 & 0.95 & 1.54 & 1.50 & 1.50\\

 &  & 1500 & 0.01 & 0.11 & 0.02 & 0.37 & 0.38 & 0.36 & 0.37 & 0.37 & 0.36 & 0.95 & 0.93 & 0.95 & 1.41 & 1.38 & 1.37\\

\multirow{-4}{*}{\centering 1.00} & \multirow{-4}{*}{\centering -0.61} & 2000 & 0.00 & 0.10 & 0.02 & 0.36 & 0.37 & 0.34 & 0.36 & 0.35 & 0.34 & 0.94 & 0.92 & 0.94 & 1.33 & 1.31 & 1.29\\
\hline
\end{tabular}
}
\caption*{Table S4: Comparison of \texttt{DLL}, plug-in (\texttt{Plug}), oracle (\texttt{Orac}) estimators in Setting 3, across different sample sizes $n$ and evaluation points $\eval$. The column indexed with ``True'' represents the true value of $f'(\eval)$. The columns indexed with ``Bias'', ``RMSE'' and ``SE" report the absolute bias, the root mean square error, and the standard error computed by 500 estimates, respectively; the columns indexed with ``Coverage'' report the empirical coverage level and the columns indexed with ``Length'' report the average CI length.}
\label{table:S4}
\end{table}

\begin{table}[htb!]
\centering
\resizebox{\linewidth}{!}{
\begin{tabular}[t]{|c|c|c|ccc|ccc|ccc|ccc|ccc|}
\multicolumn{18}{c}{Setting 4, exactly sparse: $f(d) = 1.5\sin(d)$ and df=10}\\
\hline
\multicolumn{1}{|c}{ } & \multicolumn{1}{c}{ } & \multicolumn{1}{c|}{ } & \multicolumn{3}{c|}{Bias} & \multicolumn{3}{c|}{RMSE} & \multicolumn{3}{c|}{SE} & \multicolumn{3}{c|}{Coverage} & \multicolumn{3}{c|}{CI Length} \\
\hline
$\eval$ & True & $n$ & \texttt{DLL} & \texttt{Plug} & \texttt{Orac} & \texttt{DLL} & \texttt{Plug} & \texttt{Orac} & \texttt{DLL} & \texttt{Plug} & \texttt{Orac} & \texttt{DLL} & \texttt{Plug} & \texttt{Orac} & \texttt{DLL} & \texttt{Plug} & \texttt{Orac}\\
\hline
 &  & 500 & 0.13 & 0.27 & 0.00 & 0.50 & 0.55 & 0.49 & 0.48 & 0.48 & 0.49 & 0.94 & 0.88 & 0.94 & 1.84 & 1.76 & 1.82\\

 &  & 1000 & 0.04 & 0.13 & 0.00 & 0.41 & 0.43 & 0.41 & 0.41 & 0.41 & 0.41 & 0.95 & 0.92 & 0.92 & 1.55 & 1.47 & 1.46\\

 &  & 1500 & 0.01 & 0.04 & 0.03 & 0.33 & 0.33 & 0.31 & 0.33 & 0.33 & 0.31 & 0.97 & 0.95 & 0.97 & 1.37 & 1.31 & 1.27\\

\multirow{-4}{*}{\centering -1.25} & \multirow{-4}{*}{\centering 0.47} & 2000 & 0.02 & 0.00 & 0.04 & 0.34 & 0.33 & 0.32 & 0.33 & 0.33 & 0.32 & 0.94 & 0.92 & 0.93 & 1.28 & 1.21 & 1.17\\
\cline{1-18}
 &  & 500 & 0.16 & 0.33 & 0.01 & 0.39 & 0.49 & 0.36 & 0.36 & 0.36 & 0.36 & 0.92 & 0.82 & 0.95 & 1.41 & 1.35 & 1.41\\

 &  & 1000 & 0.08 & 0.22 & 0.01 & 0.31 & 0.36 & 0.30 & 0.30 & 0.29 & 0.30 & 0.96 & 0.86 & 0.94 & 1.19 & 1.11 & 1.11\\

 &  & 1500 & 0.05 & 0.16 & 0.01 & 0.27 & 0.31 & 0.25 & 0.27 & 0.27 & 0.25 & 0.94 & 0.88 & 0.95 & 1.05 & 0.98 & 0.97\\

\multirow{-4}{*}{\centering -0.50} & \multirow{-4}{*}{\centering 1.32} & 2000 & 0.03 & 0.13 & 0.01 & 0.23 & 0.26 & 0.22 & 0.23 & 0.23 & 0.22 & 0.96 & 0.91 & 0.95 & 0.98 & 0.92 & 0.90\\
\cline{1-18}
 &  & 500 & 0.17 & 0.36 & 0.04 & 0.38 & 0.49 & 0.33 & 0.34 & 0.33 & 0.32 & 0.92 & 0.78 & 0.95 & 1.33 & 1.27 & 1.32\\

 &  & 1000 & 0.10 & 0.28 & 0.03 & 0.30 & 0.40 & 0.28 & 0.29 & 0.28 & 0.28 & 0.94 & 0.80 & 0.95 & 1.13 & 1.06 & 1.06\\

 &  & 1500 & 0.06 & 0.23 & 0.02 & 0.25 & 0.34 & 0.23 & 0.24 & 0.24 & 0.23 & 0.96 & 0.83 & 0.97 & 1.00 & 0.94 & 0.92\\

\multirow{-4}{*}{\centering 0.10} & \multirow{-4}{*}{\centering 1.49} & 2000 & 0.04 & 0.21 & 0.01 & 0.24 & 0.31 & 0.23 & 0.23 & 0.23 & 0.23 & 0.95 & 0.83 & 0.93 & 0.93 & 0.88 & 0.85\\
\cline{1-18}
 &  & 500 & 0.18 & 0.38 & 0.05 & 0.40 & 0.51 & 0.36 & 0.35 & 0.35 & 0.35 & 0.89 & 0.76 & 0.94 & 1.34 & 1.28 & 1.34\\

 &  & 1000 & 0.10 & 0.29 & 0.01 & 0.30 & 0.41 & 0.28 & 0.29 & 0.28 & 0.28 & 0.93 & 0.78 & 0.93 & 1.14 & 1.07 & 1.08\\

 &  & 1500 & 0.06 & 0.24 & 0.01 & 0.25 & 0.35 & 0.24 & 0.24 & 0.24 & 0.24 & 0.96 & 0.83 & 0.96 & 1.01 & 0.95 & 0.94\\

\multirow{-4}{*}{\centering 0.25} & \multirow{-4}{*}{\centering 1.45} & 2000 & 0.07 & 0.24 & 0.03 & 0.24 & 0.33 & 0.21 & 0.23 & 0.23 & 0.21 & 0.95 & 0.80 & 0.97 & 0.94 & 0.89 & 0.86\\
\cline{1-18}
 &  & 500 & 0.12 & 0.30 & 0.04 & 0.43 & 0.50 & 0.43 & 0.42 & 0.41 & 0.43 & 0.96 & 0.86 & 0.94 & 1.64 & 1.54 & 1.64\\

 &  & 1000 & 0.09 & 0.26 & 0.01 & 0.35 & 0.42 & 0.33 & 0.34 & 0.33 & 0.33 & 0.96 & 0.89 & 0.94 & 1.38 & 1.31 & 1.31\\

 &  & 1500 & 0.05 & 0.22 & 0.01 & 0.32 & 0.38 & 0.30 & 0.31 & 0.31 & 0.30 & 0.94 & 0.86 & 0.95 & 1.23 & 1.16 & 1.14\\

\multirow{-4}{*}{\centering 1.00} & \multirow{-4}{*}{\centering 0.81} & 2000 & 0.05 & 0.22 & 0.00 & 0.28 & 0.35 & 0.27 & 0.28 & 0.28 & 0.27 & 0.97 & 0.87 & 0.93 & 1.15 & 1.08 & 1.05\\
\hline

\multicolumn{18}{c}{Setting 4, exactly sparse: $f(d) = 2\exp(-d/2)$ and df=10}\\
\hline
\multicolumn{1}{|c}{ } & \multicolumn{1}{c}{ } & \multicolumn{1}{c|}{ } & \multicolumn{3}{c|}{Bias} & \multicolumn{3}{c|}{RMSE} & \multicolumn{3}{c|}{SE} & \multicolumn{3}{c|}{Coverage} & \multicolumn{3}{c|}{CI Length} \\
\hline
$\eval$ & True & $n$ & \texttt{DLL} & \texttt{Plug} & \texttt{Orac} & \texttt{DLL} & \texttt{Plug} & \texttt{Orac} & \texttt{DLL} & \texttt{Plug} & \texttt{Orac} & \texttt{DLL} & \texttt{Plug} & \texttt{Orac} & \texttt{DLL} & \texttt{Plug} & \texttt{Orac}\\
\hline
 &  & 500 & 0.05 & 0.19 & 0.04 & 0.45 & 0.48 & 0.46 & 0.44 & 0.44 & 0.45 & 0.95 & 0.93 & 0.95 & 1.80 & 1.70 & 1.81\\

 &  & 1000 & 0.06 & 0.20 & 0.04 & 0.41 & 0.44 & 0.40 & 0.40 & 0.39 & 0.40 & 0.95 & 0.90 & 0.95 & 1.61 & 1.52 & 1.58\\

 &  & 1500 & 0.02 & 0.15 & 0.02 & 0.37 & 0.40 & 0.37 & 0.37 & 0.37 & 0.37 & 0.95 & 0.92 & 0.94 & 1.49 & 1.41 & 1.45\\

\multirow{-4}{*}{\centering -1.25} & \multirow{-4}{*}{\centering -1.87} & 2000 & 0.03 & 0.15 & 0.02 & 0.33 & 0.36 & 0.34 & 0.33 & 0.33 & 0.34 & 0.97 & 0.93 & 0.96 & 1.41 & 1.34 & 1.38\\
\cline{1-18}
 &  & 500 & 0.03 & 0.22 & 0.04 & 0.35 & 0.42 & 0.36 & 0.35 & 0.35 & 0.36 & 0.94 & 0.89 & 0.94 & 1.38 & 1.31 & 1.39\\

 &  & 1000 & 0.03 & 0.20 & 0.01 & 0.31 & 0.36 & 0.30 & 0.31 & 0.30 & 0.30 & 0.96 & 0.89 & 0.95 & 1.23 & 1.16 & 1.20\\

 &  & 1500 & 0.03 & 0.18 & 0.01 & 0.29 & 0.34 & 0.29 & 0.29 & 0.29 & 0.29 & 0.96 & 0.90 & 0.94 & 1.14 & 1.08 & 1.11\\

\multirow{-4}{*}{\centering -0.50} & \multirow{-4}{*}{\centering -1.28} & 2000 & 0.02 & 0.16 & 0.01 & 0.26 & 0.30 & 0.27 & 0.26 & 0.25 & 0.27 & 0.97 & 0.92 & 0.95 & 1.08 & 1.01 & 1.05\\
\cline{1-18}
 &  & 500 & 0.07 & 0.26 & 0.01 & 0.34 & 0.42 & 0.33 & 0.34 & 0.33 & 0.33 & 0.94 & 0.88 & 0.95 & 1.31 & 1.25 & 1.32\\

 &  & 1000 & 0.04 & 0.19 & 0.00 & 0.28 & 0.34 & 0.28 & 0.28 & 0.28 & 0.28 & 0.96 & 0.89 & 0.96 & 1.17 & 1.10 & 1.14\\

 &  & 1500 & 0.01 & 0.15 & 0.00 & 0.27 & 0.30 & 0.28 & 0.27 & 0.26 & 0.28 & 0.96 & 0.90 & 0.94 & 1.08 & 1.01 & 1.05\\

\multirow{-4}{*}{\centering 0.10} & \multirow{-4}{*}{\centering -0.95} & 2000 & 0.03 & 0.15 & 0.01 & 0.27 & 0.31 & 0.27 & 0.27 & 0.26 & 0.27 & 0.94 & 0.87 & 0.93 & 1.02 & 0.96 & 1.00\\
\cline{1-18}
 &  & 500 & 0.08 & 0.26 & 0.01 & 0.34 & 0.41 & 0.32 & 0.33 & 0.32 & 0.32 & 0.93 & 0.86 & 0.96 & 1.32 & 1.26 & 1.33\\

 &  & 1000 & 0.02 & 0.17 & 0.01 & 0.30 & 0.34 & 0.30 & 0.30 & 0.29 & 0.30 & 0.96 & 0.90 & 0.94 & 1.18 & 1.11 & 1.15\\

 &  & 1500 & 0.01 & 0.14 & 0.00 & 0.26 & 0.29 & 0.27 & 0.26 & 0.26 & 0.27 & 0.97 & 0.91 & 0.96 & 1.09 & 1.03 & 1.06\\

\multirow{-4}{*}{\centering 0.25} & \multirow{-4}{*}{\centering -0.88} & 2000 & 0.01 & 0.14 & 0.01 & 0.25 & 0.29 & 0.26 & 0.25 & 0.25 & 0.26 & 0.96 & 0.91 & 0.95 & 1.04 & 0.97 & 1.01\\
\cline{1-18}
 &  & 500 & 0.01 & 0.13 & 0.03 & 0.41 & 0.42 & 0.42 & 0.41 & 0.40 & 0.42 & 0.95 & 0.92 & 0.94 & 1.61 & 1.49 & 1.61\\

 &  & 1000 & 0.01 & 0.09 & 0.01 & 0.35 & 0.36 & 0.36 & 0.35 & 0.35 & 0.36 & 0.98 & 0.93 & 0.93 & 1.44 & 1.35 & 1.41\\

 &  & 1500 & 0.00 & 0.08 & 0.01 & 0.32 & 0.33 & 0.33 & 0.32 & 0.32 & 0.33 & 0.96 & 0.95 & 0.95 & 1.33 & 1.25 & 1.29\\

\multirow{-4}{*}{\centering 1.00} & \multirow{-4}{*}{\centering -0.61} & 2000 & 0.01 & 0.06 & 0.01 & 0.30 & 0.30 & 0.30 & 0.30 & 0.29 & 0.30 & 0.96 & 0.94 & 0.94 & 1.26 & 1.18 & 1.22\\
\hline
\end{tabular}
}
\caption*{Table S5: Comparison of \texttt{DLL}, plug-in (\texttt{Plug}), oracle (\texttt{Orac}) estimators in Setting 4 with df = 10, across different sample sizes $n$ and evaluation points $\eval$. The column indexed with ``True'' represents the true value of $f'(\eval)$. The columns indexed with ``Bias'', ``RMSE'' and ``SE" report the absolute bias, the root mean square error, and the standard error computed by 500 estimates, respectively; the columns indexed with ``Coverage'' report the empirical coverage level and the columns indexed with ``Length'' report the average CI length.}
\label{table:S5}
\end{table}

\begin{table}[htb!]
\centering
\resizebox{\linewidth}{!}{
\begin{tabular}[t]{|c|c|c|ccc|ccc|ccc|ccc|ccc|}
\multicolumn{18}{c}{Setting 4, exactly sparse: $f(d) = 1.5\sin(d)$ and df=15}\\
\hline
\multicolumn{1}{|c}{ } & \multicolumn{1}{c}{ } & \multicolumn{1}{c|}{ } & \multicolumn{3}{c|}{Bias} & \multicolumn{3}{c|}{RMSE} & \multicolumn{3}{c|}{SE} & \multicolumn{3}{c|}{Coverage} & \multicolumn{3}{c|}{CI Length} \\
\hline
$\eval$ & True & $n$ & \texttt{DLL} & \texttt{Plug} & \texttt{Orac} & \texttt{DLL} & \texttt{Plug} & \texttt{Orac} & \texttt{DLL} & \texttt{Plug} & \texttt{Orac} & \texttt{DLL} & \texttt{Plug} & \texttt{Orac} & \texttt{DLL} & \texttt{Plug} & \texttt{Orac}\\
\hline
 &  & 500 & 0.17 & 0.36 & 0.01 & 0.51 & 0.59 & 0.48 & 0.48 & 0.47 & 0.48 & 0.94 & 0.88 & 0.97 & 1.97 & 1.91 & 1.99\\

 &  & 1000 & 0.09 & 0.23 & 0.04 & 0.45 & 0.49 & 0.41 & 0.44 & 0.43 & 0.41 & 0.94 & 0.90 & 0.95 & 1.69 & 1.62 & 1.62\\

 &  & 1500 & 0.01 & 0.11 & 0.03 & 0.39 & 0.40 & 0.38 & 0.39 & 0.38 & 0.38 & 0.96 & 0.94 & 0.94 & 1.53 & 1.48 & 1.45\\

\multirow{-4}{*}{\centering -1.25} & \multirow{-4}{*}{\centering 0.47} & 2000 & 0.00 & 0.08 & 0.02 & 0.33 & 0.34 & 0.33 & 0.33 & 0.33 & 0.33 & 0.97 & 0.96 & 0.95 & 1.41 & 1.36 & 1.34\\
\cline{1-18}
 &  & 500 & 0.19 & 0.39 & 0.05 & 0.41 & 0.52 & 0.37 & 0.36 & 0.35 & 0.37 & 0.93 & 0.78 & 0.94 & 1.47 & 1.41 & 1.47\\

 &  & 1000 & 0.07 & 0.24 & 0.02 & 0.32 & 0.38 & 0.31 & 0.31 & 0.30 & 0.31 & 0.95 & 0.87 & 0.95 & 1.26 & 1.21 & 1.22\\

 &  & 1500 & 0.04 & 0.18 & 0.01 & 0.29 & 0.33 & 0.28 & 0.28 & 0.28 & 0.28 & 0.96 & 0.91 & 0.95 & 1.14 & 1.09 & 1.09\\

\multirow{-4}{*}{\centering -0.50} & \multirow{-4}{*}{\centering 1.32} & 2000 & 0.02 & 0.14 & 0.01 & 0.26 & 0.29 & 0.25 & 0.26 & 0.25 & 0.25 & 0.97 & 0.93 & 0.95 & 1.06 & 1.01 & 1.00\\
\cline{1-18}
 &  & 500 & 0.15 & 0.36 & 0.00 & 0.40 & 0.51 & 0.35 & 0.36 & 0.36 & 0.35 & 0.92 & 0.77 & 0.96 & 1.39 & 1.33 & 1.39\\

 &  & 1000 & 0.07 & 0.27 & 0.02 & 0.31 & 0.40 & 0.30 & 0.30 & 0.30 & 0.30 & 0.95 & 0.82 & 0.95 & 1.20 & 1.15 & 1.15\\

 &  & 1500 & 0.07 & 0.26 & 0.03 & 0.27 & 0.36 & 0.26 & 0.26 & 0.25 & 0.26 & 0.96 & 0.87 & 0.94 & 1.08 & 1.04 & 1.03\\

\multirow{-4}{*}{\centering 0.10} & \multirow{-4}{*}{\centering 1.49} & 2000 & 0.06 & 0.23 & 0.03 & 0.25 & 0.33 & 0.24 & 0.25 & 0.24 & 0.23 & 0.95 & 0.84 & 0.95 & 1.01 & 0.96 & 0.95\\
\cline{1-18}
 &  & 500 & 0.17 & 0.38 & 0.02 & 0.40 & 0.52 & 0.34 & 0.36 & 0.36 & 0.34 & 0.90 & 0.77 & 0.95 & 1.41 & 1.35 & 1.42\\

 &  & 1000 & 0.07 & 0.27 & 0.01 & 0.32 & 0.41 & 0.30 & 0.31 & 0.31 & 0.30 & 0.95 & 0.84 & 0.95 & 1.21 & 1.16 & 1.17\\

 &  & 1500 & 0.04 & 0.23 & 0.01 & 0.26 & 0.35 & 0.25 & 0.26 & 0.26 & 0.25 & 0.97 & 0.87 & 0.96 & 1.09 & 1.05 & 1.04\\

\multirow{-4}{*}{\centering 0.25} & \multirow{-4}{*}{\centering 1.45} & 2000 & 0.04 & 0.22 & 0.01 & 0.26 & 0.34 & 0.25 & 0.26 & 0.26 & 0.25 & 0.94 & 0.88 & 0.95 & 1.02 & 0.97 & 0.96\\
\cline{1-18}
 &  & 500 & 0.16 & 0.34 & 0.03 & 0.46 & 0.55 & 0.46 & 0.43 & 0.43 & 0.46 & 0.95 & 0.85 & 0.94 & 1.74 & 1.64 & 1.73\\

 &  & 1000 & 0.07 & 0.26 & 0.00 & 0.39 & 0.46 & 0.38 & 0.39 & 0.38 & 0.38 & 0.94 & 0.87 & 0.94 & 1.49 & 1.43 & 1.44\\

 &  & 1500 & 0.05 & 0.23 & 0.00 & 0.34 & 0.40 & 0.33 & 0.34 & 0.33 & 0.34 & 0.95 & 0.90 & 0.94 & 1.35 & 1.29 & 1.28\\

\multirow{-4}{*}{\centering 1.00} & \multirow{-4}{*}{\centering 0.81} & 2000 & 0.04 & 0.20 & 0.00 & 0.31 & 0.36 & 0.30 & 0.31 & 0.30 & 0.30 & 0.96 & 0.90 & 0.95 & 1.25 & 1.20 & 1.18\\
\hline
\multicolumn{18}{c}{Setting 4, exactly sparse: $f(d) = 2\exp(-d/2)$ and df=15}\\
\hline
\multicolumn{1}{|c}{ } & \multicolumn{1}{c}{ } & \multicolumn{1}{c|}{ } & \multicolumn{3}{c|}{Bias} & \multicolumn{3}{c|}{RMSE} & \multicolumn{3}{c|}{SE} & \multicolumn{3}{c|}{Coverage} & \multicolumn{3}{c|}{CI Length} \\
\hline
$\eval$ & True & $n$ & \texttt{DLL} & \texttt{Plug} & \texttt{Orac} & \texttt{DLL} & \texttt{Plug} & \texttt{Orac} & \texttt{DLL} & \texttt{Plug} & \texttt{Orac} & \texttt{DLL} & \texttt{Plug} & \texttt{Orac} & \texttt{DLL} & \texttt{Plug} & \texttt{Orac}\\
\hline
 &  & 500 & 0.04 & 0.21 & 0.02 & 0.46 & 0.50 & 0.44 & 0.46 & 0.45 & 0.44 & 0.93 & 0.89 & 0.95 & 1.78 & 1.70 & 1.73\\

 &  & 1000 & 0.01 & 0.16 & 0.00 & 0.40 & 0.42 & 0.40 & 0.40 & 0.39 & 0.40 & 0.96 & 0.92 & 0.93 & 1.55 & 1.49 & 1.50\\

 &  & 1500 & 0.01 & 0.16 & 0.01 & 0.36 & 0.39 & 0.35 & 0.36 & 0.36 & 0.35 & 0.96 & 0.91 & 0.96 & 1.44 & 1.39 & 1.38\\

\multirow{-4}{*}{\centering -1.25} & \multirow{-4}{*}{\centering -1.87} & 2000 & 0.01 & 0.14 & 0.01 & 0.32 & 0.35 & 0.31 & 0.32 & 0.32 & 0.31 & 0.96 & 0.94 & 0.96 & 1.34 & 1.30 & 1.30\\
\cline{1-18}
 &  & 500 & 0.06 & 0.27 & 0.01 & 0.33 & 0.42 & 0.33 & 0.33 & 0.33 & 0.33 & 0.95 & 0.87 & 0.94 & 1.34 & 1.31 & 1.32\\

 &  & 1000 & 0.00 & 0.18 & 0.01 & 0.30 & 0.35 & 0.29 & 0.30 & 0.30 & 0.29 & 0.95 & 0.89 & 0.95 & 1.16 & 1.11 & 1.12\\

 &  & 1500 & 0.00 & 0.16 & 0.01 & 0.29 & 0.32 & 0.28 & 0.29 & 0.28 & 0.28 & 0.93 & 0.90 & 0.95 & 1.07 & 1.04 & 1.04\\

\multirow{-4}{*}{\centering -0.50} & \multirow{-4}{*}{\centering -1.28} & 2000 & 0.01 & 0.16 & 0.01 & 0.26 & 0.31 & 0.26 & 0.26 & 0.26 & 0.26 & 0.95 & 0.92 & 0.94 & 1.01 & 0.97 & 0.97\\
\cline{1-18}
 &  & 500 & 0.07 & 0.27 & 0.01 & 0.32 & 0.41 & 0.30 & 0.31 & 0.31 & 0.30 & 0.94 & 0.85 & 0.96 & 1.27 & 1.22 & 1.25\\

 &  & 1000 & 0.02 & 0.19 & 0.00 & 0.28 & 0.33 & 0.28 & 0.28 & 0.27 & 0.28 & 0.94 & 0.90 & 0.92 & 1.10 & 1.05 & 1.06\\

 &  & 1500 & 0.01 & 0.15 & 0.01 & 0.26 & 0.30 & 0.25 & 0.26 & 0.26 & 0.25 & 0.95 & 0.91 & 0.95 & 1.01 & 0.97 & 0.98\\

\multirow{-4}{*}{\centering 0.10} & \multirow{-4}{*}{\centering -0.95} & 2000 & 0.01 & 0.12 & 0.01 & 0.24 & 0.27 & 0.24 & 0.24 & 0.24 & 0.24 & 0.95 & 0.92 & 0.94 & 0.96 & 0.92 & 0.92\\
\cline{1-18}
 &  & 500 & 0.06 & 0.26 & 0.01 & 0.33 & 0.41 & 0.32 & 0.33 & 0.32 & 0.32 & 0.95 & 0.87 & 0.93 & 1.29 & 1.23 & 1.27\\

 &  & 1000 & 0.02 & 0.18 & 0.00 & 0.28 & 0.33 & 0.27 & 0.28 & 0.28 & 0.27 & 0.95 & 0.90 & 0.95 & 1.11 & 1.07 & 1.07\\

 &  & 1500 & 0.01 & 0.14 & 0.00 & 0.24 & 0.28 & 0.24 & 0.24 & 0.24 & 0.24 & 0.97 & 0.92 & 0.96 & 1.03 & 0.98 & 0.99\\

\multirow{-4}{*}{\centering 0.25} & \multirow{-4}{*}{\centering -0.88} & 2000 & 0.02 & 0.14 & 0.02 & 0.25 & 0.28 & 0.25 & 0.25 & 0.25 & 0.25 & 0.94 & 0.91 & 0.94 & 0.97 & 0.93 & 0.93\\
\cline{1-18}
 &  & 500 & 0.00 & 0.15 & 0.05 & 0.42 & 0.44 & 0.41 & 0.42 & 0.41 & 0.41 & 0.93 & 0.92 & 0.92 & 1.59 & 1.49 & 1.56\\

 &  & 1000 & 0.01 & 0.13 & 0.00 & 0.31 & 0.34 & 0.32 & 0.32 & 0.31 & 0.32 & 0.97 & 0.95 & 0.96 & 1.38 & 1.31 & 1.33\\

 &  & 1500 & 0.01 & 0.09 & 0.01 & 0.30 & 0.32 & 0.30 & 0.30 & 0.30 & 0.30 & 0.97 & 0.95 & 0.96 & 1.27 & 1.21 & 1.22\\

\multirow{-4}{*}{\centering 1.00} & \multirow{-4}{*}{\centering -0.61} & 2000 & 0.01 & 0.08 & 0.00 & 0.29 & 0.30 & 0.29 & 0.29 & 0.29 & 0.29 & 0.95 & 0.93 & 0.94 & 1.19 & 1.14 & 1.15\\
\hline
\end{tabular}
}
\caption*{Table S6: Comparison of \texttt{DLL}, plug-in (\texttt{Plug}), oracle (\texttt{Orac}) estimators in Setting 4 with df = 15, across different sample sizes $n$ and evaluation points $\eval$. The column indexed with ``True'' represents the true value of $f'(\eval)$. The columns indexed with ``Bias'', ``RMSE'' and ``SE" report the absolute bias, the root mean square error, and the standard error computed by 500 estimates, respectively; the columns indexed with ``Coverage'' report the empirical coverage level and the columns indexed with ``Length'' report the average CI length.}
\label{table:S6}
\end{table}

\begin{table}[htb!]
\centering
\scalebox{0.75}{
\begin{tabular}[t]{|c|c|c|cccc|cccc|cccc|}
\multicolumn{15}{c}{Non-linear Treatment Model, exactly sparse: $f(d) = 1.5\sin(d)$}\\
\hline
\multicolumn{1}{|c}{ } & \multicolumn{1}{c}{ } & \multicolumn{1}{c|}{ } & \multicolumn{4}{c|}{Bias} & \multicolumn{4}{c|}{Coverage} & \multicolumn{4}{c|}{Length} \\
\hline
$\eval$ & True & $n$ & \texttt{DLL} & \texttt{DLL-S} & \texttt{Plug} & \texttt{Orac} & \texttt{DLL} & \texttt{DLL-S} & \texttt{Plug} & \texttt{Orac} & \texttt{DLL} & \texttt{DLL-S} & \texttt{Plug} & \texttt{Orac}\\
\hline
 &  & 500 & 0.29 & 0.12 & 0.33 & 0.01 & 0.94 & 0.95 & 0.90 & 0.90 & 4.36 & 4.27 & 3.41 & 3.60\\

 &  & 1000 & 0.21 & 0.13 & 0.22 & 0.06 & 0.94 & 0.96 & 0.91 & 0.93 & 3.55 & 3.61 & 3.31 & 3.51\\

 &  & 1500 & 0.06 & 0.11 & 0.07 & 0.15 & 0.94 & 0.96 & 0.91 & 0.92 & 3.24 & 3.24 & 3.10 & 3.25\\

\multirow{-4}{*}{\centering -1.25} & \multirow{-4}{*}{\centering 0.47} & 2000 & 0.18 & 0.12 & 0.19 & 0.01 & 0.96 & 0.95 & 0.93 & 0.94 & 2.98 & 3.02 & 2.84 & 2.95\\
\cline{1-15}
 &  & 500 & 0.32 & 0.26 & 0.36 & 0.01 & 0.90 & 0.93 & 0.87 & 0.95 & 2.04 & 2.08 & 1.97 & 2.20\\

 &  & 1000 & 0.24 & 0.15 & 0.27 & 0.00 & 0.93 & 0.94 & 0.90 & 0.94 & 1.78 & 1.80 & 1.71 & 1.81\\

 &  & 1500 & 0.18 & 0.10 & 0.21 & 0.02 & 0.94 & 0.94 & 0.91 & 0.96 & 1.61 & 1.65 & 1.55 & 1.61\\

\multirow{-4}{*}{\centering -0.50} & \multirow{-4}{*}{\centering 1.32} & 2000 & 0.14 & 0.08 & 0.17 & 0.01 & 0.92 & 0.95 & 0.90 & 0.95 & 1.50 & 1.52 & 1.45 & 1.51\\
\cline{1-15}
 &  & 500 & 0.31 & 0.23 & 0.37 & 0.01 & 0.84 & 0.91 & 0.78 & 0.95 & 1.40 & 1.44 & 1.34 & 1.50\\

 &  & 1000 & 0.20 & 0.12 & 0.25 & 0.00 & 0.87 & 0.93 & 0.83 & 0.95 & 1.24 & 1.24 & 1.18 & 1.25\\

 &  & 1500 & 0.18 & 0.10 & 0.23 & 0.03 & 0.92 & 0.95 & 0.88 & 0.96 & 1.12 & 1.14 & 1.07 & 1.12\\

\multirow{-4}{*}{\centering 0.10} & \multirow{-4}{*}{\centering 1.49} & 2000 & 0.13 & 0.08 & 0.17 & 0.00 & 0.91 & 0.94 & 0.87 & 0.94 & 1.05 & 1.05 & 1.01 & 1.05\\
\cline{1-15}
 &  & 500 & 0.33 & 0.23 & 0.39 & 0.06 & 0.80 & 0.87 & 0.74 & 0.94 & 1.34 & 1.36 & 1.28 & 1.42\\

 &  & 1000 & 0.19 & 0.12 & 0.24 & 0.01 & 0.90 & 0.92 & 0.85 & 0.95 & 1.17 & 1.18 & 1.12 & 1.19\\

 &  & 1500 & 0.16 & 0.10 & 0.21 & 0.02 & 0.92 & 0.95 & 0.87 & 0.95 & 1.06 & 1.07 & 1.02 & 1.06\\

\multirow{-4}{*}{\centering 0.25} & \multirow{-4}{*}{\centering 1.45} & 2000 & 0.13 & 0.08 & 0.18 & 0.00 & 0.91 & 0.95 & 0.87 & 0.95 & 1.00 & 0.99 & 0.95 & 0.99\\
\cline{1-15}
 &  & 500 & 0.20 & 0.19 & 0.26 & 0.02 & 0.90 & 0.92 & 0.86 & 0.95 & 1.31 & 1.32 & 1.23 & 1.38\\

 &  & 1000 & 0.16 & 0.13 & 0.21 & 0.01 & 0.92 & 0.93 & 0.87 & 0.94 & 1.15 & 1.14 & 1.10 & 1.16\\

 &  & 1500 & 0.13 & 0.11 & 0.17 & 0.01 & 0.93 & 0.93 & 0.86 & 0.96 & 1.04 & 1.05 & 0.99 & 1.04\\

\multirow{-4}{*}{\centering 1.00} & \multirow{-4}{*}{\centering 0.81} & 2000 & 0.11 & 0.09 & 0.15 & 0.01 & 0.93 & 0.93 & 0.88 & 0.94 & 0.98 & 0.97 & 0.92 & 0.97\\
\hline

\multicolumn{15}{c}{Non-linear Treatment Model, exactly sparse: $f(d) = 2\exp(-d/2)$}\\
\hline
\multicolumn{1}{|c}{ } & \multicolumn{1}{c}{ } & \multicolumn{1}{c|}{ } & \multicolumn{4}{c|}{Bias} & \multicolumn{4}{c|}{Coverage} & \multicolumn{4}{c|}{Length} \\
\hline
$\eval$ & True & $n$ & \texttt{DLL} & \texttt{DLL-S} & \texttt{Plug} & \texttt{Orac} & \texttt{DLL} & \texttt{DLL-S} & \texttt{Plug} & \texttt{Orac} & \texttt{DLL} & \texttt{DLL-S} & \texttt{Plug} & \texttt{Orac}\\
\hline
 &  & 500 & 0.11 & 0.12 & 0.11 & 0.09 & 0.94 & 0.94 & 0.91 & 0.91 & 2.88 & 2.89 & 2.68 & 2.81\\

 &  & 1000 & 0.03 & 0.09 & 0.02 & 0.04 & 0.95 & 0.94 & 0.93 & 0.92 & 2.39 & 2.42 & 2.31 & 2.35\\

 &  & 1500 & 0.01 & 0.03 & 0.00 & 0.02 & 0.97 & 0.95 & 0.95 & 0.96 & 2.17 & 2.22 & 2.11 & 2.13\\

\multirow{-4}{*}{\centering -1.25} & \multirow{-4}{*}{\centering -1.87} & 2000 & 0.06 & 0.02 & 0.04 & 0.05 & 0.94 & 0.96 & 0.92 & 0.93 & 2.05 & 2.08 & 2.01 & 2.04\\
\cline{1-15}
 &  & 500 & 0.03 & 0.03 & 0.07 & 0.01 & 0.96 & 0.95 & 0.94 & 0.94 & 1.43 & 1.46 & 1.38 & 1.42\\

 &  & 1000 & 0.01 & 0.00 & 0.04 & 0.04 & 0.96 & 0.94 & 0.95 & 0.95 & 1.22 & 1.25 & 1.19 & 1.21\\

 &  & 1500 & 0.02 & 0.01 & 0.05 & 0.02 & 0.95 & 0.94 & 0.94 & 0.95 & 1.11 & 1.14 & 1.08 & 1.09\\

\multirow{-4}{*}{\centering -0.50} & \multirow{-4}{*}{\centering -1.28} & 2000 & 0.01 & 0.03 & 0.04 & 0.03 & 0.93 & 0.94 & 0.92 & 0.93 & 1.04 & 1.07 & 1.02 & 1.03\\
\cline{1-15}
 &  & 500 & 0.10 & 0.08 & 0.14 & 0.01 & 0.94 & 0.92 & 0.92 & 0.94 & 1.00 & 1.02 & 0.96 & 1.00\\

 &  & 1000 & 0.08 & 0.06 & 0.12 & 0.00 & 0.92 & 0.95 & 0.90 & 0.94 & 0.85 & 0.87 & 0.83 & 0.84\\

 &  & 1500 & 0.05 & 0.05 & 0.10 & 0.00 & 0.95 & 0.95 & 0.93 & 0.94 & 0.78 & 0.79 & 0.76 & 0.77\\

\multirow{-4}{*}{\centering 0.10} & \multirow{-4}{*}{\centering -0.95} & 2000 & 0.05 & 0.03 & 0.10 & 0.01 & 0.95 & 0.94 & 0.92 & 0.95 & 0.73 & 0.74 & 0.71 & 0.72\\
\cline{1-15}
 &  & 500 & 0.11 & 0.10 & 0.15 & 0.00 & 0.92 & 0.92 & 0.88 & 0.95 & 0.94 & 0.96 & 0.91 & 0.94\\

 &  & 1000 & 0.09 & 0.07 & 0.13 & 0.02 & 0.93 & 0.92 & 0.90 & 0.94 & 0.81 & 0.82 & 0.79 & 0.80\\

 &  & 1500 & 0.04 & 0.06 & 0.09 & 0.01 & 0.94 & 0.93 & 0.91 & 0.93 & 0.74 & 0.75 & 0.72 & 0.73\\

\multirow{-4}{*}{\centering 0.25} & \multirow{-4}{*}{\centering -0.88} & 2000 & 0.04 & 0.05 & 0.09 & 0.01 & 0.95 & 0.93 & 0.93 & 0.95 & 0.69 & 0.70 & 0.68 & 0.68\\
\cline{1-15}
 &  & 500 & 0.11 & 0.16 & 0.16 & 0.01 & 0.92 & 0.88 & 0.86 & 0.95 & 0.93 & 0.93 & 0.90 & 0.93\\

 &  & 1000 & 0.07 & 0.12 & 0.12 & 0.01 & 0.93 & 0.90 & 0.90 & 0.95 & 0.80 & 0.79 & 0.77 & 0.78\\

 &  & 1500 & 0.08 & 0.07 & 0.13 & 0.01 & 0.92 & 0.93 & 0.87 & 0.95 & 0.73 & 0.73 & 0.71 & 0.71\\

\multirow{-4}{*}{\centering 1.00} & \multirow{-4}{*}{\centering -0.61} & 2000 & 0.06 & 0.07 & 0.11 & 0.00 & 0.92 & 0.93 & 0.87 & 0.94 & 0.68 & 0.68 & 0.67 & 0.67\\
\hline
\end{tabular}}
\caption*{Table S7: Comparison of \texttt{DLL}, \texttt{DLL-S}, plug-in (\texttt{Plug}), oracle (\texttt{Orac}) estimators for the non-linear treatment model, across different sample sizes $n$ and evaluation points $\eval$. The column indexed with ``True'' represents the true value of $f'(\eval)$. The columns indexed with ``Bias'' report the absolute bias; the columns indexed with ``Coverage'' report the empirical coverage level and the columns indexed with ``Length'' report the average CI length.}
\label{table:S7}
\end{table}

\subsection{Other Bandwidth Selection Methods}
\label{sec: other bw}
We also investigate the performances of \texttt{DLL}, plug-in (\texttt{Plug}), and oracle (\texttt{Orac}) estimators using other bandwidth selection methods: the methods \texttt{regCVBwSelC()} implemented in \cite{locpol} and \texttt{npregbw()} implemented in \cite{nppackage}. We generate $X_i, D_i$ as in Setting 2, and generate the outcome model as the exactly sparse model. The results are summarised in Table \hyperref[table:S8]{S8} and Table \hyperref[table:S9]{S9}. We observe that using these two bandwidth selection methods might lead to a bad coverage for \texttt{DLL}, or a wide confidence interval. For the undercoverage settings for \texttt{DLL}, the oracle CI (the benchmark) does not attain the desired coverage level. This indicates that these bandwidth selections are not stable for our simulation studies. Hence, we select the bandwidth by the function \texttt{thumbBw()} in \texttt{locpol} as mentioned in the main paper.

\begin{table}[htb!]
\centering
\resizebox{\linewidth}{!}{
\begin{tabular}[t]{|c|c|c|ccc|ccc|ccc|ccc|ccc|}
\multicolumn{18}{c}{Setting 2, exactly sparse: $f(d) = 2\exp(-d/2)$ with \texttt{regCVBwSelC()} in \texttt{locpol}}\\
\hline
\multicolumn{1}{|c}{ } & \multicolumn{1}{c}{ } & \multicolumn{1}{c|}{ } & \multicolumn{3}{c|}{Bias} & \multicolumn{3}{c|}{RMSE} & \multicolumn{3}{c|}{SE} & \multicolumn{3}{c|}{Coverage} & \multicolumn{3}{c|}{Length} \\
\hline
$\eval$ & True & $n$ & \texttt{DLL} & \texttt{Plug} & \texttt{Orac} & \texttt{DLL} & \texttt{Plug} & \texttt{Orac} & \texttt{DLL} & \texttt{Plug} & \texttt{Orac} & \texttt{DLL} & \texttt{Plug} & \texttt{Orac} & \texttt{DLL} & \texttt{Plug} & \texttt{Orac}\\
\hline
 &  & 500 & 0.04 & 0.11 & 0.01 & 0.24 & 0.19 & 0.14 & 0.24 & 0.15 & 0.14 & 0.95 & 0.82 & 0.92 & 0.90 & 0.52 & 0.49\\

 &  & 1000 & 0.01 & 0.10 & 0.01 & 0.20 & 0.19 & 0.11 & 0.20 & 0.16 & 0.11 & 0.96 & 0.78 & 0.94 & 0.74 & 0.41 & 0.36\\

 &  & 1500 & 0.02 & 0.10 & 0.01 & 0.18 & 0.15 & 0.09 & 0.18 & 0.11 & 0.09 & 0.95 & 0.74 & 0.94 & 0.67 & 0.34 & 0.31\\

\multirow{-4}{*}{\centering -1.25} & \multirow{-4}{*}{\centering -0.41} & 2000 & 0.01 & 0.09 & 0.01 & 0.17 & 0.13 & 0.08 & 0.17 & 0.09 & 0.08 & 0.96 & 0.71 & 0.96 & 0.62 & 0.31 & 0.28\\
\cline{1-18}
 &  & 500 & 0.01 & 0.05 & 0.06 & 0.19 & 0.15 & 0.13 & 0.19 & 0.14 & 0.11 & 0.94 & 0.84 & 0.85 & 0.71 & 0.42 & 0.39\\

 &  & 1000 & 0.01 & 0.03 & 0.06 & 0.16 & 0.12 & 0.10 & 0.16 & 0.11 & 0.09 & 0.95 & 0.81 & 0.79 & 0.59 & 0.33 & 0.29\\

 &  & 1500 & 0.02 & 0.04 & 0.06 & 0.15 & 0.11 & 0.09 & 0.15 & 0.10 & 0.07 & 0.94 & 0.78 & 0.74 & 0.53 & 0.27 & 0.25\\

\multirow{-4}{*}{\centering -0.50} & \multirow{-4}{*}{\centering -0.52} & 2000 & 0.01 & 0.03 & 0.05 & 0.13 & 0.08 & 0.08 & 0.13 & 0.08 & 0.07 & 0.96 & 0.81 & 0.73 & 0.49 & 0.24 & 0.22\\
\cline{1-18}
 &  & 500 & 0.00 & 0.10 & 0.01 & 0.20 & 0.16 & 0.11 & 0.20 & 0.13 & 0.11 & 0.93 & 0.80 & 0.95 & 0.72 & 0.43 & 0.40\\

 &  & 1000 & 0.01 & 0.08 & 0.00 & 0.17 & 0.13 & 0.08 & 0.17 & 0.10 & 0.08 & 0.95 & 0.79 & 0.94 & 0.59 & 0.33 & 0.29\\

 &  & 1500 & 0.01 & 0.08 & 0.01 & 0.16 & 0.12 & 0.07 & 0.16 & 0.09 & 0.07 & 0.93 & 0.73 & 0.94 & 0.54 & 0.28 & 0.25\\

\multirow{-4}{*}{\centering 0.10} & \multirow{-4}{*}{\centering -0.40} & 2000 & 0.01 & 0.07 & 0.00 & 0.14 & 0.10 & 0.06 & 0.14 & 0.07 & 0.06 & 0.95 & 0.76 & 0.96 & 0.49 & 0.25 & 0.22\\
\cline{1-18}
 &  & 500 & 0.01 & 0.11 & 0.01 & 0.19 & 0.17 & 0.12 & 0.19 & 0.13 & 0.12 & 0.94 & 0.82 & 0.94 & 0.74 & 0.44 & 0.41\\

 &  & 1000 & 0.03 & 0.09 & 0.02 & 0.17 & 0.13 & 0.09 & 0.16 & 0.10 & 0.09 & 0.95 & 0.80 & 0.94 & 0.61 & 0.34 & 0.30\\

 &  & 1500 & 0.02 & 0.08 & 0.01 & 0.15 & 0.13 & 0.08 & 0.15 & 0.10 & 0.07 & 0.95 & 0.73 & 0.95 & 0.55 & 0.29 & 0.26\\

\multirow{-4}{*}{\centering 0.25} & \multirow{-4}{*}{\centering -0.35} & 2000 & 0.02 & 0.08 & 0.01 & 0.14 & 0.11 & 0.07 & 0.14 & 0.08 & 0.07 & 0.95 & 0.76 & 0.94 & 0.51 & 0.25 & 0.23\\
\cline{1-18}
 &  & 500 & 0.04 & 0.12 & 0.06 & 0.26 & 0.22 & 0.17 & 0.26 & 0.18 & 0.16 & 0.96 & 0.87 & 0.87 & 1.04 & 0.61 & 0.56\\

 &  & 1000 & 0.03 & 0.09 & 0.07 & 0.22 & 0.17 & 0.14 & 0.22 & 0.15 & 0.13 & 0.96 & 0.85 & 0.81 & 0.86 & 0.48 & 0.42\\

 &  & 1500 & 0.01 & 0.09 & 0.06 & 0.20 & 0.16 & 0.13 & 0.20 & 0.14 & 0.12 & 0.94 & 0.85 & 0.80 & 0.77 & 0.40 & 0.36\\

\multirow{-4}{*}{\centering 1.00} & \multirow{-4}{*}{\centering -0.15} & 2000 & 0.01 & 0.08 & 0.05 & 0.20 & 0.14 & 0.11 & 0.20 & 0.11 & 0.10 & 0.94 & 0.84 & 0.82 & 0.71 & 0.35 & 0.32\\
\hline

\multicolumn{18}{c}{Setting 2, exactly sparse: $f(d) = 1.5\sin(d)$ with \texttt{regCVBwSelC()} in \texttt{locpol}}\\
\hline
\multicolumn{1}{|c}{ } & \multicolumn{1}{c}{ } & \multicolumn{1}{c|}{ } & \multicolumn{3}{c|}{Bias} & \multicolumn{3}{c|}{RMSE} & \multicolumn{3}{c|}{SE} & \multicolumn{3}{c|}{Coverage} & \multicolumn{3}{c|}{Length} \\
\hline
$\eval$ & True & $n$ & \texttt{DLL} & \texttt{Plug} & \texttt{Orac} & \texttt{DLL} & \texttt{Plug} & \texttt{Orac} & \texttt{DLL} & \texttt{Plug} & \texttt{Orac} & \texttt{DLL} & \texttt{Plug} & \texttt{Orac} & \texttt{DLL} & \texttt{Plug} & \texttt{Orac}\\
\hline
 &  & 500 & 0.00 & 0.06 & 0.11 & 0.27 & 0.25 & 0.25 & 0.27 & 0.24 & 0.22 & 0.94 & 0.85 & 0.76 & 0.97 & 0.73 & 0.73\\

 &  & 1000 & 0.03 & 0.04 & 0.11 & 0.20 & 0.20 & 0.21 & 0.20 & 0.19 & 0.18 & 0.96 & 0.86 & 0.74 & 0.77 & 0.58 & 0.57\\

 &  & 1500 & 0.02 & 0.04 & 0.08 & 0.18 & 0.18 & 0.18 & 0.18 & 0.18 & 0.16 & 0.95 & 0.87 & 0.78 & 0.70 & 0.54 & 0.53\\

\multirow{-4}{*}{\centering -1.25} & \multirow{-4}{*}{\centering -1.00} & 2000 & 0.03 & 0.02 & 0.09 & 0.17 & 0.16 & 0.18 & 0.16 & 0.16 & 0.16 & 0.94 & 0.83 & 0.72 & 0.63 & 0.48 & 0.47\\
\cline{1-18}
 &  & 500 & 0.10 & 0.00 & 0.20 & 0.23 & 0.19 & 0.28 & 0.20 & 0.19 & 0.19 & 0.89 & 0.86 & 0.55 & 0.77 & 0.58 & 0.58\\

 &  & 1000 & 0.10 & 0.02 & 0.18 & 0.19 & 0.17 & 0.24 & 0.17 & 0.16 & 0.15 & 0.87 & 0.80 & 0.51 & 0.61 & 0.46 & 0.45\\

 &  & 1500 & 0.09 & 0.02 & 0.16 & 0.19 & 0.16 & 0.22 & 0.16 & 0.16 & 0.15 & 0.87 & 0.77 & 0.53 & 0.55 & 0.43 & 0.42\\

\multirow{-4}{*}{\centering -0.50} & \multirow{-4}{*}{\centering -0.56} & 2000 & 0.09 & 0.03 & 0.16 & 0.16 & 0.14 & 0.21 & 0.13 & 0.14 & 0.14 & 0.85 & 0.74 & 0.54 & 0.50 & 0.38 & 0.37\\
\cline{1-18}
 &  & 500 & 0.27 & 0.47 & 0.26 & 0.34 & 0.52 & 0.33 & 0.21 & 0.22 & 0.21 & 0.63 & 0.23 & 0.49 & 0.78 & 0.59 & 0.59\\

 &  & 1000 & 0.20 & 0.41 & 0.23 & 0.26 & 0.44 & 0.29 & 0.17 & 0.17 & 0.17 & 0.68 & 0.21 & 0.44 & 0.62 & 0.47 & 0.45\\

 &  & 1500 & 0.16 & 0.35 & 0.20 & 0.23 & 0.39 & 0.27 & 0.17 & 0.18 & 0.18 & 0.72 & 0.24 & 0.46 & 0.56 & 0.43 & 0.43\\

\multirow{-4}{*}{\centering 0.10} & \multirow{-4}{*}{\centering 0.74} & 2000 & 0.16 & 0.33 & 0.20 & 0.21 & 0.37 & 0.26 & 0.14 & 0.16 & 0.16 & 0.68 & 0.25 & 0.47 & 0.50 & 0.38 & 0.38\\
\cline{1-18}
 &  & 500 & 0.30 & 0.51 & 0.29 & 0.38 & 0.56 & 0.37 & 0.23 & 0.23 & 0.23 & 0.59 & 0.24 & 0.44 & 0.80 & 0.61 & 0.61\\

 &  & 1000 & 0.22 & 0.44 & 0.27 & 0.28 & 0.48 & 0.33 & 0.18 & 0.19 & 0.20 & 0.65 & 0.20 & 0.42 & 0.64 & 0.48 & 0.47\\

 &  & 1500 & 0.18 & 0.37 & 0.22 & 0.26 & 0.42 & 0.29 & 0.19 & 0.20 & 0.19 & 0.64 & 0.24 & 0.44 & 0.58 & 0.45 & 0.44\\

\multirow{-4}{*}{\centering 0.25} & \multirow{-4}{*}{\centering 0.94} & 2000 & 0.16 & 0.36 & 0.23 & 0.22 & 0.40 & 0.29 & 0.14 & 0.18 & 0.18 & 0.67 & 0.27 & 0.46 & 0.52 & 0.40 & 0.39\\
\cline{1-18}
 &  & 500 & 0.06 & 0.19 & 0.01 & 0.31 & 0.33 & 0.27 & 0.31 & 0.27 & 0.27 & 0.94 & 0.80 & 0.92 & 1.11 & 0.85 & 0.84\\

 &  & 1000 & 0.00 & 0.12 & 0.03 & 0.25 & 0.24 & 0.19 & 0.25 & 0.20 & 0.19 & 0.96 & 0.78 & 0.93 & 0.89 & 0.68 & 0.66\\

 &  & 1500 & 0.01 & 0.11 & 0.03 & 0.21 & 0.22 & 0.18 & 0.21 & 0.19 & 0.18 & 0.96 & 0.82 & 0.94 & 0.81 & 0.63 & 0.62\\

\multirow{-4}{*}{\centering 1.00} & \multirow{-4}{*}{\centering 0.81} & 2000 & 0.01 & 0.12 & 0.00 & 0.19 & 0.20 & 0.15 & 0.19 & 0.16 & 0.15 & 0.95 & 0.80 & 0.93 & 0.72 & 0.55 & 0.54\\
\hline
\end{tabular}}
\caption*{Table S8: Comparison of \texttt{DLL}, plug-in (\texttt{Plug}), and oracle (\texttt{Orac}) estimators using \texttt{regCVBwSelC()} for bandwidth selection, across different sample sizes $n$ and evaluation points $\eval$. The column indexed with ``True'' represents the true value of $f'(\eval)$. The columns indexed with ``Bias'', ``RMSE'' and ``SE" report the absolute bias, the root mean square error, and the standard error computed by 500 estimates, respectively; the columns indexed with ``Coverage'' report the empirical coverage level and the columns indexed with ``Length'' report the average CI length.}
\label{table:S8}
\end{table}

\begin{table}[htb!]
\centering
\resizebox{\linewidth}{!}{
\begin{tabular}[t]{|c|c|c|ccc|ccc|ccc|ccc|ccc|}
\multicolumn{18}{c}{Setting 2, exactly sparse: $f(d) = 2\exp(-d/2)$ with \texttt{npregbw()} in \texttt{np}}\\
\hline
\multicolumn{1}{|c}{ } & \multicolumn{1}{c}{ } & \multicolumn{1}{c|}{ } & \multicolumn{3}{c|}{Bias} & \multicolumn{3}{c|}{RMSE} & \multicolumn{3}{c|}{SE} & \multicolumn{3}{c|}{Coverage} & \multicolumn{3}{c|}{Length} \\
\hline
$\eval$ & True & $n$ & \texttt{DLL} & \texttt{Plug} & \texttt{Orac} & \texttt{DLL} & \texttt{Plug} & \texttt{Orac} & \texttt{DLL} & \texttt{Plug} & \texttt{Orac} & \texttt{DLL} & \texttt{Plug} & \texttt{Orac} & \texttt{DLL} & \texttt{Plug} & \texttt{Orac}\\
\hline
 &  & 500 & 0.02 & 0.12 & 0.01 & 0.88 & 0.90 & 0.71 & 0.88 & 0.89 & 0.71 & 0.94 & 0.85 & 0.90 & 1.89 & 1.76 & 1.43\\

 &  & 1000 & 0.01 & 0.08 & 0.02 & 0.56 & 0.56 & 0.41 & 0.56 & 0.56 & 0.41 & 0.94 & 0.87 & 0.92 & 1.44 & 1.36 & 1.12\\

 &  & 1500 & 0.01 & 0.08 & 0.02 & 0.46 & 0.47 & 0.43 & 0.46 & 0.46 & 0.43 & 0.96 & 0.84 & 0.94 & 1.32 & 1.27 & 1.00\\

\multirow{-4}{*}{\centering -1.25} & \multirow{-4}{*}{\centering -0.41} & 2000 & 0.02 & 0.12 & 0.00 & 0.47 & 0.49 & 0.45 & 0.47 & 0.47 & 0.45 & 0.95 & 0.88 & 0.93 & 1.10 & 1.06 & 0.91\\
\cline{1-18}
 &  & 500 & 0.02 & 0.02 & 0.02 & 0.63 & 0.62 & 0.62 & 0.63 & 0.62 & 0.62 & 0.88 & 0.83 & 0.74 & 1.48 & 1.41 & 1.16\\

 &  & 1000 & 0.03 & 0.02 & 0.03 & 0.49 & 0.50 & 0.41 & 0.49 & 0.50 & 0.41 & 0.91 & 0.86 & 0.81 & 1.14 & 1.08 & 0.89\\

 &  & 1500 & 0.01 & 0.02 & 0.03 & 0.49 & 0.49 & 0.32 & 0.49 & 0.49 & 0.32 & 0.88 & 0.86 & 0.76 & 1.05 & 1.02 & 0.78\\

\multirow{-4}{*}{\centering -0.50} & \multirow{-4}{*}{\centering -0.52} & 2000 & 0.01 & 0.00 & 0.04 & 0.36 & 0.36 & 0.33 & 0.36 & 0.36 & 0.33 & 0.88 & 0.86 & 0.82 & 0.86 & 0.84 & 0.71\\
\cline{1-18}
 &  & 500 & 0.02 & 0.11 & 0.00 & 0.66 & 0.67 & 0.49 & 0.66 & 0.66 & 0.50 & 0.92 & 0.84 & 0.91 & 1.51 & 1.45 & 1.15\\

 &  & 1000 & 0.00 & 0.10 & 0.00 & 0.46 & 0.46 & 0.39 & 0.46 & 0.45 & 0.39 & 0.94 & 0.83 & 0.93 & 1.16 & 1.10 & 0.90\\

 &  & 1500 & 0.01 & 0.11 & 0.00 & 0.59 & 0.64 & 0.34 & 0.59 & 0.63 & 0.34 & 0.95 & 0.84 & 0.94 & 1.07 & 1.03 & 0.80\\

\multirow{-4}{*}{\centering 0.10} & \multirow{-4}{*}{\centering -0.40} & 2000 & 0.03 & 0.06 & 0.00 & 0.36 & 0.36 & 0.30 & 0.36 & 0.36 & 0.30 & 0.95 & 0.87 & 0.95 & 0.87 & 0.84 & 0.71\\
\cline{1-18}
 &  & 500 & 0.05 & 0.09 & 0.02 & 0.59 & 0.61 & 0.52 & 0.59 & 0.60 & 0.52 & 0.92 & 0.85 & 0.91 & 1.55 & 1.47 & 1.20\\

 &  & 1000 & 0.05 & 0.18 & 0.03 & 0.63 & 0.65 & 0.50 & 0.63 & 0.62 & 0.50 & 0.94 & 0.87 & 0.94 & 1.19 & 1.14 & 0.94\\

 &  & 1500 & 0.05 & 0.17 & 0.00 & 0.47 & 0.49 & 0.32 & 0.47 & 0.46 & 0.32 & 0.94 & 0.84 & 0.95 & 1.10 & 1.06 & 0.81\\

\multirow{-4}{*}{\centering 0.25} & \multirow{-4}{*}{\centering -0.35} & 2000 & 0.01 & 0.12 & 0.01 & 0.33 & 0.35 & 0.34 & 0.33 & 0.33 & 0.34 & 0.93 & 0.86 & 0.95 & 0.90 & 0.88 & 0.73\\
\cline{1-18}
 &  & 500 & 0.05 & 0.11 & 0.08 & 1.13 & 1.02 & 0.89 & 1.13 & 1.02 & 0.89 & 0.84 & 0.84 & 0.71 & 2.21 & 2.07 & 1.68\\

 &  & 1000 & 0.02 & 0.14 & 0.03 & 0.81 & 0.80 & 0.54 & 0.81 & 0.79 & 0.54 & 0.85 & 0.86 & 0.77 & 1.67 & 1.58 & 1.28\\

 &  & 1500 & 0.02 & 0.16 & 0.03 & 0.73 & 0.70 & 0.54 & 0.73 & 0.68 & 0.54 & 0.92 & 0.86 & 0.81 & 1.54 & 1.49 & 1.14\\

\multirow{-4}{*}{\centering 1.00} & \multirow{-4}{*}{\centering -0.15} & 2000 & 0.03 & 0.11 & 0.02 & 0.50 & 0.49 & 0.44 & 0.50 & 0.48 & 0.44 & 0.90 & 0.88 & 0.83 & 1.26 & 1.23 & 1.05\\
\hline

\multicolumn{18}{c}{Setting 2, exactly sparse: $f(d) = 1.5\sin(d)$ with \texttt{npregbw()} in \texttt{np}}\\
\hline
\multicolumn{1}{|c}{ } & \multicolumn{1}{c}{ } & \multicolumn{1}{c|}{ } & \multicolumn{3}{c|}{Bias} & \multicolumn{3}{c|}{RMSE} & \multicolumn{3}{c|}{SE} & \multicolumn{3}{c|}{Coverage} & \multicolumn{3}{c|}{Length} \\
\hline
$\eval$ & True & $n$ & \texttt{DLL} & \texttt{Plug} & \texttt{Orac} & \texttt{DLL} & \texttt{Plug} & \texttt{Orac} & \texttt{DLL} & \texttt{Plug} & \texttt{Orac} & \texttt{DLL} & \texttt{Plug} & \texttt{Orac} & \texttt{DLL} & \texttt{Plug} & \texttt{Orac}\\
\hline
 &  & 500 & 0.05 & 0.18 & 0.03 & 0.79 & 0.76 & 0.76 & 0.79 & 0.74 & 0.76 & 0.93 & 0.90 & 0.93 & 2.12 & 1.93 & 2.15\\

 &  & 1000 & 0.02 & 0.14 & 0.02 & 0.57 & 0.58 & 0.65 & 0.57 & 0.56 & 0.65 & 0.95 & 0.93 & 0.95 & 1.67 & 1.58 & 1.76\\

 &  & 1500 & 0.03 & 0.08 & 0.07 & 0.73 & 0.73 & 0.62 & 0.73 & 0.73 & 0.61 & 0.94 & 0.94 & 0.95 & 1.58 & 1.50 & 1.62\\

\multirow{-4}{*}{\centering -1.25} & \multirow{-4}{*}{\centering -1.00} & 2000 & 0.00 & 0.10 & 0.01 & 0.51 & 0.51 & 0.52 & 0.51 & 0.50 & 0.52 & 0.95 & 0.92 & 0.94 & 1.36 & 1.28 & 1.41\\
\cline{1-18}
 &  & 500 & 0.05 & 0.09 & 0.11 & 0.57 & 0.57 & 0.63 & 0.57 & 0.56 & 0.62 & 0.93 & 0.93 & 0.90 & 1.65 & 1.56 & 1.75\\

 &  & 1000 & 0.07 & 0.07 & 0.07 & 0.51 & 0.51 & 0.50 & 0.50 & 0.50 & 0.50 & 0.93 & 0.92 & 0.91 & 1.33 & 1.26 & 1.42\\

 &  & 1500 & 0.07 & 0.05 & 0.07 & 0.50 & 0.47 & 0.47 & 0.49 & 0.47 & 0.47 & 0.94 & 0.94 & 0.92 & 1.24 & 1.19 & 1.28\\

\multirow{-4}{*}{\centering -0.50} & \multirow{-4}{*}{\centering -0.56} & 2000 & 0.05 & 0.06 & 0.07 & 0.44 & 0.43 & 0.42 & 0.43 & 0.43 & 0.41 & 0.93 & 0.90 & 0.91 & 1.07 & 1.02 & 1.12\\
\cline{1-18}
 &  & 500 & 0.22 & 0.36 & 0.11 & 0.63 & 0.70 & 0.67 & 0.60 & 0.60 & 0.66 & 0.81 & 0.66 & 0.86 & 1.66 & 1.56 & 1.77\\

 &  & 1000 & 0.13 & 0.26 & 0.08 & 0.56 & 0.59 & 0.68 & 0.54 & 0.53 & 0.67 & 0.87 & 0.70 & 0.90 & 1.34 & 1.28 & 1.43\\

 &  & 1500 & 0.13 & 0.25 & 0.10 & 0.47 & 0.54 & 0.44 & 0.45 & 0.48 & 0.43 & 0.89 & 0.72 & 0.93 & 1.26 & 1.22 & 1.31\\

\multirow{-4}{*}{\centering 0.10} & \multirow{-4}{*}{\centering 0.74} & 2000 & 0.08 & 0.18 & 0.04 & 0.33 & 0.36 & 0.38 & 0.32 & 0.31 & 0.38 & 0.92 & 0.78 & 0.93 & 1.09 & 1.04 & 1.14\\
\cline{1-18}
 &  & 500 & 0.20 & 0.35 & 0.16 & 0.69 & 0.74 & 0.70 & 0.66 & 0.65 & 0.68 & 0.81 & 0.66 & 0.84 & 1.74 & 1.64 & 1.83\\

 &  & 1000 & 0.14 & 0.27 & 0.09 & 0.41 & 0.46 & 0.52 & 0.39 & 0.38 & 0.51 & 0.86 & 0.71 & 0.89 & 1.40 & 1.33 & 1.48\\

 &  & 1500 & 0.14 & 0.25 & 0.08 & 0.43 & 0.48 & 0.51 & 0.41 & 0.40 & 0.51 & 0.88 & 0.72 & 0.89 & 1.30 & 1.25 & 1.34\\

\multirow{-4}{*}{\centering 0.25} & \multirow{-4}{*}{\centering 0.94} & 2000 & 0.07 & 0.18 & 0.05 & 0.35 & 0.38 & 0.38 & 0.34 & 0.33 & 0.38 & 0.91 & 0.77 & 0.88 & 1.13 & 1.09 & 1.18\\
\cline{1-18}
 &  & 500 & 0.02 & 0.15 & 0.07 & 0.84 & 0.76 & 0.96 & 0.84 & 0.74 & 0.96 & 0.96 & 0.91 & 0.94 & 2.44 & 2.22 & 2.48\\

 &  & 1000 & 0.07 & 0.19 & 0.02 & 0.74 & 0.74 & 0.87 & 0.74 & 0.72 & 0.87 & 0.96 & 0.93 & 0.94 & 1.98 & 1.86 & 2.07\\

 &  & 1500 & 0.06 & 0.06 & 0.05 & 0.94 & 0.87 & 0.70 & 0.94 & 0.87 & 0.69 & 0.92 & 0.90 & 0.92 & 1.82 & 1.71 & 1.85\\

\multirow{-4}{*}{\centering 1.00} & \multirow{-4}{*}{\centering 0.81} & 2000 & 0.01 & 0.11 & 0.01 & 0.51 & 0.53 & 0.52 & 0.51 & 0.52 & 0.52 & 0.96 & 0.93 & 0.95 & 1.57 & 1.49 & 1.63\\
\hline
\end{tabular}}
\caption*{Table S9: Comparison of \texttt{DLL}, plug-in (\texttt{Plug}), and oracle (\texttt{Orac}) estimators using \texttt{npregbw()} for bandwidth selection, across different sample sizes $n$ and evaluation points $\eval$. The column indexed with ``True'' represents the true value of $f'(\eval)$. The columns indexed with ``Bias'', ``RMSE'' and ``SE" report the absolute bias, the root mean square error, and the standard error computed by 500 estimates, respectively; the columns indexed with ``Coverage'' report the empirical coverage level and the columns indexed with ``Length'' report the average CI length.}
\label{table:S9}
\end{table}

\subsection{Data Swap and Quantile Transformation}
\label{sec: swap and trans}
In Table \hyperref[table:S10]{S10} and Table \hyperref[table:S11]{S11}, we compare the \texttt{DLL} estimator without data swapping and the \texttt{DLL} with data swapping (\texttt{Swap}). The data is generated as in Setting 2, with the sparse additive model being exactly sparse.
The CIs with and without data swapping attain the desired coverage level. When the sample size is relatively large, they have similar performance; for relatively small sample size, the confidence interval without data swapping can be shorter than that with data swapping. 
This happens because no data swapping uses the entire data to construct initial estimators of $g$ and $\gamma$. When the sample size is relatively small (e.g., $n=500$ and $p=1500$), the \texttt{DLL} with data swapping might be slightly noisier than the one without data swapping.

We now
investigate the performance of our proposed method with quantile transformation (\texttt{Trans}), which is detailed in Section \ref{sec: quantile}.
We report the comparison with the \texttt{Trans} estimator in Tables \hyperref[table:S10]{S10} and \hyperref[table:S11]{S11}. The data is generated as in Setting 2, with the sparse additive model being exactly sparse. As reported in Table \hyperref[table:S11]{S11}, the method using quantile transformation leads to slightly better performance for $f(d) = 2\exp(-d/2)$: the bias is slightly smaller, and the CI is shorter. Nevertheless, the regular \texttt{DLL} method still attains the desired coverage.

\begin{table}[htb!]
\centering
\resizebox{\linewidth}{!}{
\begin{tabular}[t]{|c|c|c|ccc|ccc|ccc|ccc|ccc|}
\multicolumn{18}{c}{Setting 2, exactly sparse: $f(d) = 1.5\sin(d)$}\\
\hline
\multicolumn{1}{|c}{ } & \multicolumn{1}{c}{ } & \multicolumn{1}{c|}{ } & \multicolumn{3}{c|}{Bias} & \multicolumn{3}{c|}{RMSE} & \multicolumn{3}{c|}{SE} & \multicolumn{3}{c|}{Coverage} & \multicolumn{3}{c|}{Length} \\
\hline
$\eval$ & True & $n$ & \texttt{DLL} & \texttt{Swap} & \texttt{Trans} & \texttt{DLL} & \texttt{Swap} & \texttt{Trans} & \texttt{DLL} & \texttt{Swap} & \texttt{Trans} & \texttt{DLL} & \texttt{Swap} & \texttt{Trans} & \texttt{DLL} & \texttt{Swap} & \texttt{Trans}\\
\hline
 &  & 500 & 0.00 & 0.09 & 0.04 & 0.49 & 0.55 & 0.55 & 0.49 & 0.55 & 0.55 & 0.95 & 0.94 & 0.93 & 2.01 & 2.13 & 2.02\\

 &  & 1000 & 0.02 & 0.02 & 0.04 & 0.41 & 0.43 & 0.42 & 0.41 & 0.43 & 0.42 & 0.96 & 0.94 & 0.96 & 1.68 & 1.70 & 1.63\\

 &  & 1500 & 0.02 & 0.06 & 0.02 & 0.38 & 0.40 & 0.37 & 0.38 & 0.40 & 0.37 & 0.94 & 0.94 & 0.95 & 1.51 & 1.52 & 1.47\\

\multirow{-4}{*}{\centering -1.25} & \multirow{-4}{*}{\centering -0.41} & 2000 & 0.01 & 0.03 & 0.00 & 0.37 & 0.36 & 0.33 & 0.37 & 0.36 & 0.33 & 0.96 & 0.94 & 0.96 & 1.42 & 1.42 & 1.37\\
\cline{1-18}
 &  & 500 & 0.00 & 0.07 & 0.00 & 0.40 & 0.45 & 0.40 & 0.40 & 0.45 & 0.40 & 0.95 & 0.93 & 0.94 & 1.59 & 1.66 & 1.58\\

 &  & 1000 & 0.02 & 0.05 & 0.00 & 0.32 & 0.34 & 0.33 & 0.32 & 0.33 & 0.33 & 0.96 & 0.96 & 0.95 & 1.32 & 1.34 & 1.29\\

 &  & 1500 & 0.01 & 0.03 & 0.02 & 0.29 & 0.31 & 0.30 & 0.29 & 0.31 & 0.30 & 0.97 & 0.95 & 0.94 & 1.19 & 1.20 & 1.16\\

\multirow{-4}{*}{\centering -0.50} & \multirow{-4}{*}{\centering -0.52} & 2000 & 0.03 & 0.01 & 0.03 & 0.28 & 0.29 & 0.27 & 0.28 & 0.30 & 0.26 & 0.95 & 0.94 & 0.96 & 1.12 & 1.12 & 1.08\\
\cline{1-18}
 &  & 500 & 0.17 & 0.19 & 0.11 & 0.45 & 0.49 & 0.43 & 0.42 & 0.45 & 0.41 & 0.92 & 0.90 & 0.93 & 1.61 & 1.69 & 1.60\\

 &  & 1000 & 0.08 & 0.10 & 0.08 & 0.33 & 0.36 & 0.33 & 0.32 & 0.35 & 0.32 & 0.95 & 0.94 & 0.94 & 1.34 & 1.36 & 1.31\\

 &  & 1500 & 0.07 & 0.10 & 0.05 & 0.31 & 0.34 & 0.28 & 0.30 & 0.32 & 0.28 & 0.95 & 0.93 & 0.96 & 1.21 & 1.22 & 1.18\\

\multirow{-4}{*}{\centering 0.10} & \multirow{-4}{*}{\centering -0.40} & 2000 & 0.07 & 0.07 & 0.06 & 0.29 & 0.30 & 0.27 & 0.29 & 0.29 & 0.27 & 0.96 & 0.93 & 0.96 & 1.13 & 1.13 & 1.09\\
\cline{1-18}
 &  & 500 & 0.17 & 0.18 & 0.14 & 0.46 & 0.50 & 0.45 & 0.43 & 0.47 & 0.42 & 0.92 & 0.92 & 0.95 & 1.67 & 1.75 & 1.66\\

 &  & 1000 & 0.08 & 0.12 & 0.09 & 0.36 & 0.38 & 0.35 & 0.35 & 0.37 & 0.34 & 0.93 & 0.93 & 0.94 & 1.38 & 1.40 & 1.35\\

 &  & 1500 & 0.08 & 0.11 & 0.07 & 0.31 & 0.35 & 0.31 & 0.30 & 0.34 & 0.30 & 0.96 & 0.93 & 0.95 & 1.25 & 1.26 & 1.22\\

\multirow{-4}{*}{\centering 0.25} & \multirow{-4}{*}{\centering -0.35} & 2000 & 0.07 & 0.06 & 0.08 & 0.30 & 0.32 & 0.30 & 0.29 & 0.32 & 0.29 & 0.95 & 0.91 & 0.95 & 1.17 & 1.17 & 1.13\\
\cline{1-18}
 &  & 500 & 0.08 & 0.10 & 0.03 & 0.61 & 0.64 & 0.59 & 0.61 & 0.64 & 0.59 & 0.93 & 0.93 & 0.95 & 2.33 & 2.45 & 2.33\\

 &  & 1000 & 0.00 & 0.09 & 0.02 & 0.48 & 0.51 & 0.48 & 0.48 & 0.50 & 0.48 & 0.96 & 0.94 & 0.94 & 1.94 & 1.97 & 1.89\\

 &  & 1500 & 0.03 & 0.03 & 0.02 & 0.44 & 0.45 & 0.45 & 0.44 & 0.45 & 0.45 & 0.95 & 0.95 & 0.94 & 1.75 & 1.76 & 1.70\\
 
\multirow{-4}{*}{\centering 1.00} & \multirow{-4}{*}{\centering -0.15} & 2000 & 0.01 & 0.00 & 0.00 & 0.40 & 0.44 & 0.40 & 0.40 & 0.44 & 0.40 & 0.95 & 0.94 & 0.97 & 1.63 & 1.64 & 1.58\\
\hline

\multicolumn{18}{c}{Setting 2, approximately sparse: $f(d) = 1.5\sin(d)$}\\
\hline
\multicolumn{1}{|c}{ } & \multicolumn{1}{c}{ } & \multicolumn{1}{c|}{ } & \multicolumn{3}{c|}{Bias} & \multicolumn{3}{c|}{RMSE} & \multicolumn{3}{c|}{SE} & \multicolumn{3}{c|}{Coverage} & \multicolumn{3}{c|}{Length} \\
\hline
$\eval$ & True & $n$ & \texttt{DLL} & \texttt{Swap} & \texttt{Trans} & \texttt{DLL} & \texttt{Swap} & \texttt{Trans} & \texttt{DLL} & \texttt{Swap} & \texttt{Trans} & \texttt{DLL} & \texttt{Swap} & \texttt{Trans} & \texttt{DLL} & \texttt{Swap} & \texttt{Trans}\\
\hline
 &  & 500 & 0.02 & 0.07 & 0.00 & 0.47 & 0.61 & 0.48 & 0.47 & 0.61 & 0.48 & 0.96 & 0.94 & 0.96 & 1.89 & 2.23 & 1.88\\

 &  & 1000 & 0.04 & 0.02 & 0.02 & 0.39 & 0.44 & 0.44 & 0.39 & 0.44 & 0.44 & 0.96 & 0.96 & 0.92 & 1.62 & 1.77 & 1.61\\

 &  & 1500 & 0.02 & 0.02 & 0.03 & 0.35 & 0.41 & 0.39 & 0.35 & 0.41 & 0.39 & 0.98 & 0.96 & 0.94 & 1.47 & 1.57 & 1.44\\

\multirow{-4}{*}{\centering -1.25} & \multirow{-4}{*}{\centering -0.41} & 2000 & 0.03 & 0.03 & 0.03 & 0.34 & 0.37 & 0.35 & 0.34 & 0.37 & 0.35 & 0.95 & 0.95 & 0.95 & 1.37 & 1.45 & 1.35\\
\cline{1-18}
 &  & 500 & 0.03 & 0.01 & 0.02 & 0.38 & 0.44 & 0.38 & 0.38 & 0.44 & 0.38 & 0.95 & 0.97 & 0.94 & 1.48 & 1.75 & 1.47\\

 &  & 1000 & 0.03 & 0.00 & 0.02 & 0.33 & 0.38 & 0.35 & 0.33 & 0.38 & 0.35 & 0.94 & 0.95 & 0.94 & 1.28 & 1.41 & 1.26\\

 &  & 1500 & 0.03 & 0.02 & 0.02 & 0.28 & 0.32 & 0.29 & 0.28 & 0.32 & 0.29 & 0.96 & 0.95 & 0.95 & 1.16 & 1.24 & 1.13\\

\multirow{-4}{*}{\centering -0.50} & \multirow{-4}{*}{\centering -0.52} & 2000 & 0.02 & 0.01 & 0.01 & 0.28 & 0.29 & 0.26 & 0.28 & 0.29 & 0.26 & 0.96 & 0.95 & 0.96 & 1.08 & 1.15 & 1.07\\
\cline{1-18}
 &  & 500 & 0.19 & 0.18 & 0.16 & 0.42 & 0.51 & 0.44 & 0.38 & 0.47 & 0.41 & 0.92 & 0.93 & 0.91 & 1.51 & 1.78 & 1.49\\

 &  & 1000 & 0.15 & 0.09 & 0.11 & 0.36 & 0.37 & 0.34 & 0.32 & 0.36 & 0.32 & 0.93 & 0.94 & 0.93 & 1.30 & 1.42 & 1.28\\

 &  & 1500 & 0.11 & 0.08 & 0.09 & 0.33 & 0.34 & 0.32 & 0.31 & 0.34 & 0.31 & 0.94 & 0.92 & 0.93 & 1.17 & 1.25 & 1.15\\

\multirow{-4}{*}{\centering 0.10} & \multirow{-4}{*}{\centering -0.40} & 2000 & 0.10 & 0.09 & 0.06 & 0.29 & 0.30 & 0.29 & 0.27 & 0.29 & 0.28 & 0.95 & 0.95 & 0.94 & 1.10 & 1.17 & 1.08\\
\cline{1-18}
 &  & 500 & 0.21 & 0.18 & 0.22 & 0.45 & 0.47 & 0.45 & 0.40 & 0.44 & 0.39 & 0.92 & 0.96 & 0.92 & 1.55 & 1.84 & 1.53\\

 &  & 1000 & 0.18 & 0.11 & 0.14 & 0.37 & 0.41 & 0.36 & 0.32 & 0.39 & 0.33 & 0.92 & 0.94 & 0.93 & 1.34 & 1.47 & 1.32\\

 &  & 1500 & 0.12 & 0.08 & 0.09 & 0.34 & 0.34 & 0.32 & 0.32 & 0.33 & 0.30 & 0.92 & 0.94 & 0.94 & 1.22 & 1.30 & 1.18\\

\multirow{-4}{*}{\centering 0.25} & \multirow{-4}{*}{\centering -0.35} & 2000 & 0.12 & 0.08 & 0.08 & 0.30 & 0.32 & 0.30 & 0.28 & 0.31 & 0.30 & 0.95 & 0.96 & 0.93 & 1.13 & 1.21 & 1.12\\
\cline{1-18}
 &  & 500 & 0.12 & 0.01 & 0.07 & 0.59 & 0.64 & 0.57 & 0.58 & 0.65 & 0.57 & 0.92 & 0.95 & 0.94 & 2.19 & 2.58 & 2.17\\

 &  & 1000 & 0.04 & 0.08 & 0.01 & 0.46 & 0.54 & 0.48 & 0.46 & 0.54 & 0.48 & 0.96 & 0.93 & 0.94 & 1.88 & 2.06 & 1.85\\

 &  & 1500 & 0.01 & 0.01 & 0.02 & 0.42 & 0.48 & 0.43 & 0.42 & 0.48 & 0.43 & 0.95 & 0.95 & 0.95 & 1.70 & 1.82 & 1.66\\

\multirow{-4}{*}{\centering 1.00} & \multirow{-4}{*}{\centering -0.15} & 2000 & 0.05 & 0.02 & 0.03 & 0.40 & 0.46 & 0.41 & 0.40 & 0.46 & 0.41 & 0.94 & 0.93 & 0.94 & 1.58 & 1.68 & 1.57\\
\hline
\end{tabular}}
\caption*{Table S10: Comparison of \texttt{DLL}, \texttt{DLL} with data swapping (\texttt{Swap}), \texttt{DLL} with quantile transformation (\texttt{Trans})) in Setting 2 when $f(d) = 1.5\sin(d)$, across different sample sizes $n$ and evaluation points $\eval$. The column indexed with ``True'' represents the true value of $f'(\eval)$. The columns indexed with ``Bias'', ``RMSE'' and ``SE" report the absolute bias, the root mean square error, and the standard error computed by 500 estimates, respectively; the columns indexed with ``Coverage'' report the empirical coverage level and the columns indexed with ``Length'' report the average CI length.}
\label{table:S10}
\end{table}

\begin{table}[htb!]
\centering
\resizebox{\linewidth}{!}{
\begin{tabular}[t]{|c|c|c|ccc|ccc|ccc|ccc|ccc|}
\multicolumn{18}{c}{Setting 2, exactly sparse: $f(d) = 2\exp(-d/2)$}\\
\hline
\multicolumn{1}{|c}{ } & \multicolumn{1}{c}{ } & \multicolumn{1}{c|}{ } & \multicolumn{3}{c|}{Bias} & \multicolumn{3}{c|}{RMSE} & \multicolumn{3}{c|}{SE} & \multicolumn{3}{c|}{Coverage} & \multicolumn{3}{c|}{Length} \\
\hline
$\eval$ & True & $n$ & \texttt{DLL} & \texttt{Swap} & \texttt{Trans} & \texttt{DLL} & \texttt{Swap} & \texttt{Trans} & \texttt{DLL} & \texttt{Swap} & \texttt{Trans} & \texttt{DLL} & \texttt{Swap} & \texttt{Trans} & \texttt{DLL} & \texttt{Swap} & \texttt{Trans}\\
\hline
 &  & 500 & 0.03 & 0.12 & 0.04 & 0.30 & 0.37 & 0.29 & 0.30 & 0.35 & 0.29 & 0.96 & 0.94 & 0.94 & 1.19 & 1.33 & 1.07\\

 &  & 1000 & 0.05 & 0.06 & 0.02 & 0.25 & 0.27 & 0.20 & 0.24 & 0.26 & 0.20 & 0.94 & 0.95 & 0.97 & 0.92 & 0.96 & 0.82\\

 &  & 1500 & 0.02 & 0.04 & 0.01 & 0.21 & 0.21 & 0.19 & 0.21 & 0.21 & 0.19 & 0.94 & 0.94 & 0.95 & 0.82 & 0.83 & 0.71\\

\multirow{-4}{*}{\centering -1.25} & \multirow{-4}{*}{\centering -0.41} & 2000 & 0.02 & 0.02 & 0.00 & 0.20 & 0.19 & 0.16 & 0.20 & 0.19 & 0.16 & 0.95 & 0.97 & 0.96 & 0.74 & 0.76 & 0.64\\
\cline{1-18}
 &  & 500 & 0.00 & 0.05 & 0.01 & 0.24 & 0.27 & 0.23 & 0.24 & 0.27 & 0.23 & 0.95 & 0.96 & 0.94 & 0.93 & 1.04 & 0.85\\

 &  & 1000 & 0.01 & 0.00 & 0.00 & 0.19 & 0.19 & 0.17 & 0.19 & 0.19 & 0.17 & 0.93 & 0.96 & 0.94 & 0.73 & 0.75 & 0.64\\

 &  & 1500 & 0.02 & 0.00 & 0.00 & 0.17 & 0.16 & 0.15 & 0.16 & 0.16 & 0.15 & 0.95 & 0.96 & 0.93 & 0.65 & 0.66 & 0.56\\

\multirow{-4}{*}{\centering -0.50} & \multirow{-4}{*}{\centering -0.52} & 2000 & 0.02 & 0.01 & 0.01 & 0.15 & 0.16 & 0.13 & 0.15 & 0.16 & 0.13 & 0.95 & 0.95 & 0.94 & 0.59 & 0.60 & 0.51\\
\cline{1-18}
 &  & 500 & 0.01 & 0.03 & 0.02 & 0.25 & 0.28 & 0.22 & 0.25 & 0.27 & 0.22 & 0.94 & 0.95 & 0.96 & 0.94 & 1.06 & 0.86\\

 &  & 1000 & 0.01 & 0.01 & 0.02 & 0.19 & 0.20 & 0.18 & 0.19 & 0.20 & 0.18 & 0.95 & 0.95 & 0.96 & 0.74 & 0.77 & 0.65\\

 &  & 1500 & 0.01 & 0.01 & 0.00 & 0.16 & 0.17 & 0.15 & 0.16 & 0.17 & 0.15 & 0.96 & 0.97 & 0.95 & 0.66 & 0.67 & 0.57\\

\multirow{-4}{*}{\centering 0.10} & \multirow{-4}{*}{\centering -0.40} & 2000 & 0.02 & 0.01 & 0.01 & 0.15 & 0.16 & 0.14 & 0.15 & 0.16 & 0.14 & 0.96 & 0.96 & 0.94 & 0.60 & 0.61 & 0.52\\
\cline{1-18}
 &  & 500 & 0.01 & 0.02 & 0.01 & 0.26 & 0.29 & 0.22 & 0.26 & 0.29 & 0.22 & 0.95 & 0.94 & 0.96 & 0.97 & 1.10 & 0.89\\

 &  & 1000 & 0.01 & 0.03 & 0.03 & 0.20 & 0.20 & 0.17 & 0.20 & 0.20 & 0.17 & 0.95 & 0.95 & 0.96 & 0.76 & 0.79 & 0.68\\

 &  & 1500 & 0.02 & 0.02 & 0.00 & 0.17 & 0.17 & 0.16 & 0.17 & 0.17 & 0.16 & 0.95 & 0.96 & 0.95 & 0.68 & 0.69 & 0.59\\

\multirow{-4}{*}{\centering 0.25} & \multirow{-4}{*}{\centering -0.35} & 2000 & 0.02 & 0.02 & 0.01 & 0.16 & 0.17 & 0.14 & 0.16 & 0.17 & 0.13 & 0.95 & 0.95 & 0.96 & 0.62 & 0.63 & 0.53\\
\cline{1-18}
 &  & 500 & 0.01 & 0.02 & 0.02 & 0.37 & 0.44 & 0.34 & 0.37 & 0.44 & 0.34 & 0.94 & 0.93 & 0.92 & 1.37 & 1.54 & 1.23\\

 &  & 1000 & 0.04 & 0.00 & 0.01 & 0.28 & 0.28 & 0.25 & 0.28 & 0.28 & 0.25 & 0.95 & 0.95 & 0.91 & 1.06 & 1.11 & 0.94\\

 &  & 1500 & 0.02 & 0.02 & 0.04 & 0.25 & 0.26 & 0.22 & 0.25 & 0.26 & 0.21 & 0.96 & 0.94 & 0.94 & 0.94 & 0.96 & 0.82\\

\multirow{-4}{*}{\centering 1.00} & \multirow{-4}{*}{\centering -0.15} & 2000 & 0.01 & 0.01 & 0.02 & 0.23 & 0.23 & 0.20 & 0.23 & 0.23 & 0.20 & 0.93 & 0.96 & 0.93 & 0.86 & 0.88 & 0.74\\
\hline

\multicolumn{18}{c}{Setting 2, approximately sparse: $f(d) = 2\exp(-d/2)$}\\
\hline
\multicolumn{1}{|c}{ } & \multicolumn{1}{c}{ } & \multicolumn{1}{c|}{ } & \multicolumn{3}{c|}{Bias} & \multicolumn{3}{c|}{RMSE} & \multicolumn{3}{c|}{SE} & \multicolumn{3}{c|}{Coverage} & \multicolumn{3}{c|}{Length} \\
\hline
$\eval$ & True & $n$ & \texttt{DLL} & \texttt{Swap} & \texttt{Trans} & \texttt{DLL} & \texttt{Swap} & \texttt{Trans} & \texttt{DLL} & \texttt{Swap} & \texttt{Trans} & \texttt{DLL} & \texttt{Swap} & \texttt{Trans} & \texttt{DLL} & \texttt{Swap} & \texttt{Trans}\\
\hline
 &  & 500 & 0.01 & 0.05 & 0.01 & 0.29 & 0.36 & 0.28 & 0.29 & 0.36 & 0.28 & 0.95 & 0.94 & 0.96 & 1.11 & 1.34 & 1.03\\

 &  & 1000 & 0.01 & 0.03 & 0.00 & 0.22 & 0.25 & 0.21 & 0.22 & 0.25 & 0.21 & 0.96 & 0.95 & 0.94 & 0.89 & 1.00 & 0.79\\

 &  & 1500 & 0.02 & 0.03 & 0.01 & 0.21 & 0.23 & 0.18 & 0.21 & 0.23 & 0.18 & 0.95 & 0.95 & 0.94 & 0.80 & 0.87 & 0.69\\

\multirow{-4}{*}{\centering -1.25} & \multirow{-4}{*}{\centering -0.41} & 2000 & 0.01 & 0.04 & 0.01 & 0.19 & 0.21 & 0.17 & 0.19 & 0.20 & 0.17 & 0.94 & 0.95 & 0.94 & 0.73 & 0.79 & 0.62\\
\cline{1-18}
 &  & 500 & 0.05 & 0.01 & 0.04 & 0.24 & 0.27 & 0.23 & 0.23 & 0.27 & 0.22 & 0.92 & 0.94 & 0.93 & 0.87 & 1.05 & 0.81\\

 &  & 1000 & 0.03 & 0.00 & 0.02 & 0.19 & 0.20 & 0.17 & 0.19 & 0.20 & 0.17 & 0.93 & 0.95 & 0.92 & 0.70 & 0.79 & 0.62\\

 &  & 1500 & 0.02 & 0.01 & 0.02 & 0.17 & 0.18 & 0.14 & 0.17 & 0.18 & 0.14 & 0.92 & 0.94 & 0.95 & 0.63 & 0.68 & 0.54\\

\multirow{-4}{*}{\centering -0.50} & \multirow{-4}{*}{\centering -0.52} & 2000 & 0.01 & 0.00 & 0.02 & 0.15 & 0.16 & 0.13 & 0.15 & 0.16 & 0.13 & 0.95 & 0.96 & 0.92 & 0.57 & 0.62 & 0.49\\
\cline{1-18}
 &  & 500 & 0.01 & 0.02 & 0.05 & 0.23 & 0.27 & 0.23 & 0.23 & 0.27 & 0.22 & 0.95 & 0.96 & 0.93 & 0.89 & 1.07 & 0.82\\

 &  & 1000 & 0.00 & 0.01 & 0.02 & 0.19 & 0.23 & 0.18 & 0.19 & 0.23 & 0.18 & 0.93 & 0.93 & 0.92 & 0.71 & 0.79 & 0.63\\

 &  & 1500 & 0.01 & 0.00 & 0.01 & 0.16 & 0.18 & 0.14 & 0.16 & 0.18 & 0.14 & 0.95 & 0.95 & 0.95 & 0.64 & 0.70 & 0.55\\

\multirow{-4}{*}{\centering 0.10} & \multirow{-4}{*}{\centering -0.40} & 2000 & 0.01 & 0.00 & 0.02 & 0.15 & 0.15 & 0.13 & 0.15 & 0.16 & 0.12 & 0.95 & 0.97 & 0.95 & 0.58 & 0.63 & 0.50\\
\cline{1-18}
 &  & 500 & 0.00 & 0.03 & 0.05 & 0.23 & 0.28 & 0.24 & 0.23 & 0.28 & 0.23 & 0.96 & 0.96 & 0.94 & 0.91 & 1.11 & 0.84\\

 &  & 1000 & 0.00 & 0.01 & 0.00 & 0.20 & 0.23 & 0.16 & 0.20 & 0.23 & 0.16 & 0.92 & 0.94 & 0.95 & 0.74 & 0.82 & 0.65\\

 &  & 1500 & 0.01 & 0.01 & 0.00 & 0.16 & 0.18 & 0.15 & 0.16 & 0.18 & 0.15 & 0.96 & 0.95 & 0.94 & 0.66 & 0.72 & 0.57\\

\multirow{-4}{*}{\centering 0.25} & \multirow{-4}{*}{\centering -0.35} & 2000 & 0.02 & 0.01 & 0.00 & 0.15 & 0.16 & 0.13 & 0.15 & 0.16 & 0.13 & 0.96 & 0.96 & 0.95 & 0.60 & 0.65 & 0.52\\
\cline{1-18}
 &  & 500 & 0.01 & 0.07 & 0.02 & 0.35 & 0.41 & 0.34 & 0.35 & 0.40 & 0.34 & 0.94 & 0.95 & 0.94 & 1.28 & 1.55 & 1.19\\

 &  & 1000 & 0.02 & 0.01 & 0.00 & 0.29 & 0.32 & 0.23 & 0.29 & 0.32 & 0.23 & 0.94 & 0.94 & 0.95 & 1.03 & 1.15 & 0.90\\

 &  & 1500 & 0.01 & 0.01 & 0.01 & 0.23 & 0.25 & 0.20 & 0.23 & 0.25 & 0.20 & 0.93 & 0.96 & 0.95 & 0.92 & 1.00 & 0.79\\

\multirow{-4}{*}{\centering 1.00} & \multirow{-4}{*}{\centering -0.15} & 2000 & 0.00 & 0.01 & 0.02 & 0.22 & 0.23 & 0.18 & 0.22 & 0.23 & 0.18 & 0.95 & 0.95 & 0.96 & 0.84 & 0.90 & 0.72\\
\hline
\end{tabular}}
\caption*{Table S11: Comparison of \texttt{DLL}, \texttt{DLL} with data swapping (\texttt{Swap}), \texttt{DLL} with quantile transformation (\texttt{Trans}) in Setting 2 when $f(d) = 2\exp(-d/2)$, across different sample sizes $n$ and evaluation points $\eval$. The column indexed with ``True'' represents the true value of $f'(\eval)$. The columns indexed with ``Bias'', ``RMSE'' and ``SE" report the absolute bias, the root mean square error, and the standard error computed by 500 estimates, respectively; the columns indexed with ``Coverage'' report the empirical coverage level and the columns indexed with ``Length'' report the average CI length.}
\label{table:S11}
\end{table}

\subsection{Comparison with ReSmoothing Method}
\label{sec: compare RS}
For comparison with \texttt{RS} method, we generate the data as in Setting 1 and we present the full results with $f(d) = 1.5\sin(d), g_1(x) = 2\exp(-x/2)$ or $f(d) = 2\exp(-d/2), g_1(x) = 1.5\sin(x)$ in Table \hyperref[table:S12]{S12}. Our \texttt{DLL} method achieves the desired coverage, but \texttt{RS} and \texttt{OraRS} do not.  Even though the standard error for \texttt{OraRS} is computed in an oracle way, the confidence interval is still undercoverage since the \texttt{RS} estimator suffers from a large bias. In addition, the CI length for \texttt{OraRS} is large and the length of our \texttt{DLL} method is similar to the oracle confidence interval.

In addition, we generate the data following the simulation setting of \cite{gregory2016optimal} and compare \texttt{DLL}, \texttt{RS}, and \texttt{OraRS}. Specifically, the outcome model is $Y_i = \sum_{j=1}^pf_j(X_{ij}) + \epsilon_i$ with four non-zero functions:
\begin{equation*}
    \begin{aligned}
         f_1(x) &= -\sin(2x); & f_2(x) &= x^2-25/12; & f_3(x) &= x; & f_4(x) = e^{-x}-2/5\sinh(5/2).
         \end{aligned}
\end{equation*}
and $f_j(x)=0$ for $j \geq 5$. The sample size $n$ is varied across \{100,1000\} and the dimension $p$ is varied across \{50,150\}. For $1\leq i\leq n$, we generate $\{X_{i,j}\}_{1\leq j\leq p}$ as follows: the marginal distribution of $X_{i,j}$ is Uniform$(-2.5,2.5)$ and the correlation between $X_{i,j}$ and $X_{i,l}$ is $r^{|j-l|}$ for $1 \leq j \neq l \leq p$, where $r$ is a correlation parameter varied across \{0,0.1,0.3,0.5\}. The error term $\epsilon_i \sim N(0,1)$. We estimate the derivatives of different functions $\{f_j(x_0)\}_{j \in \mathcal{J}}$ where $\mathcal{J} = \{1,2,3,4,5\}$ and $x_0$ is varied across \{-1, 0.5\}. 

The results for the empirical coverage are reported in Table \hyperref[table:S13]{S13}. Our \texttt{DLL} method achieves the desired coverage across all sample sizes, dimensions, and correlation parameters. \texttt{OraRS} does not attain the expected coverage when the sample size is small, even if its standard error is computed in an oracle way. The \texttt{RS} method does not have coverage in almost all settings. For a small sample size, the \texttt{RS} estimator suffers from a large bias while our proposed \texttt{DLL} effectively corrects the bias in these settings. We report the results for the absolute bias in Table \hyperref[table:S14]{S14}.

%%%%%%%%%%%%%%%%%%%%%%%%%%%%%%%
% comparison settings
%%%%%%%%%%%%%%%%%%%%%%%%%%%%%%%

\begin{table}[htb!]
\centering
\resizebox{\linewidth}{!}{
\begin{tabular}[t]{|c|c|c|ccc|ccc|cccc|cccc|}
\multicolumn{17}{c}{Setting 1, exactly sparse: Comparison with ReSmoothing}\\
\hline
\multicolumn{1}{|c}{ } & \multicolumn{1}{c}{ } & \multicolumn{1}{c|}{ } & \multicolumn{3}{c|}{Bias} & \multicolumn{3}{c|}{SE} & \multicolumn{4}{c|}{Coverage} & \multicolumn{4}{c|}{Length} \\
\hline
$\eval$ & $f$ & $n$ & \texttt{DLL} & \texttt{RS} & \texttt{Orac} & \texttt{DLL} & \texttt{RS} & \texttt{Orac} & \texttt{DLL} & \texttt{RS} & \texttt{OraRS} & \texttt{Orac} & \texttt{DLL} & \texttt{RS} & \texttt{OraRS} & \texttt{Orac}\\
\hline
 &  & 500 & 0.05 & 0.42 & 0.03 & 0.39 & 0.94 & 0.35 & 0.95 & 0.01 & 0.91 & 0.95 & 1.50 & 0.05 & 3.68 & 1.42\\

 &  & 750 & 0.02 & 0.38 & 0.02 & 0.36 & 0.91 & 0.34 & 0.94 & 0.01 & 0.92 & 0.94 & 1.37 & 0.03 & 3.58 & 1.29\\

 & \multirow{-3}{*}{\centering $exp$} & 1000 & 0.01 & 0.22 & 0.01 & 0.34 & 0.87 & 0.32 & 0.94 & 0.01 & 0.95 & 0.93 & 1.27 & 0.02 & 3.41 & 1.20\\

 &  & 500 & 0.20 & 0.91 & 0.01 & 0.42 & 0.83 & 0.41 & 0.94 & 0.01 & 0.80 & 0.96 & 1.69 & 0.04 & 3.24 & 1.68\\

 &  & 750 & 0.07 & 0.69 & 0.01 & 0.40 & 0.90 & 0.39 & 0.93 & 0.01 & 0.88 & 0.95 & 1.55 & 0.02 & 3.53 & 1.51\\

\multirow{-6}{*}{\centering -1.0} & \multirow{-3}{*}{\centering $sin$} & 1000 & 0.05 & 0.53 & 0.01 & 0.35 & 0.89 & 0.35 & 0.96 & 0.01 & 0.92 & 0.95 & 1.46 & 0.02 & 3.48 & 1.42\\
\cline{1-17}
 &  & 500 & 0.00 & 0.49 & 0.01 & 0.38 & 1.01 & 0.37 & 0.96 & 0.01 & 0.93 & 0.94 & 1.51 & 0.04 & 3.97 & 1.41\\

 &  & 750 & 0.00 & 0.37 & 0.01 & 0.35 & 0.90 & 0.33 & 0.96 & 0.01 & 0.92 & 0.96 & 1.36 & 0.02 & 3.53 & 1.27\\

 & \multirow{-3}{*}{\centering $exp$} & 1000 & 0.02 & 0.26 & 0.01 & 0.31 & 0.84 & 0.31 & 0.97 & 0.00 & 0.93 & 0.93 & 1.27 & 0.01 & 3.30 & 1.20\\

 &  & 500 & 0.20 & 0.97 & 0.02 & 0.42 & 0.78 & 0.41 & 0.92 & 0.01 & 0.78 & 0.95 & 1.69 & 0.04 & 3.06 & 1.68\\

 &  & 750 & 0.09 & 0.66 & 0.01 & 0.39 & 0.91 & 0.40 & 0.95 & 0.01 & 0.89 & 0.93 & 1.57 & 0.02 & 3.57 & 1.53\\

\multirow{-6}{*}{\centering 0.5} & \multirow{-3}{*}{\centering $sin$} & 1000 & 0.09 & 0.51 & 0.03 & 0.37 & 0.87 & 0.37 & 0.94 & 0.01 & 0.92 & 0.95 & 1.46 & 0.02 & 3.41 & 1.41\\
\hline

\multicolumn{17}{c}{Setting 1, approximately sparse: Comparison with ReSmoothing}\\
\hline
\multicolumn{1}{|c}{ } & \multicolumn{1}{c}{ } & \multicolumn{1}{c|}{ } & \multicolumn{3}{c|}{Bias} & \multicolumn{3}{c|}{SE} & \multicolumn{4}{c|}{Coverage} & \multicolumn{4}{c|}{Length} \\
\hline
$\eval$ & $f$ & $n$ & \texttt{DLL} & \texttt{RS} & \texttt{Orac} & \texttt{DLL} & \texttt{RS} & \texttt{Orac} & \texttt{DLL} & \texttt{RS} & \texttt{OraRS} & \texttt{Orac} & \texttt{DLL} & \texttt{RS} & \texttt{OraRS} & \texttt{Orac}\\
\hline
 &  & 500 & 0.09 & 0.73 & 0.00 & 0.35 & 0.85 & 0.35 & 0.93 & 0.01 & 0.86 & 0.95 & 1.33 & 0.04 & 3.34 & 1.41\\

 &  & 750 & 0.07 & 0.59 & 0.03 & 0.32 & 0.88 & 0.33 & 0.95 & 0.01 & 0.90 & 0.95 & 1.26 & 0.02 & 3.46 & 1.28\\

 & \multirow{-3}{*}{\centering $exp$} & 1000 & 0.04 & 0.61 & 0.01 & 0.31 & 0.84 & 0.31 & 0.94 & 0.00 & 0.88 & 0.95 & 1.20 & 0.02 & 3.28 & 1.20\\

 &  & 500 & 0.26 & 0.80 & 0.02 & 0.39 & 0.46 & 0.46 & 0.86 & 0.00 & 0.65 & 0.91 & 1.46 & 0.02 & 1.81 & 1.67\\

 &  & 750 & 0.18 & 0.73 & 0.01 & 0.35 & 0.63 & 0.38 & 0.94 & 0.00 & 0.80 & 0.95 & 1.42 & 0.01 & 2.46 & 1.51\\

\multirow{-6}{*}{\centering -1.0} & \multirow{-3}{*}{\centering $sin$} & 1000 & 0.12 & 0.66 & 0.01 & 0.35 & 0.73 & 0.38 & 0.94 & 0.00 & 0.86 & 0.93 & 1.37 & 0.01 & 2.87 & 1.42\\
\cline{1-17}
 &  & 500 & 0.15 & 0.62 & 0.01 & 0.34 & 0.81 & 0.36 & 0.94 & 0.01 & 0.88 & 0.94 & 1.33 & 0.03 & 3.17 & 1.42\\

 &  & 750 & 0.06 & 0.53 & 0.04 & 0.33 & 0.89 & 0.33 & 0.94 & 0.01 & 0.90 & 0.95 & 1.26 & 0.02 & 3.47 & 1.28\\

 & \multirow{-3}{*}{\centering $exp$} & 1000 & 0.08 & 0.48 & 0.00 & 0.31 & 0.87 & 0.31 & 0.95 & 0.00 & 0.91 & 0.96 & 1.20 & 0.01 & 3.39 & 1.20\\

 &  & 500 & 0.29 & 1.08 & 0.02 & 0.40 & 0.44 & 0.46 & 0.86 & 0.00 & 0.26 & 0.93 & 1.47 & 0.02 & 1.72 & 1.68\\

 &  & 750 & 0.23 & 0.96 & 0.03 & 0.38 & 0.60 & 0.40 & 0.89 & 0.00 & 0.63 & 0.94 & 1.43 & 0.01 & 2.34 & 1.53\\

\multirow{-6}{*}{\centering 0.5} & \multirow{-3}{*}{\centering $sin$} & 1000 & 0.14 & 0.78 & 0.00 & 0.33 & 0.73 & 0.35 & 0.94 & 0.01 & 0.81 & 0.95 & 1.37 & 0.01 & 2.88 & 1.41\\
\hline
\end{tabular}
}
\caption*{Table S12: Comparison of \texttt{DLL}, ReSmoothing (\texttt{RS}), \texttt{OraRS}, and oracle (\texttt{Orac}) estimators in Setting 1 with $p=750$, across different sample sizes $n$, evaluation points $\eval$, and function of interest $f$. The columns indexed with ``Bias'', and ``SE" report the absolute bias, and the standard error computed by 500 estimates, respectively; the columns indexed with ``Coverage'' report the empirical coverage level and the columns indexed with ``Length'' report the average CI length.}
\label{table:S12}
\end{table}

\begin{table}[htb!]
\centering
\resizebox{\linewidth}{!}{
\begin{tabular}[t]{|c|c|c|ccc|ccc|ccc|ccc|ccc|}
\multicolumn{18}{c}{Coverage in the Setting of \cite{gregory2016optimal}: $\eval$ = -1}\\
\hline
\multicolumn{1}{|c}{ } & \multicolumn{1}{c}{ } & \multicolumn{1}{c|}{ } & \multicolumn{3}{c|}{$f_1$} & \multicolumn{3}{c|}{$f_2$} & \multicolumn{3}{c|}{$f_3$} & \multicolumn{3}{c|}{$f_4$} & \multicolumn{3}{c|}{$f_5$} \\
\hline
$n$ & $p$ & $r$ & \texttt{DLL} & \texttt{RS} & \texttt{OraRS} & \texttt{DLL} & \texttt{RS} & \texttt{OraRS} & \texttt{DLL} & \texttt{RS} & \texttt{OraRS} & \texttt{DLL} & \texttt{RS} & \texttt{OraRS} & \texttt{DLL} & \texttt{RS} & \texttt{OraRS}\\
\hline
 &  & 0.0 & 0.94 & 0.09 & 0.84 & 0.93 & 0.16 & 0.83 & 0.95 & 0.17 & 0.90 & 0.96 & 0.21 & 0.89 & 0.96 & 0.25 & 0.95\\

 &  & 0.1 & 0.96 & 0.09 & 0.79 & 0.93 & 0.13 & 0.83 & 0.96 & 0.15 & 0.88 & 0.94 & 0.19 & 0.90 & 0.93 & 0.22 & 0.94\\

 &  & 0.3 & 0.94 & 0.07 & 0.85 & 0.95 & 0.11 & 0.89 & 0.95 & 0.09 & 0.87 & 0.95 & 0.17 & 0.87 & 0.94 & 0.21 & 0.92\\

 & \multirow{-4}{*}{\centering 50} & 0.5 & 0.94 & 0.09 & 0.90 & 0.94 & 0.11 & 0.88 & 0.93 & 0.04 & 0.76 & 0.94 & 0.14 & 0.83 & 0.94 & 0.24 & 0.93\\
 \cline{2-18}

 &  & 0.0 & 0.94 & 0.06 & 0.43 & 0.94 & 0.11 & 0.82 & 0.95 & 0.12 & 0.91 & 0.94 & 0.17 & 0.87 & 0.95 & 0.23 & 0.94\\

 &  & 0.1 & 0.92 & 0.03 & 0.87 & 0.93 & 0.09 & 0.74 & 0.94 & 0.11 & 0.88 & 0.95 & 0.15 & 0.87 & 0.96 & 0.22 & 0.94\\

 &  & 0.3 & 0.91 & 0.07 & 0.92 & 0.93 & 0.11 & 0.84 & 0.93 & 0.07 & 0.82 & 0.95 & 0.13 & 0.85 & 0.95 & 0.25 & 0.93\\

\multirow{-8}{*}{\centering 100} & \multirow{-4}{*}{\centering 150} & 0.5 & 0.93 & 0.08 & 0.92 & 0.92 & 0.15 & 0.84 & 0.92 & 0.05 & 0.73 & 0.94 & 0.10 & 0.79 & 0.93 & 0.22 & 0.95\\
\cline{1-18}
 &  & 0.0 & 0.96 & 0.01 & 0.95 & 0.95 & 0.06 & 0.95 & 0.97 & 0.00 & 0.94 & 0.96 & 0.07 & 0.95 & 0.94 & 0.01 & 0.95\\

 &  & 0.1 & 0.96 & 0.01 & 0.95 & 0.97 & 0.05 & 0.95 & 0.96 & 0.00 & 0.94 & 0.97 & 0.06 & 0.95 & 0.96 & 0.00 & 0.95\\

 &  & 0.3 & 0.96 & 0.01 & 0.95 & 0.96 & 0.05 & 0.94 & 0.94 & 0.01 & 0.95 & 0.94 & 0.07 & 0.94 & 0.95 & 0.00 & 0.95\\

 & \multirow{-4}{*}{\centering 50} & 0.5 & 0.96 & 0.01 & 0.96 & 0.95 & 0.05 & 0.96 & 0.95 & 0.00 & 0.93 & 0.95 & 0.07 & 0.95 & 0.94 & 0.00 & 0.94\\
  \cline{2-18}

 &  & 0.0 & 0.96 & 0.01 & 0.94 & 0.95 & 0.05 & 0.95 & 0.93 & 0.01 & 0.95 & 0.95 & 0.09 & 0.95 & 0.94 & 0.00 & 0.94\\

 &  & 0.1 & 0.94 & 0.01 & 0.95 & 0.97 & 0.04 & 0.95 & 0.95 & 0.01 & 0.94 & 0.96 & 0.09 & 0.93 & 0.94 & 0.00 & 0.95\\

 &  & 0.3 & 0.95 & 0.02 & 0.94 & 0.96 & 0.04 & 0.95 & 0.95 & 0.01 & 0.95 & 0.96 & 0.06 & 0.94 & 0.95 & 0.00 & 0.94\\

\multirow{-8}{*}{\centering 1000} & \multirow{-4}{*}{\centering 150} & 0.5 & 0.96 & 0.01 & 0.95 & 0.94 & 0.05 & 0.95 & 0.95 & 0.01 & 0.95 & 0.96 & 0.08 & 0.95 & 0.94 & 0.02 & 0.94\\
\hline

\multicolumn{18}{c}{Coverage in the Setting of \cite{gregory2016optimal}: $\eval$ = 0.5}\\
\hline
\multicolumn{1}{|c}{ } & \multicolumn{1}{c}{ } & \multicolumn{1}{c|}{ } & \multicolumn{3}{c|}{$f_1$} & \multicolumn{3}{c|}{$f_2$} & \multicolumn{3}{c|}{$f_3$} & \multicolumn{3}{c|}{$f_4$} & \multicolumn{3}{c|}{$f_5$} \\
\hline
$n$ & $p$ & $r$ & \texttt{DLL} & \texttt{RS} & \texttt{OraRS} & \texttt{DLL} & \texttt{RS} & \texttt{OraRS} & \texttt{DLL} & \texttt{RS} & \texttt{OraRS} & \texttt{DLL} & \texttt{RS} & \texttt{OraRS} & \texttt{DLL} & \texttt{RS} & \texttt{OraRS}\\
\hline
 &  & 0.0 & 0.95 & 0.06 & 0.44 & 0.94 & 0.21 & 0.92 & 0.93 & 0.14 & 0.88 & 0.97 & 0.22 & 0.94 & 0.96 & 0.21 & 0.95\\

 &  & 0.1 & 0.95 & 0.07 & 0.56 & 0.95 & 0.21 & 0.93 & 0.96 & 0.14 & 0.90 & 0.94 & 0.23 & 0.94 & 0.95 & 0.25 & 0.94\\

 &  & 0.3 & 0.93 & 0.07 & 0.81 & 0.94 & 0.22 & 0.94 & 0.95 & 0.11 & 0.85 & 0.95 & 0.23 & 0.94 & 0.96 & 0.18 & 0.92\\

 & \multirow{-4}{*}{\centering 50} & 0.5 & 0.96 & 0.07 & 0.87 & 0.95 & 0.21 & 0.95 & 0.96 & 0.07 & 0.75 & 0.94 & 0.20 & 0.92 & 0.95 & 0.23 & 0.94\\
  \cline{2-18}

 &  & 0.0 & 0.92 & 0.02 & 0.30 & 0.94 & 0.21 & 0.92 & 0.95 & 0.13 & 0.87 & 0.97 & 0.25 & 0.95 & 0.95 & 0.23 & 0.94\\

 &  & 0.1 & 0.90 & 0.03 & 0.23 & 0.95 & 0.18 & 0.92 & 0.95 & 0.09 & 0.83 & 0.96 & 0.18 & 0.95 & 0.94 & 0.23 & 0.94\\

 &  & 0.3 & 0.90 & 0.06 & 0.80 & 0.95 & 0.20 & 0.95 & 0.93 & 0.06 & 0.83 & 0.93 & 0.19 & 0.94 & 0.95 & 0.23 & 0.95\\

\multirow{-8}{*}{\centering 100} & \multirow{-4}{*}{\centering 150} & 0.5 & 0.92 & 0.06 & 0.51 & 0.95 & 0.21 & 0.93 & 0.93 & 0.06 & 0.82 & 0.93 & 0.20 & 0.93 & 0.96 & 0.19 & 0.93\\
\cline{1-18}
 &  & 0.0 & 0.96 & 0.03 & 0.94 & 0.94 & 0.03 & 0.94 & 0.95 & 0.00 & 0.95 & 0.95 & 0.03 & 0.95 & 0.96 & 0.01 & 0.94\\

 &  & 0.1 & 0.97 & 0.01 & 0.95 & 0.95 & 0.03 & 0.95 & 0.94 & 0.00 & 0.94 & 0.96 & 0.02 & 0.95 & 0.95 & 0.00 & 0.94\\

 &  & 0.3 & 0.96 & 0.02 & 0.96 & 0.96 & 0.04 & 0.95 & 0.96 & 0.02 & 0.94 & 0.98 & 0.03 & 0.95 & 0.94 & 0.01 & 0.94\\

 & \multirow{-4}{*}{\centering 50} & 0.5 & 0.95 & 0.02 & 0.94 & 0.95 & 0.03 & 0.95 & 0.96 & 0.00 & 0.94 & 0.96 & 0.03 & 0.95 & 0.96 & 0.01 & 0.94\\
  \cline{2-18}

 &  & 0.0 & 0.97 & 0.02 & 0.95 & 0.97 & 0.03 & 0.94 & 0.96 & 0.01 & 0.95 & 0.97 & 0.02 & 0.95 & 0.96 & 0.00 & 0.95\\

 &  & 0.1 & 0.93 & 0.01 & 0.97 & 0.95 & 0.03 & 0.94 & 0.96 & 0.01 & 0.96 & 0.96 & 0.03 & 0.94 & 0.95 & 0.01 & 0.94\\

 &  & 0.3 & 0.95 & 0.01 & 0.95 & 0.95 & 0.03 & 0.94 & 0.96 & 0.01 & 0.95 & 0.97 & 0.04 & 0.94 & 0.95 & 0.01 & 0.94\\

\multirow{-8}{*}{\centering 1000} & \multirow{-4}{*}{\centering 150} & 0.5 & 0.94 & 0.01 & 0.95 & 0.94 & 0.03 & 0.95 & 0.96 & 0.01 & 0.94 & 0.96 & 0.03 & 0.96 & 0.96 & 0.01 & 0.95\\
\hline
\end{tabular}}
\caption*{Table S13: Comparison of coverage of \texttt{DLL}, ReSmoothing (\texttt{RS}), and \texttt{OraRS} in the setting of \cite{gregory2016optimal}, across different sample sizes $n$, dimension of covariates $p$, and the correlation parameter $r$. The $f_1$, $f_2$, $f_3$, $f_4$ and $f_5$ represent the functions of interest to estimate the derivatives. The entries of the table represents the empirical coverage across 500 simulations.}
\label{table:S13}
\end{table}

\begin{table}[htb!]
\centering
\scalebox{0.85}{
\begin{tabular}[t]{|c|c|c|cc|cc|cc|cc|cc|}
\multicolumn{13}{c}{Bias in the Setting of \cite{gregory2016optimal}: $\eval$ = -1}\\
\hline
\multicolumn{1}{|c}{ } & \multicolumn{1}{c}{ } & \multicolumn{1}{c|}{ } & \multicolumn{2}{c|}{$f_1$} & \multicolumn{2}{c|}{$f_2$} & \multicolumn{2}{c|}{$f_3$} & \multicolumn{2}{c|}{$f_4$} & \multicolumn{2}{c|}{$f_5$} \\
\hline
$n$ & $p$ & $r$ & DLL & RS & DLL & RS & DLL & RS & DLL & RS & DLL & RS\\
\hline
 &  & 0.0 & 0.18 & 0.66 & 0.00 & 0.90 & 0.05 & 0.53 & 0.14 & 0.71 & 0.06 & 0.04\\

 &  & 0.1 & 0.20 & 0.64 & 0.13 & 0.94 & 0.03 & 0.57 & 0.14 & 0.80 & 0.00 & 0.04\\

 &  & 0.3 & 0.22 & 0.66 & 0.40 & 0.95 & 0.18 & 0.71 & 0.15 & 0.90 & 0.01 & 0.03\\

 & \multirow{-4}{*}{\centering 50} & 0.5 & 0.14 & 0.61 & 0.13 & 0.91 & 0.05 & 0.94 & 0.13 & 1.06 & 0.01 & 0.11\\
 \cline{2-13}

 &  & 0.0 & 0.34 & 0.70 & 0.07 & 1.03 & 0.04 & 0.61 & 0.16 & 0.84 & 0.08 & 0.00\\

 &  & 0.1 & 0.22 & 0.78 & 0.20 & 1.13 & 0.01 & 0.65 & 0.23 & 0.94 & 0.01 & 0.00\\

 &  & 0.3 & 0.29 & 0.69 & 0.06 & 1.02 & 0.27 & 0.84 & 0.13 & 1.02 & 0.01 & 0.07\\

\multirow{-8}{*}{\centering 100} & \multirow{-4}{*}{\centering 150} & 0.5 & 0.30 & 0.61 & 0.08 & 0.96 & 0.22 & 0.91 & 0.01 & 1.24 & 0.08 & 0.10\\
\cline{1-13}
 &  & 0.0 & 0.01 & 0.04 & 0.01 & 0.01 & 0.01 & 0.01 & 0.02 & 0.04 & 0.01 & 0.01\\

 &  & 0.1 & 0.02 & 0.03 & 0.03 & 0.08 & 0.05 & 0.06 & 0.01 & 0.05 & 0.02 & 0.03\\

 &  & 0.3 & 0.04 & 0.05 & 0.04 & 0.14 & 0.03 & 0.06 & 0.02 & 0.06 & 0.03 & 0.02\\

 & \multirow{-4}{*}{\centering 50} & 0.5 & 0.05 & 0.08 & 0.02 & 0.09 & 0.01 & 0.16 & 0.09 & 0.12 & 0.03 & 0.04\\
  \cline{2-13}

 &  & 0.0 & 0.07 & 0.05 & 0.00 & 0.05 & 0.03 & 0.04 & 0.10 & 0.01 & 0.00 & 0.00\\

 &  & 0.1 & 0.05 & 0.03 & 0.05 & 0.09 & 0.02 & 0.02 & 0.11 & 0.10 & 0.02 & 0.04\\

 &  & 0.3 & 0.01 & 0.03 & 0.03 & 0.05 & 0.04 & 0.10 & 0.01 & 0.11 & 0.03 & 0.02\\

\multirow{-8}{*}{\centering 1000} & \multirow{-4}{*}{\centering 150} & 0.5 & 0.02 & 0.01 & 0.01 & 0.16 & 0.03 & 0.18 & 0.03 & 0.18 & 0.00 & 0.04\\
\hline

\multicolumn{13}{c}{Bias in the Setting of \cite{gregory2016optimal}: $\eval$ = 0.5}\\
\hline
\multicolumn{1}{|c}{ } & \multicolumn{1}{c}{ } & \multicolumn{1}{c|}{ } & \multicolumn{2}{c|}{$f_1$} & \multicolumn{2}{c|}{$f_2$} & \multicolumn{2}{c|}{$f_3$} & \multicolumn{2}{c|}{$f_4$} & \multicolumn{2}{c|}{$f_5$} \\
\hline
$n$ & $p$ & $r$ & \texttt{DLL} & \texttt{RS} & \texttt{DLL} & \texttt{RS} & \texttt{DLL} & \texttt{RS} & \texttt{DLL} & \texttt{RS} & \texttt{DLL} & \texttt{RS}\\
\hline
 &  & 0.0 & 0.12 & 0.82 & 0.04 & 0.40 & 0.08 & 0.55 & 0.12 & 0.19 & 0.04 & 0.00\\

 &  & 0.1 & 0.25 & 0.85 & 0.05 & 0.41 & 0.07 & 0.57 & 0.09 & 0.17 & 0.02 & 0.01\\

 &  & 0.3 & 0.27 & 0.83 & 0.04 & 0.28 & 0.09 & 0.70 & 0.11 & 0.28 & 0.01 & 0.01\\

 & \multirow{-4}{*}{\centering 50} & 0.5 & 0.23 & 0.82 & 0.00 & 0.25 & 0.12 & 0.81 & 0.22 & 0.51 & 0.04 & 0.08\\
  \cline{2-13}

 &  & 0.0 & 0.16 & 0.93 & 0.09 & 0.52 & 0.11 & 0.63 & 0.17 & 0.20 & 0.08 & 0.00\\

 &  & 0.1 & 0.29 & 0.94 & 0.02 & 0.46 & 0.10 & 0.72 & 0.17 & 0.20 & 0.08 & 0.01\\

 &  & 0.3 & 0.29 & 0.92 & 0.13 & 0.38 & 0.19 & 0.82 & 0.27 & 0.40 & 0.06 & 0.02\\

\multirow{-8}{*}{\centering 100} & \multirow{-4}{*}{\centering 150} & 0.5 & 0.33 & 0.87 & 0.08 & 0.28 & 0.16 & 0.84 & 0.31 & 0.49 & 0.03 & 0.07\\
\cline{1-13}
 &  & 0.0 & 0.01 & 0.01 & 0.00 & 0.05 & 0.01 & 0.03 & 0.02 & 0.02 & 0.05 & 0.03\\

 &  & 0.1 & 0.06 & 0.04 & 0.01 & 0.03 & 0.01 & 0.03 & 0.04 & 0.03 & 0.02 & 0.01\\

 &  & 0.3 & 0.01 & 0.02 & 0.00 & 0.04 & 0.05 & 0.09 & 0.08 & 0.08 & 0.04 & 0.02\\

 & \multirow{-4}{*}{\centering 50} & 0.5 & 0.10 & 0.05 & 0.04 & 0.01 & 0.04 & 0.12 & 0.04 & 0.10 & 0.02 & 0.02\\
  \cline{2-13}

 &  & 0.0 & 0.00 & 0.00 & 0.01 & 0.01 & 0.02 & 0.02 & 0.02 & 0.02 & 0.00 & 0.02\\

 &  & 0.1 & 0.05 & 0.06 & 0.01 & 0.01 & 0.02 & 0.03 & 0.08 & 0.08 & 0.05 & 0.01\\

 &  & 0.3 & 0.07 & 0.07 & 0.05 & 0.03 & 0.03 & 0.01 & 0.09 & 0.13 & 0.02 & 0.01\\

\multirow{-8}{*}{\centering 1000} & \multirow{-4}{*}{\centering 150} & 0.5 & 0.04 & 0.01 & 0.04 & 0.02 & 0.00 & 0.11 & 0.06 & 0.16 & 0.04 & 0.06\\
\hline
\end{tabular}
}
\caption*{Table S14: Comparison of bias of \texttt{DLL}, ReSmoothing (\texttt{RS}) in the setting of \cite{gregory2016optimal}, across different sample sizes $n$, dimension of covariates $p$, and the correlation parameter $r$. The $f_1$, $f_2$, $f_3$, $f_4$ and $f_5$ represent the functions of interest to estimate the derivatives. The entries of the table represents the bias across 500 simulations.}
\label{table:S14}
\end{table}

\end{document}